\documentclass[a4paper,11pt]{article}
\usepackage[letterpaper,margin=0.80in,bottom=0.90in,top=0.85in]{geometry}
\usepackage[latin,english]{babel}
\usepackage[utf8]{inputenc}
\usepackage[T1]{fontenc}
\usepackage{eufrak, color, mathrsfs,dsfont}
\usepackage{geometry}
\usepackage{amsmath}
\usepackage{amssymb}
\usepackage{amsthm}
\usepackage{mathtools}
\usepackage{mathabx}
\usepackage{cases}
\usepackage{braket}
\usepackage[toc,page]{appendix}
\usepackage{pdfsync}
\usepackage{hyperref}
\usepackage{psfrag}
\usepackage[mathscr]{euscript}

\newcommand{\perd}{\mathtt{d}}
\everymath{\displaystyle}
\usepackage{amssymb,amsmath,amsthm}
\usepackage{graphicx}

\newcommand{\tOm}{\mathtt \Omega}
 \newcommand{\bw}{\boldsymbol{\omega} }
\def\N{\mathbb N}
\def\R{\mathbb R}
\usepackage{amsfonts}
\usepackage{latexsym}
\usepackage{amsopn}
\usepackage{amscd}
\usepackage{mathrsfs}
\usepackage{dsfont}
\usepackage[english]{babel}
\usepackage{caption}

\numberwithin{equation}{section}

\font\teneufm=eufm10
\font\seveneufm=eufm7
\font\fiveeufm=eufm5
\newfam\eufmfam
\textfont\eufmfam=\teneufm
\scriptfont\eufmfam=\seveneufm
\scriptscriptfont\eufmfam=\fiveeufm

\newcommand{\zug}{{k_0,\upsilon}}
\newcommand{\Ke}{M}
\newcommand{\normk}[2]{\| #1 \|_{#2}^{k_0,\upsilon}}

\newcommand{\tK}{{\mathtt K}}
\newcommand{\cK}{ {\mathcal K}}
 \newcommand{\sK}{{\mathscr K}}
  \newcommand{\sN}{{\mathscr N}}
  \newcommand{\sP}{{\mathscr P}}
\newcommand{\W}{ W}

\newcommand{\acca}{\frak H}
\newcommand{\we}{{\mathtt w}_0}

\newcommand{\Om}{\Omega}
\newcommand{\la}{\langle}
\newcommand{\ra}{\rangle}
\newcommand{\tn}{\mathtt n}
\newcommand{\tb}{\mathtt b}
 \newcommand{\mG}{\mathcal G}
\newcommand{\ta}{\mathtt a}
\newcommand{\bl}{[\![}
\newcommand{\br}{]\!]}

\newcommand{\tc}{\mathtt c}
\newcommand{\ts}{\mathtt s}
\newcommand{\ka}{\kappa}

\newcommand{\tg}{{\mathtt g}}
\newcommand{\IK}{\mathfrak G}
\newcommand{\fB}{\mathfrak B}
\newcommand{\symm}{M}

\newenvironment{pf}{\noindent{\sc Proof}.\enspace}{\rule{2mm}{2mm}\smallskip}
\newenvironment{pfn}{\noindent{\sc Proof} \enspace}{\rule{2mm}{2mm}\smallskip}
\newtheorem{theorem}{Theorem}[section]
\newtheorem{proposition}{Proposition}[section]
\newtheorem{lemma}[proposition]{Lemma}
\newtheorem{corollary}[proposition]{Corollary}

\newtheorem{remark}[proposition]{Remark}
\newtheorem{remarks}[proposition]{Remark}
\newtheorem{definition}[proposition]{Definition}
\newcommand{\be}{\begin{equation}}
\newcommand{\ee}{\end{equation}}
\newcommand{\teta}{\theta}

\newcommand{\tLm}{{\mathtt \Lambda} }

\newcommand{\om}{\omega}
\newcommand{\ph}{\varphi}
\newcommand{\e}{\varepsilon}
\newcommand{\mF}{\mathcal{F}}

\newcommand{\bphi}{{\mathfrak S}}
\newcommand{\mb}{{\mathfrak b}}

\newcommand{\mB}{\mathcal{B}}

\newcommand{\Dom}{\om \cdot \pa_\vphi}
\newcommand{\pa}{\partial}
\newcommand{\ov}{\overline}
\newcommand{\wtilde}{\widetilde}

\newcommand{\sign}{{\rm sign}}
\renewcommand{\a }{\alpha }
\renewcommand{\b }{\beta }
\newcommand{\s }{\sigma }
\newcommand{\ii }{{\rm i} }
\renewcommand{\d }{\delta }

\newcommand{\g }{\gamma}

\renewcommand{\l }{\lambda }

\newcommand{\barka}{k_0^*}
\newcommand{\vphi}{\varphi }

\renewcommand{\t }{\tau }

\newcommand{\ST}{\mathbb{S}}
\newcommand{\fracchi}{{\mathfrak{I}}}
\newcommand{\C}{\mathbb{C}}
\newcommand{\Z}{\mathbb{Z}}
\newcommand{\T}{\mathbb{T}}

\usepackage{color}

\usepackage{todonotes}

\begin{document}

\title{\bf Time quasi-periodic vortex patches \\
of Euler equation in the plane}
\date{}

\author{M. Berti\footnote{
		% International School for Advanced Studies (
		SISSA, Via Bonomea 265, 34136, Trieste, Italy. 
		\textit{Email:} \texttt{berti@sissa.it}; 
	},  \; 
 Z. Hassainia\footnote{
NYUAD Research Institute, New York University Abu Dhabi, PO Box 129188, Abu Dhabi,   United Arab Emirates.
 	\textit{Email:} \texttt{zh14@nyu.edu}; 
 },\;
 N. Masmoudi\footnote{
NYUAD Research Institute, New York University Abu Dhabi, PO Box 129188, Abu Dhabi,   United Arab Emirates.  Courant Institute of Mathematical Sciences, New York University, 251 Mercer Street, New York, NY 10012, USA.
 \textit{Email:} \texttt{masmoudi@cims.nyu.edu} 
 }}

\maketitle

\noindent
{\bf Asbtract.}  
 We prove the existence  of 
time quasi-periodic vortex patch solutions of the 2$d$-Euler equations in $ \R^2 $, close to 
 uniformly rotating Kirchhoff elliptical vortices, with aspect ratios belonging to  a 
set of asymptotically full Lebesgue measure. The problem is reformulated into a 
quasi-linear Hamiltonian equation for a radial displacement from the ellipse.
 A major  difficulty of the KAM proof  is the presence of a zero normal mode frequency, which is due to the conservation of the angular momentum. The key novelty to overcome this degeneracy is to perform a perturbative 
 symplectic reduction of the angular momentum, introducing it as a symplectic variable
 in the spirit of the Darboux-Carath\'eodory 
  theorem of symplectic rectification, valid in finite dimension. 
 This approach is particularly delicate in a infinite dimensional phase space: 
our symplectic change of variables  is a  nonlinear 
modification of the transport flow generated by the angular momentum itself.
This is the first time such an 
idea is implemented in KAM for  PDEs. 
Other difficulties are 
the lack of rotational symmetry of the equation
and the presence of hyperbolic/elliptic normal modes.
The latter difficulties -as well as the degeneracy of a normal frequency-
 are absent in other vortex patches problems which  have been recently studied  using  the formulation introduced in this paper. 

\smallskip

\noindent
{\small

\smallskip

}

\tableofcontents

\section{Introduction and main result}

The Euler equations, which date back to 1757,  
are the fundamental equations describing the time evolution 
of  inviscid and incompressible fluids. 
They are quasi-linear  partial differential equations of transport-type. 
In view of their physical and mathematical importance, 
the possible development of KAM theory
for the Euler equations is 
one of 
its main motivations and most ambitious 
goals. 

In this paper we prove the  bifurcation 
of an abundance of time quasi-periodic vortex patch solutions 
of the 2$d$-Euler equations close to the uniformly rotating ellipses discovered by Kirchhoff  \cite{Kir} in 1874, 
see Theorems \ref{thm:VP0} and \ref{thm:VP}. 
This is  the first existence result of quasi-periodic vortex patches for  Euler equations. 
Vortex patches of constant vorticity 
play an important role in modeling from hurricanes and cyclones 
to Jupiter's Great Red Spot, and 
their dynamical behavior has been widely investigated 
by means of either analytical, numerical and experimental approaches. 
Let us now introduce our main result  in detail. 

 \smallskip 

We consider the 
incompressible Euler equations in the plane $ \R^2 $,  
written in terms of the scalar vorticity $\bw$ and the stream function $\psi$ as
\begin{equation}
\label{eq:euler}
\partial_t \bw+ ( \nabla^\bot\psi )\cdot \nabla \bw = 0\, ,  \qquad
 \Delta \psi:=  \bw\,  , 
\end{equation}
where $\nabla^\bot\psi := (- \pa_y \psi,\pa_x \psi) $ is the fluid velocity. 

Global existence and uniqueness of weak solutions  of \eqref{eq:euler} for vorticity in  $L^1(\mathbb{R}^2)\cap L^\infty(\mathbb{R}^2)$ is a classical result of Yudovitch \cite{Y}. As a by-product, if the initial datum has a patch form, that is the characteristic function 
$ \bw(0,\cdot)=  {\bf 1}_{D_0} (\cdot) $  of a simply connected bounded domain  $D_0\subset \mathbb{R}^2$,  
then the  
flow $ X^t : \R^2 \mapsto \R^2 $  generated by the velocity field $ \nabla^\bot\psi $ 
is well defined, symplectic 
and carries the vorticity, 
\be\label{VPt}
\bw (t,\cdot) =  {\bf 1}_{D(t)} (\cdot)\, ,  \quad  \quad D(t) :=  X ^t(D_0)\, , \quad  \quad   \textnormal{det} \, D_z X^t (z) =1  \, , \quad \forall t\in \R\, . 
\ee
The dynamics of a vortex patch $ {\bf 1}_{D(t)} (\cdot) $ is reduced 
to the motion of its  boundary $ \pa D(t)  $, which 
evolves according to the so called {\it contour dynamics} equation.  
It has been proved by Chemin
\cite{Ch}  and  Bertozzi-Constantin \cite{BC}
that its regularity is preserved along the evolution.

The simplest example of a vortex patch is the "Rankine vortex" which is the 
circular steady solution with 
$ D(t) = D_0 = \{ |z| \leq 1 \} $ at any time $ t $. 
Another remarkable family of exact vortex patch solutions are the 
ellipses  discovered by Kirchhoff  \cite{Kir} in 1874. Such solutions initially occupy the 
elliptic region $ D_0 $ given by 
 \be\label{def:ell}
 D_\gamma := 
\Big\{ (x,y) \in \R^2 \, : \, \frac{x^2}{\gamma} + \gamma y^2 = 1  \Big\} \, , \quad
\gamma > 1 \, ,  
\ee
and then 
evolve by rotation
at a uniform angular velocity,  namely
$ D(t) = e^{\ii \Omega_\g t}  D_\gamma  $ with 
\be\label{choiceOmega}
\Omega_\g := \frac{\gamma}{(1+\g)^2} \, .
\ee
 For a proof, see for instance \cite[p.304]{BM}  
 (Lemma \ref{lem:equilibrium} provides an alternative proof).
  The constant $ \gamma $ is the aspect ratio of the ellipse, namely the ratio between
 its major and  minor axis.   
Love \cite{Love} established in 1893 the linear stability of the ellipses \eqref{def:ell} of aspect ratio $\gamma<3$ and their linear instability if $\gamma>3$.  
The complete solutions to the linearized Euler equations at the  Kirchhoff ellipses 
 have been   calculated much later by Guo, Hallstrom and Spirn \cite{GHS}, who proved also
 nonlinear instability for $ \gamma > 3 $.
 The nonlinear  stability    for $\gamma<3 $ 
 was proved by Wan \cite{Wan} and  Tang \cite{Tang}  in certain $ L^p $-norms, 
 for the circular case by Wang-Pulvirenti \cite{WP}.

In addition to the ellipses, other uniformly rotating vortex patch solutions 
close to 
the circular Rankine vortices, termed $ V $-states,  
where  numerically computed
by Deem and Zabusky  \cite{DZ}. These solutions look steady in a rotating frame. 
The local bifurcation of $ V $-states was 
analytically proved by Burbea  \cite{Bur},  
recently extended to global branches by Hassaina-Masmoudi-Wheeler \cite{Hass-Mass-Wheel}. In the last years several existence results of 
uniformly rotating patches  have been  proved 
by 
Hmidi, Mateu \cite{HM,HM2,HM3}, with Verdera \cite{HMV},  de la Hoz \cite{HFMV}, Hassaina  
\cite{HR}, \cite{HH2}, and 
Castro-Cordoda-Serrano \cite{CCS}, for $ 2d$-Euler and other  active scalar equations. 
Non-trivial compactly supported stationary solutions have been very recently 
constructed by  Serrano-Park-Shi  \cite{GPS}.

On the other hand,  a particularly rich and complex dynamics 
 is expected for the Euler equations, as also recently pointed out in \cite{TdL}.

The goal of this paper   is to prove the existence of time {\it quasi-periodic} vortex patch solutions for the Euler equations  	\eqref{eq:euler}  close to the rotating Kirchhoff ellipses $ {\bf  1}_{e^{\ii \Omega_\g t}  D_\gamma}$, with an arbitrary  number of frequencies. These solutions are {\it not} steady in any rotating reference frame.

Before stating our main result, we recall that a function 
$  f(t) $ taking values in a Banach space $ E $  is 
quasi-periodic  if 
$ f(t)= F(\omega t) $
where 
$ F:\mathbb{T}^{\nu}\rightarrow E $  is a continuous function 
defined on 
$\mathbb{T}^{\nu} := (\R / 2 \pi \Z)^\nu $,  
$ \nu\geqslant 1$,  
 and the frequency vector $\omega\in\mathbb{R}^{\nu}$ is  non-resonant, namely 
$ \omega\cdot \ell \neq 0 $ for any $  \ell \in\mathbb{Z}^{\nu}\backslash\{ 0\} $.  
Informally stated, our main result is the following;  we refer to Theorem \ref{thm:VP}  for a precise version since some preparation is required.   
\begin{theorem}\label{thm:main}\label{thm:VP0}
{\bf (Time quasi-periodic vortex patches)}  
Consider a compact  interval of aspect ratios 
$ [ \gamma_1, \gamma_2] \subset (1, + \infty) $. 
Then, for any $ \nu \in \N $, there exists a 
 set $ {\cal G} \subset [\g_1, \g_2] $ with asymptotically full  Lebesgue measure 
such that, for any $ \gamma \in {\cal G} $,  there exist 
$ \Omega $ close to $ \Omega_\gamma $ 
and, in the uniformly rotating frame with angular velocity  $ \Omega $, 
 a time quasi-periodic vortex patch solution 
$\bw(t)={\bf 1}_{D(t)} $ 
of the Euler equations \eqref{eq:euler},  with a
diophantine frequency vector $\wtilde \omega  \in \R^{\nu}  $, 
close to the Kirchhoff elliptical patch $ {\bf  1}_{e^{\ii \Omega_\g t}  D_\gamma}$
described 
in \eqref{def:ell}-\eqref{choiceOmega}. 
\end{theorem}

Theorem \ref{thm:VP0}, i.e. \ref{thm:VP}, 
is a KAM  perturbative result which, 
to the best of our knowledge, is the  first  existence result   of quasi-periodic solutions for 
the Euler equation \eqref{eq:euler} in the vortex patches setting. 
This is a difficult  small divisor problem.   
Postponing its detailed description 
in Section \ref{sub:ideas}, we  anticipate some of the difficulties and key ideas of the proof.
After formulating  the contour dynamics equation in terms of 
 a suitable radial deformation of the Kirchhoff ellipses (Lemma \ref{lem:eq-rad-def}), 
obtaining a 
Hamiltonian PDE (Proposition \ref{prop:Ham}),  a major difficulty is the following:
\begin{itemize} \label{diff1}
\item ($i$) the mode $ 2 $ normal frequency of the linearized contour dynamics equation at the Kirchhoff ellipses is {\it degenerate} for {\it any} aspect ratio  $ \gamma $, see \eqref{anbn2}. 
This degeneracy is a consequence of  the  conservation of the 
angular momentum, as we explain in Remark \ref{rem:Ham-dege}.  
\end{itemize}
We remark   that the Hamiltonian formulation provided by 
Proposition \ref{prop:Ham} is an original contribution in the study of the vortex patches problem for Euler. 
Other 
difficulties are that 
($ii$) the nonlinearity  of the contour dynamics equation \eqref{Evera} is quasi-linear and it is 
expressed as an integral operator with singular kernel; ($iii$) the aspect ratio parameter 
$ \gamma $ modifies effectively  the other normal mode 
frequencies by just an exponentially small term;  
($iv$) 
the contour dynamics equation 
is not rotationally invariant 
(usually referred in KAM language as being 
``not momentum preserving"); 
($v$)   hyperbolic 
 and elliptic 
normal mode frequencies coexist.

\smallskip

The first difficulty ($i$) is reminiscent of the well known KAM problem 
arising in Celestial Mechanics  for proving the stability 
of the solar system,  considered by Arnold in \cite{Ar} and completely solved by Chierchia-Pinzari \cite{CP}. 
 In order to eliminate the degeneracy of the second mode a key novelty  
 of this paper is to implement a {\it perturbative} symplectic 
reduction of the angular momentum, introducing it as a symplectic variable, 
 in the spirit of the Darboux-Carath\'eodory 
  ``theorem of symplectic rectification", valid in finite dimension. 
  As far as we know,
this idea was not used previously in KAM theory for PDEs. 
In an infinite dimensional phase space this procedure is particularly delicate:
 our symplectic change of variables $ \Phi $  is  a nonlinear 
modification of the 
transport flow generated by the angular momentum itself 
(Section \ref{sec:SR}) and,  in view of the KAM iteration, 
 we need strong quantitative estimates for $ \Phi $ and its inverse $ \Phi^{-1} $ 
 (Theorem \ref{thm:FM}), as well as 
$ d \Phi $ and $ d \Phi^{-1} $ (Lemma \ref{lem:dphi}). 
In this infinite dimensional context it is not even trivial 
to show the well-posedness and invertibility of $ \Phi $.  
We shall explain in detail the ideas, difficulties  and techniques of the symplectic reduction 
 in Section \ref{sub:ideas}. 
 We expect that  our procedure can be effectively  implemented to other equations.  
Other major differences between our approach and
the symplectic reduction in \cite{CP}, 
in addition to the fact that 
\begin{itemize}
\item[(I)] in this paper the phase space is 
{\it infinite} dimensional, 
\end{itemize}
are the following: 
\begin{itemize}
\item[(II)] the reduction of the angular momentum in \cite{CP}  is {\it not}  
obtained by a perturbative argument,  but provided by an exact formula (given by the special `Deprit coordinates"). 
It is interesting to remark  that it was Arnold's aim \cite{Ar} 
to construct  perturbatively the  symplectic variables 
for the  reduction of the prime integral.
\item[(III)] 
The angular momentum in \cite{CP} vanishes 
quadratically at the equilibrium, 
whereas, in this paper, it is linear, see \eqref{moangJ2} and Remark~\ref{rem:Ham-dege}.
\end{itemize}
For these reasons, the construction of the reducing variables here and in \cite{CP} is 
completely different.  

Concerning other KAM papers for PDEs,  we also 
point out that, thanks to the  explicit
structure of the linearized operator 
for the 2$d$-Euler equation 
 (computed in Section \ref{sec:Lin},  see \eqref{QPLINintr},
 and 
 preserved under conjugation by $ \Phi $, see \eqref{linKOKintr}), 
 the reduction of the linearized operator into a constant coefficient one up to smoothing remainders, which is the longer part  in other KAM works, 
  is done here  in only one step, see  
Section \ref{sec:diffeo} (we 
exploit the Hamiltonian nature of the equation to avoid 
intermediate remainders).
Thus the majority of the effort revolves around the symplectic reduction of the angular momentum to overcome the difficulty ($i$), and its impact on the linearized operator.
\smallskip

In order to formulate the precise statement of Theorem \ref{thm:VP},
we first introduce the contour dynamics  equation 
in convenient symplectic variables near the Kirchhoff ellipses. 
\\[1mm]
{\bf Notation.} Along the paper  we identify $\mathbb{C}$ with $\mathbb{R}^{2}$
and the Euclidean structure of $\mathbb{R}^{2}$ is seen 
through the usual inner  product defined, for any $z=z_{1}+\ii z_{2} $ and $  w=w_{1}+\ii w_{2}\in\mathbb{C} $,  by
$ z\cdot w:= $ $ \langle z,w\rangle_{\mathbb{R}^{2}}= $ $ \mbox{Re}\left(z\overline{w}\right)= $ $ z_{1}w_{1}+z_{2}w_{2} $. 

\subsection{Contour dynamics equation}  \label{sec:counter-dynamic}
Given a vortex patch $ \bw (t)  ={\bf 1}_{D(t)}$ as in \eqref{VPt}, we are interested in the  motion of its boundary $\pa D(t)$.
Since the  tangential component of the velocity does not change the dynamics, but just the parametrization, 
 then any smooth parametrization of the boundary $z(t,\cdot ): \mathbb{T}\mapsto \partial D(t)$ is subject to the {\it contour dynamics} equation
\be\label{NDyn}
 \pa_t z (t, \theta)  \cdot \vec n (t) =-\pa_\theta \psi(t, z(t, \theta)) \, ,
  \ee
where $ \vec n (t)$  is the unit outer normal to the contour $ \pa D(t) $ at the point $z (t, \theta)$ and $ \psi $ is the stream function of the vortex patch.
For a detailed proof of \eqref{NDyn}  see for instance \cite[p.174]{HMV}. Up to a real constant of renormalization, 
$ \vec {n}(t)=-\ii {\partial_{\theta}z(t,\theta)} $,  
and 
 the  equation \eqref{NDyn} can be  written as
\be\label{NDynC-theta}
{\rm Im} \big[ \pa_t z (t, \theta) \overline{\pa_\theta z(t, \theta) } \big] =
\pa_\theta \, \psi (t,z (t, \theta))  \, .
\ee
We  recall that the stream function of the vortex patch is 
\begin{equation}\label{stream-patch}
\psi (t, z) = 
\frac{1}{4\pi} \int_{D(t)} \ln |z-\zeta|^2  dA(\zeta)  
\end{equation}
where $ d A(\zeta) $ denotes the Lebesgue measure in $ \R^2 $.   
\\[1mm]
{\bf Rotating frame.} 
Then we look for solutions of \eqref{NDynC-theta}   of the form 
\be\label{rot-frame}
z(t, \theta) = e^{ \ii \Omega t} w(t,\theta) \, , 
\quad \text{for some} \quad  \Omega \in \R \, . 
\ee
Differentiating \eqref{rot-frame},  using that 
$ \partial_\theta \psi (t,e^{ \ii \Omega t} w(t, \theta ) ) = 
\partial_\theta \psi (t,  w(t, \theta ) ) $, see  \eqref{formpapsi}, 
and 
$ {\rm Im} \big[ \ii w(\theta) \overline{ \pa_\theta w(\theta)} \big]= \tfrac12 \pa_\theta | w(\theta)|^2 $,  the equation \eqref{NDynC-theta} becomes 
\be\label{sub-HS}
{\rm Im} \big[ \pa_t w (t, \theta) \overline{\pa_\theta w(t, \theta) } \big]   =
\pa_\theta \psi_\Omega  (t,w (t, \theta) )\, , \quad
 \psi_\Omega  (t,w (t, \theta) ) :=\psi (t,w (t, \theta))-   \tfrac{\Omega}{2} | w(t, \theta)|^2  \, . 
 \ee
\noindent 
{\bf Equation for the radial deformation.}
In order to study 
vortex patches 
close to the Kirchhoff elliptices
we parametrize the boundary   patch 
as
\be\label{eq:AACD}
\T \ni  \theta\mapsto w(t, \theta)= (1 + 2 \xi (t, \theta))^{\frac12} ( \g^{\frac12} \cos (\theta) +  
\ii \g^{-\frac12}  \sin (\theta) ) \, . 
\ee
The real variable $ \xi (t, \theta ) $
 describes the radial deviation of the contour of the vortex patch from the ellipse.
We emphasize that the  particular form  \eqref{eq:AACD}
is required to obtain a Hamiltonian PDE 
(it is reminiscent of action-angle variables). 
The following lemma shows how the 
equation 
 \eqref{sub-HS} transforms in the  variable $ \xi (t, \theta ) $. 

\begin{lemma}\label{lem:eq-rad-def} {\bf (Equation for the radial deformation)}
 Let $\Omega\in\mathbb{R}$, $\gamma \geq 1$  and  $ D(t) $ 
be a bounded simply connected region with smooth boundary $ \pa D(t) $ parametrized as 
 \be\label{time-para}
 \theta\mapsto z(t, \theta)= e^{\ii \Omega t } w(t, \theta) \quad \text{ with} \quad \begin{aligned}
 w(t, \theta) &= 
 (1+2\xi(t, \theta))^\frac12
 \we(\theta) \, , \\
 \we(\theta) & := \g^{\frac12} \cos (\theta) +  
\ii \g^{-\frac12}  \sin (\theta) \, . 
 \end{aligned}
 \ee
If $ w(t, \theta)$ solves \eqref{sub-HS},
then the radial deformation $\xi (t, \theta)$ solves the equation
\be\label{Evera}
\begin{aligned}
\pa_t \xi (t, \theta)
& = \frac{\Omega}{2} \pa_\theta 
\big(  g_\gamma (\theta) \big(1 + 2 \xi (t, \theta) \big)  \big) \\
& \quad+ \frac{1}{4 \pi} 
\int_{\T} \ln ( \Ke(\xi)(\theta, \theta') ) \pa_{\theta \theta'}^2 
\big[ (1 + 2 \xi (t, \theta))^\frac12
(1 + 2 \xi (t, \theta'))^\frac12 \sin (\theta' - \theta )\big] d \theta' 
\end{aligned}
\ee
where 
\be\label{expR}
\begin{aligned}
\Ke(\xi)(\theta, \theta') 
& := 
 \g \Big( \sqrt{1 + 2 \xi (t, \theta)} \cos \theta - 
 \sqrt{1 + 2 \xi (t, \theta') } \cos \theta' \Big)^2 \\ 
 & 
 \quad \, + \g^{-1}
 \Big( \sqrt{1 + 2 \xi (t, \theta) } \sin \theta - 
 \sqrt{1 + 2 \xi (t, \theta') } \sin \theta' \Big)^2 \, , 
\end{aligned}
\ee
and $ g_\gamma(\theta) $ is the even, $ \pi $-periodic  function
\be\label{fg0}
 g_\gamma(\theta) :=\gamma \cos^2 (\theta ) + \gamma^{-1} \sin^2 (\theta)= \frac{\gamma + \gamma^{-1}}{2} +  \frac{\gamma  - \gamma^{-1}}{2}  \cos (2 \theta) \, .  
\ee
\end{lemma}

The proof of Lemma \ref{lem:eq-rad-def} is given in Appendix \ref{app:CDE}. 
We also anticipate that, as proved in Proposition~\ref{prop:Ham}, 
the evolutionary equation \eqref{Evera} has the Hamiltonian form 
$$
\pa_t \xi = \pa_\theta \nabla H_\Omega (\xi) \, . 
$$
In addition, the  vector field in the right hand side of \eqref{Evera} 
vanishes 
at $ \xi (\theta) = 0 $ if and only if the angular velocity $ \Omega $ is equal to 
$ \Omega_\g $ defined in \eqref{choiceOmega}, confirming that the uniformly 
rotating  Kirchhoff ellipses $ {\bf  1}_{e^{\ii \Omega_\g t}  D_\gamma}$  are 
solutions of 2$d$-Euler.

\begin{lemma}\label{lem:equilibrium} {\bf (Equilibrium solution)}
For any  value of $ \gamma \geq 1 $, the equation \eqref{Evera} with $\Omega=\Omega_\gamma=\tfrac{\gamma}{(1+\g)^2} $  has the equilibrium solution 
$ \xi = 0 $.
\end{lemma}

\noindent
Also this lemma is proved in  Appendix \ref{app:CDE}
 to not interrupt the formulation of the main result.

\subsection{Main result: bifurcation of  quasi-periodic vortex patches}

The goal of this work is to prove that, close to the equilibrium solution 
$ \xi = 0 $, there exist Cantor-like families 
of  
geometrically distinct 
time-quasi periodic solutions of the contour dynamics equation \eqref{Evera}, with an 
arbitrary number of frequencies. 

We need to  anticipate the
following informations about the normal mode frequencies 
of the linearized equation \eqref{Evera} at  
$ \xi = 0 $, 
computed in \eqref{linunp}, 
that will be explained with more detail in the subsequent Section \ref{sub:ideas}. 
 According to the analysis of Section \ref{sec:3}, in accordance with \cite{Love,GHS}, the  linear frequency of oscillations close to Kirchhoff  ellipses is given by \be\label{Omn}
 \Om_n (\gamma) :=  
\Big| \Big(\frac{n \gamma}{(1+\gamma)^2}  -\frac12\Big)^2   -
   \frac{1}{4} \Big( \frac{\gamma-1}{\gamma + 1} \Big)^{2n}\Big|^{\frac12} \,  , \quad \forall \gamma \geq 1  \, .
\ee
Notice that  the mode $n=1$ oscillates with a frequency  equals to the Kirchhoff ellipse angular velocity $ \Omega_\gamma $ and the mode  $n=2$   is {\it degenerate} for all the values of $\gamma$, since $\Omega_2 (\gamma) \equiv 0$. In Remarks \ref{rem:1fre} and \ref{rem:Ham-dege} we provide an explanation of these facts. Moreover,  
for $\gamma\in [1,3)$ 
all  the frequencies of oscillations are elliptic, whereas for 
$ \gamma > 3 $  finitely many hyperbolic directions appear. 
 More precisely, there exists a sequence of aspect ratios $(\underline{\gamma}_{n})_{n\geq 3 }$ such that for  $\bar n\geq 3$ and any $\gamma\in\big( \, \underline{\gamma}_{\, \bar n}, \underline{\gamma}_{ \, \bar n+1} \big)$  the 
 normal modes frequencies \eqref{Omn} are 
 elliptic  for all $n\geq \bar n+1$ and hyperbolic  for $ 3 \leq n\leq \bar n $, see  Lemma 	\ref{lem:eige}. We refer to 
  \eqref{Hanbn} where the decomposition in a sum of harmonic oscillators and finitely many `hyperbolic repulsors" appears  clearly. 
  Fixed $\bar n\geq 2$ we denote by $\IK$ the  interval of aspect ratios
\begin{equation}\label{defgammaset0} 
\IK  := [ \gamma_1, \gamma_2] \subset 
\begin{cases}
\big( 1, 3 \big)   \;\, \quad \qquad \text{if} \ \bar n = 2\, , \\
\big( \, \underline{\gamma}_{\, \bar n}, \underline{\gamma}_{ \, \bar n+1} \big) 
 \quad \text{if} \  \bar n \geq 3\, .   
\end{cases}
\end{equation}
For any $ \gamma \in \IK $
 the linearized contour dynamics equation \eqref{linunp} 
 possesses 
the reversible, i.e. even in $ ( t, \theta ) $, oscillating in time solutions 
\begin{align}\label{q-r-l} 
 q(t, \theta) & =
 {\mathop \sum}_{n \geq \bar n + 1 }  
a_n   \symm_n  \cos (\Omega_n (\gamma) t) \cos ( n\theta) 
+  a_n \symm_n^{-1} \sin (\Omega_n(\gamma) t) \sin (n \theta)   
\end{align}
where  $ a_n \in \R $ and 
 \be\label{symn}
 \symm_n :={\Big| n \Omega_\gamma   -   \frac12
+ \frac12 \Big( \frac{\gamma-1}{\gamma + 1} \Big)^n \Big|}^{\frac14}{\Big| n \Omega_\gamma   -   \frac12
-\frac12 \Big( \frac{\gamma-1}{\gamma + 1} \Big)^n \Big|}^{-\frac14}.
\ee
 The main result of this work proves that  that these structures persist,
 slightly deformed, 
 for most values of the aspect ratio 
 $ \gamma $, at the non-linear level. 
 In order state it we also  introduce the Sobolev spaces
\begin{equation} \label{Sobonorm}
\begin{aligned}
 H^s  := H^s(\T^{|\ST|} \times \T, \R) 
 & = \Big\{ \xi (\vphi, \theta) 
:= {\mathop \sum}_{(\ell,j) \in \Z^{|\ST|+1}} \xi_{\ell,j} \, e^{\ii (\ell \cdot \vphi + j\theta )}  \, : \, \overline{\xi_{\ell,j}} = \xi_{-\ell,-j} \, , \\ 
& \qquad\quad \text{with} \quad
 \| \xi \|_s^2 := {\mathop \sum}_{(\ell,j) \in \Z^{|\ST|+1}} | \xi_{\ell, j}|^2 \langle \ell,j \rangle^{2s} < \infty 
\Big\}
\end{aligned}
\end{equation}
where  $\langle \ell,j \rangle := \max \{ 1, |\ell|, |j| \} $. For
$s >   (|\ST|+1)/2 $ we have $ H^s  \subset 
C^0 (\T^{|\ST|} \times \T, \R) $.

\begin{theorem}\label{KAM finite gap}\label{thm:VP}
{\bf (Time quasi-periodic vortex patches)}
For any $\bar n\geq 2$, consider an interval of aspect ratios 
$\IK  := [ \gamma_1, \gamma_2] $ as in \eqref{defgammaset0}.
Let $ {\mathbb S} $ be any finite subset of distinct  integers 
in 
$ \{\bar n+1,\bar n+2, \ldots \} $ (tangential sites).
Then there exists  
$ \bar s > (|{\mathbb S}| +1) / 2  $, $ \e_0 \in (0, 1) $  such that,
for any $ \vec a = (a_n)_{n \in {\mathbb S}} $ with 
$ a_n > 0 $, for any $ n \in {\mathbb S} $, and $ |\vec a| \leq \e_0 $, 
 the following holds: 
\\[1mm]
1. there exists a subset $ {\cal G} := {\cal G}_{\vec a} \subset [\g_1, \g_2] $
with asymptotically full Lebesgue measure, i.e.
\begin{equation}\label{stima di misura teoremone}
\lim_{\vec a  \to 0} | {\cal G}| = \g_2 - \g_1  \, ; 
\end{equation}
2. for any $ \gamma \in {\cal G} $, there exist 
\begin{enumerate}
\item[(a)]
a constant $\mu(\gamma):= \mu (\gamma, \vec a) \in \R $  satisfying 
$  \mu (\gamma) \to 0 $ as $ \vec a \to 0 $;
\item[(b)] 
a diophantine 
frequency vector $ \wtilde \omega(\gamma) := \wtilde \omega (\g,\vec a) 
:= ( \wtilde \omega_{n} )_{n\in {\mathbb S}} \in \R^{|\ST|}$
with  $ \wtilde \omega (\g) 
\to \vec \omega (\gamma) := (\Omega_n(\gamma))_{n \in {\mathbb S}} $ as $ \vec a \to 0 $ where
$  \Omega_n(\gamma) $ 
are the linear normal mode frequencies of oscillations of the Kirchhoff ellipses, given by \eqref{Omn};
\item[(c)] a quasi-periodic solution 
$ t \mapsto \xi ( \wtilde \omega  t, \theta)$ of 
the contour dynamics equation 
\eqref{Evera} with 
$ \Omega =  \Omega_\gamma +  \mu  (\g) $, $ \Omega_\g = \tfrac{\gamma}{(1+\g)^2} \, , $
 with
Sobolev regularity $  \xi (\vphi, \theta) \in H^{\bar s} (\T^{|\ST|} \times \T, \R) $, 
even in $ (\vphi, \theta) $, 
of the form 
\be\label{xiexpa}
\xi  ( \wtilde \omega t , \theta) 
  =   
 {\mathop \sum}_{n \in \ST}  
 a_n   \symm_n  \cos (\wtilde \omega_n t) \cos ( n \theta) 
+   a_n \symm_n^{-1} \sin (\wtilde \omega_n t) \sin ( n \theta)  + 
{\rm r} (\wtilde \omega t , \theta) 
\ee
where $\symm_n$ is introduced in \eqref{symn} and the remainder
 $ {\rm r}(\vphi, \theta) \in H^{\bar s} $ 
 satisfies 
 $ \| {\rm r} \|_{\bar s}/ |\vec a | \to 0 $ as $ \vec a  \to 0 $. 
 \end{enumerate}
\end{theorem}

Finally, recalling \eqref{time-para}  we 
have proved the existence of a quasi-periodic solution 
$$
z(t, \theta)  
= e^{ \ii \Omega t} (1 + 2 \xi  ( \wtilde \omega t , \theta)  )^{\frac12} \big( \g^{\frac12} \cos (\theta) +  
\ii \g^{-\frac12}  \sin (\theta) \big)  
$$
of the contour dynamics equation \eqref{NDyn}. 
 Theorem \ref{thm:VP}
thus proves the existence 
of  quasi-periodic vortex patch  
solutions of the {\it parameter independent}  
$ 2d $-Euler equation \eqref{eq:euler} on $ \R^2 $, 
bifurcating from  Cantor  families of 
 Kirkhoff ellipses 
with aspect ratio $ \gamma $,
that may be regarded as an  {\it internal} 
parameter. 

\smallskip

Let us make some further comments on the result:
\begin{enumerate}
\item {{\sc Linear Instability}.}
The quasi-periodic solutions \eqref{xiexpa} are linearly {\it unstable}  for $ \gamma > 3 $, 
i.e. the linearized equation \eqref{Evera} at the quasi-periodic solution \eqref{xiexpa}
possesses  real eigenvalues, 
see \eqref{anbn-rin}. 
  This is 
in accordance with the presence of hyperbolic directions of the 
linearized  contour dynamics 
equation 
at $ \xi = 0 $, see 
 \eqref{anbn-ri0}, {{and could 
 account for the  
 chaotic dynamics of the vortex patches}}. 
As $ \gamma \to + \infty $ the number  of hyperbolic directions $ \bar n (\gamma ) - 2 \to +  \infty $ 
 increases to infinity, whereas for 
$ \gamma \in (1,3) $ all the normal mode frequencies are elliptic
(except the degenerate mode $ n = 2 $).  
\item 
  {\sc 
  Set of ``good" aspect-ratio.} The fact that the quasi-periodic 
solutions \eqref{xiexpa} exist  only for  
$ \gamma \in {\cal G} $, is not a technical issue.  In the complementary set we expect 
that these solutions may break down, due to resonances or near resonances among the frequencies.  Actually numerical simulations reveal a wild  dynamical 
behavior of the vortex patches.  
We also note that the values of the aspect ratio 
$ \underline{\gamma}_{\, \bar n} $ in Lemma \ref{lem:eige}, 
that we exclude in \eqref{defgammaset0}, are those $ \gamma $'s 
where a normal frequency $ \Omega_n ( \underline{\gamma}_{\, \bar n}   )$ 
vanishes (and an eigenvalue of the linearized  
equation  \eqref{Evera} at $ \xi = 0 $ changes  nature, passing 
from being 
purely imaginary  
to real).   
\item 
{\sc V-states.}
For $ \gamma $ close to the aspect-ratio 
$ \underline{\gamma}_{\, \bar n} $'s, 
there  exist  uniformly  rotating vortex patches which are deformations of the 
Kirchhoff ellipses,
as proved in \cite{HM2} and \cite{CCS}.
In the setting of this paper, these  solutions can be obtained 
as stationary  solutions $ \xi (\theta )$ 
of the contour dynamics equation \eqref{Evera} with $ \Omega = \Omega_\gamma $ 
and $ \gamma $ close to  $\underline{\gamma}_{\, \bar n} $. 
We also remark that the solutions \eqref{xiexpa} with $ |\ST| = 1 $
are time periodic solutions, but not 
V-states. 
\item 
 {\sc Absence of rotational invariance.} 
For any  $ \gamma > 1 $ the equation \eqref{Evera} is {\it not} 
rotationally $ \theta $-invariant (it does not preserve momentum), unlike for $ \gamma = 1 $, which 
corresponds to study 
vortex patches close to the circular Rankine vortex.  
\item 
{\sc Parameter $ \gamma$}.
We are able to use  the aspect ratio   $ \gamma $  to verify all the 
Melnikov non-resonance conditions along the KAM iteration 
(Theorem \ref{Teorema stima in misura}) thanks to the 
term 
$ \big( \tfrac{\gamma-1}{\gamma + 1} \big)^{2n} $  in  \eqref{Omn}, 
despite the fact that it is 
exponentially small as $ n \to + \infty $, cfr. Section \ref{sec:ND}.
\end{enumerate}

\smallskip

Before discussing in the next section the main difficulties and ideas 
in the proof of Theorem \ref{thm:VP}, we shortly describe some 
 literature regarding KAM theory for PDEs, mainly focusing on 
 quasi-linear PDEs. 
\\[1mm]
{\it Related
KAM Literature.} The first KAM  results for 
 infinite dimensional Hamiltonian PDEs have been proved by Kuksin \cite{Kuksin} and 
Wayne \cite{Wayne} 
for 1$d$ semilinear perturbations of the linear wave (NLW)  and Schr\"odinger  (NLS) equations, 
 using the potential as a parameter to verify Diophantine non-resonance conditions among the linear frequencies. KAM results for parameter independent  NLS  have then been proved in Kuksin-P\"oschel \cite{KP} and, for NLW, by  P\"oschel \cite{Po3},
exploiting, via Birkhoff normal form,   
 the genericity of the initial data to generate large families of small amplitude quasi-periodic solutions. 

For PDEs in higher dimension the first breakthrough results 
are due to Bourgain  \cite{B3,B5} for NLS and NLW,   
and Eliasson-Kuksin \cite{Eliasson-Kuksin} for NLS, 
using the convolution potential  as a parameter.
The presence of external parameters improves significantly the structure of the resonances of the system, which are particularly complex on $ \T^d $, $ d \geq 2 $. Later on   
Procesi-Procesi \cite{PP,PP3} proved KAM results 
for the 
cubic NLS, building on the Birkhoff 
normal form analysis in \cite{PP1},
Eliasson-Gr\'ebert-Kuksin \cite{Eliasson-Greber-Kuksin} for the beam equation
and Berti-Bolle \cite{BB20,BBo1,BB10} for NLS and NLW equations with multiplicative potential (\cite{BBo1,BB10} deal with  quasi-periodically  forced nonlinearities with 
 frequencies used as external parameters).

Now we quote  KAM results  with unbounded nonlinearities, which hold 
in $ 1 d $. In the semilinear case, 
 results 
are proved in Kuksin and Kappeler-P\"oschel
\cite{K2,Kuksin2,Kappeler-Poschel} for 
KdV,  in 
Berti-Biasco-Procesi \cite{BBP2} for derivative NLW,  in Liu-Yuan \cite{Liu-Yuan}  for derivative NLS. 
In  the quasi-linear case, small amplitude 
quasi-periodic solutions were constructed  
in Baldi-Berti-Montalto \cite{BBM-Airy}  for
quasi-periodically forced perturbations of Airy 
equations,  in  \cite{BBM-auto,BKM1}  for  autonomous 
perturbed KdV equations --using a weak Birkhoff normal form analysis to 
modulate the frequencies via the initial data--, in Giuliani \cite{Giu} for  gKdV,  
and  Feola-Procesi \cite{Feola-Procesi} for 
quasi-periodically forced NLS. We also mention 
 Berti-Kappeler-Montalto \cite{BKM1} that proves the persistence of 
Cantor families of finite gap solutions of KdV of arbitrary size,  under quasi-linear 
Hamiltonian perturbations. 
The first bifurcation results of time quasi-periodic  
standing solutions of the 
water waves equations were proved in 
Berti-Montalto \cite{BertiMontalto}, in the gravity-capillary case (using the 
surface tension as a parameter), 
and in Baldi-Berti-Haus-Montalto
\cite{BBHM} for pure gravity waves (using the depth as a parameter). 
The proof that the linear normal mode frequencies 
satisfy Diophantine non-resonance conditions
relies on a generalization of the ``degenerate KAM approach" in 
Bambusi-Berti-Magistrelli \cite{BaBM}. 
Previous  results for periodic solutions were obtained by Iooss, Plotnikov, Toland
\cite{PlTo,Ioo-Plo-Tol,IP-Mem-2009} and Alazard-Baldi  \cite{Alaz-Bal}.
Traveling quasi-periodic water waves with constant vorticity --which may be regarded as
a linear superposition of simple Stokes waves with  
Diophantine speeds-- have been 
recently obtained in 
\cite{BFM1} for most surface tension coefficients, 
and in  \cite{BFM} for pure gravity waves, using
the vorticity as a physical parameter. For pure gravity 
irrotational water waves  in infinite depth,  Feola-Giuliani 
\cite{FG}  were able to use the integrable  Zakharov-Dyachenko 
 normal form 
to prove existence of 
traveling quasi-periodic solutions for
`generic"  initial conditions.  
We also mention Feola-Giuliani-Procesi \cite{FGP} for the existence of 
small amplitude quasi-periodic 
solutions of quasi-linear perturbations of 
the Degasperis-Procesi equation via a weak Birkhoff normal form analysis.
We wish also to quote  
Baldi-Montalto  \cite{BM20}, where an external  quasi-periodic forcing term 
is added to generate quasi-periodic solutions of 3$d$ Euler. 
 
\smallskip

  We finally mention that the formulation introduced in this paper was very recently  used to prove  the existence of quasi-periodic patches  close to the circular Rankine vortices  
 for  $ (SQG)_\alpha$ in \cite{HHM21}, 
 for the  
$ (QGSW)_\lambda $  equations in \cite{HR21} and 
   for  2$d$-Euler set in the {\it disk} in \cite{HaR} 
(not in $ \R^2 $).   
We point out that in these cases the difficulty ($i$) of the degeneracy of the 
normal mode frequency for $ n = 2 $ is {\it not} present.
In addition,  these works 
concern perturbations of  circular Rankine vortices, 
and thus the corresponding  equations  are rotationally invariant  
(i.e. momentum preserving) 
 and 
the normal mode frequencies are  
{\it all} elliptic (in this paper
we deal with elliptic {\it and} hyperbolic directions in Section \ref{sec:redu}).  We also emphasize that in \cite{HHM21,HR21} the existence of quasi-periodic patches is proved
varying an exterior parameter ($ \alpha $ or $ \lambda$) of the equations.
Instead, in this work, we are able use the natural {\it inner} geometrical 
parameter $\gamma$ of the Kirchhoff ellipses to prove the existence 
of solutions of the $ 2d$-Euler equation \eqref{eq:euler} on $ \R^2 $.

\subsection{Ideas of proof and  plan of the paper}
\label{sub:ideas}

As we shall prove  in Section \ref{sec:sym},  the contour dynamics equation 
\eqref{Evera}  enjoys the Hamiltonian structure 
\be\label{PDEintr}
\pa_t \xi = \pa_\theta \nabla H_\Omega (\xi) \, , \quad H_\Omega (\xi) := - \tfrac12 E (\xi) +  \tfrac{\Omega}{2} J  (\xi) \, , 
\ee
where $ E  $ is   the pseudo-energy, $ J $ is the angular momentum, $ \nabla H_\Omega $ denotes the $ L^2 $-gradient of the  Hamiltonian and the rotating frequency $ \Omega $ 
is a free parameter. The angular momentum 
$ J $ is a prime integral of \eqref{PDEintr}.

A significant part of the work is 
 to show how the equation \eqref{Evera}, i.e. \eqref{PDEintr}, 
falls into a 
 framework convenient for KAM techniques.
We now  explain it, focusing 
on the difficulties 
mentioned after the statement of  Theorem \ref{thm:main}, especially the first ($i$).  

As already explained, by Lemma  \ref{lem:equilibrium}  it results $ \nabla H_{\Omega_\gamma} (0)  = 0 $ for any $ \gamma $, i.e. $ \xi = 0 $ is an {\it equilibrium}.

In  Section \ref{sec:3} we compute the 
 {\it linearized  equation} of \eqref{Evera} at the equilibrium state 
$ \xi = 0 $, see \eqref{linunp},  and we diagonalize it. 
Actually \eqref{linunp} follows as a particular case of 
the results of Section \ref{sec:Lin}
where -in view of a Nash-Moser iterative scheme- 
we compute the linearized 
equation of \eqref{Evera} at any 
$ \xi (\theta ) $ in a neighborhood of  zero.
Expanding  in Fourier series
\be\label{fouexpa}
 \xi (\theta)  
 =  
  {\mathop \sum}_{n \in \N} \, \a_n \tc_n (\theta) 
  +  \b_n \ts_n ( \theta)  \quad \text{where} \quad  
  \tc_n  = \tfrac{1}{\sqrt{\pi}} 
  \cos ( n\theta) \, , \quad \ts_n  = \tfrac{1}{\sqrt{\pi}} \sin (n \theta) \, , 
\ee
(for $ n \neq 2 $ we  choose in \eqref{defcnsn} 
a different normalization constant)
we obtain  that \eqref{linunp} 
is equivalent,  in the variables $(\alpha_n,\beta_n)_{n\geq 1}$,  to the infinitely many decoupled linear systems 
\be\label{anbn}
\begin{pmatrix}
\dot \a_n \\
\dot \b_n
\end{pmatrix} = \begin{pmatrix}
0 & -\mu_n^{+}  \\
 \mu_n^{-} & 0 
\end{pmatrix}\begin{pmatrix}
\a_n \\
\b_n
\end{pmatrix}  \, , \quad \textnormal{where}\qquad
\begin{aligned}
 \mu_n^+ &:= \mu_n^+ (\gamma) := n \Omega_\gamma   -    \tfrac{1}{2}
+ \tfrac{1}{2} \Big( \frac{\gamma-1}{\gamma + 1} \Big)^{n} \, , 
\\   \mu_n^-  &:= \mu_n^- (\gamma) :=    n \Omega_\gamma      
-      \tfrac{1}{2}     -   \tfrac{ 1}{2}\Big( \frac{\gamma-1}{\gamma + 1} \Big)^{n}  \, . 
\end{aligned}
\ee
It turns out that, for any 
aspect ratio $ \gamma \in \IK $ defined in  \eqref{defgammaset0} 
each  system \eqref{anbn} 
is a harmonic oscillator 
 for $ n = 1 $ and  $n\geq \bar n+1$ 
 and a hyperbolic  repulsor  for $ 3 \leq n\leq \bar n $, cfr.  \eqref{anbn-ri0}. 
 For $ n = 2 $  it results 
 $ \mu_2^+ (\gamma) = 0$  and the dynamics of $(\alpha_2,\beta_2)$ is described by a Jordan block with eigenvalue $ 0 $,  
 namely the {\it degenerate} linear system
\be\label{anbn2}
\begin{pmatrix}
\dot \a_2 \\
\dot \b_2
\end{pmatrix} = \begin{pmatrix}
0 & 0  \\
 \mu_2^{-} & 0 
\end{pmatrix}\begin{pmatrix}
\a_2 \\
\b_2
\end{pmatrix}  \, . 
\ee
This is  the serious difficulty  ($i$)  for the KAM proof  
mentioned after the statement of  Theorem \ref{thm:main},  
because the higher order 
nonlinear terms in \eqref{Evera} actually depend 
on $ (\alpha_2, \beta_2 )$.
The corresponding Hamiltonian looks, in the symmetrized 
variables
$  \breve \a_n = \symm_n^{-1} \a_n $, $   \breve  \b_n = \symm_n \b_n $, $ \forall n \neq 2 $, 
cfr. \eqref{def:Lan}, 
\be\label{Hanbn}
{ H}_L =  \frac{\Omega_1}{2} (\breve \a_1^2 +  \breve \b_1^2) - 
\frac{\Omega_2}{2}  \a_2^2  
+ 
\sum_{ 3 \leq n \leq \bar n} \frac{\Omega_n}{2n} ( \breve \a_n^2 - \breve\b_n^2)
- 
\sum_{n \geq \bar n + 1} \frac{\Omega_n}{2n} ( \breve \a_n^2 + \breve \b_n^2) \, .
\ee 
The degeneracy of the mode $ 2 $ --note that $ H_L $ does not depend on $ \beta_2 $--
is actually due to the conservation of the  angular momentum
that has a linear $ \a_2 $-component, see Remark  \ref{rem:Ham-dege}.  
In order to 
overcome this problem, we  have to fully exploit the
conservation law of   $ J $ along the KAM procedure. We proceed as follows.
\\[1mm]
{\bf The symplectic rectification of the angular momentum.} We implement the 
normal form approach of Darboux-Caratheodory, also termed, for finite dimensional systems,  the theorem of "symplectic rectification", stated e.g. in  Theorem 10.20 in  Fasano-Marmi \cite{Marmi-Fasano}.

In order to explain our approach we identify below 
 $ \xi (\theta ) $ with the set of coordinates $ (\alpha, \beta) := 
(\alpha_n, \beta_n)_{n \in \N} $ defined in \eqref{fouexpa}. 
Moreover in the sequel  $ Q(\alpha, \beta)  $
denotes a function which is quadratic in $ (\alpha, \beta) $, i.e. $ \xi $, 
without specifying the norms.   

The goal is  to  construct a symplectic diffeomorphism  
$ \wtilde 	\xi = \Phi (\xi) $, acting in Sobolev spaces $ H^s (\T) $,  
which introduces the prime integral (which is   a multiple of the angular momentum $J$, cfr. \eqref{Jsvi})
$$ 
{\cal J } = \frac12 (\xi, \tc_2) + \aleph \int_\T \xi^2 (\theta) \, g_\gamma (\theta) \, d \theta \, ,
\quad \text{in coordinates} \ {\cal J} = \a_2 + Q(\alpha, \beta)  \, ,  
$$ 
as a symplectic coordinate, see Theorem \ref{thm:FM}. Here $ \aleph \in \R $ 
is a normalization constant. 
Note that  
the associated Hamiltonian vector field  (cfr. 
Lemma \ref{lem:AM})
\be\label{HAMJ}
X_{\cal J} (\xi) =  - \ts_2  (\theta) + 
2 \aleph \big( \pa_\teta \circ g_\gamma (\theta)\big) (\xi)
\ee
is an affine transport  operator, which 
does not  vanish at $ \xi = 0 $, 
\be\label{XJVF}
 X_{\cal J} (0) = - \ts_2  (\theta) \, , \quad \text{in coordinates} \   X_{\cal J} (0)  \equiv 
 (\underbrace{0}_{\alpha_2},\underbrace{-1}_{\beta_2}, \underbrace{0,0,\ldots}_{\a_n,\b_n, n \neq 2}) \, . 
\ee
In other words 
we look for real functions $   \bar t (\xi) $, 
$  (\wtilde \a_n (\xi), \widetilde \beta_n (\xi))_{n \neq 2} $, such that the map 
$$
 \Phi(\xi)  = {\cal J}(\xi) \tc_2 (\theta) 
  +    \bar t(\xi ) \ts_2 ( \theta) + 
{\mathop \sum}_{n \neq 2} \, \wtilde \a_n (\xi) \tc_n (\theta) 
  +   \widetilde \beta_n (\xi) \ts_n ( \theta) 
$$ 
is a symplectic diffeomorphism of $ H^s $.
This procedure is called a symplectic 
  rectification of $ {\cal J } $ because,  in the new symplectic variables 
\be\label{coordinew}
\wtilde \a_2 := {\cal J } (\xi) \, , \quad \widetilde \beta_2 :=  \bar t 
(\xi) \, , \quad 
\wtilde \a_n := \wtilde \a_n (\xi) \, , \quad \widetilde \beta_n := \widetilde \beta_n (\xi), \ \forall n \neq 2 \, , 
\ee
  the angular momentum $ {\cal J } $ is  just $ \widetilde\alpha_2 $ 
  and the vector field $ X_{\cal J } (\xi) $ 
  is transformed, in a {\it full} neighborhood of $ \xi = 0  $,  into the {\it straightened} field 
  $ (0, -1, 0, 0, \ldots) $, 
  wheres \eqref{XJVF} holds just at $ \xi = 0 $.

The new coordinates \eqref{coordinew} are constructed as follows.   
  The symplectic variable $ \wtilde \beta_2  = \bar t (\xi)  $ 
  conjugated to $ \wtilde \a_2 = {\cal J } (\xi)  $ is  the 
  {\it time of impact} of the flow $ \Phi^t_{\cal J} (\xi) $ on the manifold $ \{ \beta_2 = 0 \} $, which is  transverse to
  the flow, namely $\bar t (\xi)  $ is defined as the unique local solution of 
  $$
\big( \ts_2, \Phi^{\bar t (\xi)}_{\cal J} (\xi) \big)_{L^2} = 0 \, . 
  $$
  In Section \ref{flowJ} we prove that the flow generated by  
  $ X_{\cal J} $ 
  is well posed in any $ H^s (\T) $.
In particular, in view of \eqref{XJVF}, the Hamiltonian system generated by $  \cal J $ has  the form 
  $ \dot \beta_2 = - 1 + \ldots $ up to linear terms (cfr. \eqref{xj1coo}),
and thus 
$
    \widetilde \beta_2 := \bar t(\xi )  = \beta_{2}  + Q(\alpha, \beta) $.  
     The other symplectic variables  
$ (\wtilde \a_n (\xi), \wtilde \b_n (\xi))_{n \neq 2}$ 
  are defined by 
  $$ 
 {\mathop \sum}_{n \neq 2} \, \wtilde \a_n (\xi) \tc_n (\theta)  
  +   \widetilde \beta_n (\xi) \ts_n ( \theta)    =   \Pi_2^\bot \Phi^{\bar t (\xi)}_{\cal J} (\xi) 
  $$  
  where $ \Pi_2^\bot $ is the projection on the subspace 
  supplementary to the mode $ 2 $.
 This perturbative construction 
  is local around $ 0 $, which is sufficient for our purposes.  
  The proof of the invertibility of $ \Phi $ is delicate and does not follow by the 
  implicit function theorem, since $ d \Phi (\xi) - {\rm Id} $ is small 
  with loss of $ 1 $ derivative (the vector field \eqref{HAMJ} in unbounded). 
  Indeed it relies on the transport structure 
  of the equation. 
   In view of the KAM-Nash-Moser iteration to be performed 
  after this change of variables, 
  very    strong quantitative  estimates for $ \Phi $ and $ \Phi^{-1} $, as well as  
  for $ d \Phi $ and $ d \Phi^{-1} $,
  are actually required, see Theorem \ref{thm:FM} and Lemma \ref{phi-1}. The map  $ \Phi $ is also reversibility preserving. 
  
  \begin{remark}
   Interestingly the finite dimensional 
   symplectic rectification theorem of Darboux-Caratheodory has a quantum analogue which is the Duistermaat-H\"ormander \cite{DH}
  theorem of (microlocal) rectification  of a non-degenerate 
  pseudo-differential   operator via 
  Fourier integral operators. The results of 
  Section \ref{flowJ} can be seen as a result of this kind.
\end{remark}
\begin{remark}
The dynamics on mode $ 1$ is related to other prime integrals of \eqref{Evera}, 
namely the center in \eqref{forCZJ}, originating by the translation invariance of Euler equations in the plane, see Remarks \ref{HWS} and \ref{rem:1fre}.
\end{remark}

\noindent
{\bf The symplectic reduction of the angular momentum.} 
Since $ {\cal J} $ is a prime integral of $ H_\Omega $, the Hamiltonian 
$ K := H_\Omega \circ \Phi^{-1} $ in the new symplectic  variables $ \widetilde \xi $
is,  being $ \widetilde\alpha_2 = {\cal J} $ a prime integral of $ K $,  
{\it independent} of $  \widetilde \beta_2  $, namely it 
 has  the form  $ K({\cal J}, \wtilde u ) $ where we  
 denote 
$  \wtilde u   \equiv ( \wtilde \alpha_n, \wtilde \beta_n)_{n \neq 2} $, cfr. \eqref{coordinew}. 
 The Hamiltonian system has been reduced to the Hamiltonian PDE 
 at  fixed $ \underline{\cal J} $,  in the variable 
$ \widetilde u $, 
 \be\label{equtilde}
  \pa_t \wtilde  u = \pa_\theta 
  \nabla_{\wtilde u}  \cK (\, \underline{\cal J} , \wtilde u \, ) 
  \quad \text{where} \quad \nabla_{\wtilde u}  \cK ( \underline{\cal J}, \wtilde u \, ) = \Pi_2^\bot  
  (\nabla_{\wtilde \xi} K) ( \underline{\cal J} \tc_2  + \wtilde u ) \, , \quad \wtilde u = \Pi_2^\bot \wtilde \xi \, , 
  \ee
 where the degenerate mode $ 2 $ has been removed,
see Corollary \ref{cor:new} and subsequent discussion. 
This natural construction would 
enable to  conclude, 
in finite dimension, 
   the KAM proof, reducing it to the known non-degenerate case.  
   
For infinite dimensional systems however  there is a further serious difficulty.
 In the $ \xi $ variable the quasi-periodic linearized Hamiltonian operator  
$ \omega \cdot \pa_\varphi - d X_{H_\Omega} (\xi ) $ at an approximate solution  $ \xi (\varphi, \theta) $
has the 
PDE structure, 
derived by the computations of Section \ref{sec:Lin},  
\be\label{QPLINintr}
\omega \cdot \pa_\varphi - d X_{H_\Omega} (\xi )  = 
\om \cdot \pa_\vphi - \underbrace{\pa_\theta\circ  V}_{transport}
  + \underbrace{\pa_\theta  \circ \W_0}_{unperturbed \ Hamiltonian \ operator} + 
 \underbrace{\pa_\theta 	\circ { R}_\e}_{smoothing} 
\ee
where 
 $ V  $ is a  real valued function close to $ - \Omega_\gamma $,  
the unperturbed Hamiltonian operator $\pa_\theta \circ \W_0 $ has the form 
$  \pa_\theta \circ \W_0 = \tfrac12 {\cal H} +  {\cal Q}_\infty $ where 
${\cal H}  $ is the Hilbert transform and 
$ {\cal Q}_\infty \in  {\rm OPS}^{-\infty } $ is 
 defined in \eqref{defQinfty}, 
 and 
$ { R}_\e $ is a small smoothing remainder in $ {\rm OPS}^{-\infty }$, see Lemma \ref{lem11.4}.
A linear operator of the form \eqref{QPLINintr} has the nice dynamical feature that 
may be reduced to constant coefficients up to smoothing remainders in only one step (cfr.   
Section \ref{sec:diffeo}).

But, for implementing a 
Nash-Moser iteration
for \eqref{equtilde} in the variable $ \widetilde u = \Pi_2^\bot  \widetilde \xi  $,
where $ \widetilde \xi $ is the variable obtained after the symplectic rectification,
we have to
\begin{itemize} 
\item   
invert the linearized operators 
$ \Pi_2^\bot \big( \omega \cdot \pa_\varphi - d X_{K} (\wtilde \xi ) ) \Pi_2^\bot $ in the new 
 coordinate $ \widetilde u = \Pi_2^\bot  \widetilde \xi  $. 
\end{itemize}

This is not an easy task since the rectification 
map $ \Phi $ is somehow implicit and might significantly affect
the PDE structure \eqref{QPLINintr} of 
$ \omega \cdot \pa_\varphi - d X_{H_\Omega} (\xi ) $. Nevertheless 
we observe that $ \Pi_2^\bot d \Phi (\xi  )$ is the restriction of the 
linear transport $ \Phi^{\bar t (\xi)}_{{\cal J}_2} $ flow plus a regularizing operator,
see Lemma \ref{phi-1},  and similarly  
$ \big[ d \Phi (\xi) \big]^{-1}  \Pi_{2}^\bot 
=   \Phi^{- \bar t (\xi)}_{{\cal J}_2} \Pi_{2}^\bot+  \ldots $ up to a regularizing operator.
Then in Lemma \ref{lem:AI} we implement an approximate-inverse 
argument \`a la Zehnder, noting 
that  
the linearized operator 
$ \omega \cdot \pa_\varphi - d X_{K} (\wtilde \xi )$ 
 is obtained by conjugating the linearized operator
 $ \omega \cdot \pa_\varphi - d X_{H_\Omega} (\xi ) $ 
in the $ \xi $ variables via the linearized map $ d \Phi (\xi ) $,  
 plus a term which vanishes at a solution, 
see \eqref{terzoquasi}.
In this way, 
we are able to deduce the conjugation formula \eqref{bellaconiug} 
(where    a new small transport operator appears). Finally,  
using an Egorov type argument and the representation \eqref{reprePhi2}
of the flow $ \Phi^t_{{\cal J}_2}$, we conclude
that  the linearized operator
$ \Pi_2^\bot ( \omega \cdot \pa_\varphi - d X_{K} (\wtilde \xi )){|\Pi_2^\bot} $
 has still a 
 PDE structure
similar to $ \omega \cdot \pa_\varphi - d X_{H_\Omega} (\xi ) $ in \eqref{QPLINintr}, 
 namely 
 \be\label{linKOKintr}
\Pi_2^\bot ( \omega \cdot \pa_\varphi - d X_{K} (\wtilde \xi )){|\Pi_2^\bot} = 
  \Pi_{2}^\bot ( 
\omega \cdot \partial_\vphi  -  
\pa_\theta \, \circ \,   {\cal V}+ \pa_\theta \, \circ \,   W_0  +{\rm R}_\e
 ){|\Pi_2^\bot}
+{\cal R}_{Z} \, , 
\ee
where   ${\cal V} $ is a real valued function close to $  -  \Omega_\gamma $,
 ${\rm R_\e}$ 
is a small smoothing operator in $ {\rm OPS}^{-\infty} $ 
 and the remainder $ {\cal R}_{Z} $  vanishes at a solution,  
  see Proposition  \ref{lem:Lorto}. 
The operator   \eqref{linKOKintr}  is 
then reduced in Section \ref{sec:diffeo}  
to a constant coefficient one up to smoothing remainders,    
and  finally 
completely diagonalized via a KAM algorithm  in
Section \ref{sec:redu}.

  We think that this  approach is a 
 non trivial step of the proof, which could be applicable 
 in other PDE contexts. 

\smallskip

The paper is organized as follows. In Section  \ref{sec:sym} we present some conserved quantities of the contour dynamics equation and its Hamiltonian formulation. 
In Section \ref{sec:Lin} we compute the linearized contour dynamics equation  at any state $\xi (\theta) $, not necessarily the equilibrium. 
At  $\xi=0$, we diagonalize the linearized equations
and  compute  their linear frequencies.   
 This is done  in Section \ref{sec:3}, where we also establish the aspect ratio 
 $ \gamma $ for which 
 hyperbolic  frequencies 
 appear.   Section \ref{sec:ND} is devoted to prove transversality 
 properties  of the  linear frequencies with respect to
 $ \gamma $. 
  In Section \ref{sec:SR} we implement the  symplectic reduction of the angular momentum  
  to eliminate the degenerate 
  mode $2$. 
  In Section  \ref{sec:NMT}  we start the proof of the main Theorem \ref{KAM finite gap}. 
 First, we reformulate it 
 in a more dynamical system language, 
  that is Theorem \ref{main theorem}. 
In Section \ref{sec:measure} we prove the measure estimates
for the set of $ \gamma $'s where all the Melnikov non-resonance conditions are verified.   
The key result for proving Theorem \ref{main theorem}
is the existence 
of an almost approximate inverse
as stated in Theorem~\ref{thm:stima inverso approssimato}, which is proved 
 in Sections~\ref{reduction}-\ref{sec:redu}.   
As preparation, in Section~\ref{sec:regularity} we prove tame estimates for the symplectic rectification map constructed in Section~\ref{sec:SR} and the transformed 
Hamiltonian vector field. 
In Section~\ref{reduction}  
we reduce the linearized operator  to constant coefficients up to smoothing remainders, 
and in Section~\ref{sec:redu} we diagonalize it 
by a KAM scheme  with 
hyperbolic and elliptic directions. 
Finally Section~\ref{sec:NM} is devoted to the Nash-Moser result. 
Appendix \ref{app:CDE} contains the derivation of the contour dynamics equation
and Appendix \ref{sec:B} technical lemmata. 

We underline that the parts which are close to the KAM strategy 
implemented  
for 
PDEs  in previous works \cite{BertiMontalto,BBHM,BBM-auto,BFM1,BFM,BB20} are only Sections \ref{sezione approximate inverse} and \ref{sec:NM}, which, for this reason, are short. A significantly long portion of this paper 
is due to the symplectic change of variables constructed to reduce the angular momentum, cfr. Sections \ref{sec:SR} and 
\ref{Sec:KAMnew}-\ref{sec:LH}.
\\[1mm]
{\bf Notations.} We denote by $\mathbb{N}$ the set of natural numbers, that is 
$ \mathbb{N}:=\{1,2,\cdots\} $ and $ \N_0 := \{0\} \cup \N $.  
For $a\lesssim_s b$ means that $a\leq C(s) b$ for some positive constant $C(s)$.
For $ \upsilon, \tau > 0 $ and an integer $|\ST| \in \N $, we denote the set of Diophantine vectors as 
\be\label{DC tau0 gamma0}
\mathtt {DC} (\upsilon, \tau) := \Big\{ \omega \in \R^{|\ST|} \, : \, 
|\omega \cdot \ell | \geq \upsilon \langle \ell \rangle^{-\tau} \, , \ \forall \ell \in \Z^{|\ST|} \setminus \{0\} \Big\} \, .
\ee
{\bf Acknowledgments.}
M. Berti thanks M. Procesi and G. Pinzari for  useful discussions and L. Chierchia and J-M. Delort for suggestions which improved the presentation. M. Berti was partially supported 
by PRIN 2020 (2020XB3EFL001) 
``Hamiltonian and dispersive PDEs". 
The work of Z. Hassainia and  N. Masmoudi  is supported   by Tamkeen under the NYU Abu Dhabi Research Institute grant of the center SITE.    The work of N. Masmoudi is supported by  NSF grant DMS-1716466. 

\section{Algebraic properties and symmetries}\label{sec:sym}

The main concern of this section is to describe some conserved quantities of
the contour dynamic equation \eqref{Evera} 
as well as its  Hamiltonian and reversible 
structure. 

\subsection{Prime integrals}

For vortex patches data, the four fundamental conserved quantities are   the circulation $C$,  the center of mass $Z$, the angular momentum $J$ and  the pseudo-energy $E$,
 \begin{equation}\label{pr-vp}
 C:=\int_{D(t)}dA(z) \, , \quad
 Z:=\int_{D(t)}zdA(z) \, , \quad 
J:=\int_{D(t)}|z|^2dA(z) \, , \quad  E := \int_{D(t)} \psi (t, z) d A(z) \, .
 \end{equation}
In order to recover the Hamiltonian structure of  the contour dynamics  equation  \eqref{Evera},  we  need to describe  the prime integrals \eqref{pr-vp}  in the patch setting using the radial  deformation  $\xi$ defined by \eqref{time-para}. 
Using the complex form of Green's formula \eqref{green}, the above prime integrals 
can be written as 
 \begin{equation}\label{CZJ}
 C=\frac{1}{2 \ii}\int_{\partial D(t)}\overline{z} dz \, , \quad 
 Z=\frac{1}{2 \ii}\int_{\partial D(t)}|z|^2dz \, , \quad 
J=\frac{1}{4 \ii}\int_{\partial D(t)}|z|^2\overline{z} dz \, .
 \end{equation}
 From \eqref{time-para} and \eqref{CZJ}  we obtain the following formulae.

{\begin{lemma}	\label{lem:prin}
Let $ D(t) $ 
be a bounded simply connected region with smooth boundary $ \pa D(t) $ parametrized as in \eqref{time-para}.
 Assume that $\bw(t)={\bf 1}_{D(t)} $. Then the circulation $C$, the center of mass $Z$,  the angular momentum $ J $  
  in \eqref{CZJ}, i.e. \eqref{pr-vp}, 
  are expressed as
  \begin{align}\label{forCZJ} 
 C=\pi+\int_{\T}\xi(\theta) d\theta\, , \quad
 Z= e^{ \ii \Omega t}\frac{1}{3}\int_{\T}\big(1+2\xi(\theta)\big)^\frac{3}{2}\we(\theta)d\theta\, , \quad
J=\frac{1}{4}\int_{\T}\big(1+2\xi(\theta)\big)^2g_\gamma(\theta)d\theta \, , 
 \end{align}
 where   $ \we(\theta)$ and $g_\gamma(\theta)$ are defined in \eqref{time-para} and \eqref{fg0} respectively.
Note that  $ C $, $ J $ 
 are autonomous prime integrals
 of the 
 equation \eqref{Evera}, 
 whereas $Z$  is a time-dependent  one. 
 \end{lemma}

\begin{pf}
By \eqref{time-para} and  \eqref{CZJ} we have 
\begin{align*}
C 
& = 
\frac{1}{2\ii}\int_{\T} \big(1+2\xi(\theta)\big)^\frac12\overline{ \we(\theta)} \partial_\theta\big[\big(1+2\xi(\theta)\big)^\frac12 \we(\theta)\big] d\theta \, , 
\\
Z &= e^{\ii \Omega t}\frac{1}{2 \ii }\int_{\T}\big(1+2\xi(\theta)\big)| \we(\theta)|^2 
 \partial_\theta\big[\big(1+2\xi(\theta)\big)^\frac12 \we(\theta)\big] d\theta \, ,
\\
J &= \frac{1}{4 \ii }\int_{\T}\big(1+2\xi(\theta)\big)^{\frac32}| \we(\theta)|^2 \overline{ \we(\theta)} 
\partial_\theta\big[\big(1+2\xi(\theta)\big)^\frac12\we(\theta)\big] d\theta \, , 
\end{align*}
and, since
$ \partial_\theta\big[\big(1+2\xi(\theta)\big)^\frac12\we(\theta)\big]=\partial_\theta\xi(\theta)\big(1+2\xi(\theta)\big)^{-\frac12}\we(\theta)+\big(1+2\xi(\theta)\big)^\frac12\partial_\theta \we(\theta) $, 
we get  
\begin{align*}
C &=\frac{1}{4\ii}\int_{\T}\partial_\theta\big(1+2\xi(\theta)\big)| \we(\theta)|^2  d\theta
+\frac{1}{2\ii}\int_{\T}\big(1+2\xi(\theta)\big)\overline{ \we(\theta)}\partial_\theta \we(\theta)  d\theta \, , 
\\
Z&= e^{\ii \Omega t}\frac{1}{6\ii}\int_{\T}
 \partial_\theta \big(1+2\xi(\theta)\big)^\frac32  | \we(\theta)|^2 \we(\theta)d\theta
+ e^{\ii \Omega t} \frac{1}{2\ii}\int_{\T}\big(1+2\xi(\theta)\big)^{\frac32}| \we(\theta)|^2 \partial_\theta \we(\theta) d\theta,
\\
J&=\frac{1}{16\ii}\int_{\T}
 \partial_\theta\big(1+2\xi(\theta)\big)^2 | \we(\theta)|^4d\theta
 +\frac{1}{4\ii}\int_{\T}\big(1+2\xi(\theta)\big)^{2}| \we(\theta)|^2\overline{ \we(\theta)}\partial_\theta \we(\theta) d\theta \, . 
\end{align*}
Integrating  by parts the first term of these expressions, 
by \eqref{id:we} and  $ |\we(\theta)|^2 = g_\gamma (\theta) $,
 we deduce  \eqref{forCZJ}.  
 \end{pf}
 
 \smallskip
 
 The next lemma provides  the formula for  the pseudo-energy $E$ in term of the radial deformation~$\xi$. 
\begin{lemma}
Let $ D(t) $ 
be a bounded simply connected region with smooth boundary $ \pa D(t) $ 
parametrized as
in \eqref{time-para}. 
Then
the pseudo-energy $E$, defined in   \eqref{pr-vp},   is 
  \be
 E (\xi) =
 \frac{ 1}{32\pi} \int_{\T^2}\big[\ln \big(\Ke(\xi)(\theta, \theta')\big)-2\big] \Ke(\xi)(\theta, \theta') \, \partial^2_{\theta\theta'} \Ke(\xi)(\theta, \theta')d\theta'd\theta   \label{forE}
\ee
 where $\Ke(\xi)(\theta, \theta') $ is given by \eqref{expR}.
Note that $ E $ 
is an autonomous prime integral of \eqref{Evera}. 
\end{lemma}
\begin{pf}
By  \eqref{pr-vp}, \eqref{stream-c} and Green's formula \eqref{green},
the pseudo-energy $E$ can be written as 
$$
E=\frac{ 1}{32\pi}\int_{\partial D(t)}\int_{\partial D(t)}\Big[\ln\big(\vert \zeta-z \vert^2\big)-\frac32\Big](\overline{\zeta}-\overline{z})^2dz  d\zeta \, .
$$
Using the parametrization \eqref{time-para} of the boundary $ \pa D(t) $
we deduce that 
\begin{align}
E
 &=\frac{1}{32\pi}\int_{\T^2}\Big[\ln \big(|w(\theta')- w(\theta)|^2\big)-\frac32\Big]\big( \overline{w(\theta')}- \overline{w(\theta)}\big)^2\partial_{\theta} w(\theta)\partial_{\theta'} w(\theta')d\theta'd\theta\, ,
\label{E}
\end{align}
where, to simplify notations, we have denoted $w(\theta):=w(t,\theta)$.
Inserting the identity 
$$
\big( \overline{w(\theta')}- \overline{w(\theta)}\big)^2\partial_{\theta} w(\theta)=-|w(\theta')- w(\theta)|^2\partial_{\theta}\overline{ w(\theta)}-\big( \overline{w(\theta')}- \overline{w(\theta)}\big)\partial_{\theta}|w(\theta')- w(\theta)|^2  
$$
into \eqref{E} we get
\begin{align*}
E&=-\frac{1}{32\pi}\int_{\T^2}\Big[\ln\big( |w(\theta')- w(\theta)|^2\big)-\frac32\Big]|w(\theta')- w(\theta)|^2\partial_{\theta}\overline{ w(\theta)}\partial_{\theta'} w(\theta')d\theta'd\theta
\\&\quad-\frac{ 1}{32\pi}\int_{\T^2}\Big[\ln( |w(\theta')- w(\theta)|^2\big)-\frac32\Big]
\big( \partial_{\theta}|w(\theta')- w(\theta)|^2 \big) \big( \overline{w(\theta')}- \overline{w(\theta)}\big)\partial_{\theta'} w(\theta')d\theta'd\theta \\
&=-\frac{1}{32\pi}\int_{\T^2}\Big[\ln\big( |w(\theta')- w(\theta)|^2\big)-\frac32\Big]|w(\theta')- w(\theta)|^2\partial_{\theta}\overline{ w(\theta)}\partial_{\theta'} w(\theta')d\theta'd\theta
\\&\quad-\frac{ 1}{32\pi}\int_{\T^2}\partial_{\theta}\Big(\Big[\ln\big( |w(\theta')- w(\theta)|^2\big)-\frac52\Big]|w(\theta')- w(\theta)|^2\Big)\big( \overline{w(\theta')}- \overline{w(\theta)}\big)\partial_{\theta'} w(\theta')d\theta'd\theta \, , 
\end{align*}
having used the identity 
$$
\partial_{\theta}\Big(\Big[\ln \big(|w(\theta')- w(\theta)|^2\big)-\frac52\Big]|w(\theta')- w(\theta)|^2\Big) 
=\Big[\ln |w(\theta')- w(\theta)|^2-\frac32\Big]\partial_{\theta}|w(\theta')- w(\theta)|^2 \, . 
$$
Integrating the last term by parts we obtain
\begin{align}
E  & 
=-\frac{1}{16\pi}\int_{\T^2}\big[\ln \big(|w(\theta')- w(\theta)|^2\big)-2\big]|w(\theta')- w(\theta)|^2\partial_{\theta}\overline{ w(\theta)}\partial_{\theta'} w(\theta')d\theta'd\theta \notag \\
& =-\frac{1}{16\pi }\int_{\T^2}\big[\ln \big(|w(\theta')- w(\theta)|^2\big)-2\big]|w(\theta')- w(\theta)|^2\textnormal{Re}\big[\partial_{\theta} \overline{ w(\theta)}\partial_{\theta'} w(\theta')\big]d\theta'd\theta \label{ultE}
\end{align}
since the pseudo-energy $E$ is real valued. Finally,
by the identity 
$ -2\textnormal{Re}\big[\partial_{\theta} \overline{ w(\theta)}\partial_{\theta'} w(\theta') \big] = \partial^2_{\theta\theta'}|w(\theta')- w(\theta)|^2 $, 
we conclude by \eqref{ultE} that
\begin{align}\label{forEw}
E&=\frac{1}{32\pi }\int_{\T^2}\big[\ln\big( |w(\theta')- w(\theta)|^2\big)-2\big]|w(\theta')- w(\theta)|^2\partial^2_{\theta\theta'}|w(\theta')- w(\theta)|^2d\theta'd\theta\, .
\end{align}
Then \eqref{forE} follows by \eqref{forEw} and \eqref{zq}.
\end{pf}

\subsection{Hamiltonian and reversible structure}\label{sec:HAM}
The main result of this section is the following: 
\begin{proposition} \label{prop:Ham}
{\bf (Hamiltonian formulation of the contour dynamics equation)}
The PDE \eqref{Evera} is the Hamiltonian equation
\be\label{eq:main}
\pa_t \xi = \pa_\theta \nabla H_\Omega (\xi) \, , 
\ee
where $ \nabla H_\Omega $ denotes the $ L^2 $-gradient of the  Hamiltonian
\be 
\label{Ham:main}
H_\Omega (\xi) := -\tfrac12 E (\xi) +  \tfrac{\Omega}{2} J  (\xi)\, , 
\ee
where the pseudo-energy $ E  $ and the angular momentum $ J $ are  expressed as  a
function of $ \xi $ respectively in \eqref{forE}  
and  \eqref{forCZJ}.  
\end{proposition}
Before proving Proposition \ref{prop:Ham} we describe 
the phase space and the Hamiltonian aspects  of \eqref{eq:main}.
 The average 
$$
\langle \xi \rangle_{\theta} := \frac{1}{2 \pi} \int_{\T}
\xi(\theta) d \theta 
$$ 
is a prime integral of \eqref{eq:main}, coherently with 
the existence of the prime integral  $ C $ 
in Lemma~\ref{lem:prin}.
Hence we  consider the phase space $ H^s_0 := H^s_0 (\T) := H^s_0 (\T, \R) $ 
of periodic real functions with zero average, 
\be\label{phase-space}
H^s_0 (\T, \R) := \Big\{ \xi(\theta) := \sum_{j \in \Z \setminus \{0\}} \xi_j e^{\ii j \theta} \, , 
\
\overline{\xi_j} = \xi_{-j} \, , \
\| \xi \|_s^2 = \sum_{j \in \Z \setminus \{0\}}|\xi_j|^2 |j|^{2s} < + \infty \Big\} \,.
\ee
We  denote by $ \Pi_0^\bot $ the orthogonal projection on $ H_0^s$. 
Note that the Hamiltonian vector field $X_{H_\Omega} (\xi) :=\partial_\theta\nabla H_\Omega (\xi) $, associated to the Hamiltonian $H_\Omega(\xi ) $, is determined by the identity
\be\label{Hamide}
dH_\Omega(\xi)[ \cdot ] =\mathcal{W} ( X_{H_\Omega}(\xi),  \, \cdot\, ) \, , 
\ee
where ${\cal W}$ is the non-degenerate symplectic form defined,
for any $ \xi_1,\xi_2\in L_0^2(\mathbb{T})  $, as
\be\label{sy2form}
{\cal W} (\xi_1, \xi_2) := \int_{\T}(\pa_\theta^{-1} \xi_1) (\theta) \, \xi_2 (\theta) d \theta \, , 
\quad\partial_\theta^{-1} \xi= {\mathop \sum}_{j \in \Z \setminus \{0\}}\frac{1}{\ii j} \xi_j e^{\ii j \theta}  \, .  
\ee
The corresponding Poisson tensor is $\partial_\theta$ and the  Poisson bracket  is 
\be\label{Poib}
\{ F(\xi), G(\xi)\} =\mathcal{W}(X_F,X_G)= \int_{\T} \nabla F (\xi)  \pa_\theta\nabla G (\xi)  \, d \theta = d F(\xi) [ X_G (\xi)]
\ee
where $ \nabla F $, $ \nabla G $  denote the $ L^2 $-gradients  of $ F, G  $.

In order to prove Proposition \ref{prop:Ham} we need  the following two lemmata.
\begin{lemma}\label{lem:gradE}{\bf (Pseudo-energy)} \label{lemma:expression-Energy}
The $ L^2$-gradient of the pseudo energy $E$, given by \eqref{forE},  is 
\begin{align}
\nabla E(\xi)(\theta) 
 & =   \frac{1}{4\pi}\int_{\T}\big[\ln\big( |w(\theta)- w(\theta')|^2\big)-1\big] 
 {\rm Im}\big[ \big( \overline{w(\theta')}- \overline{w(\theta)}\, \big)\partial_{\theta'}  w(\theta') \big] d\theta' \label{nablaE} \\
& = 2 \psi(w( \theta))  \label{deltaE}
\end{align}
where $ \psi $ is the stream function given by \eqref{formapsi} and   $w(\theta)$ is introduced in \eqref{time-para}. 
Moreover 
\begin{align}\label{grad-pse}
\pa_\theta \nabla E (\xi)&=\frac{1}{2 \pi} 
\int_{\T} \ln ( \Ke(\xi)(\theta, \theta') ) \pa_{\theta \theta'}^2 
\big[ (1 + 2 \xi (\theta))^\frac12
(1 + 2 \xi (\theta'))^\frac12 \sin (\theta' - \theta )\big] d \theta' \, .
\end{align}
\end{lemma}
\begin{pf}
Differentiating  \eqref{forEw} with respect to $\xi$ 
we get that $ dE(\xi)[q] $ is equal to 
\begin{align*}
&\frac{ 1}{32\pi}\int_{\T^2 }\big[\ln  \big( |w( \theta')-w( \theta)|^2\big)-1\big]\big(A(\xi)(\theta,\theta')q(\theta)+A(\xi)(\theta',\theta)q(\theta')\big) \partial^2_{\theta\theta'}   |w( \theta')-w( \theta)|^2d\theta'd\theta
\\&+\frac{ 1}{32\pi}\int_{\T^2} \big[\ln  \big( |w( \theta')-w( \theta)|^2\big)-2\big] |w( \theta')-w( \theta)|^2\partial^2_{\theta\theta'} \big(A(\xi)(\theta,\theta')q(\theta)+A(\xi)(\theta',\theta)q(\theta')\big) d\theta'd\theta 
\end{align*}
with
\be
A(\xi)(\theta,\theta') :=2\big(1+2\xi(\theta)\big)^{-\frac12}\textnormal{Re}\big[\we(\theta)\big(\overline{w(\theta)}- \overline{w(\theta')}\big)\big]
 \, .\label{B0}
\ee
By symmetry, we may write
\be
\begin{aligned}
dE(\xi)[q]&=  \frac{ 1}{16\pi}\int_{\T^2} \big[\ln  \big( |w( \theta')-w( \theta)|^2\big)-1\big] \partial^2_{\theta\theta'}  |w( \theta')-w( \theta)|^2A(\xi)(\theta,\theta')q(\theta)d\theta'd\theta
\\&\quad+\frac{ 1}{16\pi}\int_{\T^2} \big[\ln \big(  |w( \theta')-w( \theta)|^2\big)-2\big]  |w( \theta')-w( \theta)|^2\partial^2_{\theta\theta'}  \big(A(\xi)(\theta,\theta')q(\theta)\big)d\theta'd\theta \, .  \label{L3}
\end{aligned}
\ee
Integrating  by parts the second term in   \eqref{L3}  gives
\begin{align}
&dE(\xi)[q]
=\frac{ 1}{16\pi}\int_{\T^2} \big[\ln  \big( |w( \theta')-w( \theta)|^2\big)-1\big] B(\xi)(\theta, \theta')q(\theta)d\theta'd\theta \label{L33}
\,  
\end{align}
with 
$ B(\xi)(\theta, \theta')  := $ $ -\partial_{\theta}   |w( \theta')-w( \theta)|^2\partial_{\theta'} A(\xi)(\theta, \theta')+ $ $ \partial^2_{\theta\theta'}  |w( \theta')-w( \theta)|^2 A(\xi)(\theta, \theta') $. 
By \eqref{time-para} and \eqref{B0}  one has
\begin{align*}
\partial_{\theta}  |w( \theta')-w( \theta)|^2 & 
= \partial_{\theta}\xi(\theta)A(\xi)(\theta, \theta')+2\big(1+2\xi(\theta)\big)^{\frac12}\textnormal{Re}\big[\partial_{\theta} \we( \theta)\big(\overline{w(\theta)}- \overline{w(\theta')}\big)\big]\, ,\\
\partial^2_{\theta\theta'}  |w( \theta')-w( \theta)|^2 & 
= \partial_{\theta}\xi(\theta) \partial_{\theta'} A(\xi)(\theta, \theta')-2\big(1+2\xi(\theta)\big)^{\frac12}\textnormal{Re}\big[\partial_{\theta} \we( \theta)\partial_{\theta'} \overline{w(\theta')}\big)\big]\, . 
\end{align*}
It follows that
\begin{align*}
 B(\xi)(\theta, \theta')
 &=-2\big(1+2\xi(\theta)\big)^{\frac12}\textnormal{Re}\big[\partial_\theta \we( \theta)\big(\overline{w(\theta)}- \overline{w(\theta')}\big)\big] \partial_{\theta'} A(\xi)(\theta, \theta')
 \\ &\quad-2\big(1+2\xi(\theta)\big)^{\frac12}\textnormal{Re}\big[\partial_\theta \we( \theta)\partial_{\theta'} \overline{w(\theta')}\big]A(\xi)(\theta, \theta') 
 \stackrel{\eqref{B0}} =  4\,\textnormal{Re}\big[\partial_{\theta'} \overline{w(\theta')}D(\xi)(\theta, \theta')\big]
\end{align*}
with
$ D(\xi)(\theta, \theta'):= $ $ \we( \theta)\textnormal{Re}\big[\partial_\theta \we( \theta)\big(\overline{w(\theta)}- \overline{w(\theta')}\big)\big]-\partial_\theta \we( \theta)\textnormal{Re}\big[\we( \theta)\big(\overline{w(\theta)}- \overline{w(\theta')}\big)\big] $. 
Straightforward computations, using the identity
\be\label{id:we}
{\rm Im}\big[\overline{ \we(\theta)}\partial_\theta \we(\theta)\big]=1\, ,
\ee
leads to
$$
 D(\xi)(\theta, \theta')
 =\ii\,\textnormal{Im}\big[\we(\theta)\partial_{\theta}\overline{\we(\theta)}\big]\big(w(\theta)- w(\theta')\big)   =-\ii\big(w(\theta)- w(\theta')\big) \, . 
$$
Therefore
$ B(\xi)(\theta, \theta') = 4\,\textnormal{Im}\big[\partial_{\theta'} \overline{w(\theta')}\big(w(\theta)- w(\theta')\big)\big] $. 
Inserting the last identity into \eqref{L33} we get
\begin{align*}
dE(\xi)[q]
 &=\frac{1}{4\pi}\int_{\T^2}\big[\big(\ln |w(\theta')- w(\theta)|^2\big)-1\big]\textnormal{Im}\big[\partial_{\theta'} \overline{w(\theta')}\big(w(\theta)- w(\theta')\big)\big]q(\theta)d\theta'd\theta \, , 
\end{align*}
implying \eqref{nablaE}.  
Comparing \eqref{nablaE} and \eqref{formapsi} we deduce \eqref{deltaE}.   Then,   by \eqref{deltaE} and \eqref{parpsi-new} we get \eqref{grad-pse}. This ends the proof of Lemma \ref{lem:gradE}.
\end{pf}

\noindent
{\bf Angular momentum.} 
The angular momentum $J$ in \eqref{forCZJ}  can be expanded as
 \begin{equation} 
J  =\frac{\pi}{4}(\gamma + \gamma^{-1})+ J_1 + J_2\, ; \quad
J_1  := \int_{\T}  \xi (\theta)  g_\gamma (\theta) 
\, d \theta \, , \quad 
J_2  :=  
\int_{\T}  \xi^2 (\theta) g_\gamma (\theta)
\, d \theta \, . \label{moangJ2}
\end{equation}
\begin{lemma}\label{lem:AM}
The Hamiltonian vector field $ X_J (\xi)  $ 
associated to the angular momentum $ J $ in \eqref{forCZJ}  is  
\be\label{Mom1}
\pa_\theta \nabla J (\xi) = 
\pa_\theta 
\big[\big(1+2\xi(\theta)\big)g_\gamma(\theta)\big] \, ,
\ee
where $ g_\gamma(\theta) $ is given by \eqref{fg0}.
\end{lemma}

\begin{pf}
Differentiating  \eqref{moangJ2} 
we get 
$ dJ(\xi)[q] = \int_{\T}q(\theta)\big(1+2\xi(\theta)\big)g_\gamma(\theta)d\theta $
and therefore 
$ \nabla J(\xi) 
= \Pi_0^\bot \big(1+2\xi(\theta)\big)g_\gamma(\theta)  $
proving \eqref{Mom1}. 
\end{pf}

\smallskip

\noindent{\sc Proof of Proposition \ref{prop:Ham}}.\enspace  Follows comparing  \eqref{Evera} with \eqref{eq:main}, \eqref{Ham:main}, \eqref{grad-pse}, \eqref{Mom1}. \rule{2mm}{2mm}
\\[1mm]
{\bf Reversible structure.} 
 We finally  point out that the  2$d$-Euler equation  \eqref{eq:euler} is a reversible 
system and this property persists at the level  of the contour dynamics equation \eqref{eq:main}. 
We introduce the involution
\be\label{S-invo}
( {\cal S} \xi )(\theta) :=  \xi (- \theta ) \, . 
\ee
Notice that $ {\cal S} $ satisfies 
\be\label{Stra}
{\cal S}^2 = {\rm Id} \, , \quad 
{\cal S}^\top = {\cal S} \, , \quad \pa_\theta \circ {\cal S} = - {\cal S} \circ \pa_\theta \, , 
\ee
where the transpose is taken with respect to the $ L^2 $-scalar product. 
\begin{lemma} {\bf (Reversibility)}\label{lem:rev}
The Hamiltonian vector field 
$ X_{H_{\Omega}} = \pa_\theta \nabla H_\Omega $ is reversible with respect to the involution ${\cal S}$ defined in \eqref{S-invo}, namely 
\be\label{XrevS}
X_{H_{\Omega}} \circ {\cal S} = - {\cal S} \circ X_{H_{\Omega}} \, .
\ee
Equivalently the  Hamiltonian $ H_{\Omega} $ satisfies
\be\label{XrevSH}
H_\Omega \circ {\cal S} = H_\Omega \, . 
\ee
\end{lemma}

\begin{pf}
In view of the expression \eqref{forCZJ} of 
$ J $, and since $ g_\gamma (\theta )$ is even, it follows that
$ J \circ {\cal S} = J $. Moreover the
 pseudo-energy $ E $ given in \eqref{forE} satisfies 
 $ E \circ {\cal S} = E $, since   the function $ \Ke(\xi)  (\theta, \theta') $ 
in \eqref{expR} satisfies 
$ \Ke(\xi)(-\theta,-\theta') = \Ke({\cal S}\xi)(\theta,\theta') $. 
In view of \eqref{Ham:main} the identity \eqref{XrevSH} is proved. 
\end{pf}

\begin{remark} {\bf (Half-wave symmetry)} \label{HWS} 
The Hamiltonian PDE  \eqref{eq:main},  i.e. \eqref{Evera},
leaves invariant the subspace  
$ \big\{ \xi (\theta) \, : \, \xi (\theta + \pi) = \xi (\theta) \} $ of $ \pi $-periodic functions, which, 
in the Fourier expansion 
\eqref{qthanbn}, have only 
the even harmonics $ (\alpha_n,\beta_n)$. 
This symmetry is related to the center $ Z  $ in \eqref{forCZJ}:
if $ \xi (\theta) $ is $ \pi $-periodic then $ Z = 0 $. We  
do not restrict to this subspace,   considering more general solutions.
\end{remark}

\subsection{Hamiltonian and reversible linear operators}

The Hamiltonian and reversible properties of the linearized equation that we now describe
play a  key  role. 

Along the paper we encounter $\vphi$-dependent (possibly also constant in $ \vphi $)  linear operators 
$A:\T^{|\ST|}\mapsto {\cal L}(L^2(\T_\theta))$, $\vphi\mapsto A(\vphi)$, acting on subspaces of $L^2(\T_\theta)$. We also regard $A$ as an operator (denoted for simplicity by $A$ as well) that acts on functions $ \xi (\vphi,\theta)$, 
that is $ (A \xi )(\vphi,\theta) := ( A(\vphi) \xi (\vphi,\,\cdot \,) )(\theta) $.

\begin{definition} {\bf (Hamiltonian)}\label{def:HL}
A  linear real 
 operator of the form $ \pa_\theta \circ A (\vphi) $
where $  A (\vphi) $ is self-adjoint, is called Hamiltonian.  
We also say that
$ \om \cdot \pa_\vphi + \pa_\theta A(\vphi )   $
is a Hamiltonian operator.
\end{definition}

\begin{definition} {\bf (Symplectic)}\label{def:sympl}
A $ \vphi $-dependent family $ \Phi(\vphi ) $ of linear transformations of the phase space 
 is symplectic, if, for any $ \vphi \in \T^{|\ST|} $, each $ \Phi(\vphi ) $ 
preserves the symplectic $ 2 $-form defined in \eqref{sy2form}, i.e. 
$ {\cal W}( \Phi(\vphi )  u, \Phi(\vphi ) v ) = {\cal W}(  u,  v ) $. 
Equivalently
$ \Phi(\vphi )^\top \circ \pa_\theta^{-1} \circ \Phi(\vphi ) = \pa_\theta^{-1} $. 
\end{definition}

A Hamiltonian operator transforms under conjugation with a symplectic  
family of $ \vphi $-dependent  transformations $ \Phi (\vphi ) $
 into another Hamiltonian operator, see e.g. \cite{BBM-Airy}.

\begin{definition}\label{def:R-AR}
{\bf (Reversible)}. 
A linear  operator
$  A (\vphi) $
is reversible, resp. reversibility preserving, if
$ A( - \vphi) \circ {\cal S} = - {\cal S} \circ A(  \vphi) $,  resp. 
$ A( - \vphi) \circ {\cal S} =  {\cal S} \circ A(  \vphi) $, 
where $ {\cal S} $ is the involution defined in \eqref{S-invo}. 
\end{definition}

Composition of reversible operators with reversibility preserving operators is reversible.  
Reversible multiplication and integral operators are characterized as follows:

\begin{lemma}\label{lem:IOK}
A $ \vphi$-dependent family of  multiplication
 operators  for the function $ V (\vphi, \theta )$ is 
reversible if and only if 
$ V( \vphi, \theta) $ is odd in $ (\vphi, \theta ) $; 
reversibility preserving 
if and only if  $ V( \vphi, \theta) $ is even in $ (\vphi, \theta ) $.
A $ \vphi$-dependent family of  integral operators 
\be\label{def:int-op}
\big( {\cal K}(\varphi) \xi \big) (\theta) := \int_\T K( \varphi, \theta, \theta') \xi (\theta')\,d \theta' 
\ee
 is
reversible iff $ K(-\vphi, -\theta, -\theta') =  - K(\vphi, \theta, \theta') $; 
reversibility-preserving iff $ K(-\vphi, -\theta, -\theta') =  K(\vphi, \theta, \theta') $.
\end{lemma}

\begin{definition} \label{defRAR} {\bf (Reversible and anti-reversible function)}
A function $ \xi (\vphi, \cdot) $ is called 
{\sc Reversible} if $ {\cal S} \xi (\vphi, \cdot ) = \xi(-\vphi, \cdot ) $;  
{\sc Anti-reversible} if $ - {\cal S} \xi (\vphi, \cdot ) = \xi(-\vphi, \cdot) $. 
\end{definition}

A reversibility preserving operator maps reversible, respectively anti-reversible, functions into 
reversible, respectively anti-reversible functions. 

\begin{lemma}\label{lem:REV}
Let $ X $ be a reversible vector field, cfr. \eqref{XrevS}.
Let $  \xi(\vphi, \theta)  $ be a  reversible quasi-periodic function, according to Def. 
\ref{defRAR}.
	Then the linearized operator  $ d_\xi X( \xi(\vphi, \cdot) ) $  is reversible, cfr. 
	Def. 
	\ref{def:R-AR}.
\end{lemma}

We finally note that the involution $ {\cal S} $ defined in \eqref{S-invo} is reversible, 
by  \eqref{Stra}, and therefore 
a Hamiltonian operator $ \pa_\theta A(\vphi) $ 
is   reversible if and only if  $ A(\vphi) $ is reversibility preserving.

\section{Linearized vector field of $ X_{H_{\Omega}} $ at any $ \xi $}\label{sec:Lin}

The goal  of this section is to compute the linear Hamiltonian PDE obtained 
 linearizing  the equation \eqref{eq:main}, i.e. \eqref{Evera}, at any state $ \xi (\theta )$. 
 The main result of this section is the following proposition: 
\begin{proposition}\label{lin-eq}
The linearized equation of \eqref{eq:main}, i.e. \eqref{Evera}, 
at  a small state $ \xi (  \theta ) $ is 
\begin{align}
\pa_t q (\theta) & = 
  \pa_\theta
 \Big(   \big(\Omega\,  g_\gamma(\theta)+v(\xi) (\theta ) \big) q ( \theta)-\W (\xi)\, [ q] (\theta)\Big)
 \label{linVF1}
\end{align}
where $v (\xi)(\theta)  $ is the real  function
\be\label{defV}
v ( \xi )(\theta) := \frac{1}{4 \pi } 
\int_{\T} \ln ( \Ke(\xi)(\theta, \theta')) \pa_{\theta'} 
\Big[  \Big(\frac{1 + 2 \xi (\theta')}{1 + 2 \xi (\theta)}\Big)^{\frac12} \sin (\theta' - \theta) \Big] d \theta'\, ,
\ee
$\W (\xi)$ is the self-adjoint integral  operator
\be\label{int-op}
\W( \xi) [ q] ( \theta):= 
\frac{1}{4  \pi} \int_{\T} \ln \big( \Ke(\xi) (\theta, \theta') \big) q(\theta') \, d \theta' \, ,  
\ee
the real  function $g_\gamma(\theta)$ is given by \eqref{fg0} and   $ \Ke(\xi)(\theta, \theta')  $ is defined in \eqref{expR}.
\end{proposition}

\begin{pf}
In view of  Proposition \ref{prop:Ham},  one has that 
\be\label{lin:nablaH}
d X_{H_\Omega} (\xi ) [q] (\theta ) = \partial_\theta \big(-\tfrac12 d \nabla E (\xi)[q] (\theta )+\tfrac{\Omega}{2} d \nabla J (\xi)[q] (\theta )\big) \, . 
\ee
By \eqref{Mom1} we have  
\be\label{lin:nablaJ}
d \nabla J (\xi)[q] (\theta )=2g_\gamma(\theta)q(\theta ) \, . 
\ee
Thus, it remains to compute $d \nabla E (\xi)[q] $. Differentiating \eqref{nablaE} with respect to $\xi$ % in the direction $q$ 
we get
 \begin{align*}
d \nabla E (\xi)[q](\theta&) =(1 + 2 \xi(\theta))^{ -\frac12} q(\theta)\frac{1}{4\pi}\int_{\T}\tfrac{ 2{\rm Re}\big[ \big( \overline{w(\theta)}- \overline{w(\theta')}\, \big) \we(\theta) \big]}{|w(\theta)- w(\theta')|^2}
 {\rm Im}\big[ \big( \overline{w(\theta')}- \overline{w(\theta)}\, \big)\partial_{\theta'}  w(\theta') \big]d\theta'
 \\
 &
 +\frac{1}{4\pi}\int_{\T}\tfrac{2 {\rm Re}\big[ \big( \overline{w(\theta')}- \overline{w(\theta)}\, \big) \we(\theta') \big]q(\theta')}{|w(\theta)- w(\theta')|^2(1 + 2 \xi(\theta'))^{\frac12}}
 {\rm Im}\big[ \big( \overline{w(\theta')}- \overline{w(\theta)}\, \big)\partial_{\theta'}  w(\theta') \big]d\theta'
 \\
 &-(1 + 2 \xi(\theta))^{ -\frac12} q(\theta)\frac{1}{4\pi}\int_{\T}\ln\big( |w(\theta)- w(\theta')|^2\big) 
\partial_{\theta'}   {\rm Im}\big[ \overline{\we(\theta)} w(\theta') \big]d\theta'\\
  &+\frac{1}{4\pi}\int_{\T}\big[\ln \big(|w(\theta)- w(\theta')|^2\big)-1\big] 
 {\rm Im}\big[\overline{\we(\theta')}\, \partial_{\theta'}  w(\theta') \big] (1 + 2 \xi(\theta'))^{ -\frac12} q(\theta')d\theta'
 \\
 &+\frac{1}{4\pi}\int_{\T}\big[\ln \big(|w(\theta)- w(\theta')|^2\big)-1\big] 
 {\rm Im}\big[ \big( \overline{w(\theta')}- \overline{w(\theta)}\, \big)\partial_{\theta'}\big( (1 + 2 \xi(\theta'))^{ -\frac12} q(\theta')\we(\theta')\big)  \big]  d\theta' \, . 
\end{align*}
The last term  writes
 \begin{align*}
I_5&:=\frac{1}{4\pi}\int_{\T}\big[\ln \big(|w(\theta)- w(\theta')|^2\big)-1\big] 
 {\rm Im}\big[ \big( \overline{w(\theta')}- \overline{w(\theta)}\, \big)\partial_{\theta'}\big( (1 + 2 \xi(\theta'))^{ -\frac12} q(\theta')\we(\theta')\big)  \big]  d\theta'
 \\ &=\frac{1}{4\pi}\int_{\T}\big[\ln\big( |w(\theta)- w(\theta')|^2\big)-1\big] 
 {\rm Im}\big[ \big( \overline{w(\theta')}- \overline{w(\theta)}\, \big)\partial_{\theta'}\we(\theta')  \big] (1 + 2 \xi(\theta'))^{ -\frac12} q(\theta')\we(\theta')  d\theta'
 \\ &\quad+\frac{1}{4\pi}\int_{\T}\big[\ln\big( |w(\theta)- w(\theta')|^2\big)-1\big] 
 {\rm Im}\big[ \big( \overline{w(\theta')}- \overline{w(\theta)}\, \big)\we(\theta') \big] \partial_{\theta'}\big( (1 + 2 \xi(\theta'))^{ -\frac12} q(\theta') \big) d\theta'
 \,  .
\end{align*}
Integrating by parts the second term in $I_5$ we get
 \begin{align*}
I_5 &=-\frac{1}{4\pi}\int_{\T}\frac{\partial_{\theta'} |w(\theta)- w(\theta')|^2}{|w(\theta)- w(\theta')|^2}
 {\rm Im}\big[ \big( \overline{w(\theta')}- \overline{w(\theta)}\, \big)\we(\theta') \big] (1 + 2 \xi(\theta'))^{ -\frac12} q(\theta') d\theta'
 \\ &\quad-\frac{1}{4\pi}\int_{\T}\big[\ln\big( |w(\theta)- w(\theta')|^2\big)-1\big] 
 {\rm Im}\big[ \partial_{\theta'} \overline{w(\theta')}\, \we(\theta') \big](1 + 2 \xi(\theta'))^{ -\frac12} q(\theta')d\theta'
 \,  .
\end{align*}
Inserting $I_5$ into the expression of $  d\nabla E (\xi) $, we obtain
\begin{align}\label{eq0}
\nonumber d &\nabla  E (\xi)[q](\theta) =q(\theta)(1 + 2 \xi(\theta))^{-\frac12}\frac{1}{4\pi}\int_{\T}\tfrac{ 2{\rm Re}\big[ \big( \overline{w(\theta)}- \overline{w(\theta')}\, \big) \we(\theta) \big]}{|w(\theta)- w(\theta')|^2}
 {\rm Im}\big[ \big( \overline{w(\theta')}- \overline{w(\theta)}\, \big)\partial_{\theta'}  w(\theta') \big]d\theta'\\
\nonumber  &-q(\theta)(1 + 2 \xi(\theta))^{ -\frac12} \frac{1}{4\pi}\int_{\T}\ln \big(|w(\theta)- w(\theta')|^2\big)\partial_{\theta'}  {\rm Im}\big[ \overline{\we(\theta)} w(\theta') \big]d\theta'+\frac{1}{4\pi}\int_{\T}\tfrac{ D(\theta',\theta) q(\theta')}{|w(\theta)- w(\theta')|^2}  
d\theta' \\
 & +\frac{1}{2\pi}\int_{\T}\big[\ln \big(|w(\theta)- w(\theta')|^2\big)-1\big] 
 {\rm Im}\big[\overline{\we(\theta')}\, \partial_{\theta'}  \we(\theta') \big]  q(\theta')d\theta'
\, ,
\end{align}
where
\begin{align*}
D(\theta',\theta)&:=2{(1 + 2 \xi(\theta'))^{-\frac12}} {{\rm Re}\big[ \big( \overline{w(\theta')}- \overline{w(\theta)}\, \big) \we(\theta') \big]} {\rm Im}\big[ \big( \overline{w(\theta')}- \overline{w(\theta)}\, \big)\partial_{\theta'}  w(\theta') \big]\\
&\qquad -{(1 + 2 \xi(\theta'))^{-\frac12}}{\partial_{\theta'}  |w( \theta')-w( \theta)|^2} {\rm Im}\big[ \big( \overline{w(\theta')}- \overline{w(\theta)}\, \big)\we(\theta') \big] 
\\ 
&=  2 |w(\theta)- w(\theta')|^2 {\rm Im}\big[   \overline{\we(\theta')}\partial_{\theta'}  \we(\theta')\big]    \stackrel{\eqref{id:we}} =2 |w(\theta)- w(\theta')|^2 \, . 
\end{align*}
Moreover, straightforward computations, using  \eqref{time-para} and \eqref{id:we},  lead to
\begin{align*}
2{{\rm Re}\big[ \big( \overline{w(\theta)}- \overline{w(\theta')}\, \big) \we(\theta) \big]{\rm Im}\big[ \big( \overline{w(\theta')}- \overline{w(\theta)}\, \big)\partial_{\theta'}  w(\theta') \big]}&={\partial_{\theta'}  |w( \theta')-w( \theta)|^2}{{\rm Im}\big[ \overline{\we(\theta)}w(\theta')\big]}\\ &\quad -2{|w(\theta)- w(\theta')|^2} {\partial_{\theta'}  {\rm Im}\big[\overline{\we(\theta)} w(\theta') \big]}\, .
\end{align*}
Inserting  the last two identity    into \eqref{eq0} and using \eqref{id:we}  we get 
\begin{align}\label{eq111} 
d \nabla E (\xi)[q](\theta) &=q(\theta)(1 + 2 \xi(\theta))^{-\frac12}\frac{1}{4\pi}\int_{\T}\tfrac{\partial_{\theta'}  |w( \theta')-w( \theta)|^2}{|w(\theta)- w(\theta')|^2}\,{\rm Im}\big[ \overline{\we(\theta)}w(\theta') \big]d\theta'\nonumber
  \\ \nonumber
 &\quad- q(\theta)(1 + 2 \xi(\theta))^{-\frac12}\frac{1}{4\pi}\int_{\T}\ln \big(|w(\theta)- w(\theta')|^2\big)
\partial_{\theta'} {\rm Im}\big[ \overline{\we(\theta)}  w(\theta') \big]d\theta'   \\
 &\quad+\frac{1}{2\pi}\int_{\T}\ln\big( |w(\theta)- w(\theta')|^2\big)q(\theta')d\theta'  
\end{align}
having used  that $q(\theta)$  has zero average.
Integrating the first term by parts we conclude that 
\begin{align*}
 d\nabla E (\xi) [q](\theta)&=-q(\theta)(1 + 2 \xi(\theta))^{-\frac12}\frac{1}{2\pi}\int_{\T}\ln \big(|w(\theta)- w(\theta')|^2\big)\,
 \partial_{\theta'} {\rm Im}\big[ \overline{\we(\theta)} w(\theta')\big]d\theta'
\\ &\quad+\frac{1}{2\pi}\int_{\T}\ln \big(|w(\theta)- w(\theta')|^2\big)q(\theta')d\theta'
 \, . 
\end{align*}
This with \eqref{zq} and \eqref{time-para} yields
\be\label{lin:nablaE}
d\nabla E(\xi)[q] (\theta) =  -2 v(\xi) (\theta ) q( \theta)+2 \W (\xi)\, [ q] ( \theta) \, ,
\ee
where  $ v(\xi) (\theta )$  and $ \W (\xi) [ q] ( \theta)$ are  defined in \eqref{defV} and \eqref{int-op}, respectively. Putting together \eqref{lin:nablaH},  \eqref{lin:nablaJ} and \eqref{lin:nablaE} completes the proof of Proposition \ref{lin-eq}.
The integral operator $\W (\xi) $ is  self-adjoint because 
$  \Ke(\xi) (\theta, \theta') $ is real and symmetric, i.e. 
$  \Ke(\xi) (\theta, \theta')  =  \Ke(\xi) (\theta', \theta)   $. 
\end{pf}

\section{Linearized vector field at $ \xi = 0 $}\label{sec:3}
 
In this section we compute 
the linear equation in \eqref{linVF1} at the equilibrium $\xi=0 $ and find 
 its solutions.
By Proposition  \ref{lin-eq}  the linearized equation  of \eqref{eq:main} (that is \eqref{Evera}) with  $ \Omega=\Omega_\gamma  $  at  $ \xi = 0 $,  is 
\begin{equation}\label{linunp}
\pa_t q (\theta) =   \pa_\theta 
 \Big( \big(\Omega_\gamma g_\gamma(\theta)+v_0 ( \theta ) \big) q(\theta)- \W_0\, [q ](\theta)\Big) 
\end{equation}
where, by \eqref{defV} and \eqref{int-op},  
\begin{align}\label{defV0}
v_0 (\theta) & := v(0) (\theta) := \frac{1}{4 \pi} 
\int_{\T} \ln(\Ke(0)(\theta, \theta'))  \cos (\theta' - \theta) d \theta'\, ,  
\\
\W_0 [q] & := \W (0) [q]  := \frac{1}{4 \pi}  \int_{\T} \ln ( \Ke(0)(\theta, \theta')) q (\theta' ) d \theta' \label{defK0y}
\end{align}
with $\Ke(0)(\theta, \theta') $ defined in \eqref{expR0}.  We first compute $ v_0 (	\theta)$. 

\begin{lemma}\label{lem:v0}
The function $ v_0(\theta )$ in \eqref{defV0} is  
\be\label{cv0}
v_0 (\theta) = \frac12 \Big( - 1 + \frac{1-\g}{1+\g} \cos (2 \theta) \Big)  = 
- \frac{1}{1+ \gamma} 
 \big( \gamma \cos^2 (\theta) +  \sin^2 (\theta) \big) \, .
\ee
\end{lemma}

\begin{pf}
By the  expression of  $ \Ke(0)(\theta, \theta') $ 
in \eqref{rel-dis0} we get 
\be
v_0 (\theta) 
=  
\int_{\T} \ln \Big(  \sin^2 \Big( \frac{\theta'- \theta}{2} \Big)  \Big) 
\cos  (\theta' - \theta) d \theta'  + \int_{\T} \ln  \Big(  
\frac{\g^2+1}{\g^2-1} -
  \cos (\theta + \theta' )  \Big) 
\cos  (\theta' - \theta) d \theta' \, . \label{ll2} 
\ee
The first integral in \eqref{ll2} is
\be \label{L02}
\int_{\T} \ln \Big(  \sin^2 \Big( \frac{\theta'- \theta}{2} \Big)  \Big) 
\cos  (\theta' - \theta) d \theta'  = 
\int_{\T} \ln \Big(  \sin^2 \Big( \frac{ \theta}{2} \Big)  \Big) 
\cos  ( \theta) d \theta \stackrel{\eqref{Iccs0}} = -  2 \pi \, .
\ee
For 
the second one 
in \eqref{ll2}  we use 
$ \cos ( \theta' - \theta) = \cos (\theta + \theta') \cos (2 \theta) + 
\sin (\theta ' + \theta) \sin (2 \theta) $ and \eqref{Iccs0}, 
obtaining
 \begin{align}
\int_{\T} \ln  \Big(  
\frac{\g^2+1}{\g^2-1} -
  \cos (\theta + \theta' )  \Big) 
\cos  (\theta' - \theta) d \theta'  
= \cos (2 \theta ) 2 \pi \frac{1- \g}{1+\g} \, .  \label{L03} 
\end{align}
By \eqref{ll2},  \eqref{L02}, \eqref{L03} we deduce   \eqref{cv0}.  
\end{pf}

We now compute the integral operator in \eqref{linunp}.

\begin{lemma}\label{intop0}
Let $ q(\theta )$ be a $ 2 \pi $-periodic real function 
with Fourier expansion 
\be\label{qFou}
 q (\theta) = {\mathop \sum}_{j \in \Z \setminus \{0\}} q_j e^{\ii j \theta}\, , \quad q_{-j} = \ov{q_j}\, .  
\ee
 Then the  integral operator $\W_0 [q]$ in \eqref{defK0y}
is 
\be\label{defK0}
\W_0\,  [q] (\theta)= - \sum_{j \in \Z \setminus \{0\}} \frac{q_j}{2|j|} \Big( e^{\ii j \theta} + \ka_j e^{- \ii j \theta} \Big)
=  - \sum_{j \in \Z \setminus \{0\}} \frac{1}{2|j|}   \Big( q_j  + \ka_j q_{-j}  \Big) e^{ \ii j \theta}  
\ee
where $ \ka_j := \ka_j (\gamma) :=  \Big( \frac{\gamma-1}{\gamma + 1} \Big)^{|j|} $. 
Moreover the operator $ \W_0$  is  self-adjoint and reversibility-preserving.  
\end{lemma}

\begin{pf}
 Formula    \eqref{defK0} follows by a computation using   
\eqref{rel-dis0}, \eqref{Iccs0}   and the fact that $ q (\theta )$ has 
zero average.
The integral operator $ \W_0 $ in \eqref{defK0y} is self-adjoint 
because the kernel $ \ln ( \Ke(0)(\theta, \theta')) $ is real and symmetric in 
$ (\theta, \theta') $. 
Moreover it 
is  reversibility preserving 
since, by \eqref{expR0}, it is even in $ (\theta, \theta') $, see  Lemma \ref{lem:IOK}. Lemma \ref{intop0} is proved.
\end{pf}

From \eqref{choiceOmega}, \eqref{fg0} and  \eqref{cv0} we get the identity
\be\label{id-g-v0}
\Omega_\gamma g_\gamma(\theta)+v_0 ( \theta )=-\Omega_\gamma \, .
\ee
Then combining \eqref{linunp} with \eqref{id-g-v0} and  Lemma \ref{intop0}  we obtain the 
  following lemma:
\begin{lemma}\label{lem:Linearized system at 0}{\bf (Linearized system at $ \xi = 0 $)}
The linearization of the 
Hamiltonian PDE \eqref{Evera}, with $\Omega=\Omega_\g = \tfrac{\gamma}{(1+\g)^2} $, at  
$ \xi = 0  $,  is  the linear Hamiltonian equation   
\begin{align}
\pa_t q  & =  - \Omega_\gamma \pa_\theta q  - 
\pa_\theta \W_0\, [q]    
\label{Lin1}
\end{align}
 generated by the quadratic Hamiltonian
\be\label{defQH}
H_L (q) := \tfrac{1}{2} ( {\bf \Omega}(\gamma) q, q )_{L^2}   \quad
\text{where} \quad
{\bf \Omega}(\gamma) := - \Omega_\gamma - \W_0   \, .   
\ee
\end{lemma}
We  now compute explicitly the  solutions of \eqref{Lin1}.
\\[1mm]
{\bf Unperturbed normal mode frequencies.}\label{sec:UL}
Fourier expanding $ q(\theta) $ as in \eqref{qFou}
the equation \eqref{Lin1} becomes, by \eqref{defK0}, 
$ \dot q_j =  
- \ii j \Omega_\g q_j +   \ii 
\tfrac{j}{2|j|}   \big( q_j  + \ka_j q_{-j}  \big) $,   
for an $  j \in \Z \setminus \{0\} $. 
Since the function $ q (\theta )$ 
is real valued, we have $ q_{-j} = \ov{q_j} $, for any 
$ j \in \Z \setminus \{0\} $.
Therefore the independent variables  
$ (q_n)_{n \geq 1 } $ satisfy 
\be\label{Lin1F}
\dot q_n =  -
\ii n \Omega_\g  q_n +   
\tfrac{\ii }{2}   \big( q_n  + \ka_n \ov{q_{n}}  \big)    \, , \quad \forall n \in \N \, . 
\ee
We conveniently describe \eqref{Lin1F}  in the real variables $ (\a_n, \b_n) \in \R^2 $, $ n \geq 1 $,  defined by
\be\label{qnanbn}
q_n =  \frac{\a_n -\ii \b_n}{2 \sqrt{\pi}} \, , \  \forall n \neq 2 \, , \quad 
q_2 = \frac{\a_2 -\ii \b_2}{\sqrt{2 \pi}} \, , \  n = 2 \, , 
\ee
obtaining, for any $ n \geq 1 $, 
\be\label{caseanbn}
\dot \a_n      = 
-\big( n \Omega_\g    +
   \tfrac{1}{2}   (  \ka_n  -   1  )  \big) \b_n    \, , \quad
\dot  \b_n      = 
\big(  n \Omega_\g   
-    \tfrac{1}{2}   ( 1   +     \ka_n  ) \big) \a_n  \, . 
\ee
The  normalization constants in \eqref{qnanbn} are such 
that  $ q (\theta )$ is expressed, in  
the coordinates 
$ (\a_n, \b_n )_{n \geq 1}  $, as  
\be\label{qthanbn}
 q(\theta)  
 =  
  {\mathop \sum}_{n \in \N} \a_n \tc_n ( \theta) + 
 \b_n \ts_n ( \theta) \, ,  
\ee
where 
\be\label{defcnsn}
\tc_n (\theta) := 
\begin{cases}  \tfrac{\cos (2 \theta)}{\sqrt{\pi}}\sqrt{2}  \quad \text{if} \ n = 2 \\ 
\tfrac{\cos (n \theta)}{\sqrt{\pi}}  \, \qquad  \text{if} \ n \neq 2 \, , 
\end{cases} \quad
\ts_n (\theta) := 
\begin{cases}  \tfrac{\sin (2 \theta)}{\sqrt{\pi}}\sqrt{2}  \quad \text{if} \ n = 2 \\ 
\tfrac{\sin (n \theta)}{\sqrt{\pi}}  \, \qquad  \text{if} \ n \neq 2 \, , 
\end{cases}
\ee
and, denoting by $ ( \ , \ ) = (\ , \ )_{L^2}$ the $ L^2 (\T) $-scalar product, it results   
  \begin{equation}\label{newa2b2}
  \a_2 = \tfrac{1}{2} ( q, \tc_2)  \, , 
  \ 
    \beta_2 = \tfrac{1}{2} (q,  \ts_2) \, , \quad 
   \a_n = (q, \tc_n)   \, , \ 
  \  
    \beta_n =  (q,  \ts_n)   \, ,  \quad \forall n \neq 2 \, .
\end{equation}
In the variables $ (\a,\beta) :=
 (\a_n, \b_n )_{n \geq 1} $ the involution $ {\cal S} $ defined in \eqref{S-invo} reads
\be\label{Sanbn}
(\a_n, \b_n) \mapsto (\a_n, - \b_n ) \, , \quad \forall n \geq 1 \,  ,
\ee
the symplectic form $ {\cal W} $ in \eqref{sy2form}  reads
\be\label{sy2form0}
 d \alpha_2 \wedge d \beta_2 +  {\mathop \sum}_{n \neq 2 } 
\tfrac{1}{n}d \a_n\wedge  d\b_n \, , 
\ee
and the Hamiltonian vector field generated by a Hamiltonian $ H(\a,\b)$ is  
\begin{equation}\label{XVS}
[X_H]_2=  \begin{psmallmatrix}
\partial_{\beta_2} H  \\
-\partial_{\alpha_2} H
\end{psmallmatrix} \, ,  \quad
[X_H]_n= n \begin{psmallmatrix}
\partial_{\beta_n} H  \\
-\partial_{\alpha_n} H
\end{psmallmatrix}\, , \quad \forall n \neq 2 \, .
\ee
The  linear Hamiltonian system \eqref{caseanbn} writes as in \eqref{anbn}
which is a harmonic oscillators with frequency
\be\label{defomegan}
\Omega_n  (\gamma)  =  \sqrt{\mu_n^+ \mu_n^-} =
\Big[ \Big(\frac{n \gamma}{(1+\gamma)^2}  - \frac{1}{2}\Big)^2   -
   \frac{\ka_n^2}{4} \Big]^{1/2}       
\ee
if and only if  $ \mu_n^+ \mu_n^- > 0 $, namely $ \mu_n^+ $, $ \mu_n^-  $ have the same sign. According if $ \mu_n^+ , \mu_n^-  $ are both positive 
or both negative they rotate in opposite direction. 
For $ n = 1 $ it results 
\be\label{freq1}
\mu_1^- = - \frac{\g^2}{(1+\g)^2} < 0  \, , 
\quad \mu_1^+ = - \frac{1}{(1+\g)^2}  < 0 \, , \quad 
\Omega_1^2  = \frac{ \gamma^2}{(\gamma + 1)^4 } \, . 
\ee
The fact that $ \Omega_1 $ is equal to the Kirchhoff ellipse rotating frequency $ \Omega_\gamma $ in \eqref{choiceOmega}  is due to the  presence of the prime integral $ Z $ (the center of mass), see Remark \ref{rem:1fre}. 
For $ n = 2 $  it turns out that 
\be\label{umod2}
\mu_2^-  = - \Big( \frac{\g-1}{\g+1} \Big)^2 < 0 \, , \quad 
\mu_2^+ (\gamma) = 0  \, , \ \ \forall \gamma \, ,   
\ee
and so \eqref{anbn} reduces to \eqref{anbn2}. 
The degeneracy of \eqref{anbn2} 
is due to the presence of the prime integral 
$ J $ (the angular momentum) that has a linear $ \a_2 $-component, see Remark  \ref{rem:Ham-dege}. 
For $ n = 3  $  we have   
\be\label{mu3-mu3+}
\mu_3^{-} = \frac{\g^2 (3 - \g)}{(1+\g)^3} \, , \quad 
\mu_3^{+} = \frac{3\g - 1}{(1+\g)^3} \, , 
\ee
and thus, according to the values of $ \gamma $, the mode $ n = 3 $ may  
be stable or unstable. In particular, for
$ 1 < \gamma < 3 $ it results $  0 <  \mu_3^{-} < \mu_3^{+ } $ and 
 thus the mode $ n =3 $ is linearly stable.  
Actually in this range of values all the modes $ n \geq 3 $ are harmonic oscillators, 
since $ n \mapsto  \mu_n^- (\gamma )$ is increasing and thus 
$  0 <   \mu_3^{-} < \mu_n^{-} < \mu_n^{+ } $ for any $ n \geq 4 $ and 
$ \gamma \in (1,3)$.
Note that, for any $\gamma > 1 $, we have  $ \mu_n^{-} \to + \infty $ 
 for $ n \to + \infty $, and therefore  the systems 
\eqref{anbn} are harmonic oscillators for $ n \geq n_0 (\gamma ) $ large enough.  
For 
$ \gamma > 3 $  finitely many hyperbolic directions appear.  
We have the following lemma.

\begin{lemma}\label{lem:eige}
{\bf (Critical aspect ratios)} 
For any $ n \geq 1 $
we have 
$  \mu_n^{-} (\gamma ) < \mu_{n+1}^{-} (\gamma ) $, $ \forall   \gamma > 1 $,  
each function $ \gamma \mapsto \mu_n^{-} (\gamma ) $ is monotone decreasing
on $ (1,+\infty) $ and 
\begin{equation}\label{limits}
 \lim_{\gamma \to 1^+ } \mu_n^{-} (\gamma ) = \frac{n}{4} - \frac12  \, , \quad
  \lim_{\gamma \to + \infty} \mu_n^{-} (\gamma ) = - 1 \, . 
\end{equation}
For any $ n \geq 3 $ there exists a unique $ \underline{\gamma}_n \geq 3 $  such that
\begin{equation}\label{monotonia}
\mu_n^{-} (\gamma ) < 0 \, ,  \ \forall \gamma > \underline{\gamma}_n \, , \quad
 \mu_n^{-} (\underline{\gamma}_n) = 0 \, , \quad 
\mu_n^{-} (\gamma ) > 0 \, , \  \forall \gamma < \underline{\gamma}_n \, .
\end{equation}
It results that $ 3 = \underline{\gamma}_3 < 	\ldots < \underline{\gamma}_n  < \underline{\gamma}_{n+1}  < \ldots $. 

Moreover  $ \mu_n^+ (\gamma) > \mu_n^- (\gamma) $ and 
$ \mu_{n+1}^+ (\gamma) > \mu_n^+ (\gamma)  $ for any $ \gamma > 1 $.
For all $ n \geq 3 $ 
the function $ \gamma \mapsto \mu_n^{+} (\gamma ) $ is positive and monotone decreasing
on $ (1,+\infty) $ and 
\begin{equation}\label{limits+}
 \lim_{\gamma \to 1^+ } \mu_n^{+} (\gamma ) = \frac{n}{4} - \frac12  \, , \quad
  \lim_{\gamma \to + \infty} \mu_n^{+} (\gamma ) = 0 \, . 
\end{equation}
\end{lemma}

\begin{pf}
By \eqref{anbn} we immediately conclude  that
$  \mu_n^{-} (\gamma)$ is increasing  in $ n $ and \eqref{limits} holds. 
Moreover, the derivative of  $  \mu_n^{-} (\gamma)$ with respect to $\gamma$ satisfies
\be\label{dermun-}
\frac{d \mu_n^{-} (\gamma)  }{d \gamma} = - n \frac{\gamma^2-1}{(1+\gamma)^4}
-  n \Big( \frac{\g-1}{\g+1}\Big)^{n-1}  \frac{1}{(1+\g)^2} < 0 
\ee
for any $ \gamma > 1 $. Thus, the  function $ \gamma \mapsto \mu_n^{-} (\gamma ) $ is monotone decreasing
on $ (1,+\infty) $ and \eqref{monotonia} follows. 
By \eqref{mu3-mu3+} it results that $ \underline{\gamma}_3  = 3 $. 
 Since $  \mu_n^{-} (\gamma ) < \mu_{n+1}^{-} (\gamma )  $, for any $ \gamma > 1 $, 
 we deduce that the sequence of $ \underline{\gamma}_n $ where 
 $ \mu_n^{-} (\underline{\gamma}_n) = 0 $ is monotone increasing in $n$. 
  
In view of  \eqref{anbn}  we trivially  have that  $ \mu_n^+ (\gamma) > \mu_n^- (\gamma) $ and that 
$ \mu_{n+1}^+ (\gamma) > \mu_n^+ (\gamma)  $ if and only if
$  (\gamma + 1 )^{n-1} \gamma > (\g-1)^{n-1} (\g-1)$, namely for any $ \gamma > 1 $.  
Moreover, for any $ n \geq 3 $,  (cfr. \eqref{dermun-})
$$
\begin{aligned}
\frac{d \mu_n^{+} (\gamma)  }{d \gamma} 
& = 
 \frac{n(\g-1)}{(1+\g)^4} \Big( - (\g+1) + 
 \frac{(\g-1)^{n-2}}{(\g+1)^{n-3}} \Big) < 0 
\end{aligned}
$$
because $ (\g-1)^{n-2} <  (\g+1)^{n-2}$, for any $ \g > 1 $.  This ends the proof of Lemma \ref{lem:eige}.
\end{pf}

For any $ \bar n \geq 2 $ we fix a compact interval $ \IK $ 
of values of $ \gamma $ as
in \eqref{defgammaset0}  
so that 
$ \mu_n^- \mu_n^+ > 0 $, for  $ n =1 $ and any $ n = \bar n + 1, \bar n + 2, \ldots $, whereas
$ \mu_n^- \mu_n^+ < 0 $ for $ n = 3, \ldots, \bar n $.  
The finitely many directions $ n =3, \ldots, \bar n $ are hyperbolic  and will not produce resonance phenomena.

We  perform the  symplectic and reversibility preserving 
change of variable
\be\label{def:Lan}
\begin{pmatrix}
\a_n \\
 \b_n
\end{pmatrix} = 
\begin{pmatrix}
\symm_n & 0 \\
0 & \symm_n^{-1}
\end{pmatrix}
\begin{pmatrix}
\breve \a_n \\
\breve  \b_n
\end{pmatrix}  \, , \  
\symm_n := \Big( \frac{|\mu_n^+|}{|\mu_n^-|} \Big)^{\frac14}\, , 
\ {\rm if} \ n \neq 2  \, ,  \
\symm_2 := 1 \, ,    
\ee
in order to  reduce the linear Hamiltonian systems \eqref{anbn}  to the more symmetric form  
\be\label{anbn-ri0} 
\begin{aligned}
& \begin{pmatrix}
\dot {\breve \a}_1 \\
\dot {\breve\b}_1
\end{pmatrix} =
\begin{pmatrix}
0 &  \Om_1  \\
 -\Om_1 & 0 
\end{pmatrix} 
\begin{pmatrix}
\breve \a_1 \\
\breve \b_1
\end{pmatrix}  \, , \ \Omega_1 = \Omega_\gamma = \frac{\g}{(1+\g)^2} \\
& 
\begin{pmatrix}
\dot {\breve \a}_2 \\
\dot {\breve \b}_2
\end{pmatrix} =
\begin{pmatrix}
0 & 0  \\
\Omega_2 & 0 
\end{pmatrix} 
\begin{pmatrix}
\breve \a_2 \\
\breve \b_2
\end{pmatrix} \, , \ \Omega_2 :=  \mu_2^{-} = -\Big( \frac{\g-1}{\g+1} \Big)^2  \\
& 
\begin{pmatrix}
\dot {\breve \a}_n \\
\dot {\breve \b}_n
\end{pmatrix} =
\begin{pmatrix}
0 & - \Om_n  \\
- \Om_n & 0 
\end{pmatrix} 
\begin{pmatrix}
\breve \a_n \\
\breve \b_n
\end{pmatrix} \, , \ 3 \leq n \leq \bar n  \, ,  \\ 
& 
\begin{pmatrix}
\dot {\breve \a}_n \\
\dot {\breve \b}_n
\end{pmatrix} =
\begin{pmatrix}
0 & -\Om_n  \\
\Om_n & 0 
\end{pmatrix} 
\begin{pmatrix}
\breve \a_n \\
\breve \b_n
\end{pmatrix} \, , \ n \geq \bar n + 1  \, ,  \\ 
& \Om_n = \Om_n (\gamma)=  |\mu_n^+ \mu_n^- |^{\frac12} = 
\Big| \Big(\frac{n \gamma}{(1+\gamma)^2}  - \frac{1}{2}\Big)^2   -
   \frac{\ka_n^2}{4}\Big|^{\frac12}   \, . 
\end{aligned} 
\ee
In a more geometrical point of view,  the phase space decomposes into the direct sum
\be\label{Hdir}
 H^s_0 (\T) :=  \oplus_{n \geq 1} V_n  
\ee
of $ 2 $-dimensional $ L^2$-orthogonal  real Lagrangian subspaces 
\be\label{def:Vn}
V_n := \Big\{ q(\theta ) =  \symm_n \breve \a_n \tc_n(\theta) + \symm_n^{-1} 
\breve \b_n \ts_n(\theta) \, , \ (\breve \a_n, \breve \b_n ) \in \R^2   \Big\} \, ,
\ee
where the multipliers $\symm_n $ are defined in \eqref{def:Lan}. 
These subspaces are also symplectic orthogonal.
In the symplectic coordinates $ (\breve\a_n, \breve \b_n) $ in \eqref{def:Vn}, 
the restriction of the symplectic form $ {\cal W} $ in \eqref{sy2form} 
to each subspace $ V_n  $
is $ \tfrac{1 }{n} d \breve \a_n \wedge d \breve \b_n $, for $ n \neq 2 $, and 
$ d \breve\a_2 \wedge d \breve \b_2 $, for $ n = 2 $, see \eqref{sy2form0}.
The Hamiltonian associated 
to the linear system  \eqref{anbn-ri0} 
is the quadratic Hamiltonian  \eqref{Hanbn}. 
We also conclude that 
  the linearized  equation \eqref{Lin1} possesses 
the reversible 
oscillating in time solutions \eqref{q-r-l}. 
\smallskip

We end this section by  few comments about  the first and 
the second frequencies.
 \begin{remark}\label{rem:1fre}
{\bf (Explanation that the first frequency $ \Omega_1 =  \Omega_\gamma $).} 
 The first order expansion in $ \xi $ of the center of mass $ Z $ in \eqref{forCZJ} is,   
expanding $ \xi (\theta )$ as in  \eqref{qthanbn}-\eqref{defcnsn},  
\begin{align}
Z &  =
 e^{ \ii \Omega t} \int_{\T}
\xi(\theta)   \big(\gamma^{\frac12}\cos\theta + 
\ii \gamma^{-\frac12}\sin\theta\big)d\theta  + O(\| \xi \|^2) \nonumber
\\
& = 	\sqrt{\pi} e^{ \ii \Omega t} 
\Big( \sqrt{\gamma} \a_1  +
\ii  \tfrac{\b_1}{\sqrt{\gamma} } \big)   + O(\| \xi \|^2) 
 = \sqrt{\pi}  e^{ \ii \Omega t} 
\big( \breve \a_1  + \ii  \breve \b_1 \big)   + O(\| \xi \|^2) \label{Zab}
\end{align}
where
$ \breve \a_1 := \sqrt{\g} \a_1 $, $  \breve \b_1 := \b_1/  \sqrt{\g} $
are the variables defined in \eqref{def:Lan}, for $ n = 1 $. 
Since $ Z $ is a prime integral, we deduce by \eqref{Zab} that the mode-$ 1$ variables
$ (\breve \a_1, \breve \b_1 ) $ perform, at the linear level, a rotation with angular velocity
$  \Omega $. This is coherent with the 
linear system 
\eqref{anbn-ri0} for the first mode.  
\end{remark}

\begin{remark}\label{rem:Ham-dege}
{\bf (Degeneracy of the second frequency)}.
Recalling \eqref{moangJ2}, \eqref{fg0}, 
 \eqref{defcnsn} and \eqref{newa2b2}, we have that 
$$
\frac{\sqrt{2} \, J_1  (\xi)}{\sqrt{\pi}(\g- \g^{-1})}  	=  \frac{1}{2} 
\int_{\T} \xi (\theta) \tc_2 (\theta ) \, d \theta = 
\a_2  
$$
where  $ J_1 (\xi) $  is linear component in $ \xi  $ of 
the angular momentum $ J $. 
Note that the linear system in \eqref{anbn-ri0} is degenerate exactly in the 
$ \a_2 $ component. This is a general fact about the structure of a prime integral
close to an equilibrium of a dynamical system
$ \dot x = f(x) $ with $ f(0) = 0 $. Let $  A := Df (0) $. 
If $ b(x) $ is a prime integral we have
$ \nabla b(x) \cdot f(x) = 0 $, $  \forall x $. 
Hence, differentiating and using that $ f(0) = 0 $, we get 
$ \nabla b(0) \cdot A y  = 0 $, $ \forall y $. 
If $ A $ is non singular we deduce that $ \nabla b(0) = 0 $, i.e. 
$ b$ is quadratic at  $ x = 0 $. Here $ A  $ is degenerate in the $ \a_2 $ variable and indeed $ J $ has a linear term only in $ \a_2 $. 
\end{remark}

\section{Transversality properties of the linear frequencies}\label{sec:ND}

We recall  
that 
the linear frequencies $ \Omega_n (\gamma)  $, for $ n = 1 $ and 
$ n \geq  \bar n +1 $,  are
(see \eqref{anbn-ri0}, \eqref{defomegan})
\begin{equation}\label{Omg}
\Omega_1 (\gamma) = \Omega_\g =\tfrac{\gamma}{(1+\gamma)^2} \, , 
\quad
\Omega_n (\gamma) = 
\sqrt{ \big(\tfrac{n \gamma}{(1+\gamma)^2}  - \tfrac{1}{2}\big)^2 - \tfrac14 \big( \tfrac{\gamma-1}{\gamma + 1} \big)^{2n}} \, , \ \forall n \geq \bar n + 1 \, . 
\end{equation}
For  any $\gamma $ varying in  
$ \IK := [ \gamma_1, \gamma_2] $ 
defined in \eqref{defgammaset0}, 
 all the frequencies $ \Om_n (\gamma)$, $ n \geq \bar n + 1  $, in \eqref{Omg} 
are real. 
Note that the maps $ \gamma \mapsto \Omega_n (\gamma) $ are analytic, for any $ n \geq 1 $.

\smallskip

 Following the degenerate KAM theory approach as developed  
in \cite{BaBM} and \cite{BertiMontalto,BBHM},   
the key argument is to prove a non-degeneracy condition for these frequencies. 

\begin{definition}\label{def:non-deg} {\bf (Non-degeneracy)}
A function $ f := (f_1, \ldots, f_N ) :  [ \gamma_1, \gamma_2] \to \R^N $ is called non-degenerate 
if, for any vector $ c := (c_1, \ldots, c_N ) \in \R^N \setminus \{0\}$, 
the function $ f \cdot c = f_1 c_1 + \ldots + f_N c_N $ 
is not identically zero on the whole interval $  [ \gamma_1, \gamma_2]  $. 
\end{definition}

From a geometric point of view, $ f $ is non-degenerate means that the image of the curve $ f( [ \gamma_1, \gamma_2]) \subset \R^N $ is not contained in any  hyperplane of $ \R^N $. 

\begin{lemma}\label{non degenerazione frequenze imperturbate}
{\bf (Non-degeneracy)}  For any $ N \in \N $, 
$ 1 \leq n_1 < n_2 < \ldots < n_N $, $ n_i \notin \{ 2, \ldots, \bar n   \}, \forall i= 1, \ldots, N $, 
the functions 
\begin{align}\label{freND1}
& [ \gamma_1, \gamma_2] \ni \gamma \mapsto (\Omega_{n_1}(\gamma ), \ldots, \Omega_{n_N}( \gamma )) \in \R^N \\
& \label{freND2}
[ \gamma_1, \gamma_2] \ni \gamma \mapsto (1, \Omega_{n_1}(\gamma ), \ldots, \Omega_{n_N}( \gamma )) \in \R^{N+1} 
\end{align}
are non-degenerate according to Definition \ref{def:non-deg}.
\end{lemma}

\begin{pf}
We prove the non-degeneracy of the functions \eqref{freND1}-\eqref{freND2}
after performing a change of variable. 
The function $  (1, +\infty) \ni \gamma \mapsto \tfrac{\gamma-1}{\gamma+1}  $ 
 is monotone increasing and takes values in the interval $ (0,1) $.  
Then we consider its inverse function 
$ \tfrac{\gamma-1}{\gamma+1} =: y  $ $  \Leftrightarrow $
$ \gamma = \tfrac{1+y}{1-y} $, $ y \in (0,1) $.
Setting  
$ y^2 =: z $, $ z \in (0,1) $, 
we express the frequencies $ \Omega_n $ as functions of $ z $, obtaining by \eqref{Omg}, 
(that for simplicity we denote with the same letters) 
\be\label{Om1n}
\Omega_1 (z) = \tfrac14 (1- z) \, , 
\quad
\Omega_n (z) = \tfrac14
\sqrt{ (n - 2 - nz )^2 - 
4 z^{n}} \, , \quad \forall n \geq  \bar n + 1  \, . 
\ee
Note that each of these  functions is well defined and analytic in a full neighborhood of $ z = 0 $. 
Then, for $ n \geq \bar n + 1 $, we Taylor 
expand in a neighborhood of $ z = 0 $ the functions $ 4 \Omega_n(z) $ obtaining 
\begin{align}
4 \Omega_n (z) 
 =  (n - 2 - nz ) \Big( 1 - \frac12  \frac{4z^{n}}{(n-2- nz)^2} + O(z^{2n}) \Big)  
   =  (n - 2) - nz   - \frac{2z^{n}  }{n-2}   + O(z^{n+1}) \label{Taylor-fre} \, . 
\end{align}
The Taylor expansion \eqref{Taylor-fre} of $ 4 \Omega_n(z) $
at $ z = 0 $, up to order $ n $, proves that  the $ n $-th derivative 
\be\label{NDERN}
D^{(n)}_z \Omega_n (0) = - \frac{ n!}{2(n-2)} \, .
\ee 
We are now able to prove the non-degeneracy of the functions in \eqref{freND1}. 
We first prove that, if
\be \label{caso1}
\bar n + 1  \leq n_1 <  n_2 < \ldots < n_N \, ,
\ee
for any  vector $ c = (c_1, \ldots, c_N) \in \R^N \setminus \{0\} $, 
the analytic function 
$ z \mapsto c_1 \Omega_{n_1}( z) + \ldots + c_N \Omega_{n_N}(z)  $
is not identically zero on an interval $ (-\delta, \delta) $ for some $ \delta > 0 $, 
thus on the whole $ (0,1) $ (below 
we shall also consider the case  $ n_1 = 1 $). 
Suppose, by contradiction, that there exists $ c \in \R^N \setminus \{0\} $ such that 
\be\label{la-identita}
c_1 \Omega_{n_1}( z) + \ldots + c_N \Omega_{n_N}(z) = 0 \, , 
\quad \forall  |z| < \delta \,  . 
\ee
Differentiating with respect to $ z $ the identity in \eqref{la-identita}, we find
 $$
 \begin{cases}
c_1 D_z^{(n_N)} \Omega_{n_1}( z) + \ldots +
 c_N D_z^{(n_N)} \Omega_{n_N}( z)  = 0 \cr 
 \ldots  \ldots \ldots \cr
 c_1 D_z^{(n_1)} \Omega_{n_1}( z) + \ldots +
 c_N D_z^{(n_1)} \Omega_{n_N}( z)  = 0 \, .
 \end{cases}
$$
As a consequence the $ N \times N $-matrix 
$  {\cal A}(z) := \begin{psmallmatrix}
D_z^{(n_N)} \Omega_{n_1}( z)  & 
 \dots & D_z^{(n_N)} \Omega_{n_N}( z)
\\
\vdots & \ddots & \vdots \\
D_z^{(n_1)} \Omega_{n_1}( z)  & 
\dots & D_z^{(n_1)} \Omega_{n_N}( z)
\end{psmallmatrix}
$
is singular for all $ z \in \R $, $ |z| < 	\d $,  and 
in particular at $ z  = 0 $ we have  $ \det {\cal A}(0) = 0 $. 
On the other hand, by \eqref{Taylor-fre}, \eqref{NDERN}, and since $ 
\bar n + 1  \leq n_1 < n_2 < \ldots < n_N $, 
we have that, for some real constants $ c_{i,j}  $, 
$$ 
{\cal A}(0)  =
 \begin{psmallmatrix}
c_{1,1} & c_{2,1}  &  \ldots & \dots & - \frac{ n_{N}!}{2(n_{N}-2)} \\
c_{1,2} & c_{2,2}  &  \ldots &  
- \frac{n_{N-1}!}{2(n_{N-1}-2)}  & 0 \\
& \ddots & \vdots   \\
c_{1,N-1} & - \frac{n_2!}{2(n_2-2)}  &  0 & \dots & 0 \\
- \frac{ n_1!}{2(n_1-2)}  &  0 & 0 & \dots & 0
\end{psmallmatrix} 
$$
is triangular.  
Thus the determinant 
$ \det {\cal A}(0) = \prod_{i=1,\ldots, N}  (-1)^{i+1}  \tfrac{ n_{i}!}{2(n_{i}-2)}  \neq 0  $.  
This contradiction proves the non-degeneracy of the function defined in 
\eqref{freND1} in the case \eqref{caso1} holds. 
On the other hand, if $  n_1 = 1 $ and 
$\bar n + 1 \leq  n_2 < \ldots < n_N $ we can write \eqref{la-identita} as 
$ c_2 \Omega_{n_2}( z) + \ldots + c_N \Omega_{n_N}(z) + c_1 \Omega_{1}( z)   = 0 $, 
for all $  |z|  < \d $, 
obtaining, by differentiation,  arguing as before that 
the $ N \times N $-matrix 
$$
{\cal A}(z) := \begin{pmatrix}
\Omega_{n_2}( z)  & 
 \dots & \Omega_{n_N}( z) & \Omega_{1}( z)
\\
D_z^{(n_{N})} \Omega_{n_2}( z) &  
\dots & D_z^{(n_{N})} \Omega_{n_N}( z) &  D_z^{(n_{N})} \Omega_{1}( z)  \\
\vdots & \ddots & \vdots \\
D_z^{(n_2)} \Omega_{n_2}( z)  & 
\dots & D_z^{(n_2)} \Omega_{n_N}( z) &  D_z^{(n_{2})} \Omega_{1}( z) 
\end{pmatrix}
$$
is singular for all $ z \in \R $, $ |z| < 	\d $, and   $ \det {\cal A}(0) = 0 $. 
On the other hand, by \eqref{Om1n} and 
\eqref{NDERN}, and since $ \bar n + 1 \leq n_2 < \ldots < n_N $, 
we have that, for some real constants $ c_{i,j} $,   
$$
{\cal A}(0)  =
 \begin{psmallmatrix}
c_{1,1} & c_{2,1}  &  \ldots & \dots & \frac14 \\
c_{1,2} & c_{2,2}  &  \ldots &  
- \frac{n_{N}!}{2(n_{N}-2)}  & 0 \\
& \ddots & \vdots   \\
c_{1,N-1} & - \frac{n_3!}{2(n_3-2)}  &  0 & \dots & 0 \\
- \frac{ n_2!}{2(n_2-2)}  &  0 & 0 & \dots & 0
\end{psmallmatrix} 
$$
and thus $ {\rm det}{\cal A}(0) \neq 0 $. 
This contradiction proves
the non-degeneracy of the function defined in 
\eqref{freND1}.
The  non-degeneracy of the function in \eqref{freND2} follows similarly. 
\end{pf}

We shall use the following asymptotic expansions of the frequencies. 
\begin{lemma}\label{lem:asy}
{\bf (Asymptotics)} 
There exists $ c > 0 $ such that  for  $ n =1 $ and any $ n \geq \bar n + 1$ one has, for
 any $ \gamma \in \IK $ defined in  \eqref{defgammaset0},   
\be\label{ASYFR}
\Omega_n (\gamma) \geq c  n \, .
\ee
Moreover  for  any $ n \geq \bar n + 1$, 
\be\label{ASYFR1}
\Omega_n (\gamma) = n \Omega_1 (\gamma) - \frac12 +   r(n,\gamma)
 \, , 
\quad \sup_{n \geq \bar n+1, \gamma \in [\g_1,\g_2]} 
n |\pa_\gamma^k r(n,\gamma)| \leq C_k  \, , \ \forall k \in \N_0  \,  . 
\ee
\end{lemma}

\begin{pf}
For $ n \geq \bar n+1 $ we write, by  \eqref{Omg},  
\be\label{Omnex}
\Omega_n (\gamma) = 
 \Big| n \Omega_1 (\gamma)  - \tfrac{1}{2}\Big| 
\sqrt{ 1 - \tfrac{\ka_n^2}{4}
\big(n  \Omega_1 (\gamma) - \tfrac{1}{2}\big)^{-2}  } \, , 
\quad \ka_n = \big( \tfrac{\gamma-1}{\gamma + 1} \big)^{n}.
\ee
By Lemma \ref{lem:eige} and the compactness of $  \IK $ 
defined in  \eqref{defgammaset0},  there exists  $ \underline{c}  > 0 $  
such that $ \mu_n^{-} (\gamma) > 0 $ for any $ \gamma \in \IK $ and
 $ n \geq \bar n + 1 $. Thus, recalling 
 \eqref{anbn} and that  $ \Omega_1 (\gamma) = \Omega_\gamma  $, 
\be\label{nOmebarn+1}
 n \Omega_\gamma  -      \tfrac{1}{2}   >  \tfrac{ \ka_n}{2} + \underline{c} 
 \geq \underline{c}  \, , \quad 
\forall \gamma \in \IK \, , \  n \geq \bar n + 1  \, . 
\ee
Then, for any $ n \geq \bar n+1 $, we write \eqref{Omnex} as 
$$
\Omega_n (\gamma) = 
n \Omega_1 (\gamma)  - \tfrac{1}{2}  +
\big( n \Omega_1 (\gamma)  - \tfrac{1}{2} \big) 
\big[ \sqrt{ 1 - \tfrac{\ka_n^2}{4}
\big(n  \Omega_1 (\gamma) - \tfrac{1}{2}\big)^{-2}} -1 \big]\, , 
$$
which has the form  \eqref{ASYFR1} with
\be\label{rng}
 r(n,\gamma)  := 
   -
\frac{1}{4 ( n  \Omega_1 (\gamma) - \frac{1}{2})} 
\frac{ \ka_n^2 }{\sqrt{ 1 - \frac{\ka_n^2}{4}
\big(n  \Omega_1 (\gamma) - \frac{1}{2}\big)^{-2}} + 1} \, , \quad  \ka_n = \Big( \frac{\gamma-1}{\gamma + 1} \Big)^{n} \, . 
\ee
Using \eqref{nOmebarn+1} 
we get the  bounds for  $ r(n,\gamma) $ claimed in \eqref{ASYFR1}. Moreover 
the expansion \eqref{ASYFR1} and since each 
$ \Omega_n (\gamma) $, $ n 	\geq \bar n+1  $, is  positive, 
implies  \eqref{ASYFR}.
\end{pf}

\begin{remark}
Actually the $ r(n,\gamma) $ decay much faster than stated in  \eqref{ASYFR1}
because of the exponentially small terms $ \ka_n $ in \eqref{rng}. 
\end{remark}

In the next proposition 
we deduce the quantitative bounds \eqref{0 Melnikov}-\eqref{2 Melnikov+} 
from the qualitative non-degeneracy condition of Lemma 
\ref{non degenerazione frequenze imperturbate}, 
the analyticity of the linear frequencies $\Omega_n (\gamma) $, 
and their asymptotics.
For any $\bar n\geq 2$ we consider 
finitely many  `tangential" sites 
\be\label{Tsites}
{\mathbb S} := \{ n_1, \ldots, n_{|\ST|} \} \, , \quad 
\bar n + 1  \leq n_1 < n_2 < \ldots < n_{|\ST|} \, , 
\ee
and we denote the unperturbed tangential frequency vector by 
\be\label{tangential-normal-frequencies}
{\vec \om} ( \gamma ) := ( \Omega_n (\gamma ) )_{n \in \ST} \, , 
\ee
where $ \Omega_n (\gamma )  $ are the frequencies in \eqref{Omg}.

\begin{proposition}\label{Lemma: degenerate KAM}
{\bf (Transversality)}
There exist $ \barka \in \N $, $ \rho_0 > 0$ such that, 
for any $ \gamma \in [\g_1,\g_2] $,
\begin{align}\label{0 Melnikov}
& \max_{k \leq \barka} 
|\partial_\g^{k}  \{{\vec \om} (\gamma) \cdot \ell   \} |  
\geq \rho_0 \langle \ell \rangle\,, 
\quad \forall \ell  \in \Z^{|\ST|} \setminus \{ 0 \},   
\\
\label{T Melnikov}
& \max_{k \leq \barka}
|\partial_{\gamma}^{k}  \{{\vec \om} (\gamma) \cdot \ell  +  \Omega_1 (\gamma ) j \} | 
 \geq \rho_0 \langle \ell  \rangle\,, 
\quad \forall (\ell,j)  \in (\Z^{|\ST|} \times \Z) \setminus \{(0,0)\} \, ,  
\\
\label{1 Melnikov}
& \max_{k \leq \barka}
|\partial_{\gamma}^{k}  \{{\vec \om} (\gamma) \cdot \ell  +  \Omega_n (\gamma ) \} | 
 \geq \rho_0 \langle \ell  \rangle\,, 
\quad \forall \ell  \in \Z^{|\ST|}, \, n \in \N \setminus (\ST \cup \{ 2, \ldots, \bar n \}) \, ,  
\\
\label{2 Melnikov-}
& \max_{k \leq \barka}
|\partial_\gamma^{k}  \{{\vec \om} (\gamma) \cdot \ell  
 + \Omega_n (\gamma ) - \Omega_{n'}(\gamma ) \} | 
  \geq \rho_0 \langle \ell  \rangle\,, \\
&  
\forall \ell \in \Z^{|\ST|} , \,  n, n' \in \N \setminus (\ST \cup \{ 2, \ldots, \bar n \}) \, ,  
\quad (\ell,n,n') \neq (0,n,n) \, , 
\nonumber \\ 
& \max_{k \leq \barka}
 |\partial_\gamma^{k}  \{{\vec \om} (\gamma) \cdot \ell  
 + \Omega_n (\gamma ) + \Omega_{n'}(\gamma ) \} | 
 \geq \rho_0 \langle \ell  \rangle\,, \
 \forall \ell \in \Z^{|\ST|}, \,  n, n' \in \N \setminus (\ST \cup \{ 2, \ldots, \bar n \}) \, ,   \label{2 Melnikov+}
\end{align}
where $\vec\om(\gamma)$ and $\Omega_n(\gamma)$ are defined in \eqref{tangential-normal-frequencies} and \eqref{Omg}.
We call 
$ \rho_0 $  the ``\emph{amount of non-degeneracy}'' 
and $ \barka $ the ``\emph{index of non-degeneracy}''.
\end{proposition}

\begin{pf}
	We prove separately  \eqref{0 Melnikov}-\eqref{2 Melnikov+}. 
	\\
	{\bf Proof of \eqref{0 Melnikov}}. By contradiction, suppose that 
for all $ k_0^* \in \N $, $ \rho_0 > 0 $, there exists 
$\ell \in\Z^{|\ST|}\setminus\{0\}$, $ \gamma \in \IK := [\gamma_1, \gamma_2] $, such that 
$ \max_{k \leq \barka} 
|\partial_\g^{k}  \{{\vec \om} (\gamma) \cdot \ell   \} |  
< \rho_0 \langle \ell \rangle $. This implies that, taking $ k_0^* = m $,
$ \rho_0 = \braket{m}^{-1} $ 
	there exist $\ell_m\in\Z^{|\ST|}\setminus\{0\}$ and $\gamma_m\in \IK $   such that
$
	\Big| \partial_\gamma^k \vec{\omega}(\gamma_m) \cdot \ell_m  \Big| < $$ \frac{\braket{\ell_m}}{\braket{m}} $, $ \forall\,0\leq k \leq m $, 
	and therefore 
	\begin{equation}\label{eq:0_abs_m}
\forall k \in \N_0 \,  , \quad \forall m \geq k \, ,  \quad
	\Big| \partial_\gamma^k \vec{\omega}(\gamma_m) \cdot \frac{\ell_m}{\braket{\ell_m}}  \Big| < \frac{1}{\braket{m}}  \, .
	\end{equation}	
	The sequences $(\gamma_m)_{m\in\N}\subset\IK$ and $(\ell_m/\braket{\ell_m})_{m\in\N}\subset \R^{|\ST|}\setminus\{0\}$ are both bounded. By compactness, up to subsequences
	 $\gamma_m\to \bar\gamma \in \IK $ 
	 and $\ell_m/\braket{\ell_m}\rightarrow\bar c\neq 0$. 
	Therefore, in the limit for $m\rightarrow + \infty$, by  \eqref{eq:0_abs_m} we get $\partial_\gamma^k \vec{\omega}(\bar\gamma)\cdot \bar c = 0$ for any $k\in\N_0$.
	 By the analyticity of $ \vec{\omega}(\gamma)$, we deduce 
	 that the function $ \gamma \mapsto \vec{\omega}(\gamma)\cdot \bar c$ is identically zero on $\IK$, which contradicts Lemma \ref{non degenerazione frequenze imperturbate}.
	  \\[1mm]
	{ \bf Proof of \eqref{T Melnikov}}. The proof is similar to that of \eqref{1 Melnikov}, and thus we omit it.
	 \\[1mm]
	{ \bf Proof of \eqref{1 Melnikov}}. We divide the proof in $ 4 $ steps. \\ 
	{\sc Step 1. } Recalling 
	\eqref{ASYFR} 
	we have that, for any $\gamma \in\IK $, 
	$ 	| \vec{\omega}(\gamma)\cdot \ell + \Omega_n(\gamma)|
	 \geq $ $ |\Omega_n (\gamma)| - |\vec{\omega}(\gamma)\cdot \ell |  \geq $ $  c
	 n - C \langle \ell \rangle \geq $ $ \langle \ell \rangle  $
	whenever $ n \geq C_0 \langle \ell \rangle $, 
	for some $C_0>0$.  In this cases \eqref{1 Melnikov} is already fulfilled with $ k = 0 $.  Hence we  restrict  in the sequel to indexes $\ell\in\Z^{{|\ST|}}$ and 
	$ n \in \N \setminus (\ST \cup \{ 2, \ldots, \bar n \}) $ satisfying
	\begin{equation}\label{eq:1_restr}
	n < C_0 \langle \ell \rangle \,.
	\end{equation}
	{\sc Step 2.} By contradiction, 
	we assume that, for any $m\in\N$, 
	there exist $\gamma_m\in\IK$, 
	$\ell_m\in\Z^{|\ST|}$ and $ n_m\in \N \setminus (\ST \cup  \{ 2, \ldots, \bar n \}) $, 
	with ${n_m} <C_0 \langle \ell_m \rangle $, such that, for any 
	$ k \in\N_0$ with 
	$ k \leq m $,
	\begin{equation}\label{eq:1_abs_m}
	\big| \partial_\gamma^k \big( \vec{\omega}(\gamma) 
	\cdot\frac{\ell_m}{\braket{\ell_m}}+\frac{1}{\braket{\ell_m}}
	\Omega_{n_m}(\gamma) \big)_{|\gamma= \gamma_m} \big|  < \frac{1}{\braket{m}}  \, . 
	\end{equation}
Up to subsequences 
$\gamma_m\rightarrow\bar\gamma \in\IK$ and $\ell_m/\braket{\ell_m}\rightarrow \bar c\in\R^{|\ST|}$. \\
	{\sc Step 3.} We consider first the case when the sequence $(\ell_m)_{m\in\N}\subset \Z^{|\ST|}$ is bounded. Up to subsequences, we have definitively that 
	$\ell_m=\bar \ell\in\Z^{|\ST|} $. Moreover, since $ n_m$ and $\ell_m$ satisfy \eqref{eq:1_restr}, also the sequence $(n_m)_{m\in\N}$ is bounded and, up to subsequences,  definitively  $ n_m = \underline n \in \N \setminus (\ST \cup \{ 2,  \ldots, \bar n \}) $. 
	Therefore, in the limit $m\rightarrow \infty$, by \eqref{eq:1_abs_m} we obtain
$	\partial_\gamma^k \big( \vec{\omega}(\gamma) 
	\cdot \bar \ell + \Omega_{ \underline n}(\gamma) \big)_{|\gamma = \bar \gamma} = 0 $,  $ \forall k \in \N_0 $.  
By analyticity this implies 
$	\vec{\omega}(\gamma)\cdot \bar \ell + \Omega_{ \underline n}(\gamma) = 0 $, 
$ \forall \, \gamma \in \IK $,
which contradicts Lemma \ref{non degenerazione frequenze imperturbate}.	
\\[1mm]
{\sc Step 4.} We consider now the case  the sequence $(\ell_m)_{m\in\N}$ is unbounded. Up to subsequences
	$|{\ell_m}|\rightarrow \infty$ as $m\rightarrow\infty$ and  $  \lim_{m\rightarrow\infty}\ell_m/\braket{\ell_m} =: \bar c  \neq 0 $. By the asymptotic expansion \eqref{ASYFR1}, 
	for any $ k \in \N_0 $,
	\begin{equation*}
	\begin{split}
	\partial_\gamma^k \frac{1}{\braket{\ell_m}}\Omega_{n_m}(\gamma_m) & 
	=  \partial_\gamma^k \Big( \frac{n_m}{\braket{\ell_m}} 
	 \Omega_1(\gamma) - \frac{1}{2 \braket{\ell_m}} + 
	\frac{r(n_m,\gamma)}{\braket{\ell_m}}
	 \Big)_{|\gamma= \gamma_m}  \stackrel{\eqref{eq:1_restr}} \to 
\bar d (\partial_\gamma^k \Omega_1 (\gamma))_{|\gamma = \bar \gamma } 
\, , \ {\rm for}\ m \rightarrow \infty \, , 
	\end{split}
	\end{equation*}
where $\bar d:= \lim_{m\rightarrow\infty} n_m /\braket{\ell_m}  $ which is finite by
\eqref{eq:1_restr}. Therefore
	 \eqref{eq:1_abs_m}  becomes, in the limit $m\rightarrow\infty $,   
	$
	\partial_\gamma^k \big( \vec{\omega}(\gamma) \cdot \bar c + \bar d 
	\Omega_1 (\gamma) \, 
	\big)_{|\gamma= \bar \gamma } = 0 $, $ \forall \, k \in\N_0 $. 
By analyticity, this implies that $ \vec{\omega}(\gamma) \cdot \bar c + \bar d 
	\Omega_1(\gamma)  = 0 $ for any $\gamma \in\IK$. 
	This contradicts the non-degeneracy of the vector $(\vec{\omega}(\gamma), 
	\Omega_1 (\gamma) )$  
	provided in Lemma \ref{non degenerazione frequenze imperturbate}, 
	since $(\bar c, \bar d) \neq 0$.
\\[1mm]
	{\bf Proof of \eqref{2 Melnikov-}}. We split again the proof into 4 steps.\\
	{\sc Step 1.} By Lemma \ref{lem:asy}, for any $\gamma \in\IK $, 
	\begin{equation*}
	\begin{split}
	& | \vec{\omega}(\gamma)\cdot \ell  + \Omega_n(\gamma)-\Omega_{n'}(\gamma)|  
	 \geq 
	|\Omega_n(\gamma)-\Omega_{n'}(\gamma)| -
	| \vec{\omega}(\gamma)\cdot \ell|  \\
&  \qquad \qquad \qquad \geq | n- n' | \Omega_{1}(\gamma) - |r(n,\g)| - |r(n',\g)|
	- C  |\ell |   
	\geq | n- n' | c
	- C  \langle \ell \rangle \geq \langle \ell \rangle 
	\end{split}
	\end{equation*}
	whenever $ | n - n'  | \geq C_1 \langle \ell \rangle $ 
	for some $C_1>0$. In this case \eqref{2 Melnikov-} is already fulfilled with $ k = 0 $.	
	Thus we restrict to 
	indexes $\ell\in\Z^{|\ST|}$ and $ n , n '\in \N \setminus (\ST \cup \{ 2, \ldots, \bar n \}) $,
	such that
	\begin{equation}\label{eq:2_restr}
	\big| n - n' \big|  < C_1 \langle \ell \rangle \,.
	\end{equation}
	Furthermore we may assume also that $ n \neq n' $ because the case $ n = n' $ is included in \eqref{0 Melnikov}.
	\\[1mm]
	{\sc Step 2.} By contradiction, we assume that, for any $m\in\N$,  there exist $\gamma_m\in\IK$, $\ell_m\in\Z^{|\ST|}$ and $ n_m, n_m'\in \N \setminus (\ST \cup \{ 2, \ldots, \bar n \}) $, 
	satisfying \eqref{eq:2_restr} and $ n_m \neq n_m' $, such that, for any $0\leq k \leq m$,
	\begin{equation}\label{eq:2_abs_m}
	\big| \partial_\gamma^k\big( \vec{\omega}(\gamma) 
	\cdot \frac{\ell_m}{\braket{\ell_m}}+ \frac{1}{\braket{\ell_m}}\big( \Omega_{n_m}(\gamma)-\Omega_{n_m'}(\gamma) \big) \big)_{|
	\gamma= \gamma_m} \big| < \frac{1}{\braket{m}} \, . 
	\end{equation}
	Up to subsequences $\gamma_m\rightarrow\bar\gamma \in\IK$.
	We distinguish  two cases.
\\[1mm]
{\sc Step 3.} Suppose that $ (\ell_{m} )$ is unbounded. 
 Up to subsequences $\ell_m/\braket{\ell_m}\rightarrow \bar c\in\R^{|\ST|} \setminus \{0\} $ for $ m \to + \infty $. By \eqref{eq:2_restr} we deduce that
$(n_m - n_m' ) \langle \ell_m \rangle^{-1} \to \bar d \in \R $,  
and, by \eqref{ASYFR1}, 
\begin{align}
	\partial_\gamma^k 
\big(	\frac{1}{\braket{\ell_m}}\big( \Omega_{n_m}(\gamma)-\Omega_{n_m'}(\gamma) \big)  & 
	=  \partial_\gamma^k \Big( \frac{(n_m - n_m')}{\braket{\ell_m}} 
	 \Omega_1(\gamma)  + 
	\frac{r(n_m,\gamma)- r(n_m',\gamma)}{\braket{\ell_m}}
	 \Big)_{|\gamma= \gamma_m}  \nonumber \\
& \to 
\bar d (\partial_\gamma^k \Omega_1 (\gamma))_{|\gamma = \bar \gamma }  \, . \label{eq:1_expand2}
	\end{align}
Hence passing to the limit in \eqref{eq:2_abs_m} for $ m \to + \infty $, 
we deduce by 
\eqref{eq:1_expand2} that 
$ \partial_{\gamma}^k \big( {\vec \om}( \gamma) \cdot \bar c + \bar d 
 \Omega_1 (\gamma) \big)_{|\gamma 
= \bar \gamma} = 0 $,  $ \forall k \in \N_0 $. 
Therefore the analytic function $ \gamma  \mapsto  {\vec \om}( 
\gamma ) \cdot \bar c + \bar d \Omega_1 (\gamma)  $ is identically zero.  
This in contradiction with  Lemma \ref{non degenerazione frequenze imperturbate}, since
$ (\bar c, \bar d) \neq (0,0) $.
\\[1mm]
{\sc Step 4.} Suppose that $ (\ell_{m} )$ is bounded. 
Up to a subsequence, we have that definitively 
$ \ell_{m} = \bar \ell \in \Z^{|\ST|} $. Moreover, by \eqref{eq:2_restr} and 
since we restrict to indexes satisfying $ n_m \neq n_m' $,  we have that
\be\label{1nmC} 
1 \leq | n_m - n_m' | \leq C \, . 
\ee
The following  subcases are possible: 
\\[1mm]
$(i)$ $ n_m, n_m' \leq C $. 
Up to a subsequence, we have that definitively 
$ n_{m} =  \underline n  $, $ n'_{m} =  \underline n'  $ where $ \underline n ,  \underline n' $ are integers
in $ \N \setminus (\ST \cup \{2, \ldots, \bar n\} ) $, $ \underline n \neq  \underline n' $.
Hence passing to the limit in \eqref{eq:2_abs_m} we deduce that 
$ \partial_\gamma^k  
\big( {\vec \om}(  \gamma) \cdot {\bar \ell} + 
{ \Omega_{ \underline n}( \gamma) - 
\Omega_{ \underline n'}( \gamma )} \big)_{| \gamma = \bar \gamma } = 0 $, 
$ \forall k  \in \N_0  $.  
Hence the analytic function 
$ \gamma \mapsto {\vec \om}(  \gamma) \cdot
\bar \ell  
+   \Omega_{ \underline n}(\gamma) -  \Omega_{ \underline n'}( \gamma ) 
 $ is identically zero, which is a contradiction  
with Lemma \ref{non degenerazione frequenze imperturbate}, because
$ \underline n \neq  \underline n' $.
\\[1mm]
$(ii)$ $ n_m, n_m' \to + \infty $. 
By \eqref{1nmC} we have that definitively
$ n_m - n_m' = \bar p $ for some integer $ \bar p \in [1,C] $, and by 
\eqref{ASYFR1} we deduce that  
$$
\begin{aligned}
 \partial_\gamma^k 
\frac{1}{\langle \ell_m \rangle}
\big(\Omega_{n_m}(\gamma) - \Omega_{n_m'}(\gamma)\big)_{|\g = \g_m} 
& =
\partial_\gamma^k 
\frac{1}{\langle \bar \ell \rangle}
\big( (n_m - n_m') \Omega_1 (\gamma) + r(n_m,\gamma) - r(n_m', \gamma) \big)_{|\g = \g_m}  \\
& \to \bar d  \, \partial_\gamma^k 
(  \Omega_1 (\gamma)  )_{|\g = \g_m} \, , \quad {\rm for} 
\quad m \to + \infty \, , 
\end{aligned}
$$
where $ \bar d := \bar p / \langle \bar \ell \rangle \neq 0  $. 
Then by \eqref{eq:2_abs_m} we get 
$ \partial_{\gamma}^k  
\big\{ {\vec \om}( \gamma ) \cdot 
\bar \ell \langle \bar \ell \rangle^{-1}  + \bar d \Omega_1 (\gamma) \big\}_{|\gamma = \bar \gamma} = 0 $,  $ \forall k  \in \N_0  $.  
Hence the analytic function 
$ \gamma  \mapsto {\vec \om}( \gamma) \cdot \bar \ell \langle \bar \ell \rangle^{-1} 
+ \bar d  \Omega_1 (\gamma) $
is identically zero,
contradicting  Lemma \ref{non degenerazione frequenze imperturbate}, because
$ (\bar \ell \langle \bar \ell \rangle^{-1}, \bar d) \neq 0 $.	
\\[1mm]
{\bf Proof of \eqref{2 Melnikov+}}. The proof follows is similar 
to that of \eqref{1 Melnikov}
and we omit it.
\end{pf}

\section{The symplectic reduction of the angular momentum}\label{sec:SR}

  We now construct a symplectic diffeomorphism 
  which introduces  the prime integral $ {\cal J}  $ 
  defined in 
  \eqref{Jsvi} below, that is a multiple of the angular momentum $ J $  in \eqref{moangJ2},  
  as a symplectic coordinate, see Theorem \ref{thm:FM}. 
 We follow the Darboux-Caratheodory theorem of symplectic rectification, valid in finite dimension, 
 along the lines   of Theorem 10.20 in \cite{Marmi-Fasano}.  
 The proof is  much more delicate since
the phase space is infinite dimensional and the Hamiltonian vector field generated by 
$ J $ is a  first order 
transport operator.  As far as we are aware of, this idea wasn't used in previous PDE papers where KAM theory is applied and it is a major key idea to eliminate the degeneracy of the second mode. 
  
\smallskip

We decompose the phase space $ H^s_0 = H^s_0 (\T) $ as
\be\label{splipin0}
 H^s_0 = H_2 \oplus H^s_{\bot,2} \, , \quad 
H_2 := \langle \tc_2, \ts_2 \rangle \, , \quad 
H^s_{\bot,2} := 
  H^s_0 (\T) \cap H_2^\bot  
\ee
writing a real function $ \xi (\theta ) \in  H^s_0 (\T) $,
 in the basis $ \{ \tc_n, \ts_n \}_{n \geq 1 } $
 defined  in \eqref{defcnsn}, as
  \begin{equation}\label{a2b2n}
  \begin{aligned}
&  \qquad \qquad \qquad \quad  \xi (\theta) = v +  u   
\in H^s_0 (\T)  \, , \\
&  
 v =  \a_2 \, \tc_2  +  \b_2 \, \ts_2 \in H_2  \, , \quad 
u  = {\mathop \sum}_{n \in \N, n \neq 2} \a_n \, \tc_n  + \b_n \, \ts_n  \in 
H^s_{\bot,2} \, ,  
  \end{aligned}
  \end{equation}
  with coefficients  $(\a_n,\b_n)$ given in \eqref{newa2b2}.
 We denote by $ \Pi_2 $ the projector on $ H_2 $, 
and by 
$ \Pi_2^\bot $ the projector on $ H^s_{\bot,2}  $,  
 namely   $  \Pi_2^\bot \xi (\theta)  = u (\theta) $. 
We remind that in the coordinates $ (\a,\beta) := (\a_n, \beta_n)_{n \geq 1 } $  defined in \eqref{a2b2n}
  the symplectic form 
\eqref{sy2form} reads as in \eqref{sy2form0}, 
  \be\label{simpl2}
  (d \alpha_2 \wedge d \beta_2) \oplus  {\cal W}_{\bot,2} =
  d \alpha_2 \wedge d \beta_2 + 
   {\mathop \sum}_{n \neq 2} \, \tfrac{1}{n} d \alpha_n \wedge d \beta_n \, , 
  \ee
 where $ {\cal W}_{\bot,2} $ is the restriction to $H^s_{\bot,2}$
 of $ {\cal W} $, and  the Poisson bracket between two functions reads, see \eqref{Poib}, 
\be\label{Poibrak}
\{ f, g\} = {\mathop \sum}_{n \geq 1} p_n \big( \partial_{\a_n} f \, \partial_{\b_n} g - 
\partial_{\b_n} \, f \partial_{\a_n} g \big)  \quad \text{where}
\quad p_n := 
\begin{cases}
1 \quad \text{if}  \ n = 2 \cr
n \quad \text{if}  \ n \neq 2 \, . 
\end{cases}
\ee
We want to introduce the following prime integral   of the vortex patch equation, 
\begin{align}
 {\cal J} (\xi) \label{Jsvi}
 & :=  \aleph \big( J_1 (\xi) +  J_2 (\xi) 
  \big) 
   \stackrel{\eqref{moangJ2},\eqref{fg0},\eqref{newa2b2}} 
  = 
  \a_2 + \aleph J_2 (\xi)   \, \\
  & \label{formaJ}
= \frac12 (\xi, \tc_2) + \aleph \int_\T \xi^2 (\theta) \, g_\gamma (\theta) \, d \theta 
\, , \qquad \aleph := \tfrac{\sqrt{2}}{\sqrt{\pi}(\g -\g^{-1})} 
\end{align}
 with  $ J_1, J_2 $  defined in \eqref{moangJ2}, 
 as one of the  canonical coordinates in a neighborhood of zero.

\subsection{The flow of  the angular momentum} \label{flowJ}

For the construction of the symplectic set of variables which complete $ {\cal J }$
we shall use the symplectic flow  $ \Phi^{ t }_{{\cal J}}(\xi) $  generated by 
the Hamiltonian $ {\cal J} $, namely solving 
 \be\label{flow-J}
 \pa_t \Phi^{ t }_{{\cal J}}(\xi) = X_{{\cal J}} (\Phi^{ t }_{{\cal J}}(\xi)) \, , \quad 
 \Phi^{0}_{{\cal J}}(\xi) = \xi \, , 
 \ee
where, using \eqref{Mom1}, \eqref{fg0}, \eqref{Jsvi}, \eqref{defcnsn}, 
\begin{equation}\label{XJcal}
X_{\cal J} (\xi) = X_{{\cal J}_1} (\xi) + X_{{\cal J}_2} (\xi) =  - \ts_2  + 
2 \aleph \big( \pa_\teta \circ g_\gamma (\theta)\big) (\xi)
\end{equation}
is the Hamiltonian vector field generated by 
$ {\cal J} = {\cal J}_1 + {\cal J}_2 = \aleph J_1 +  \aleph J_2   $ 
and $ J_1, J_2 $ are defined in \eqref{moangJ2}.
Note that the Hamiltonian vector field  
\be \label{xj1}
X_{{\cal J}_1} (\xi) = - \ts_2 
\ee
is constant in $ \xi $ 
and  $ X_{{\cal J}_2} (\xi) = 2 \aleph \big( \pa_\teta \circ g_\gamma (\theta)\big) (\xi)  $ 
is linear. 
By \eqref{Jsvi}  and \eqref{XVS}, in the variables $ (\a_n,\b_n)_{n \geq 1} $ 
only the $ \beta_2$-component  of $ X_{{\cal J}_1} $ is non zero 
and it is equal to $ - 1 $  and the Hamiltonian system generated by the 
 vector field $ X_{\cal J}  $ reads 
 \be \label{xj1coo} 
\Big(
    \dot \a_1,  
    \dot \b_1, 
    \dot \a_2, 
    \dot \b_2, 
\dot \a_3, 
    \dot \b_3,    
    \ldots\Big)^\top =  
\Big( 
    0,  
    0, 
    0,  
    -1,  
0,
0, 
    \ldots
    \Big)^\top + B(\a,\b)
\ee 
where $B $ is a matrix which can be explicitly computed by \eqref{moangJ2}
and \eqref{fg0}. 
Nevertheless, since $ X_{\cal J}$ is a transport 
operator,  its dynamics 
is better understood 
in   terms of the variable $ \xi (\theta )$. 
Note that, differentiating \eqref{flow-J}-\eqref{XJcal},
the differential $d \Phi^{ t }_{{\cal J}}(\xi) $ solves 
the linear system  
\be\label{difJ}
 \pa_t d \Phi^{ t }_{{\cal J}}(\xi) = d X_{\cal J} (\Phi^{ t }_{{\cal J}}(\xi) ) \, d \Phi^{ t }_{{\cal J}}(\xi) 
 = 
  2 \aleph (\pa_\teta \circ {g}_\gamma (\theta)) \, d \Phi^{ t }_{{\cal J}}(\xi) \, , \quad 
 d \Phi^{0}_{{\cal J}}(\xi) = {\rm Id}  \, , 
\ee
namely, for any $ \xi $,  
\be\label{uguaJJ2}
d \Phi^{ t }_{{\cal J}} (\xi ) =  \Phi^{ t }_{{\cal J}_2} 
\ee
is the linear flow 
generated by the linear transport operator
$ X_{{\cal J}_2} = 2 \aleph \pa_\teta \circ {g}_\gamma (\theta) $  in \eqref{XJcal}. The  
 flow $ \Phi^{ t }_{{\cal J}_2}(\xi)  $  is well defined 
 in  $ H^s (\T )$ and it is symplectic since the 
 vector field  $ X_{{\cal J}_2} $ is Hamiltonian. 
It results 
$ [ \Phi^{t}_{{\cal J}_2} ]^{-1} =  \Phi^{-t}_{{\cal J}_2} =  \Phi^{t}_{-{\cal J}_2} $. 
 We now provide 
  an explicit characterization. 
The key step is to `symplectically rectify" 
 the quadratic term of the angular momentum
  $ J_2 (\xi) $ 
in \eqref{moangJ2}. 
 
 \begin{lemma}\label{lemma:adj}
 Consider the diffeomorphism 
$ y = \theta + \beta (\theta) $  of  $ \T $ 
where $ \beta (\theta )$ is the odd 
function  
\be\label{bestaes}
 \beta (\theta) := \arctan \Big( \frac{\tan (\theta) }{\gamma}\Big)  - \theta \, , \quad \forall \theta \in 
 \big(- \frac{\pi}{2}, \frac{\pi}{2}\big) \, , \quad  \beta ( \pm \pi/2) = 0 \, , 
 \ee
extended on $ \R $ to a smooth $ \pi $-periodic function, with 
 inverse 
\be\label{invediff}
y = \theta + \beta (\theta) =
 \arctan \Big( \frac{\tan (\theta) }{\gamma} \Big) 
\quad \Leftrightarrow \quad \theta = y + \breve \beta (y) = 
 \arctan \big( \gamma \tan (y) \big) \, .
\ee
Under the symplectic  linear  change of variable  
\be\label{LDS}
\xi (\theta) := \Psi (q) (\theta) := (1+ \beta_\theta (\theta)) q (\theta + \beta (\theta)) \, , 
\ee
 the quadratic part of the angular momentum $ J_2 (\xi ) $ in 
\eqref{moangJ2} transforms into 
\be\label{J2n}
J_2 (\Psi(q)) =  \int_{\T} q^2 (y)   \, d y \, .
\ee
 \end{lemma}

 \begin{pf}
 First note that the odd monotone increasing function 
 $$
 F(\theta) := \theta + \beta (\theta) = \arctan \Big( \frac{\tan (\theta) }{\gamma} \Big) \, , \quad
 \forall \theta \in  \big(- \frac{\pi}{2}, \frac{\pi}{2}\big) \, ,  
 $$
  solves 
\be\label{derdiffeodir}
F' (\theta) = 1 + \beta_\theta (\theta) = \frac{1}{\gamma \cos^2 (\theta ) + \gamma^{-1} \sin^2 (\theta)  } = \frac{1}{g_\gamma (\theta) } \, . 
\ee
Note also that the function $g_\g (\theta)   $ in \eqref{fg0} does not vanish since 
$$
\g^{-1} = \frac{\gamma + \gamma^{-1}}{2} -  \frac{\gamma  - \gamma^{-1}}{2}  
\leq 
g_\g (\theta)  
\leq
\frac{\gamma + \gamma^{-1}}{2} + \frac{\gamma  - \gamma^{-1}}{2} 
= \g \, , \quad \forall \theta \in \R \, , 
$$
and then the $ \pi $-periodic   function  
$ 1/ g_\gamma (\theta ) $ 
is $  \mathcal{C}^\infty $. Moreover by \eqref{derdiffeodir}
 its average is equal to $ 1 $ and   
$ \beta (\theta ) = F(\theta)- \theta $ extends 
to  a $  \mathcal{C}^\infty $, odd $ \pi $-periodic function on $ \R $.  
Finally
$$
\begin{aligned}
J_2 (\Psi(q))) 
& \stackrel{\eqref{moangJ2},\eqref{LDS}} =  \int_{\T} q^2 (\theta + \beta (\theta))  (1 + \beta_\theta (\theta))^2 g_\gamma (\theta)  \, d \theta \\
& \stackrel{\eqref{derdiffeodir}} =
\int_{\T} q^2 (\theta + \beta (\theta))  (1 + \beta_\theta (\theta))   \, d \theta
 =  \int_{\T} q^2 (y)   \, d y 
\end{aligned}
$$
proving \eqref{J2n}. 
\end{pf}

\begin{corollary}\label{cor:sez8}
{\bf (Straightening of the angular momentum)}
The Hamiltonian transport operator 
$  \pa_\theta \circ g_\gamma (\theta)$  generated by $ X_{J_2}$ is conjugated under the map 
$ \Psi $ defined in \eqref{LDS} into  the constant transport operator $ 2 \pa_y $, namely 
\be\label{psi-1psi}
\Psi^{-1} \circ (\pa_\theta \circ g_\gamma (\theta)) \circ  \Psi = 2 \pa_y \, . 
\ee
Thus the corresponding Hamiltonian flow $\Phi^t_{J_2}  $ is  
\be\label{psiflowpsi}
\Phi^t_{J_2} =  \Psi \circ T_{2t} \circ \Psi^{-1}
\ee
where $ (T_{2t} q) (y) := q (y + 2 t )$ is the symplectic transport flow.  
\end{corollary}

\begin{pf}
Since the map $ \xi = \Psi (q)  $ in \eqref{LDS} is symplectic, the linear 
Hamiltonian vector field 
$ X_{J_2} (\xi) = \pa_\theta \circ g_\gamma (\theta) \, \xi $ 
 generated by $ J_2  (\xi )$ is transformed into the linear 
 Hamiltonian vector field $ 2 \pa_y $
generated by 
$ J_2 ( \Psi(q)) = \int_{\T} q^2 (y) d y $. This proves \eqref{psi-1psi}
and \eqref{psiflowpsi} follows. 
\end{pf}

By Corollary \ref{cor:sez8} and recalling the definition of $ {\cal J} $ in 
\eqref{Jsvi}  
we have the following result.

\begin{lemma}\label{flowJ2old}
{\bf (Representation of the flow $ \Phi^{t}_{{\cal J}_2} $)}
The linear symplectic 
flow $ \Phi^{t}_{{\cal J}_2}  $ generated by $ X_{{\cal J}_2} $  can be written as 
\be\label{reprePhi2}
\Phi^{t}_{{\cal J}_2} =  \Psi \circ T_{2\aleph t} \circ  \Psi^{-1} = 
\underbrace{(1 + \beta_\theta (\theta)) \circ {\cal B}}_{\Psi}  \, 
 \circ \, T_{2\aleph t} \, \circ 
 \underbrace{ {\cal B}^{-1}  \circ \frac{1}{1 + \beta_\theta (\theta)}}_{\Psi^{-1}}  
\ee
where 
\\[1mm] 
1. the mappings  ${\cal B}  $, $ {\cal B}^{-1} $ are the composition operators   
\be\label{phidin}
({\cal B} q) (\theta) := q (\theta + \beta (\theta)) \, , \quad
({\cal B}^{-1} \xi)(y) := \xi (y + \breve \beta (y)) \,   
\ee
induced by the 
diffeomorphism of $ \T $ defined by 
$ y = \theta + \beta (\theta) $ 
 in 
\eqref{invediff}, with inverse $ \theta = y + \breve \beta (y) $; 
\\[1mm]
2. the map 
$ T_{2 \aleph t} $ is the symplectic translation operator  
\be\label{translation-op}
(T_{2 \aleph t} q)(y) := q (y + 2 \aleph t ) \, . 
\ee
 Moreover the adjoint is
\be\label{repreaggun}
\big(\Phi^{t}_{{\cal J}_2}\big)^* =  
{\cal B} \circ T_{-2\aleph t}  \circ  {\cal B}^{-1}\, , 
\quad \big(\Phi^{t}_{{\cal J}_2} - {\rm Id}\big)^* =  
{\cal B} \circ \big( T_{-2\aleph t}  - {\rm Id} \big) \circ  {\cal B}^{-1}\, .  
\ee
\end{lemma}

\begin{pf}
Formula \eqref{reprePhi2} follows by Corollary \ref{cor:sez8}.
The adjoint of 
$ \Psi = (1 + \beta_\theta (\theta)) \circ {\cal B}  $
is $
\Psi^{*} = {\cal B}^{-1} $. 
Indeed
$$
\begin{aligned}
(\Psi h, k) 
 =  \int_{\T} (1 + \beta_\theta (\theta)) h (\theta + \beta (\theta)) k (\theta) d \theta 
 =  \int_{\T}  h (y) k (y + \breve \beta (y)) d y = (h, {\cal B}^{-1} k)
\end{aligned}
$$
proving that $ \Psi^* = {\cal B}^{-1} $.  Thus $ (\Psi^{-1})^* = (\Psi^*)^{-1}  = {\cal B} $  
and, taking the adjoint of \eqref{reprePhi2}, 
and since $T_{2\aleph t}^*  = T_{-2\aleph t} $, one deduces \eqref{repreaggun}. 
\end{pf}

\begin{remark} 
We shall denote the adjoint 
either $ \Phi^*$ or $ \Phi^\top $ since we have a  real scalar product. 
\end{remark}

The flow $ \Phi^t_{{\cal J}}  $ generated by
 the affine vector field $ X_{ {\cal J}}$ in \eqref{XJcal}   is  affine,
 see \eqref{flowaff}.
 
\begin{lemma}\label{lem:FlowJ}
{\bf (The symplectic flow $ \Phi^t_{\cal J}$)}
The flow $ \Phi^t_{{\cal J}}  $ defined  in \eqref{flow-J} 
 generated by
 the affine vector field $ X_{\cal J}$ in \eqref{XJcal}  is  affine, and may be represented by
 \be\label{flowaff}
\Phi^t_{ \cal J} (\xi) = \big( {\rm Id} - \Phi^t_{{\cal J}_2} \big) \xi_p + \Phi^t_{{\cal J}_2} \xi  
 \ee
 where $\xi_p (\theta) $ is the  $ \pi $-periodic $ {\cal C}^\infty $
 function with zero average
\be\label{xipart}
\xi_p (\theta) :=  \frac12 \Big( \frac{1}{  g_\gamma(\theta)}  - 1 \Big) 
= \frac12 \beta_\theta (\theta) \, ,
\ee
and $ \beta (\theta) $ is defined in \eqref{bestaes}. 
\end{lemma}

\begin{pf}
All the solutions of the non-homogeneous  equation
$ \pa_t \xi = \aleph \pa_\theta ( g_\gamma (\theta ) (1+ 2 \xi )) $,  
generated by $ X_{\cal J} (\xi)  $ (see \eqref{Mom1})  
are given by the sum of a particular  solution 
plus all the solutions of the homogeneous equation. We look for a particular  solution
$ \xi_p (\theta ) $ 
which is constant in time, i.e. solves 
$ g_\gamma (\theta ) (1+ 2 \xi_p ( \theta)) = c $ for some 
$ c \in \R $. This amounts to
$
\xi_p (\theta) = \tfrac12 \big( \tfrac{c}{  g_\gamma(\theta)}  - 1 \big)
$
and $ \xi_p (\theta) $ has zero average if and only if $ c = 1 $
(since the average of $ 1/ g_\gamma (\theta ) $ 
 is equal to $ 1 $ by Lemma \ref{lemma:adj}). In conclusion 
 $\xi_p (\theta) $ is given in \eqref{xipart}
and \eqref{flowaff} follows, recalling that $ \Phi^t_{{\cal J}_2} $ is the  flow of the 
homogeneous linear equation
$ \pa_t \xi = 2 \aleph \pa_\theta ( g_\gamma (\theta )  \xi ) $.  
\end{pf}

\subsection{The symplectic rectification} 
We  now construct the claimed symplectic variables 
which include 
$ {\cal J}$ as a coordinate. 
We denote by $B_r (H^s_0) $ the ball of center $ 0 $ radius $ r > 0 $ in 
$ H^s_0 := H^s_0 (\T )$. We recall that 
$ \Phi^{ t }_{ {\cal J}}(\xi) $ denotes the Hamiltonian  
flow generated by the angular momentum.

 \begin{theorem}{\bf (Symplectic rectification)}\label{thm:FM}
There exists $ r > 0 $ and  a smooth function $  \overline t : B_r (L^{2}_0)  \to \R $ 
satisfying 
\begin{enumerate}
\item \label{probart}
$ \overline t (0) = 0 $, 
$  d \bar t (0)[\widehat \xi]  = ( \nabla \bar t (0),\widehat \xi )    
  = \tfrac12 ( \ts_2, \widehat \xi ) $ for any  $ \widehat \xi \in L^{2}_0 $, 
 thus $ \overline t (\xi) $ has  the form 
$  \overline t (\xi) = \beta_2 + O(\| \xi \|_{L^2}^2) $,  
 with gradient 
 \be\label{gradbart}
 \nabla \bar t (\xi) = - 
 \frac{\big(  \Phi^{\overline t (\xi) }_{{\cal J}_2}\big)^* \ts_2}{ \big( 
   \ts_2, X_{\cal J}(\Phi^{\overline t (\xi) }_{{\cal J}}(\xi) ) \big) } \, ;  
\ee
\item  \label{tSb}
$ \bar t \circ {\cal S} = - \bar t $ 
where $ {\cal S} $ is the involution in \eqref{S-invo};    
\item \label{timeimpact0}
$ \bar t (\xi ) = 0 $ if and only if  $  (\xi  , \ts_2 ) = 0 $;  
\item \label{invtimesi}
$ \overline t ( \Phi_{\cal J}^\tau (\xi)) =  \overline t (\xi ) - \tau $, for any $ \tau $ small so that 
$ \Phi_{\cal J}^\tau (\xi) \in B_r (L^2_0)$; 
\end{enumerate}
 such that the map 
 \be\label{Phiorto}
 \Phi (\xi) :=
 {\cal J}(\xi) \tc_2  + \overline t (\xi)  \ts_2  + \Phi_2^\bot (\xi)
  \, , \quad \Phi_2^\bot (\xi) :=  \Pi_2^\bot \Phi^{ \overline t (\xi) }_{{\cal J}}(\xi) \, , 
\ee
satisfies the following properties: 
\begin{itemize}
\item (i) 
$ \Phi $ is symplectic; 
\item (ii)  $ \Phi (0) = 0 $  and,  
for any $ \xi \in B_r (L^{2}_0) \cap H^s_0 $, $ s \geq 1 $,
the differential $d \Phi (\xi) $ is
\be\label{diffePhi}
d \Phi (\xi)[\widehat \xi] = d {\cal J} (\xi)[\widehat \xi] \, \tc_2 + 
d {\bar t} (\xi)[\widehat \xi] \, \ts_2 + 
 \Pi_2^\bot X_{{\cal J}_2} (\Phi_{\cal J}^{\overline t (\xi)} (\xi))  \, 
d \bar t (\xi) [\widehat \xi] + 
 \Pi_2^\bot \Phi^{  \overline t (\xi) }_{{\cal J}_2} [\widehat \xi]   \, ,
\ee
in particular 
\be\label{dif0}
d \Phi (0) = {\rm Id} \, . 
\ee
\item  (iii) 
The map $ \Phi $ is locally invertible close to $ 0 $: for any
$ \eta = \eta_\tc \tc_2 + \eta_\ts \ts_2 + \eta_\bot $, with 
$ \eta_\bot \in H_2^\bot $, and $ \| \eta \|_{L^2} $ small, 
\be\label{formainv}
\Phi^{-1} (\eta) =  \Phi_{\cal J}^{- \eta_\ts} \big(  v_\tc  ( \eta) \tc_2  + \eta_\bot \big) 
\ee
where
\be\label{explivc}
 v_\tc  (\eta) := 
v_\tc  (\eta_\tc, \eta_\bot) :=
\frac{1}{\alpha} \Big( 1 + \frac{(\tc_4 , \eta_\bot )}{\sqrt{\pi} }\Big) \psi 
\Big( \alpha \frac{ \eta_\tc - \aleph (\eta_\bot^2, g_\gamma )}{ \big( 1 + \frac{(\tc_4 , \eta_\bot )}{\sqrt{\pi} }\big)^2 }  \Big) 
\ee 
$ \alpha $ is the constant $ \alpha := \frac{(\gamma^2 + 1) \sqrt{2}}{
(\gamma^2 - 1) \sqrt{\pi}} $,
 and $ \psi (z) $  is the  inverse of the analytic 
diffeomorphism  $ g : (- \tfrac12, \tfrac12) \to (- \tfrac14, \tfrac34) $, 
$ y \mapsto g(y) := y + y^2 $.  It results $ \psi (0) = 0 $, $ \psi' (0) = 1 $. 

Note that  $ v_\tc(0)=0$ and 
$ \pa_\eta v_\tc (0)_{|  H_2^\bot  } = 0 $. 

\item (iv)
The map  $ \Phi $ satisfies the following properties
\begin{enumerate}
\item[(a)] \label{SRP} $ \Phi $ is 
reversibility preserving, namely
$ \Phi  \circ {\cal S} = {\cal S} \circ \Phi $; 
\item[(b)] 
$  \Phi ( \Phi^{\tau}_{\cal J} (\xi))   =  \Phi (\xi)- \tau \ts_2 $, 
for any $ \tau $ small so that 
$ \Phi_{\cal J}^\tau (\xi) \in B_r (L^2_0)$. 
\end{enumerate}
\end{itemize}
 \end{theorem}  
  \begin{pf}
  We first define the canonical variable associated to $ {\cal J}$.
\\[1mm]
  {\bf Definition of $ \bar t (\xi ) $. }
The function 
  $ \bar t (\xi)$ is defined as the time such that, for any  $ \| \xi \|_{L^2} < r $ small enough,   the flow $ \Phi^{\bar t(\xi)}_{{\cal J}} (\xi) $ intersects the section 
  $ \{ \beta_2 = 0 \} $, which is transverse to $ X_{\cal J} (\xi )$ close to zero
 (see  \eqref{xj1}-\eqref{xj1coo}),  namely  
\be\label{solulemmaco}
 (  \ts_2 ,  \Phi^{\bar t(\xi)}_{{\cal J}} (\xi))_{L^2} = 0 \, . 
  \ee
This equation can be solved by the 
  implicit function theorem. 
  The function $ F : \R \times L^2_0 (\T) \to \R $ defined by 
  $$ 
\begin{aligned}
  F(t, \xi) :=  (  \ts_2 ,  \Phi^{\bar t(\xi)}_{{\cal J}} (\xi))_{L^2} 
     &  \stackrel{\eqref{flowaff}} 
   = (  \ts_2 , \big( {\rm Id} - \Phi^t_{{\cal J}_2} \big) \xi_p )_{L^2}  +
  ( \ts_2, \Phi^t_{{\cal J}_2} \xi  )_{L^2}   \\
  &   \stackrel{\eqref{repreaggun}}   = 
  ( {\cal B} \circ \big( T_{-2\aleph t}  - {\rm Id} \big) \circ  {\cal B}^{-1} \ts_2 , \xi_p )_{L^2}  +
  ( {\cal B} \circ  T_{-2\aleph t}  \circ  {\cal B}^{-1} \ts_2,  \xi  )_{L^2} 
  \end{aligned}
  $$
  is  smooth, actually it is affine in $ \xi $. 
  We have $ F(0,0)  = 0 $ and 
  $$
  \begin{aligned}
  \pa_t F(0,0)  
   = (  \ts_2 , \pa_t \Phi^t_{ {\cal J}} (\xi)_{|t=0,\xi = 0}) \stackrel{\eqref{flow-J}}  
   = (  \ts_2 , X_{ {\cal J}} (0)) \stackrel{\eqref{XJcal}} 
   =  - (  \ts_2, \ts_2) \stackrel{\eqref{defcnsn}} = - 2  \neq 0  \, .   
  \end{aligned}
  $$
  Then by the implicit function theorem there exists a unique $ \bar t (\xi)$, smooth in $ \xi $
  in $ \{ \| \xi \|_{L^2} < r \} $, such that  \eqref{solulemmaco} holds and $\bar t (0)=0$. 
 \\[1mm]
{\bf Properties of $ \bar t(\xi ) $. } 
By differentiating  \eqref{solulemmaco} with respect to $\xi $ 
in the direction $\widehat \xi$,  we get
$$
  \big( 
   \ts_2, (\pa_t \Phi^{ t }_{{\cal J}})_{|t = \overline t (\xi)} (\xi) \big) d \bar t (\xi)[\widehat \xi]
   +
   \big(  \ts_2, d \Phi^{ t }_{{\cal J}}(\xi)_{|t = \overline t (\xi)}[\widehat \xi] \big) =0\, , 
$$
and then 
 \be\label{sympjt-1}
  d \bar t (\xi)[\widehat \xi] 
     = - 
   \frac{\big(  \ts_2, d \Phi^{ t }_{{\cal J}}(\xi)_{|t = \overline t (\xi)}[\widehat \xi] \big)}
   {\big(  \ts_2,  (\pa_t \Phi^{ t }_{{\cal J}})_{|t = \overline t (\xi)} (\xi) \big) } 
   = - 
   \frac{\big(  \ts_2, d \Phi^{ t }_{{\cal J}}(\xi)_{|t = \overline t (\xi)}[\widehat \xi] \big)}{\big( 
   \ts_2, X_{\cal J}(\Phi^{\overline t (\xi) }_{{\cal J}}(\xi) ) \big) } \, .
\ee
By  \eqref{sympjt-1} and \eqref{uguaJJ2} we deduce 
 $$
  (\nabla \bar t (\xi), \widehat \xi ) 
  = - 
   \frac{\big(  \ts_2,  \Phi^{\overline t (\xi) }_{{\cal J}_2} [\widehat \xi] \big)}{\big( 
   \ts_2, X_{\cal J}(\Phi^{\overline t (\xi) }_{{\cal J}}(\xi) ) \big) } 
    = - 
\Big( \frac{\big(  \Phi^{\overline t (\xi) }_{{\cal J}_2}\big)^* \ts_2}{ \big( 
   \ts_2, X_{\cal J}(\Phi^{\overline t (\xi) }_{{\cal J}}(\xi) ) \big) },  \widehat \xi \Big)
$$
proving \eqref{gradbart}.
In particular, evaluating \eqref{sympjt-1} at $ \xi=0 $, we get
 \begin{align}\label{sympjt1}
 d \bar t (0)[\widehat \xi]
 &=- \frac{(   \ts_2, d \Phi^{ 0 }_{{\cal J}}(0)[\widehat \xi])}{(  \ts_2, X_{{\cal J}}( \Phi^{ 0 }_{{\cal J}} (0)))} 
 \stackrel{\eqref{difJ},\eqref{XJcal}}  =  \frac{(  \ts_2, \widehat \xi)}{ (   \ts_2, \ts_2)}
  \stackrel{\eqref{defcnsn}} = \frac12 (  \ts_2, \widehat \xi)  
  \stackrel{\eqref{newa2b2}}  = \widehat\beta_2 
   \end{align}
proving Item \ref{probart}.  
Thus $  \nabla \bar t (0) = \tfrac12 \ts_2 $. 

We now prove Item \ref{tSb}. Since 
 $ \bar t (\xi)$ is defined as the unique local solution of \eqref{solulemmaco}
we have 
\be\label{tbarS}
(  \ts_2, \Phi^{\bar t({\cal S}\xi)}_{{\cal J}}({\cal S}\xi)) = 0 \, . 
\ee
The angular momentum $ {\cal J}$  satisfies the reversibility property
(cfr. Lemma \ref{lem:rev}),  
\be\label{JSc}
 {\cal J} \circ  {\cal S}  =   {\cal J}  \, , 
\ee
or equivalently   its Hamiltonian vector field $ X_{{\cal J}} $  
satisfies  $ X_{{\cal J}} \circ {\cal S} = - {\cal S} \circ X_{{\cal J}} $.
This is tantamount to say that the flow $ \Phi^t_{{\cal J}}$ fulfills the property  
\be\label{flowc}
\Phi^{t}_{{\cal J}} \circ {\cal S} = {\cal S} \circ \Phi^{-t}_{{\cal J}} \, .
\ee
Therefore,  by \eqref{tbarS}, \eqref{flowc} and since $ {\cal S}^\top= {\cal S} $ we deduce that 
$ (  \ts_2 , \Phi^{- \bar t({\cal S}\xi)}_{{\cal J}}(\xi)) = 0 $.
By the unicity of the solutions of \eqref{solulemmaco} 
this implies $ - \bar t({\cal S}\xi) = \bar t (\xi )$, which is Item \ref{tSb}.  

Item \ref{timeimpact0} follows by the definition of $ \bar t (\xi )$ 
and the uniqueness  of the solutions of \eqref{solulemmaco}: the time of impact 
$  \bar t (\xi ) = 0 $ if and only if $ (\xi, \ts_2 ) = 0 $. 
Finally Item \ref{invtimesi} follows because, by \eqref{solulemmaco}, 
$$
0 = (  \ts_2, \Phi^{\bar t(\xi)}_{{\cal J}}\xi)  =
(  \ts_2, \Phi^{\bar t(\xi)}_{{\cal J}} \Phi^{-\tau}_{\cal J} \Phi^{\tau}_{\cal J} \xi)  
= 
(  \ts_2, \Phi^{\bar t(\xi) -\tau}_{{\cal J}}  \Phi^{\tau}_{\cal J} \xi)  
$$
and therefore $\bar t(\xi) -\tau = \bar t (\Phi^{\tau}_{\cal J} \xi ) $, by 
the uniqueness of the solutions of \eqref{solulemmaco}. 
\\[1mm]
{\bf Properties of $ \Phi $. }
 We now prove the properties ($i$)-($iv$)  of $ \Phi $  in \eqref{Phiorto}.   
We first prove Item ($ii$). Since 
$ {\cal J} (0) = 0 $ and $ \bar t (0) = 0 $, we have $ \Phi (0) = 0 $. Moreover \eqref{diffePhi} directly follows differentiating \eqref{Phiorto}. 
Finally  \eqref{dif0} because $ d {\cal J} (0) [\widehat \xi] = \widehat \alpha_2  $,  
\eqref{sympjt1} and 
 $ X_{{\cal J}_2} (0) = 0 $.  

Next we prove Item ($iv$)-(a). 
The map $ \Phi $ in \eqref{Phiorto}  is reversibility-preserving
 by \eqref{JSc}, Item \ref{tSb} and 
since, using \eqref{flowc},  
$$ 
\Phi_2^\bot ({\cal S}\xi)  =  \Pi_2^\bot \Phi^{ \overline t ({\cal S} \xi) }_{{\cal J}}({\cal S}\xi) 
 = 
 \Pi_2^\bot \Phi^{ - \overline t ( \xi) }_{{\cal J}}{\cal S}\xi 
 = 
 \Pi_2^\bot {\cal S} \Phi^{ \overline t ( \xi) }_{{\cal J}}\xi = 
{\cal S}  \Pi_2^\bot  \Phi^{ \overline t ( \xi) }_{{\cal J}}\xi = 
{\cal S}  \Phi_2^\bot (\xi) \, . 
$$ 
Item ($iv$)-(b) follows by 
$ {\cal J}( \Phi^\tau_{\cal J} (\xi)) =  {\cal J}(\xi)$ for any $ \tau \in \R $, 
the invariance of $ \bar t $ in Item \ref{invtimesi} and  
$$
\Phi_2^\bot (\Phi^{\tau}_{{\cal J}}\xi) = 
\Pi_2^\bot \Phi^{ \overline t (\Phi^{\tau}_{{\cal J}}\xi) }_{{\cal J}}(\Phi^{\tau}_{{\cal J}}\xi) 
= \Pi_2^\bot \Phi_{\cal J}^{ \overline t (\Phi^{\tau}_{{\cal J}}\xi) + \tau } (\xi) =
 \Pi_2^\bot \Phi_{\cal J}^{\overline t (\xi )} (\xi) = \Phi_2^\bot (\xi) \, . 
$$
{\bf Inversion of $  \Phi(\xi) $.}  
Given $ \eta = \eta_\tc \tc_2 + \eta_\ts \ts_2  + \eta_\bot $, 
we have to solve the equation 
\be\label{eq1}
\Phi (\xi) = \eta \, , \  \ \text{i.e.} \ \
{\cal J}(\xi) = \eta_\tc \, , \ \bar t(\xi) = \eta_\ts \, , \ 
\Pi_2^\bot \Phi^{ \overline t (\xi) }_{{\cal J}}(\xi) = \eta_\bot  \, . 
\ee
We look for a solution of \eqref{eq1} of the form
\be\label{oftheform}
\xi = \Phi_{\cal J}^{- \eta_\ts} (\breve \xi) \, . 
\ee
By Item ($iv$)-(b), this is equivalent to $ \Phi (\breve \xi ) = \eta - \eta_\ts \ts_2 $,
namely 
\be\label{eqtri}
{\cal J}(\breve \xi) = \eta_\tc \, , \ \bar t(\breve \xi)  = 0 \, , \ 
\Pi_2^\bot \breve \xi  = \eta_\bot \, .
\ee
All the solutions of the last equation in \eqref{eqtri} are 
$$
\breve \xi  = v + \eta_\bot \, , \quad v = v_\tc \tc_2 +  v_\ts \ts_2 \in H_2 \, . 
$$
Moreover, by Item \ref{timeimpact0}, 
all the solutions of also the second equation in \eqref{eqtri} are
\be\label{invbot}
\breve \xi  = v_\tc \tc_2  + \eta_\bot  \, .  
\ee
Finally we have to choose $ v_\tc  $  to be the solution of  the first scalar equation 
in \eqref{eqtri}, 
\be\label{eq4}
{\cal J}( v_\tc \tc_2  + \eta_\bot ) = \eta_\tc \, . 
\ee
{\bf Expansion of $ {\cal J }$.} 
By \eqref{formaJ}, \eqref{fg0}, \eqref{defcnsn}, \eqref{Jsvi}
and using $ \int_\T \tc_2^3  = 0 $,
using $ \cos^2 (2\theta) = (1 + \cos (4 \theta)) / 2 $, the fact that $ \eta_\bot $ has zero average, and  \eqref{newa2b2}, 
we deduce that $ {\cal J} $ in \eqref{formaJ} has the expansion
\be\label{calJex}
{\cal J} (v_\tc \tc_2  + \eta_\bot ) = v_\tc 
\Big( 1 + \frac{1}{\sqrt{\pi}}  (\tc_4 , \eta_\bot ) 	\Big) 
+ \alpha v_\tc^2 + \aleph \int_{\T} \eta_\bot^2 g_\gamma (\theta ) \, d \theta\quad \textnormal{where}\quad   \alpha = \frac{(\gamma^2 + 1) \sqrt{2}}{
(\gamma^2 - 1) \sqrt{\pi}}\,.
\ee
{\bf Solution of equation \eqref{eq4}.}
In view of \eqref{calJex}, we look for a solution of  \eqref{eq4} of the form
$$
v_\tc = \mu y \, , \quad 
\mu := \mu (\eta_\bot):= 
\frac{1}{\alpha} \Big( 1 + \frac{(\tc_4 , \eta_\bot )}{\sqrt{\pi} }\Big) \, .
$$
Thus the new variable $ y $ has to solve
$$
y +  y^2 = \frac{1}{\alpha \mu^2} ( \eta_\tc - \aleph (\eta_\bot^2, g_\gamma ) )
= \alpha \frac{ \eta_\tc - \aleph (\eta_\bot^2, g_\gamma )}{ \big( 1 + \frac{(\tc_4 , \eta_\bot )}{\sqrt{\pi} }\big)^2 } \, .
$$
Denoting by  $ \psi $  the inverse of the analytic diffeomorphism
$ y \mapsto g(y) := y + y^2 $ for $ y $ in a neighborhood of $ 0 $, 
 the solution of  \eqref{eq4} is \eqref{explivc}. 
In conclusion, by \eqref{eq1}, \eqref{oftheform}, \eqref{invbot} we have proved that 
the inverse of $ \Phi $ has the form \eqref{formainv}-\eqref{explivc}. This proves Item ($iii$).
It remains to prove Item ($i$).
\\[1mm]
{\bf The map $ \Phi $ defined in \eqref{Phiorto} is  symplectic. }   
  We find convenient to read $ \Phi $ in the coordinates $ (\a_n, \b_n)_{n \geq 1}$ 
  defined by 
  $ \xi (\theta ) = \sum_{n \geq 1} \a_n \tc_n (\theta ) + \b_n \ts_n (\theta) $   in 
  \eqref{a2b2n}
  and to prove  that $ \Phi $ preserves all the fundamental Poisson 
  brackets (see \eqref{Poibrak})
  $$ 
  \{\a_n, \a_{n'} \} =  \{\b_n, \b_{n'} \}  = 0 \, , \quad  \{\a_n, \b_{n'} \} = p_n \delta_{n,n'} \, , 
  \quad \forall n, n' \in \N \, . 
  $$ 
  We remind that 
   in the  coordinates $ (\a_n, \b_n)_{n \geq 1}$ 
 the symplectic form reads as in \eqref{simpl2}. 
Writing in coordinates the flow 
\be\label{flowinco}
  \Phi^t_{\cal J} (\xi ) = {\mathop \sum}_{n \geq 1} a_n (t, \xi ) \tc_n  +  b_n (t, \xi) \ts_n 
\ee
with 
\be\label{coflowJ}
\begin{aligned} 
& a_2 (t, \xi ) = \tfrac12 
(\Phi^t_{\cal J}(\xi), \tc_2) \, , \ b_2 (t, \xi ) = \tfrac12 (\Phi^t_{\cal J}(\xi), \tc_2) \, , \notag \\
& 
a_n (t, \xi ) = (\Phi^t_{\cal J}(\xi), \tc_n) \, , \ b_n (t, \xi ) = (\Phi^t_{\cal J}(\xi), \ts_n) \, ,  \
\forall n \neq 2 \, ,    
 \end{aligned}
 \ee
 the map $ \Phi $ defined in \eqref{Phiorto} reads  as the map
\be\label{cooPhi}
(\a_n, \b_n)_{n \geq 1} \stackrel{\Phi} 
\mapsto \Big( a_1 (\bar t(\xi),\xi), b_1 (\bar t(\xi),\xi), {\cal J}(\xi), \bar t(\xi),
\ldots, a_n (\bar t(\xi),\xi), b_n (\bar t(\xi),\xi),  \ldots 
 \Big)
 \ee
 that we still denote for simplicity $ \Phi $. 
 Substituting  $\widehat \xi=X_{{\cal J}} (\xi) $ in \eqref{sympjt-1} gives
 \begin{equation}\label{pbjt0}
   d \bar t (\xi)[X_{{\cal J}} (\xi)] 
     = - 
   \frac{\big(  \ts_2, d \Phi^{ t }_{{\cal J}}(\xi)_{|t = \overline t (\xi)}[X_{{\cal J}} (\xi)] \big)}
   {\big(  \ts_2,  (\pa_t \Phi^{ t }_{{\cal J}})_{|t = \overline t (\xi)} (\xi) \big) } = - 1 
 \end{equation}
 because 
\be\label{idenfg} 
d \Phi^{ t }_{{\cal J}}(\xi) [X_{{\cal J}} (\xi)]  = 
(\pa_t \Phi^{ t }_{{\cal J}}) (\xi)   \, , \quad \forall t  \, . 
\ee
Indeed, differentiating \eqref{flow-J} with respect to $ t $, 
we see that $ g(t) := (\pa_t \Phi^{ t }_{{\cal J}}(\xi)) =  X_{\cal J} ( \Phi^{ t }_{{\cal J}} (\xi) )$,  
solves 
$$
 \pa_t g(t) = d X_{{\cal J}} (\Phi^{ t }_{{\cal J}}(\xi))  g(t) 
 \, , \quad 
g(0) = X_{\cal J}(\xi) \, . 
$$
Moreover, by \eqref{difJ},  we get that 
$ f (t) := d \Phi^{ t }_{{\cal J}}(\xi) [X_{{\cal J}} (\xi)] $ solves
$$
 \pa_t f(t)  = d X_{\cal J} (\Phi^{ t }_{{\cal J}}(\xi) ) \, f(t) \, , \quad 
f(0)= X_{{\cal J}} (\xi)  \, . 
$$
By the unicity of the solutions of the Cauchy problem 
we deduce that $ f(t) = g(t)$ which is \eqref{idenfg}.
By \eqref{pbjt0} it follows that  the Poisson bracket between $ {\cal J} (\xi)$ and  $\bar t (\xi) $ is 
$$
 \{ {\cal J} (\xi),  \bar t (\xi) \} 
   =  - \{ \bar t (\xi), {\cal J} (\xi) \}  
  \stackrel{\eqref{Poib}} =  - d\bar t (\xi)[X_{{\cal J}} (\xi)]
  \stackrel{\eqref{pbjt0}}  =1
\, . 
$$
We now prove that the infinitely many Poisson brackets
\be\label{infpb} 
\begin{aligned} 
& \big\{ a_n(\bar t (\xi), \xi), {\cal J} (\xi) \big\} = 
\big\{  b_n(\bar t (\xi), \xi), {\cal J} (\xi) \big\} =  0 \, , \\
& \big\{ a_n(\bar t (\xi), \xi), {\bar t} (\xi) \big\} = 
\big\{  b_n(\bar t (\xi), \xi), {\bar t} (\xi) \big\} =  0 \, ,
\end{aligned}
\quad \forall  n \neq 2 \, . 
\ee
Defining  the vectorial Poisson bracket 
$ \big\{ \Phi_2^\bot (\xi), F(\xi)  \big\} := d \Phi_2^\bot (\xi)[X_{F} (\xi)] $, 
where the scalar function $ F(\xi) \in \{ {\cal J}(\xi), {\bar t}(\xi) \} $, 
and
\be\label{phi2bot}
 \Phi_2^\bot (\xi ) =  \Pi_2^\bot \Phi^{\bar t(\xi)}_{\cal J} (\xi ) 
 \stackrel{\eqref{flowinco}} = {\mathop \sum}_{n \neq 2} 
 a_n (\bar t (\xi), \xi) \tc_n  +  b_n ( \bar t (\xi), \xi) \ts_n \, ,  
\ee 
the Poisson brackets \eqref{infpb}  amount in compact form  to 
\be\label{pricompf}
\begin{aligned}
& \big\{ \Phi_2^\bot (\xi), {\cal J}(\xi)  \big\}  =  
{\mathop \sum}_{n \neq 2} 
 \big\{ a_n (\bar t (\xi), \xi), {\cal J}(\xi) \big\} \tc_n  +  
 \big\{ b_n ( \bar t (\xi), \xi), {\cal J}(\xi)  \big\}  \ts_n 
 =  0\, ,  \\ 
&  \big\{ \Phi_2^\bot (\xi), {\bar t}(\xi)  \big\}  =  
{\mathop \sum}_{n \neq 2} 
 \big\{ a_n (\bar t (\xi), \xi), {\bar t}(\xi) \big\} \tc_n  +  
 \big\{ b_n ( \bar t (\xi), \xi), {\bar t}(\xi)  \big\}  \ts_n 
 =  0  \, . 
 \end{aligned}
\ee
We have
 \begin{align}
 \{ \Phi_2^\bot (\xi) , {\cal J} (\xi)\} 
& 
 = d \Phi_2^\bot (\xi)[X_{{\cal J}} (\xi)] \notag
\\
& \stackrel{\eqref{phi2bot}}=  \Pi_2^\bot (\pa_t \Phi^{ t }_{{\cal J}})_{|t = \overline t (\xi)} (\xi)  \, d \bar t (\xi) [X_{{\cal J}} (\xi)]+
 \Pi_2^\bot d \Phi^{ t }_{{\cal J}}(\xi)_{|t = \overline t (\xi)}[X_{{\cal J}} (\xi)]   \notag
\\
&  \stackrel{\eqref{pbjt0},\eqref{idenfg} } =- 
 \Pi_2^\bot (\pa_t \Phi^{ t }_{{\cal J}})_{|t = \overline t (\xi)} (\xi)  \, 
+   \Pi_2^\bot (\pa_t \Phi^{ t }_{{\cal J}})_{|t = \overline t (\xi)} (\xi) =0 \, ,	\notag 
  \end{align}
 proving  the first identity in \eqref{pricompf}. 
We now consider
 \begin{align}
 \{ \Phi_2^\bot (\xi) , {\bar t} (\xi)\} 
& 
 := d \Phi_2^\bot (\xi)[X_{{\bar t}} (\xi)] \notag
\\
& \stackrel{\eqref{phi2bot}}=  \Pi_2^\bot (\pa_t \Phi^{ t }_{{\cal J}})_{|t = \overline t (\xi)} (\xi)  \, d \bar t (\xi) [X_{{\bar t }} (\xi)]+
 \Pi_2^\bot d \Phi^{ t }_{{\cal J}}(\xi)_{|t = \overline t (\xi)}[X_{{\bar t}} (\xi)]   \notag
\\
&  \stackrel{\eqref{pbjt0}} = 
 \Pi_2^\bot (\pa_t \Phi^{ t }_{{\cal J}})_{|t = \overline t (\xi)} (\xi)  \, 
\{ \bar t, \bar t \, \} +
 \Pi_2^\bot  d \Phi^{ t }_{{\cal J}}(\xi)_{|t = \overline t (\xi)}[ \pa_\theta \nabla {\bar t} (\xi)]  \notag  \\
& \stackrel{\eqref{uguaJJ2}} =  \Pi_2^\bot   \Phi^{\overline t (\xi)}_{{\cal J}_2} \big[ 
\pa_\theta \nabla {\bar t} (\xi) \big] \, . \label{psibotbart} 
  \end{align}
Now,  by \eqref{gradbart} and \eqref{sympjt-1} 
$$
   \nabla \bar t (\xi) 
     = - 
   \frac{  \big[     \Phi^{\overline t (\xi)}_{{\cal J}_2}  \big]^*\ts_2
   }
   {\d (\xi)}   \qquad {\rm with} \qquad 
\d(\xi) :=  \big(  \ts_2,  (\pa_t \Phi^{ t }_{{\cal J}})_{|t = \overline t (\xi)} (\xi) \big)  \, , 
$$
and therefore, by \eqref{psibotbart}, 
\begin{align}\label{passsimp}
 \{ \Phi_2^\bot (\xi) , {\bar t} (\xi)\} 
&   = 
- \delta^{-1}(\xi) \,  \Pi_2^\bot  
\underbrace{ \Phi^{\overline t (\xi)}_{{\cal J}_2}  \, 
\pa_\theta \,   
 \big[  \Phi^{\overline t (\xi)}_{{\cal J}_2} }_{= \pa_\theta} \big]^*\ts_2 \\
 & = 
 - \delta^{-1}(\xi) \,  \Pi_2^\bot  
 \pa_\theta \, \ts_2   \stackrel{\eqref{defcnsn}} =  
 - 2 \delta^{-1}(\xi) \,  \Pi_2^\bot   \tc_2   = 0 \notag 
 \end{align}
having used in \eqref{passsimp} that 
$  \Phi^{\overline t (\xi)}_{{\cal J}_2}  $ is symplectic,  thus  also its inverse $  
 \Phi^{- \overline t (\xi)}_{{\cal J}_2}   $, which means
$$
 [ \Phi^{\overline t (\xi)}_{{\cal J}_2} ]^{-*} \pa_\theta^{-1} [ \Phi^{\overline t (\xi)}_{{\cal J}_2} ]^{-1} = \pa_\theta^{-1} 
$$
and then, taking the inverse\footnote{This relation
can be also checked by the explicit expressions in Lemma \ref{flowJ2old}.}, 
$   \Phi^{\overline t (\xi)}_{{\cal J}_2}  \, \pa_\theta \, [ \Phi^{\overline t (\xi)}_{{\cal J}_2} ]^* = \pa_\theta  $.
This completes the proof of \eqref{pricompf}. 
\\[1mm] \indent
We now consider the other fundamental Poisson brackets. 
We find convenient to compute the Poisson brackets \eqref{Poibrak}
of the functions $ a_n(\bar t (\xi (\a,\b)), \xi (\a,\b)) $, 
$ b_n(\bar t (\xi (\a,\b)), \xi (\a,\b)) $ in the variables 
$ (\a,\b) = (\a_n,\b_n)_{n \geq 1} $.
For any $ n, n' \neq 2 $,  we get, by \eqref{Poibrak}, 
 $$
 \begin{aligned}
 \big\{ b_n(\bar t (\xi), \xi), b_{n'}(\bar t (\xi), \xi) \big\} 
 & = 
 \sum_{k \geq 1} p_k 
\Big( \frac{\pa b_n(t, \cdot)}{\pa \a_k} + \frac{\pa b_n}{\pa t} 
 \frac{\pa \bar t}{\pa \a_k} \Big)  
 \Big( \frac{\pa b_{n'}(t, \cdot)}{\pa \b_k} + \frac{\pa b_{n'}}{\pa t} 
 \frac{\pa \bar t}{\pa \b_k} \Big)_{|t=\bar t(\xi)}  \\
 & \ \, \quad - p_k
 \Big( \frac{\pa b_n(t, \cdot)}{\pa \b_k} + \frac{\pa b_n}{\pa t} 
 \frac{\pa \bar t}{\pa \b_k} \Big)  
 \Big( \frac{\pa b_{n'}(t, \cdot)}{\pa \a_k} + \frac{\pa b_{n'}}{\pa t} 
 \frac{\pa \bar t}{\pa \a_k} \Big)_{|t=\bar t(\xi)} 
 \end{aligned}
 $$
 where 
 $ \frac{\pa b_n(t, \cdot)}{\pa \a_k} = \frac{\pa \, b_n (t, \xi(\a,\beta))}{\pa \a_k} $,
  $ \frac{\pa {\bar t}}{\pa \a_k} = \frac{\pa \, {\bar t} ( \xi(\a,\beta))}{\pa \a_k} $ 
 and similarly for the partial derivatives with respect to $ \beta_k $. 
Then
 \begin{align}
  \big\{ b_n(\bar t (\xi), \xi), b_{n'}(\bar t (\xi), \xi) \big\}  & = 
 \{ b_n(t, \cdot), b_{n'}(t, \cdot) \}_{|t= \bar t(\xi)} \notag \\
 & + 
 \sum_{k \geq 1} p_k \frac{\pa b_n}{\pa t} 
\Big(  \frac{\pa \bar t}{\pa \a_k} \frac{\pa b_{n'}(t, \cdot)}{\pa \b_k} -
 \frac{\pa \bar t}{\pa \b_k} \frac{\pa b_{n'}(t, \cdot)}{\pa \a_k} \Big)_{|t= \bar t(\xi)} \notag \\ 
 & \ \quad
+ p_k \frac{\pa b_{n'}}{\pa t} 
\Big(  \frac{\pa \bar t}{\pa \b_k} \frac{\pa b_{n}(t, \cdot)}{\pa \a_k} -
 \frac{\pa \bar t}{\pa \a_k} \frac{\pa b_{n}(t, \cdot)}{\pa \b_k} \Big)_{|t= \bar t(\xi)}  \, . \label{ladif}
 \end{align}
 By \eqref{sympjt-1} we have
\begin{align}\label{pban}
\frac{\pa \bar t}{\pa \b_k}  =  d \bar t (\xi)[\ts_k] 
 &  = - 
   \frac{d \big(  \ts_2, \Phi^{ t }_{{\cal J}}(\xi) [\ts_k ] \big)_{|t = 
   \overline t (\xi)}   }{\big( \ts_2, X_{\cal J}(\Phi^{\overline t (\xi) }_{{\cal J}}(\xi) ) \big) } \notag \\ 
&   \stackrel{\eqref{coflowJ}}  
    = - \frac12
   \frac{d \, b_2(t, \cdot)_{|t = 
   \overline t (\xi)}  [\ts_k ] }{\d (\xi)}
   = - \frac12 
    \frac{\pa_{\b_k} b_2(t, \cdot)}{\d (\xi)} 
\end{align}
where 
$ \d (\xi) := \big( \ts_2, X_{\cal J}(\Phi^{\overline t (\xi) }_{{\cal J}}(\xi) ) \big) $ which is different from zero since it is close to $ \big( \ts_2, - \ts_2 \big)$ for $ \xi $ small. 
 Similarly
\be\label{paan}
\frac{\pa \bar t}{\pa \a_k}  
   = -  
    \frac12 \frac{\pa_{\a_k} b_2(t, \cdot)}{\d (\xi)} \, . 
\ee
By \eqref{ladif}, \eqref{pban}  and \eqref{paan}, we deduce that
\begin{align}\label{quasifinpoi}
\big\{ b_n(\bar t (\xi), \xi), b_{n'}(\bar t (\xi), \xi) \big\} 
& = 
 \{ b_n(t, \cdot), b_{n'}(t, \cdot) \}_{|t= {\bar t}(\xi)} \\ 
 & -  
 \frac{1}{\delta(\xi)}  \frac{\pa b_n}{\pa t} \{ b_2(t, \cdot), b_{n'}(t, \cdot) \}_{|t= {\bar t}(\xi)} -
  \frac{1}{\delta(\xi)} \frac{\pa b_{n'}}{\pa t} \{  b_{n}(t, \cdot), b_2(t, \cdot) \}_{|t= {\bar t}(\xi)}  \, . \notag 
  \end{align}
 For any $ t $, the flow map $ \xi \mapsto \Phi^t_{\cal J} (\xi) $ in \eqref{flowinco} is symplectic, and thus  
the map 
$$ 
(\a,\beta) = (\a_n,\b_n)_{n \geq 1} \mapsto 
\Big( a_n(t, \xi (\a,\beta)), b_n(t, \xi (\a,\beta)) \Big)_{n \geq 1}
$$
is symplectic with respect to the $ 2 $-form \eqref{simpl2}. Then 
$ \{ b_n (t, \cdot), b_{n'} (t, \cdot) \} = 0 $, for any $ n, n' $, 
and we deduce by \eqref{quasifinpoi} that
$ \big\{ b_n(\bar t (\xi), \xi), b_{n'}(\bar t (\xi), \xi) \big\}  = 0 $, $  \forall n,n' \neq 2 $. 
With similar arguments one checks that also the following identities the between Poisson brackets 
$$
\begin{aligned}
& \big\{ a_n(\bar t (\xi), \xi), a_{n'}(\bar t (\xi), \xi) \big\}  = 0 \, , \quad \forall n,n' \neq 2 \, , \\
& \big\{ a_n(\bar t (\xi), \xi), b_{n'}(\bar t (\xi), \xi) \big\}  = n \delta_{n,n'} 
\, , \quad \forall n,n' \neq 2 \, .
\end{aligned}
$$
This concludes the proof that the map $ \Phi $   in \eqref{Phiorto} 
is symplectic,
and of Theorem \ref{thm:FM}.
 \end{pf}

 We now write  the Hamiltonian 
\be\label{HEJ}
H(\xi) := H_{\Omega_\gamma}(\xi) = -\tfrac12 E ( \xi ) +  \tfrac{\Omega_\gamma}{2} J  (\xi) =  H_L (\xi) + H_{\geq 3} (\xi) \, , 
\ee
where $ \Omega_\gamma $ is  defined in \eqref{choiceOmega},  the
quadratic Hamiltonian $ H_L $ is in \eqref{defQH}
 and the Hamiltonian  $ H_{\geq 3} $ comprises the higher order cubic terms,
in the symplectic coordinates  introduced in Theorem \ref{thm:FM}.

 \begin{corollary}  {\bf (Hamiltonian in new symplectic variable)}\label{cor:new}
In the  symplectic variable
 \be\label{phisymdiffeo}
 \begin{aligned}
 \wtilde \xi 
 & :=   \wtilde \a_2 \tc_2  +\wtilde \b_2   \ts_2  + 
  \wtilde u  = \Phi (\xi) =
 {\cal J}(\xi) \tc_2  + \overline t (\xi)  \ts_2  +  \Phi_2^\bot (\xi) \, , 
 \end{aligned} 
\ee
 the Hamiltonian
\be\label{defKPhi}
  K ( \wtilde \xi)   := H (\Phi^{-1}( \wtilde \xi))
\ee
  is independent of  $ \wtilde \b_2 $, i.e.  $  \pa_{\wtilde \beta_2}  K ( \wtilde \xi) =0 $.
The Hamiltonian
  \be\label{tildeKeq} 
K ( \wtilde \xi) = K (\wtilde \a_2 \tc_2  + \wtilde u  )  =:   \cK ( \wtilde \a_2, \wtilde u )  
 \ee
can be expanded as  
  \begin{align}  \label{fre:uncha} 
 K ( \wtilde \xi)  
 & 
  =   \frac{1}{2} ( {\bf \Omega} (\gamma) 
{\wtilde \xi}, {\wtilde \xi} ) 
+K_{\geq 3}(\wtilde \a_2 \tc_2  + \wtilde u  ) 
  =  - 
\frac{\Omega_2}{2} {\wtilde \a}_2^2  + \frac{1}{2} ( {\bf \Omega} (\gamma) 
{\wtilde u} , {\wtilde u} )+   \cK_{\geq 3}( \wtilde \a_2, \wtilde u )
\end{align}
 where $ K_{\geq 3} (\wtilde \xi )$, resp. $\cK_{\geq 3} ( \wtilde \a_2, \wtilde u ) $, 
 comprises the higher order cubic terms of the Hamiltonian $ K (\wtilde \xi ) $, resp. 
$   \cK( \wtilde \a_2, \wtilde u ) $.

  Moreover the Hamiltonian $ K $ is reversible, i.e.
 $ K \circ {\cal S} = K $ where $ {\cal S} $ is the involution  defined in \eqref{S-invo}, and
  $  \cK( \wtilde \a_2, {\cal S} \wtilde u ) =  \cK( \wtilde \a_2, \wtilde u ) $.  
  The  Hamiltonian system 
      \be\label{HSK}
\pa_t \wtilde \xi = X_{K} (\wtilde \xi) = \pa_\theta \nabla_{\wtilde \xi} K (\wtilde \xi \, ) \, ,
\ee
reads, in the variables $ ( \wtilde \a_2, \wtilde \beta_2,  \wtilde u ) $,  as  
    \be\label{HSKt0}
 \pa_t \wtilde \alpha_2 = 0 \, , \quad   \pa_t \wtilde \beta_2 =  - \pa_{\wtilde \alpha_2} 
   {\cal K} ( \wtilde \a_2,  \wtilde u ) \, , \quad 
 \pa_t \wtilde  u = \pa_\theta \nabla_{\wtilde u} {\cal K} ( \wtilde \a_2,  \wtilde u  \, ) \, ,
 \ee
which possesses $ \wtilde \alpha_2  $ as a prime integral. 
   \end{corollary}

\begin{pf}
Since $ {\cal J}(\xi) $ is a prime integral of $ H (\xi ) $ then  
$ \wtilde \a_2 = {\cal J} (\Phi^{-1} (\wtilde \xi) ) $ 
is a prime integral of $  K ( \wtilde \xi)  $ and 
$ 0 = \big\{ \wtilde \a_2,  K (\wtilde \xi ) \big\}  = 
 \pa_{\wtilde \b_2}  K ( \wtilde \xi ) $, by \eqref{Poibrak} and since 
$ \Phi $ is symplectic. 
The expansion \eqref{fre:uncha}   follows since  $ \Phi $ 
satisfies \eqref{dif0}, and recalling \eqref{defQH}, \eqref{Hanbn}. 
The map  $ \Phi $ is reversibility preserving, 
as well as 
$ \Phi^{-1} $, and thus $ K = H \circ \Phi^{-1} $ is reversible as $ H $ is. 
Finally \eqref{HSKt0} follows recalling that, in the symplectic variables 
$ (\wtilde \a_2, \wtilde \b_2, \wtilde u ) $, 
the symplectic $ 2 $-form  is given  in \eqref{simpl2}, and 
\eqref{XVS}. 
\end{pf}

\noindent 
{\bf Symplectic reduction of the angular momentum.}
By \eqref{HSKt0} we have  $ \wtilde \a_2 (t) = \underline{\cal J} $ is constant in time, 
the third equation reduces to  
\be\label{HSredu}
  \pa_t \wtilde  u = \pa_\theta 
  \nabla_{\wtilde u}  \cK (\, \underline{\cal J} , \wtilde u \, ) \, , 
\ee
where, by \eqref{tildeKeq},   
\be\label{nablatuK}  
\nabla_{\wtilde u}  \cK ( \wtilde \a_2, \wtilde u \, ) = \Pi_2^\bot  
  (\nabla_{\wtilde \xi} K) ( \wtilde \a_2 \tc_2  + \wtilde u ) \, ,  
 \ee   
and the evolution of the decoupled  $ \wtilde \beta_2 (t) $ coordinate is obtained integrating
the scalar equation 
$   \dot {\wtilde \beta}_2 (t) =  
   - (\pa_{\wtilde \alpha_2}  \cK) ( \, \underline{\cal J}, \wtilde u ( t ) ) $. 
   
In the next sections we look for 
quasi-periodic solutions of the reduced equation \eqref{HSredu}.
Let us describe in detail which quasi-periodic solutions of  the 
original equation \eqref{Evera}
we have constructed. 
 If $ \wtilde {\underline{u}} (\omega t) $ 
  with $ \wtilde {\underline{u}} (\vphi) $, $ \vphi \in \T^{|\ST|} $,   is a quasi-periodic solution of \eqref{HSredu},
 then $ (\, \underline{\cal J}, \underline{\widetilde\beta_2} (\omega t), 
 \wtilde {\underline{u}} (\omega t))  $,  where
\be\label{defmubeta2}
\begin{aligned}
&  \underline{\widetilde\beta_2} (\varphi) := 
- (\omega \cdot \pa_\vphi)^{-1}
\big[ (\pa_{\wtilde \alpha_2}  \cK) 
(\, \underline{\cal J},   {\wtilde {\underline{u}}(\vphi)})
- \underline{\mu}  \,  \big] \, , \\
&    \underline{\mu}   
:= 
\big\langle (\pa_{\wtilde \alpha_2}  \cK) (\, \underline{\cal J}, \wtilde {\underline{u}}(\cdot )) 
\big\rangle :=   \frac{1}{(2 \pi)^{|\ST|}} \int_{\T^{|\ST|}} \, 
(\pa_{\wtilde \alpha_2}  \cK) 
(\, \underline{\cal J},   {\wtilde {\underline{u}}(\vphi)}) \, d \vphi\, , 
\end{aligned} 
\ee
is a quasi-periodic solution of the Hamiltonian system
$$
 \pa_t \wtilde \alpha_2 = 0 \, , \quad   \pa_t \wtilde \beta_2 =  
 - \pa_{\wtilde \alpha_2} 
  \big( \cK (\wtilde \a_2, \wtilde u )  - \underline{\mu} \wtilde \a_2 \big) \, , \quad 
 \pa_t \wtilde  u = \pa_\theta \nabla_{\wtilde u} \cK (\wtilde \a_2, \wtilde u \, ) \, ,
$$
generated by the Hamiltonian 
$  \cK (\wtilde \a_2, \wtilde u )  - \underline{\mu} \wtilde \a_2  $
(here $ \underline{\mu}$ is a fixed constant).
Equivalently 
\be\label{solpr1}
\wtilde  {\underline{\xi}}(\omega t ) 
:= 
 \underline{\cal J} \tc_2 +  \underline{ \widetilde\beta_2} (\omega t) \ts_2 + 
\wtilde  {\underline{u}} (\omega t)
\ee
solves the Hamiltonian system 
\be\label{linxi1}
\pa_t \wtilde \xi = X_{K -  \underline{\mu} \wtilde \alpha_2} (\wtilde \xi) 
= X_K (\wtilde \xi) + \underline{\mu}  \ts_2 \, , \quad \wtilde \alpha_2 \stackrel{\eqref{newa2b2}} = \tfrac12 ( \wtilde \xi , \tc_2) \, .
\ee 
Going back in the variable $ \xi $ via the symplectic  diffeomorphism 
defined in \eqref{phisymdiffeo}, 
$ \xi = \Phi^{-1} (\wtilde \xi ) $, 
we have, since  $ \wtilde \alpha_2 = {\cal J}( \Phi^{-1} (\wtilde \xi)) $,  that
$ H  -  {\underline \mu} \, {\cal J}  = 
(K -  {\underline \mu} \, \wtilde \alpha_2) \circ \Phi  $, and therefore 
\be\label{solpr2}
\underline {\xi}(\omega t ) := \Phi^{-1} (\wtilde{\underline {\xi} }(\omega t) ) 
\ee
is a  quasi-periodic solution  of the Hamiltonian system 
\be\label{linxi2}
\pa_t \xi = X_{H -  {\underline \mu} {\cal J}} ( \xi )   \, . 
\ee
By \eqref{HEJ} and \eqref{Jsvi} we have that, up to a constant,  
\be\label{Homfin}
H - {\underline \mu} {\cal J} =-  \frac12 E  +  \frac12
\Big( \Omega_\gamma - {\underline \mu}  \frac{2\sqrt{2}}{\sqrt{\pi}(\g -\g^{-1})} \Big) J 
\stackrel{\eqref{Ham:main}} = H_{\Omega} \, , \quad
\Omega = \Omega_\gamma -    {\underline \mu} \frac{2\sqrt{2}}{\sqrt{\pi}(\g -\g^{-1})} 
 \, ,
\ee
and therefore, in light of Proposition \ref{prop:Ham}, 
 $ \underline {\xi} ( \omega t ) $ is a time quasi-periodic vortex patch 
solution of the equation \eqref{Evera},  in a rotating frame with angular velocity 
$ \Omega $.

\smallskip

 For the
Nash-Moser construction of quasi-periodic solutions,  we shall provide estimates
of the transformed Hamiltonian vector field 
  \be\label{transformedHS}
  X_K (\wtilde \xi ) 
   = (d \Phi) ( \Phi^{-1} (\wtilde \xi)) X_{H} ( \Phi^{-1}(\wtilde \xi) )  
  \ee
  in Section \ref{sec:regularity}.
In view of the perturbative construction of the quasi-periodic solutions we decompose
the Hamiltonian vector field generated by the Hamiltonian 
$ H(\xi) $  in   \eqref{HEJ} as
\be\label{XHsv}
X_H (\xi) = \pa_\theta {\bf \Omega}(\gamma) \xi   + X_{H_{\geq 3}} (\xi ) \, . 
\ee
The transformed vector field  \eqref{transformedHS}  can be written as 
$ X_K (\wtilde \xi) = \pa_\theta {\bf \Omega}(\gamma) \wtilde \xi   + X_{K_{\geq 3}} (\wtilde \xi ) $
where 
 \begin{align}
  X_{ K_{\geq 3}}( \wtilde \xi )  
&  \stackrel{\eqref{transformedHS}} = d\Phi(\Phi^{-1}( \wtilde \xi ))X_{H}(\Phi^{-1}( \wtilde \xi ) )-\partial_\theta {\bf \Omega}(\gamma)   \wtilde \xi \notag \\
 \nonumber  & 
 \stackrel{\eqref{XHsv},\eqref{dif0}} = 
 d\Phi(\Phi^{-1}( \wtilde \xi))X_{H_{\geq 3}}(\Phi^{-1}( \wtilde \xi) )
+ d\Phi(\Phi^{-1}( \wtilde \xi))\partial_\theta {\bf \Omega}(\gamma) \big[\Phi^{-1}( \wtilde \xi)-  \wtilde \xi \big]\\ &\qquad \qquad + \big[d\Phi(\Phi^{-1}( \wtilde \xi))-d\Phi(0) \big]\partial_\theta {\bf \Omega}(\gamma) \wtilde \xi \, . \label{XKXP}
\end{align}
{\bf Formulas for the differential  $ d \Phi $.} 
For the sequel we provide some expressions of $ d \Phi (\xi) $ and 
\be\label{dphi-1}
\big[ d \Phi (\xi)  \big]^{-1} = d \Phi^{-1} (\eta)_{|\eta = \Phi (\xi)} \, . 
\ee
\begin{lemma} \label{lem:dphi}
We have that 
\begin{align}\label{laprim}
& \Pi_2^\bot d \Phi (\xi) =  
 \Pi_2^\bot \Phi^{  \overline t (\xi) }_{{\cal J}_2}   + {\mathscr R}_1 \\
 & \big[ d \Phi (\xi) \big]^{-1}  \Pi_{2}^\bot 
=   \Phi^{- \bar t (\xi)}_{{\cal J}_2} \Pi_{2}^\bot+   {\mathscr R}_2 \label{dphi-1bot}
\end{align}
 where $  {\mathscr R}_1 $ is  the finite rank operators 
 \be\label{defR1}
 {\mathscr R}_1 := 
(g_1 , \cdot )  \chi_1   \qquad \text{with} \qquad 
g_1 := (\nabla \bar t \, )  (\xi)  \, , \ 
\chi_1 := \Pi_2^\bot X_{{\cal J}_2} (\Phi_{\cal J}^{\overline t (\xi)} (\xi))  
\ee
with  $(\nabla \bar t \, )  (\xi) $ is in \eqref{gradbart},  
and  ${\mathscr R}_2 $ is the finite rank operator,  acting on $ H_2^\bot $, 
\be\label{pritra21}
{\mathscr R}_2   = 
( g_2, \cdot )  \chi_2
\quad \text{with} \quad g_2 := \nabla_{\eta_\bot} v_\tc (\eta)_{|\eta = \Phi (\xi)}
 \, , \quad \chi_2 := \Phi_{\cal J}^{- \bar t (\xi)}   \tc_2 \, , 
\ee
where $ \nabla v_\tc (\eta) $ is computed by \eqref{explivc}.
\end{lemma}

\begin{pf}
The identity \eqref{laprim}-\eqref{defR1} follows by 
 \eqref{diffePhi}. 
 Let us prove \eqref{dphi-1bot}. Let 
 $  \eta = \eta_\tc \tc_2 + \eta_\ts \ts_2 + \eta_\bot $.
Differentiating formula \eqref{formainv}  along the subspace $ H_2^\bot $
and recalling the expression of   $ \Phi_{\cal J}^t $  in \eqref{flowaff}, we obtain
for any $ \widehat \eta_\perp \in H_2^\bot $, 
\be\label{lasopsi}
d \Phi^{-1} (\eta)[\widehat \eta_\bot] =  \Phi_{{\cal J}_2}^{- \eta_\ts} 
( \tc_2 \, d v_\tc (\eta) [\widehat \eta_\bot] + \widehat \eta_\bot) \, . 
\ee
By \eqref{dphi-1}, \eqref{lasopsi}, and since $ \eta = \Phi(\xi )$ implies that 
$ \eta_\ts  = \bar t (\xi) $,  we get, for any $ \widehat \eta \in H_2^\bot $, 
$$
\big[ d \Phi (\xi)  \big]^{-1} 
 \widehat \eta_\bot = 
d \Phi^{-1} (\eta)_{|\eta = \Phi (\xi)} 
\widehat \eta_\bot  = 
\Phi_{{\cal J}_2}^{- \bar t (\xi)} 
 \big( (\nabla v_\tc (\eta)_{|\eta = \Phi (\xi)}, \widehat \eta_\bot)  \tc_2  + \widehat \eta_\bot \big) 
$$
which proves 
\eqref{dphi-1bot}, \eqref{pritra21}. 
\end{pf}

We finally consider the  adjoint of the operator 
$  \big[ d \Phi (\xi)\big]^{-1}  \Pi_{2}^\bot  $ in  \eqref{dphi-1bot}.

\begin{lemma}\label{lem:adj}
The  adjoint of the operator
$ L (\xi) := \big[ d \Phi (\xi)\big]^{-1}  \Pi_{2}^\bot 
=   \Phi^{- \bar t (\xi)}_{{\cal J}_2} \Pi_{2}^\bot+   {\mathscr R}_2   $ 
in  \eqref{dphi-1bot} 
is  
$$
L (\xi)^* = \Pi_{2}^\bot (\Phi^{- \bar t (\xi)}_{{\cal J}_2})^*  +  {\mathscr R}_3
$$
where $ {\mathscr R}_3 :=  {\mathscr R}_2^* =  (  \chi_2, \cdot ) g_2 $ 
is finite rank. 
\end{lemma} 

\begin{pf}
The adjoint of a finite rank operator 
$ R = (  \cdot ,  g) \chi $ is $ R^* = (  \chi ,  \cdot ) g $ because
$ (Rh, k) = $ $ ( h ,  g) (\chi, k) = $ 
$ \big( h ,   g ( \chi, k) \big) = $ $ (h, R^* (k) ) $, 
and the lemma follows taking the adjoint of \eqref{dphi-1bot}.  
\end{pf}

\section{Construction of quasi-periodic vortex patches}\label{sec:NMT}

Instead of looking for solutions of the Hamiltonian PDE 
\eqref{HSK} 
in a shrinking neighborhood of $ \wtilde \xi  = 0 $ 
 it is a convenient devise 
to perform the rescaling
$ \wtilde \xi  = \e \breve \xi  $, for $ \e > 0 $ small, 
and consider the new variable  $ \breve \xi $ to have size $  O(1) $ in $ \e $. 
In  the rescaled variables
\be\label{rescavar}
\breve \xi = \e^{-1} \wtilde \xi   \, , \quad {\rm i.e.} \quad
(\breve \alpha_2, \breve \beta_2, \breve u) = \e^{-1} 
( \wtilde \alpha_2, \wtilde \beta_2,  \wtilde u) \, , 
\ee
 the Hamiltonian PDE 
\eqref{HSK}, i.e. \eqref{HSKt0}, 
transforms  into
\be\label{HSKt0res}
\pa_t \breve \xi  
=  \pa_\theta (\nabla_{\breve \xi} \breve K) ( \breve \xi \, ) \, , \ 
\text{i.e.}  \ 
 \pa_t \breve \alpha_2 = 0 \, , \  
 \begin{aligned} 
  & \pa_t \breve \beta_2 =  - \e^{-1} 
 (\pa_{\wtilde \alpha_2} 
   {\cal K}) ( \e \breve \a_2, \e \breve u ) = - 
   (\pa_{\breve \alpha_2} 
   \breve {\cal K}) (  \breve \a_2, \breve u ) 
   \, , \\
 & \pa_t \breve  u = 
 \e^{-1} \pa_\theta \, (\nabla_{\wtilde u}  {\cal K}) ( \e \breve \a_2, \e \breve u  \, ) = 
  \pa_\theta \, (\nabla_{\breve u}  \breve {\cal K}) (  \breve \a_2, \breve u  \, )
 \, ,
 \end{aligned}
 \ee
generated by the rescaled Hamiltonian
\be\label{HS:mu-breve}
\breve K (\breve \xi) := \e^{-2} K( \e  \breve \xi) \, , \quad i.e. \ \ 
\breve {\cal K} (\breve \a_2, \breve u) := 
\e^{-2} {\cal K}( \e \breve \a_2, \e \breve u) \, .
\ee
We fix the value of the prime integral 
$ \breve \alpha_2 = {\cal J}_0 $ and then the last equation in \eqref{HSKt0res}
reads 
\be\label{HS:mu}
\pa_t \breve u 
  = \pa_\theta \nabla_{\breve  u} \breve {\cal K}  ({\cal J}_0 , \breve u \, ) 
 =
  \pa_\theta {\bf \Omega}(\g) \breve  u + 
 \e \pa_\theta \nabla_{\breve  u} {\cal P} ({\cal J}_0 , \breve u) \, ,
  \ee  
where 
\be\label{formaHep}
\breve {\cal K} ({\cal J}_0, \breve u) =  
 H_L (\breve u) + 
\e {\cal P} ({\cal J}_0 , \breve u) \, , \quad 
{\cal P} ({\cal J}_0 , \breve u) :=
\e^{-3} \cK_{\geq 3}  (\e  {\cal J}_0, \e \breve u)  \, , 
\ee
and 
$ H_L (\breve u)  = \tfrac{1}{2} ( {\bf \Omega} (\gamma)  \breve u, \breve u ) $.
If $  \underline{\breve u}(\omega t) $ is a quasi-periodic solution of \eqref{HS:mu}
then $  \underline{\wtilde u}(\omega t) = \e 
 \underline{\breve u}(\omega t)  $ is a quasi-periodic solution of \eqref{HSredu}
with $ \underline{\cal J} = \e {\cal J}_0 $. 
\\[1mm]
{\bf Tangential and normal variables.} 
For any $\bar n\geq 2$ we  
fix finitely many distinct `tangential" sites 
$ {\mathbb S} := \{ n_1, \ldots, n_{|\ST|} \} $ satisfying
\eqref{Tsites}, 
and we look for time quasi periodic reversible solutions of the equation
\eqref{HS:mu} 
close to the 
 reversible solutions  of the linear equation \eqref{Lin1}, cfr. \eqref{q-r-l}, 
\be\label{qlines}
 q(t, \theta)  =
 {\mathop \sum}_{n \in {\mathbb S}}  
a_n   \symm_n  \cos (\Omega_n (\gamma) t) \tc_n ( \theta) 
+   a_n \symm_n^{-1} \sin (\Omega_n (\gamma) t) \ts_n ( \theta)  \, .  
\ee
We now write the Hamiltonian equation \eqref{HS:mu} 
 in a convinient set of symplectic coordinates.
 In view of \eqref{Hdir} 
we  decompose the phase space $ H^s_{\bot,2} (\T) $  
(cfr. \eqref{splipin0}) of \eqref{HS:mu} as the direct sum 
\be\label{decoacca}
 H^s_{\bot,2} (\T) :=   
\acca_{\ST} 
\oplus  \acca_{\ST,2}^{\bot} 
\ee
of the symplectic tangential and normal subspaces
\begin{align}
\label{Htang}
& \acca_{\ST} := \Big\{ v = 
{\mathop \sum}_{n \in \ST}  \a_n \tc_n ( \theta) + 
   \b_n \ts_n ( \theta)  \, , \  (\a_n, \b_n ) \in \R^2 \Big\} \, , \\
&\acca_{\ST,2}^\bot := \acca_{\ST,2}^{\bot,s} := \Big\{ 
z = 
{\mathop \sum}_{n \in \N \setminus \ST, n \neq 2}  \a_n \tc_n ( \theta) + 
 \b_n \ts_n ( \theta)  \in H^s (\T)  \Big\} \, . \label{Hbot}
\end{align}
Note that these subspaces are both pairwise symplectic orthogonal 
and $ L^2 $-orthogonal.
We denote by $ \Pi_{\ST,2}^\bot $ the $ L^2 $-projector on
 $  \acca_{\ST,2}^\bot $. 
 We call $ v$  the ``tangential" variable and  $ z $ the `normal" one. 

Next, in view of \eqref{anbn-ri0}, 
we introduce, on the finite dimensional 
tangential subspace $ \acca_{\ST} $, action-angle coordinates
$ (I, \vartheta) := (I_n, \vartheta_n)_{n \in {\mathbb S}} $  by setting 
\be\label{AA-xi}
\begin{aligned}
& \a_n = \symm_n \sqrt{2n}   \sqrt{\zeta_n + I_n} \cos (\vartheta_n)  \, , \quad
 \b_n =\symm_n^{-1}\sqrt{2n}   \sqrt{\zeta_n + I_n} \sin (\vartheta_n) \, , 
\end{aligned}
\ee
where $  \zeta_n > 0  $ and  $ | I_n| < \zeta_n $, for any $  n \in {\mathbb S} $, and 
$ M_n $ are defined in \eqref{def:Lan}. 
Then we represent any function of the phase space $ H^s_{\bot,2} (\T) $ in \eqref{decoacca} as
\begin{align}  \label{aacoordinates} 
& \breve u( \theta) = \Xi  (   \vartheta,I, z)(\theta) := 
{\mathtt v}^\intercal (\vartheta,I)+ z  \quad {\rm where} 
 \\
& 
{\mathtt v}^\intercal (\vartheta, I) := {\mathop \sum}_{n \in \ST}  \symm_n  
 \sqrt{2n}  \sqrt{\zeta_n+I_n} \cos (\vartheta_n) \tc_n (\theta)  + 
 \symm_n^{-1} \sqrt{2n}   \sqrt{\zeta_n+I_n} \sin (\vartheta_n) \ts_n (\theta) \, . \nonumber
\end{align}
In  the coordinates $ 
(\vartheta, I, z) \in  \T^{|\ST|} \times \R^{|\ST|} \times \acca_{\ST,2}^\bot 
$ 
the involution $ {\cal S} $ defined in \eqref{S-invo} reads (see \eqref{Sanbn})
\begin{equation}\label{rev_aa}
	\vec {\cal S} : (   \vartheta,I, z)\mapsto ( -\vartheta,I, {\cal S} z )\, , 
\end{equation}
and the symplectic $2$-form in \eqref{sy2form0}, i.e. \eqref{sy2form}, becomes
\be\label{sympl_form}
{\cal W} = 
(d I \wedge d \vartheta) \oplus {\cal W}_{|{\frak H}_{\ST,2}^\bot}  \qquad 
{\rm where} \qquad 
d I \wedge d \vartheta := {\mathop \sum}_{n\in  \ST} d  I_n \wedge d\vartheta_n  \, . 
\ee
Note that ${\cal W}  $ is an exact $ 2 $-form since 
$ {\cal W} = d \Lambda  $
where $ \Lambda $ is the Liouville $ 1$-form
\begin{equation}\label{Lambda 1 form}
\Lambda_{(\vartheta, I, z)}[  \widehat \vartheta, \widehat I, \widehat z] := 
  I \cdot \widehat \vartheta + \tfrac12 (  \partial_\theta^{-1} z, \widehat z )_{L^2 (\T)} \, .  
\end{equation}
Hence the Hamiltonian vector field  associated to  the Hamiltonian
\begin{equation}\label{Hepsilon}
\begin{aligned}
\sK_{\e}  (\vartheta, I, z) 
& :=  {\breve \cK} ( {\cal J}_0,   \Xi  (\vartheta, I, z)) 
 \stackrel{ \eqref{formaHep}} = 
H_L (\Xi  (   \vartheta,I, z)) + 
\e {\cal P} ({\cal J}_0 ,  \Xi  ( \vartheta,I, z)) 
\end{aligned}
\end{equation}
is given by
\be\label{XHemu}
X_{\sK_{\e}} := 
( - \pa_I \sK_{\e} , \pa_\vartheta \sK_{\e} , 
\pa_\theta \nabla_{z} \sK_{\e} ) \, .
\ee 
In view of \eqref{aacoordinates}, the quadratic Hamiltonian $ H_L $ 
defined in \eqref{defQH} (see also \eqref{Hanbn})  simply reads in
the coordinates $ ( \vartheta , I , z ) $, up to a constant,
\be\label{QHAM}
(H_L \circ \Xi ) ( \vartheta , I , z )=   -\vec{\omega}(\g)\cdot I
+ \tfrac12  ( {\bf \Omega}(\g) z, z )_{L^2}  
\ee
where $ \vec{\omega}(\gamma) \in \R^{|\ST|} $ is the unperturbed 
tangential frequency vector defined in \eqref{tangential-normal-frequencies}
 and $ {\bf \Omega}(\gamma) $ is the self-adjoint operator in \eqref{defQH}.
By 
 \eqref{QHAM}, 
the Hamiltonian $ \sK_{\e} $ in \eqref{Hepsilon} reads
\begin{equation}\label{cNP}
\sK_{\e} =  \sN +  \e \sP   
\end{equation}
where
\be\label{cNP2}
\begin{aligned}
&  \sN :=   -\vec{\omega}(\g)\cdot I  + \tfrac12  ( {\bf \Omega}(\g) z, z )_{L^2} 
 \qquad {\rm and} \qquad 
  \sP ( \vartheta,I, z ) := 
{\cal P} ({\cal J}_0 ,  \Xi  ( \vartheta,I, z)) \, .  
\end{aligned}
\ee
We look for quasi-periodic solutions
of the Hamiltonian system generated by the Hamiltonian $\sK_{\e}$ in  \eqref{cNP}.
More precisely we look for an embedded invariant torus
\be\label{rev-torus}
i :\T^{|\ST|} \rightarrow
\T^{|\ST|}  \times \R^{|\ST|} \times \acca_{\ST,2}^\bot 
\,, \quad \varphi \mapsto i(\varphi):= (  \vartheta(\varphi), I(\varphi), z(\varphi)) \, , 
\ee
where $ \Theta(\vphi) := \vartheta (\vphi) - \vphi $ is a  $ (2 \pi)^{|\ST|} $-periodic function, 
invariant under the Hamiltonian vector field $ X_{\sK_{\e}} $ in \eqref{XHemu},  
filled by quasi-periodic solutions  
with Diophantine frequency 
vector $ \omega \in \R^{|\ST|}$ 
and which satisfies also other 
non-resonance conditions. 

\subsection{Nash-Moser theorem of hypothetical conjugation}

We first consider a relaxed problem 
where we introduce additional parameters. 
For $ \tg \in \R^{|\ST|} $, we consider the modified Hamiltonian
\begin{equation}\label{H alpha}
 \sK_{\tg} := \sN_\tg + \e \sP \, , \quad  \sN_\tg :=  
-\tg \cdot I 
+ \tfrac12 ({\bf \Omega}(\gamma) z, z)_{L^2}\, . 
\end{equation}
Let  $ {\mathtt \Omega} :=  $ ${\mathtt \Omega}(\delta):= $
$ \big\{  \om \in \R^{|\ST|} \, : \, {\rm dist} \big(\om,  {\vec \omega}[\g_1, \g_2]\big) < \d \big\}$,
 $ \d > 0 $  be a   $ \d $-neighborhood (independent of $ \e $) of the unperturbed linear frequencies $ {\vec \omega}[\g_1, \g_2] $ defined  in \eqref{tangential-normal-frequencies}. Given $(\om, \gamma,  \e )\in {\mathtt \Omega}\times [ \gamma_1, \gamma_2] 
\times [0,1)$,  
we look for zeros of the nonlinear operator 
\begin{align}
 {\cal F} (i, \tg ) 
& :=   {\cal F} (i, \tg, \om, \gamma,  \e )  := \Dom i (\vphi) - X_{\sK_{\tg}} ( i (\vphi))
=  \Dom i (\vphi) -  (X_{\sN_\tg}  +  \e X_{\sP})  (i(\vphi) ) \nonumber \\
&  :=  \left(
\begin{array}{c}
\Dom \vartheta (\vphi) - \tg +  \e \partial_I 	\sP ( i(\vphi)   )  \\
\Dom I (\vphi)  -  \e \partial_\vartheta \sP ( i(\vphi)  )  \\
\Dom z (\vphi) 
-  \pa_\theta {\bf \Om}(\g) z(\vphi) - \e \pa_\theta \nabla_z \sP ( i(\vphi) ))   
\end{array}
\right) \, .  \label{operatorF} 
\end{align} 
Each Hamiltonian $ \sK_{\tg} $ in \eqref{H alpha} is reversible, i.e. 
$ \sK_{\tg} \circ \vec {\cal S} = \sK_{\tg} $, 
where the involution $ \vec {\cal S} $ is defined in \eqref{rev_aa}. 
We look for reversible solutions of $ {\cal F}(i, \tg) = 0 $,  namely 
satisfying 
\be\label{parity solution0}
{\vec {\cal S}} i (\vphi ) = i (- \vphi) \, , 
\ee
that component-wise, writing $ i (\vphi) $ as in \eqref{rev-torus}, 
reads 
\begin{equation}\label{parity solution}
\vartheta(-\vphi) = - \vartheta (\vphi) \, , \,  \
I(-\vphi) = I(\vphi) \, , \, \ 
z (- \vphi ) = ( {\cal S} z)(\vphi) \, . 
\end{equation}
The norm of the periodic component of the embedded torus 
\begin{equation}\label{componente periodica}
{\mathfrak I}(\vphi)  := i (\vphi) - (\vphi,0,0) := 
( {\Theta} (\vphi), I(\vphi), z(\vphi))\,, \quad \Theta(\vphi) := \vartheta (\vphi) - \vphi \, , 
\end{equation}
is 
\be\label{def:norma-cp}
\|  {\mathfrak I}  \|_s^{k_0,\upsilon} 
:= \| \Theta \|_{H^s_\vphi}^{k_0,\upsilon} +  \| I  \|_{H^s_\vphi}^{k_0,\upsilon} 
+  \| z \|_s^{k_0,\upsilon}\, . 
\ee
{\bf Notation.}
{\it Along the paper we consider 
 Sobolev functions $\lambda\mapsto u(\lambda)\in H^s  $, $k_0$-times differentiable  with respect to 
$ \lambda :=(\omega,\gamma) \in  \R^{|\ST|} \times [\g_1,\g_2] $,
and, for  $\upsilon\in(0,1)$, we define the weighted Sobolev norm
\begin{equation}\label{weinorm}
	\| u \|_{s}^{k_0,\upsilon}:= 
	{\mathop \sum}_{|k|\leq k_0} \upsilon^{|k|} \sup_{\lambda\in \R^{|\ST|} \times [\g_1,\g_2]}\|\partial_\lambda^k u(\lambda)\|_{s} \,.
\end{equation}
If $ u $ is real valued we also denote $ | u |^{k_0,\upsilon}:=  \| u \|_{s}^{k_0,\upsilon}$.}

We define  $ k_0 := \barka + 2 	 $
where $\barka$ is the index of non-degeneracy provided by Proposition \ref{Lemma: degenerate KAM}, 
which only depends on the linear unperturbed frequencies. 
Thus $k_0$ is considered as an absolute constant,
and we will often omit to explicitly write the dependence of the various constants 
with respect to $ k_0 $.

\begin{theorem} \label{main theorem}
For any $\bar n\geq 2$, consider an interval
of aspect ratio  $\IK := [ \gamma_1, \gamma_2] $ 
as in \eqref{defgammaset0}. 
Let $ {\mathbb S} $ be any finite subset of distinct integers in 
$ \{\bar n+1,\bar n+2, \ldots \} $. Let $ {\cal J}_0 \in \R $ and  
 $\t   \geq 1 $. 
Then there exist positive constants $ a_0, \e_0, C $ 
depending on $\ST,  k_0, \t$ such that, 
for all $ \upsilon = \e^a $, $ 0 < a < a_0 $, for all $ \e \in (0, \e_0) $,
there exist
\begin{enumerate}
\item 
a $ k_0 $-times differentiable function 
$\tg:= \tg_\infty : \R^{|\ST|} \times [\g_1, \g_2] \mapsto \R^{|\ST|} $, 
\be\label{mappa aep}
 \tg_\infty (\om, \g ) = \om + r_\e (\om, \g) \, , 
\quad {with} \quad 
|r_\e|^{k_0, \upsilon}  \leq C \e \upsilon^{-1} \, ; 
\ee
\item
a family of embedded tori  $ i_\infty (\vphi) = 
(\vphi + \Theta_\infty (\vphi  ) ,I_\infty (\vphi  ), z_\infty (\vphi, \theta ))  $, defined for all $ (\om, \g) 
 \in \R^{|\ST|} \times  [\g_1, \g_2 ] $,  
satisfying the reversibility property \eqref{parity solution0} and 
\be\label{stima toro finale}
\|  i_\infty (\vphi) -  (\vphi,0,0) \|_{s_0}^{k_0, \upsilon} \leq C \e \upsilon^{-1}  \, ;
\ee
\item 
a constant $ {\mathtt m}^\infty : \R^{|\ST|} \times [\g_1, \g_2 ]  \to \R $ of the form 
\be\label{form1infty}
 {\mathtt m}^\infty (\omega, \gamma) = 
 \Omega_1 (\gamma) + {\mathtt r}_{\e}^\infty (\omega, \gamma) \, , \qquad 
\Omega_1 (\gamma)  = \Omega_\gamma = \frac{\gamma}{(1+\gamma)^2} \, , 
\ee
with 
\be\label{rnrs} 
|{\mathtt r}_{\e}^\infty |^{k_0, \upsilon} \lesssim_{k_0} \e  \, ;
\ee 
\item 
a sequence of $ {k_0} $-times differentiable functions
$ \Omega_n^\infty : \R^{|\ST|} \times [\g_1, \g_2 ]  \to \R $, for any
$  n \in \N \setminus (\ST \cup \{ 2, \ldots, \bar n \}) $, of the form
\be
\begin{aligned}
&  \Omega_n^\infty (\omega, \gamma) = 
\big| \big(\mu_n^+ (\gamma) + n 
{\mathtt r}_{\e}^\infty (\omega, \gamma) \big) 
\big( \mu_n^- (\gamma) +  n {\mathtt r}_{\e}^\infty (\omega, \gamma) \big) |^{\frac12} + {\frak r}_{n}^\infty (\omega, \g)
\label{autovalori infiniti-KAM}
\\
& \text{with} \quad \mu_n^+ (\gamma ) +  n {\mathtt r}_{\e}^\infty (\omega, \gamma)  = 
n  \, {\mathtt m}^\infty(\omega, \gamma) 
-    \tfrac{1}{2} + \tfrac{1}{2} \ka_n  \, , \\
&  \ \ \quad \quad 
 \mu_n^{-} (\gamma ) +  n {\mathtt r}_{\e}^\infty (\omega, \gamma) = 
n \, {\mathtt m}^\infty (\omega, \gamma)  -  \tfrac{1}{2}  -   \tfrac{1}{2} \ka_n  \, ,   
\quad
\ka_n =  \Big( \frac{\gamma-1}{\gamma + 1} \Big)^{n} \, ,  
\end{aligned}
\ee
and, for any $ M \geq 1  $, 
\begin{equation}\label{stime autovalori infiniti}
 \sup_{n \in \N \setminus (\ST \cup \{ 2, \ldots, \bar n \})} n^M |{\frak r}_n^\infty |^{k_0, \upsilon} 
 \lesssim_M \e \upsilon^{-1} \, , 
\end{equation}
\end{enumerate}
such that,  for all $ (\om, \g) $ in the  Cantor like set
\begin{align}\label{Cantor set infinito riccardo}
{\cal C}_\infty^{\upsilon} & := \Big\{ ( \omega, \g ) \in {\mathtt \Omega} \times [ \g_1, \g_2 ] \, : 
 \, |\om \cdot \ell  | \geq 8 \upsilon \langle \ell \rangle^{-\tau}, \, \forall \ell 
 \in \Z^{{|\ST|}} \setminus \{ 0 \} \, ,  
   \\
& \quad 
| \om  \cdot \ell + {\mathtt m}^\infty 
(\omega, \gamma) j | \geq \frac{8 \upsilon \langle j \rangle}{\langle \ell \rangle^{\tau}}, \, 
\ \forall (\ell, j ) \in (\Z^{|\ST|} \times \Z) \setminus \{(0,0)\} \, , \ 
  \label{trasp1M} \\
&  \quad  |\omega \cdot \ell  + \Omega_n^\infty (\omega, \g)  | \geq 4 n \upsilon 
\langle \ell  \rangle^{- \tau}, \, 
\forall \ell   \in \Z^{|\ST|}, \, n \in \N \setminus (\ST \cup \{ 2, \ldots, \bar n \}), 
  \nonumber \\
& \quad |\omega \cdot \ell  + 
 \Omega_n^\infty (\omega, \g)  -  \Omega_{n'}^\infty (\omega, \g)  | \geq 
 4 \upsilon \langle n-n' \rangle \langle \ell  \rangle^{-\tau} \, , \label{trasp4M} \\
& \quad \forall \ell   \in \Z^{|\ST|},\ n, n' \in \N \setminus (\ST \cup \{ 2, \ldots, \bar n \}) \, , \  
\ (\ell,n,n') \neq (0,n,n) \, , \nonumber \\ \nonumber
 & \quad |\omega \cdot \ell  + 
 \Omega_n^\infty (\omega, \g)  + \Omega_{n'}^\infty (\omega, \g)  | \geq 
 4 \upsilon (n + n' )\langle \ell \rangle^{-\tau} \, , \
 \forall \ell   \in \Z^{|\ST|}, \ n, n' \in \N \setminus (\ST \cup \{ 2, \ldots, \bar n \})  \nonumber   \Big\} \, , 
\end{align}
the function $( i_\infty (\vphi), \tg_\infty) :=( i_\infty (\omega, \g, \e)(\vphi) , \tg_\infty(\omega, \g, \e))$ solves
$ {\cal F}( i_\infty, \tg_\infty (\om, \g) , \om, \g, \e)  = 0 $. As a consequence the embedded torus 
$ \vphi \mapsto i_\infty (\vphi) $ is invariant for the Hamiltonian vector field 
$ X_{\sK_{\tg_\infty (\om, \g)}  } $  
and it is filled by quasi-periodic solutions with frequency $ \om $.
\end{theorem}

\subsection{Measure estimates and proof of Theorem \ref{thm:VP}}\label{sec:measure}

We now prove the existence of quasi-periodic solutions of the original Hamiltonian 
 $\sK_{\e}$ in \eqref{cNP}.
We proceed as follows. By  \eqref{mappa aep}, for any $\gamma \in [\gamma_1, \gamma_2]$,  the function 
$ \tg_\infty (\cdot, \gamma) $ from $ \tOm $ into the image $\tg_\infty( \tOm, \gamma)  $ is invertible:
\be\label{a-b}
\begin{aligned}
& g = \tg_\infty (\om, \gamma) = \om + r_\e (\om, \gamma)  \quad \Longleftrightarrow \quad \\
& 
\om  = \tg_\infty^{-1}(g , \gamma ) = g + {\breve r}_\e (g, \gamma) \quad {\rm with }  \quad
 |{\breve r}_\e |^{k_0, \upsilon} \leq C \e \upsilon^{-1} \, . 
 \end{aligned}
\ee 
For any $ g \in \tg_\infty ({\cal C}_\infty^\upsilon) $, Theorem \ref{main theorem} proves the existence of
 an embedded invariant torus  filled by quasi-periodic solutions with diophantine frequency 
$ \om =  \tg_\infty^{-1}(g , \gamma )  $ for the  Hamiltonian 
$ \sK_g = \sN_g +  \e \sP $. 
Consider the curve of the unperturbed linear frequencies 
$  \vec \om (\gamma )  $ in \eqref{tangential-normal-frequencies}.
For any  $ \gamma \in [\gamma_1, \gamma_2] $ such that 
the vector $(\tg_\infty^{- 1}(\vec \om (\gamma ), \gamma), \gamma)$ belongs to 
${\cal C}^\upsilon_\infty$, we thus obtain 
an embedded invariant torus for the Hamiltonian 
$ \sK_{\e} $ defined in \eqref{cNP}, 
filled by quasi-periodic solutions with diophantine frequency $ \om =  \tg_\infty^{-1}( \vec \om (\gamma ), \gamma ) $.
In Theorem \ref{Teorema stima in misura} below, we prove that such 
set of ``good" parameters 
\be\label{defG-ep}
{\cal G}_\e := 
\Big\{ \gamma \in [\gamma_1, \gamma_2] \, :  \big( \tg_\infty^{-1} ( {\vec \om} (\gamma ), \gamma), \gamma \big) \in  {\cal C}^\upsilon_\infty   \Big\} \, , 
\ee
has a large measure. 

\begin{theorem}\label{Teorema stima in misura} {\bf (Measure estimates)} 
Let
\begin{equation}\label{relazione tau k0}
 \upsilon = \e^a \, , \quad 0 < a < \frac{1}{k_0(1+6 k_0)} \, , \quad 
 \tau > k_0 (1+2k_0)( |\ST| + 1 )    \, . 
\end{equation} 
Then the measure of the set ${\cal G}_\e $ defined in \eqref{defG-ep}
satisfies  $ |{\cal G}_\e| \to  \gamma_2 - \gamma_1 $ as $ \e \to 0 $.
\end{theorem}

The rest of this subsection is devoted to the proof of Theorem \ref{Teorema stima in misura}. 
By \eqref{a-b}
the vector
\begin{equation}\label{omega epsilon kappa}
\vec \om_\e (\gamma) :=    \tg_\infty^{-1}( {\vec \om} (\gamma ), \gamma  ) = {\vec \om} (\gamma) + \vec {\mathtt r}_\e ( \gamma ) \, , \quad 
\vec {\mathtt r}_\e ( \gamma ) := {\breve r}_\e ({\vec \om } (\gamma), \gamma ) \, , 
\end{equation}
satisfies, for any $ \g \in \IK $,  
\begin{equation}\label{stima omega epsilon kappa}
| \pa_\gamma^k \vec {\mathtt r}_\e  (\gamma) | \lesssim \e \upsilon^{-(1  + k)} \, , \ \forall 0 \leq k \leq k_0 \, .
\end{equation}
In view of \eqref{form1infty}, we also denote (with a small abuse of notation)
\be\label{defm1ep}
{\mathtt m}^\infty (\gamma) := 
{\mathtt m}^\infty (\vec \om_\e(\gamma), \gamma) =
 \Omega_1 (\gamma) + {\mathtt r}_{\e}^\infty (\gamma) \, , \quad
 {\mathtt r}_{\e}^\infty (\gamma) := 
{\mathtt r}_{\e}^\infty ( \vec \om_\e(\gamma),\gamma) \, ,  
\ee
and, for all $  n \in \N \setminus ({\mathbb S} \cup \{ 2, \ldots, \bar n \}) $
(see \eqref{autovalori infiniti-KAM})
\begin{equation}\label{mu j infty kappa} 
\Omega_n^\infty (\gamma) :=
\Omega_n^\infty ( \vec \om_\e(\gamma) , \gamma) 
 = \big| \big(\mu_n^+ + n 
{\mathtt r}_{\e}^\infty (\gamma)\big) 
\big( \mu_n^- +  n {\mathtt r}_{\e}^\infty (\gamma) \big)\big|^{\frac12} + 
{\frak r}_{n}^\infty (\gamma)
\, ,
\end{equation}
 where 
$ {\frak r}_{n}^\infty (\gamma) := {\frak r}_{n}^\infty ( \vec \om_\e(\gamma),\gamma) $.
 
By \eqref{rnrs}, \eqref{stime autovalori infiniti},  \eqref{stima omega epsilon kappa},
for all $ 0 \leq k \leq k_0 $, for any $ \g \in \IK $,  
\be\label{derremu}
|\pa_\gamma^k {\mathtt r}_{\e}^\infty (\gamma)| \lesssim \e \upsilon^{-k} \, , \
\quad  
 \sup_{n \in \N \setminus (\ST \cup \{ 2, \ldots, \bar n \})}  n^M 
 |\pa_\gamma^k  {\frak r}_{n}^\infty (\gamma) | \lesssim_M \e \upsilon^{-k-1}  \, . 
\ee
\begin{lemma} {\bf (Perturbed normal frequencies)} Assume 
$ \e \upsilon^{-k_0-1}  \leq 1 $. 
Then the perturbed normal frequencies 
$\Omega_n^\infty (\gamma) $ 
defined in \eqref{mu j infty kappa} satisfy the expansion 
\be\label{Om-per}
\begin{aligned}
& \Omega_n^\infty (\gamma) =  
n \, {\mathtt m}^\infty (\gamma)  - \tfrac12  + r_\e (n,\g)  = 
n (\Omega_1 (\gamma) + {\mathtt r}_{\e}^\infty (\gamma)) 
- \tfrac12  + r_\e (n,\g) \, , \\  
& \text{with} \quad 
\sup_{n \geq \bar n + 1, \gamma \in \IK} n |\pa_\gamma^k r_\e(n,\gamma)| \leq C_k \, , \ \forall 0 \leq k \leq k_0 \, . 
\end{aligned}
\ee 
Moreover  
\be\label{Om-per-diff}
\begin{aligned}
& \ \ \Omega_n^\infty (\gamma) - \Omega_n (\gamma)
 =    n {\mathtt r}_{\e}^\infty (\gamma) 
 + {\frak r}_\e (n,  \g) \quad \text{with} \\
&   \sup_{n \geq  \bar n + 1 , \g \in \IK} 
n |\pa_\gamma^k {\frak r}_\e(n,  \gamma)| \leq C_k \e \upsilon^{-k-1} \, , 
\  \forall 0 \leq k \leq k_0 \, . 
\end{aligned}
\ee
\end{lemma}

\begin{pf}
By \eqref{mu j infty kappa}, \eqref{autovalori infiniti-KAM}, \eqref{defm1ep}  
we may write, for any $n\geq \bar n + 1 $, 
$$
 \Omega_n^\infty (\gamma) =  \Big| n 
 \big(\Omega_1 (\gamma) +    
{\mathtt r}_{\e}^\infty (\g) \big) -\frac{1}{2} \Big|\sqrt{ 
1 -        \frac{ \ka_n^2}{4} \Big( n\big( \Omega_1 (\gamma) +    {\mathtt r}_{\e}^\infty (\g) \big) - \frac{1}{2}\Big)^{-2} } + {\frak r}_{n}^\infty ( \g)  \, . 
$$
By \eqref{nOmebarn+1},  $ \Omega_1 (\gamma) = \frac{ \gamma}{(1+\gamma)^2}$,
 \eqref{rnrs},    
we deduce that, for $ \e $ small, for any $ n \geq \bar n + 1 $, $ \gamma \in \IK $, 
\be\label{nbep}
n\big(\Omega_1 (\gamma) +    
{\mathtt r}_{\e}^\infty (\g) \big) -\tfrac{1}{2} \geq 
(\bar n + 1) \big(\Omega_1 (\gamma) +    
{\mathtt r}_{\e}^\infty (\g) \big) -\tfrac{1}{2} 
>  \tfrac12 \underline{c} \, ,  
\ee  
and, since $ {\mathtt m}^\infty (\gamma) = \Omega_1 (\gamma) + {\mathtt r}_{\e}^\infty (\gamma) $  the last identity proves \eqref{Om-per} 
 with 
\be\label{renep}
\begin{aligned}
r_\e ( n, \gamma)
& := \Big( n \, {\mathtt m}^\infty (\gamma)
 -\frac{1}{2}\Big)\bigg(  
\sqrt{ 1 -        \frac{ \ka_n^2}{4} \Big( n \, {\mathtt m}^\infty (\gamma)
 - \frac{1}{2}\Big)^{-2}}-1\bigg) + {\frak r}_{n}^\infty ( \g)\\ & = -
\frac{\ka_n^2 }{4 ( n \, {\mathtt m}^\infty (\gamma)
 - \frac{1}{2})}  \bigg(\sqrt{ 1 - \frac{\ka_n^2}{4}
\Big(n \, {\mathtt m}^\infty (\gamma)
  - \frac{1}{2}\Big)^{-2}} + 1\bigg)^{-1}
+ {\frak r}_{n}^\infty (\gamma) \, .
\end{aligned}
\ee
By  \eqref{nbep} and \eqref{derremu} we deduce the estimate in \eqref{Om-per}. Then subtracting \eqref{ASYFR1} from \eqref{Om-per} gives
\eqref{Om-per-diff}
with $ {\frak r}_\e (n,  \g):=r_\e (n,\g)-r(n,\gamma) $.
By \eqref{rng}, \eqref{renep} where 
$ {\mathtt m}^\infty (\gamma) = \Omega_1 (\gamma) + {\mathtt r}_{\e}^\infty (\gamma) $ 
and  \eqref{derremu} we deduce that ${\frak r}_\e (n,  \g) $ satisfies 
 the bounds  \eqref{Om-per-diff}.
\end{pf}

By \eqref{Cantor set infinito riccardo}, 
\eqref{omega epsilon kappa},  \eqref{defm1ep}, \eqref{mu j infty kappa} the set
$ {\cal G}_\e $  in \eqref{defG-ep} writes 
\begin{align}\label{CGep}
{\cal G}_\e & := \Big\{ (\g  \in  [ \gamma_1, \gamma_2 ] \, : 
 \, | \vec \om_\e (\gamma) \cdot \ell  | 
 \geq 8 \upsilon \langle \ell \rangle^{-\tau}, \, \forall \ell \in \Z^{{|\ST|}} \setminus \{ 0 \} \, ,  
   \\
& \qquad 
| \vec \om_\e (\gamma)  \cdot \ell + {\mathtt m}^\infty 
(\gamma) j | \geq 8 \upsilon \langle j \rangle \langle \ell \rangle^{-\tau } \, 
\quad \forall (\ell, j ) \in \Z^{|\ST|} \times \Z \setminus \{(0,0)\} \, , \nonumber \\
&  \qquad | \vec \om_\e (\gamma) \cdot \ell  + \Omega_n^\infty (\g)  | \geq 4 n \upsilon 
\langle \ell  \rangle^{- \tau}, \, 
\forall \ell   \in \Z^{|\ST|}, \, n \in \N \setminus (\ST \cup \{ 2, \ldots, \bar n \})\, , 
  \nonumber \\
  & \qquad | \vec \om_\e (\gamma) \cdot \ell  + 
 \Omega_n^\infty ( \g) -  \Omega_{n'}^\infty ( \g) | \geq 
 4 \upsilon \langle n-n' \rangle \langle \ell  \rangle^{-\tau}  \, , \nonumber \\
& \qquad \forall \ell   \in \Z^{|\ST|} ,\,\,n, n' \in \N \setminus 
(\ST \cup \{ 2, \ldots, \bar n \}), 
\, (\ell,n,n') \neq (0,n,n)  \, ,  \nonumber  \\
& \qquad | \vec \om_\e (\gamma) \cdot \ell  + 
 \Omega_n^\infty (\g) + \Omega_{n'}^\infty ( \g ) | \geq 
 4 \upsilon (n + n' )\langle \ell \rangle^{-\tau}, \,\,
 \forall \ell   \in \Z^{|\ST|} ,\,\,n, n' \in \N \setminus (\ST \cup \{ 2, \ldots, \bar n \})   \nonumber  
 \Big\} \, .  \nonumber
\end{align}
We estimate the measure of the complementary set  
\begin{align} \label{complementare insieme di cantor}
 {\cal G}_\e^c  & := [\gamma_1, \gamma_2] \setminus {\cal G}_\e \\
& = \Big(\bigcup_{\ell \neq 0} R_{\ell}^{(0)} \Big) \bigcup \Big(\bigcup_{(\ell , j) \neq (0,0)} R_{\ell, j}^{(T)}  \Big)
\bigcup \Big(\bigcup_{\ell, n} R_{\ell,n}^{(I)} \Big) \bigcup
 \Big(\bigcup_{(\ell,n,n') \neq (0,n,n)} R_{\ell,n,n'}^{(-)} \Big) \bigcup  \Big(\bigcup_{\ell, n, n'} 
 R_{\ell,n,n'}^{(+)} \Big) \nonumber
\end{align}
where $ n, n' \in \N \setminus (\ST \cup \{ 2, \ldots, \bar n \})  $ and the ``resonant sets" are 
\begin{align}
 R_\ell^{(0)} &:= R_\ell^{(0)}(\upsilon,\tau)  := \big\{ \gamma \in [\gamma_1, \gamma_2] : 
| \vec \omega_\e (\gamma) \cdot \ell| <  8 \upsilon \langle \ell \rangle^{- \tau} \big\} 
\label{reso1} \\
 R_{\ell, j}^{(T)} &:=  R_{\ell, j}^{(T)}(\upsilon,\tau) := \big\{ \gamma \in [\gamma_1, \gamma_2] :  
 | \vec \om_\e (\gamma)  \cdot \ell + {\mathtt m}^\infty 
(\gamma) j | <  
8 \upsilon \langle j \rangle \langle \ell \rangle^{-\tau}   \big\}\, , \label{reso2new} 
\\
 R_{\ell, n}^{(I)} &:= R_{\ell, n}^{(I)}(\upsilon,\tau) 
:= \big\{ \gamma \in [\gamma_1, \gamma_2] :  
| \vec \om_\e (\gamma) \cdot \ell  + \Omega_n^\infty (\g)  | < 4 n \upsilon 
\langle \ell  \rangle^{- \tau}    \big\}\, , \label{reso2} 
\\
 R_{\ell,n, n'}^{(-)} &:= R_{\ell,n, n'}^{(-)}(\upsilon,\tau)  
:= \big\{ \gamma \in [\gamma_1, \gamma_2] : 
| \vec \om_\e (\gamma) \cdot \ell  + 
 \Omega_n^\infty ( \g) -  \Omega_{n'}^\infty ( \g) | < 
 4 \upsilon { \langle n-n' \rangle}{ \langle \ell  \rangle^{-\tau}} 
 \big\}\, ,  \label{reso4} 
 \\
 R_{\ell, n,  n'}^{(+)}  &:= R_{\ell, n,  n'}^{(+)}(\upsilon,\tau)  := \big\{ \gamma \in [\gamma_1, \gamma_2] : 
|\vec \om_\e (\gamma) \cdot \ell  + 
 \Omega_n^\infty (\g) + \Omega_{n'}^\infty ( \g ) | < 
 4 \upsilon {(n+n')}{\langle \ell \rangle^{-\tau}}   \big\} \, .  \label{reso3} 
 \end{align}
Some of the above sets are empty. 

\begin{lemma}\label{restrizione indici risonanti}
If $R_{\ell,j}^{(T)} \neq \emptyset$  then $ |j| \leq C \langle \ell \rangle $. 
If $R_{\ell,n}^{(I)} \neq \emptyset$ then $ n \leq C \langle \ell \rangle$. 
If $ R_{\ell,n,n'}^{(-)} \neq \emptyset$ then $ | n - n' | \leq C \langle \ell \rangle$.
 If $ R_{\ell,n,n'}^{(+)} \neq \emptyset$ then $ n + n'  \leq C \langle \ell \rangle$.
\end{lemma}

\begin{pf}
We prove the lemma for $ R_{\ell,n,n'}^{(-)} $ . 
The other cases follow similarly. We can suppose $ n \neq n' $ otherwise 
$ |n-n'| \leq C \langle \ell \rangle $ is trivial. 
If $\gamma \in R^{(-)}_{\ell,n,n'}$ then
$ |\Omega_n^\infty(\gamma) - \Omega_{n'}^\infty(\gamma)| < $
$ 4 \upsilon  \langle n-n' \rangle \langle \ell \rangle^{- \tau}  + 
| \vec \omega_\e(\gamma)| |\ell|  
 \leq  $ $ 4 \upsilon  |n-n' |  + C \langle \ell \rangle $. 
Moreover, by  \eqref{Om-per} and \eqref{derremu}, 
$ | \Omega_n^\infty - \Omega_{n'}^\infty|  \geq $ $ |n-n'| 
(\Omega_1 (\gamma) + {\mathtt r}_{\e}^\infty (\gamma) ) - C 
 \geq $ $ |n-n'| \tfrac12  \Omega_1 (\gamma)  - C $
for  $\e$ small enough  
and then $ |n- n'| \leq C_1 \langle \ell \rangle $.
\end{pf}

The key point  to estimate the measure of the resonant sets is that 
the perturbed frequencies satisfy transversality properties 
similar to the ones \eqref{0 Melnikov}-\eqref{2 Melnikov+} 
satisfied by the unperturbed frequencies.

\begin{lemma}\label{Lemma: degenerate KAM perturbato}
{\bf (Perturbed transversality)} 
For $ \e $ small enough, for all $ \gamma \in [\gamma_1, \gamma_2] $, 
\begin{align}\label{0 Melnikov-pert}
& \max_{k \leq \barka} 
|\partial_\g^{k}  \{{\vec \om}_\e (\gamma) \cdot \ell   \} |  
\geq \frac{\rho_0}{2} \langle \ell \rangle\,, 
\quad \forall \ell  \in \Z^{|\ST|} \setminus \{ 0 \},   
\\
\label{T Melnikov-pert}
& \max_{k \leq \barka}
|\partial_{\gamma}^{k}  \{{\vec \om}_\e (\gamma) \cdot \ell  
+  {\mathtt m}^\infty  (\gamma) j \} | 
 \geq \frac{\rho_0}{2} \langle \ell  \rangle\,, 
\quad \forall (\ell,j)  \in (\Z^{|\ST|} \times \Z) \setminus \{(0,0)\} \, ,  
\\
\label{1 Melnikov-pert}
& \max_{k \leq \barka}
|\partial_{\gamma}^{k}  \{{\vec \om}_\e (\gamma) \cdot \ell  +  \Omega_n^\infty (\gamma ) \} | 
 \geq \frac{\rho_0}{2} \langle \ell  \rangle\,, 
\quad \forall \ell  \in \Z^{|\ST|}, \, n \in \N \setminus (\ST \cup \{ 2, \ldots, \bar n \}) \, ,  
\\
\label{2 Melnikov-pert}
& \max_{k \leq \barka}
|\partial_\gamma^{k}  \{{\vec \om}_\e (\gamma) \cdot \ell  
 + \Omega_n^\infty (\gamma ) - \Omega_{n'}^\infty(\gamma ) \} | 
  \geq \frac{\rho_0}{2} \langle \ell  \rangle\,, \\
&  
\forall \ell \in \Z^{|\ST|} , \,  n, n' \in \N \setminus (\ST \cup \{ 2, \ldots, \bar n \}) \, ,  
\quad (\ell,n,n') \neq (0,n,n) \, , 
\nonumber \\ 
& \max_{k \leq \barka}
 |\partial_\gamma^{k}  \{{\vec \om}_\e (\gamma) \cdot \ell  
 + \Omega_n^\infty (\gamma ) + \Omega_{n'}^\infty(\gamma ) \} | 
 \geq \frac{\rho_0}{2} \langle \ell  \rangle\,, \
 \forall \ell \in \Z^{|\ST|}, \,  n, n' \in \N \setminus (\ST \cup \{ 2, \ldots, \bar n \}) \,   \label{2 Melnikov+pert}
\end{align} 
where $ k_0^* $ is the index of non-degeneracy 
defined in  Proposition \ref{Lemma: degenerate KAM}. 
\end{lemma}

\begin{pf}
We prove  \eqref{2 Melnikov-pert}. The other estimates are similar. 
By \eqref{omega epsilon kappa}, \eqref{stima omega epsilon kappa},
\eqref{Om-per-diff}, \eqref{2 Melnikov-} and \eqref{derremu} we get
\begin{align*}
& \max_{k \leq \barka}
|\partial_\gamma^{k}  \{{\vec \om}_\e (\gamma) \cdot \ell  
 + \Omega_n^\infty (\gamma ) - \Omega_{n'}^\infty(\gamma ) \} | 
  \geq \max_{k \leq \barka}
|\partial_\gamma^{k}  \{{\vec \om} (\gamma) \cdot \ell  
 + \Omega_n (\gamma ) - \Omega_{n'}(\gamma ) \} | \\
 & -
 C  \e \upsilon^{-(1  + k)}  \langle \ell \rangle -
 |n - n'| |\pa_\gamma^k {\mathtt r}_{\e}^\infty (\gamma)|  - |\pa_\gamma^k {\frak r}_\e (n,\g)| -  
 |\pa_\gamma^k {\frak r}_\e (n',\g)| 
 \geq
 \rho_0 \langle \ell \rangle / 2 
\end{align*}
and since we consider only sets such that $ |n-n'| \leq C \langle \ell \rangle $ by Lemma 
\ref{restrizione indici risonanti}. 
\end{pf}

\begin{lemma}[{\bf Estimates of resonant sets}]\label{stima risonanti Russman}
The measures of the sets in \eqref{complementare insieme di cantor}, 
\eqref{reso1}-\eqref{reso3} satisfy 
\be
\begin{aligned}
& |R_{\ell}^{(0)}| \lesssim  \big(  \upsilon \langle \ell \rangle^{- (\tau + 1)}\big)^{\frac{1}{k_0}} \,, \quad 
|R_{\ell,j}^{(T)}| \lesssim  \big(  \upsilon \langle j \rangle \langle \ell \rangle^{- (\tau + 1)}\big)^{\frac{1}{k_0}} \,, \quad
|R_{\ell,n}^{(I)}| \lesssim   \big(  \upsilon n \langle \ell \rangle^{- (\tau + 1)}\big)^{\frac{1}{k_0}} \, ,  \\
&  
|R_{\ell,n,n'}^{(-)}| \lesssim 
 \big(  \upsilon  \langle n - n' \rangle \langle \ell \rangle^{- (\tau + 1)}  \big)^{\frac{1}{k_0}}, \quad
|R_{\ell,n,n'}^{(+)}| \lesssim 
\big(  \upsilon (n+n') \langle \ell \rangle^{- (\tau + 1)}\big)^{\frac{1}{k_0}} \, . 
\label{secoriga}
\end{aligned}
\ee
\end{lemma}
\begin{pf}
We prove the estimate of $ R_{\ell,n,n'}^{(-)} $. 
The other cases are simpler. 
We write 
$$
R_{\ell,n,n'}^{(-)} = \big\{ \gamma \in [\gamma_1, \gamma_2] : |g_{\ell,n,n'}(\gamma)| < 4 \upsilon \langle n - n' \rangle \langle \ell \rangle^{- (\tau + 1)}  \big\}
$$
where
$ g_{\ell,n,n'}(\gamma) := ( \omega_\e(\gamma) \cdot \ell + \Omega_n^\infty(\gamma) - \Omega_{n'}^\infty(\gamma) ) 
\langle \ell\rangle^{-1} $. We apply Theorem 17.1 in \cite{Ru1}.  
By \eqref{2 Melnikov-pert} we derive that  
$ {\rm max}_{k \leq k_0} | \pa_\gamma^k 
g_{\ell,n,n'}(\gamma)| \geq \rho_0 / 2 $, for any 
$ \gamma \in [\gamma_1, \gamma_2] $. 
In addition,  by \eqref{omega epsilon kappa}-\eqref{stima omega epsilon kappa},
\eqref{Om-per} and since by  Lemma \ref{restrizione indici risonanti} we consider
$|n-n'| \leq C \langle \ell \rangle $, 
we deduce that 
$ \max_{k \leq k_0} | \pa_\gamma^k g_{\ell,n,n'}(\gamma) | \leq C_1 $, 
$  \forall  \gamma \in [\gamma_1, \gamma_2] $, 
provided $\e \upsilon^{- (1 + k_0)}$ is small enough.
By Theorem 17.1 in \cite{Ru1}  the bound \eqref{secoriga} for 
$ R_{\ell,n,n'}^{(-)} $ follows. 
\end{pf}

We now estimate the measure of  the sets in \eqref{complementare insieme di cantor}.
By Lemmata \ref{restrizione indici risonanti} and  \ref{stima risonanti Russman}  and the 
condition on $ \tau $ in \eqref{relazione tau k0}, we get
\begin{align}\label{meas1}
 \Big| \bigcup_{\ell \neq 0} R_{\ell}^{(0)} \Big| 
 & \leq
\sum_{\ell \neq 0} |R_\ell^{(0)}|  
\lesssim 
\sum_{\ell} \big( \upsilon \langle \ell \rangle^{- (\tau + 1)}\big)^{\frac{1}{k_0}} \lesssim 
\upsilon^{\frac{1}{k_0}} \, , \\
\label{meas2}
\Big|  \Big(\bigcup_{(\ell , j) \neq (0,0)} R_{\ell, j}^{(T)}  \Big)
\bigcup \Big(\bigcup_{\ell, n} R_{\ell,n}^{(I)} \Big)\Big|
&  \leq
\sum_{|j| \leq C \langle \ell \rangle} |R_{\ell,j}^{(T)}|
+  \sum_{n \leq C \langle \ell \rangle} |R_{\ell,n}^{(I)}| \lesssim 
\sum_{\ell } \ell \big(  \upsilon  \langle \ell \rangle^{- \tau}\big)^{\frac{1}{k_0}} 
\lesssim 
\upsilon^{\frac{1}{k_0}}\, ,\\
\label{meas3}
\Big|\bigcup_{\ell, n, n'}  R_{\ell,n,n'}^{(+)} \Big|
&\leq \sum_{n, n' \leq C \langle \ell \rangle} |R_{\ell,n,n'}^{(+)}|
\lesssim \sum_{n, n' \leq C \langle \ell \rangle}  \ell^2 
\big(  \upsilon  \langle \ell \rangle^{- \tau}\big)^{\frac{1}{k_0}}  
\lesssim \upsilon^{\frac{1}{k_0}} \, . 
\end{align}
It remains to estimate
$ \cup_{(\ell,n,n') \neq (0,n,n)} R_{\ell,n,n'}^{(-)} (\upsilon, \tau) $. 
We  need the following inclusion lemma. 

\begin{lemma}\label{lemma:inWave}
Let $ \upsilon_0 \geq   \upsilon $ and $ \tau \geq \tau_0 \geq 1 $.  
Then, for $ \e $ small, for all $ (\ell,n,n') \neq (0,n,n) $, 
\be\label{incluwave}
\min \{n,n'\} \geq \upsilon_0^{-2} \langle \ell \rangle^{\tau_0} \quad  \Longrightarrow \quad 
R_{\ell,n, n'}^{(-)}(\upsilon,\tau) \subset 
\bigcup_{(\ell , j) \neq (0,0)} R_{\ell, j}^{(T)}(\upsilon_0,\tau_0) \, . 
\ee
\end{lemma}

\begin{pf}
Let $  \gamma \in [ \g_1, \g_2] \setminus 
\cup_{(\ell , j) \neq (0,0)} R_{\ell, j}^{(T)}(\upsilon_0,\tau_0) $ (cf. \eqref{reso2new}), then   
\be\label{1Melpg}
| \vec \om_\e (\gamma)  \cdot \ell + {\mathtt m}^\infty 
(\gamma) j | \geq 
8 \upsilon_0 \langle j \rangle \langle \ell \rangle^{-\tau_0} \, , 
\quad  \forall  (\ell, j) \in (\Z^{|\ST|} \times \Z) \setminus \{ (0,0)\} \, . 
\ee
Then, by the expansion \eqref{Om-per}, 
we have, for all $ (\ell,n,n') \neq (0,n,n)  $, 
\begin{align}
| \vec \om_\e (\gamma) \cdot \ell  +  \Omega_n^\infty ( \g) -  \Omega_{n'}^\infty ( \g) |  
&  \geq  | \vec \om_\e (\gamma) \cdot \ell + (n-n') 
{\mathtt m}^\infty (\gamma) | -   \frac{C }{n} -   \frac{C }{n'} \nonumber \\
& \stackrel{\eqref{1Melpg}} 
\geq \frac{ 8 \upsilon_0 \langle n-n' \rangle}{\langle \ell \rangle^{\tau_0}} -   \frac{2 C}{\min \{n,n'\}}
\geq \frac{4 \upsilon_0 \langle n-n' \rangle }{ \langle \ell \rangle^{\tau_0}}    \nonumber \\
& 
\geq 4 \upsilon \langle n-n' \rangle  \langle \ell \rangle^{-\tau} 
 \label{2dav} 
\end{align}
for $ \min \{n,n'\} \geq \upsilon_0^{-2} \langle \ell \rangle^{\tau_0} $, 
taking $ \e $ small, 
and since $ \tau \geq \tau_0 $ and $ \upsilon_0 \geq \upsilon $.
Recalling the definition \eqref{reso4} of $ R_{\ell,n, n'}^{(-)}(\upsilon,\tau) $,
the estimate \eqref{2dav} 
proves that $ \g \in [\g_1, \g_2] \setminus R_{\ell,n, n'}^{(-)}(\upsilon,\tau)  $, thus 
 \eqref{incluwave}. 
\end{pf}

Note that the set of indices $ (\ell,n,n') \neq (0,n,n) $ such that 
$ \min \{ n , n' \}  < \upsilon_0^{-2} \langle \ell \rangle^{\tau_0}$ and
$ |n - n' | \leq C \langle \ell \rangle $ is included, for $ \upsilon_0 $ small enough, into 
 the set 
 \be\label{defsetIl}
 {\cal I}_\ell := \Big\{
 (\ell,n,n') \neq (0,n,n) \, , \ 
 n, n'  \leq \upsilon_0^{-3}  \langle  \ell \rangle^{\tau_0} \Big\} 
 \ee
because
$ \max \{ n , n' \} \leq \min \{ n , n' \} + |n - n' |  \leq  $ $
\upsilon_0^{-2} \langle \ell \rangle^{\tau_0} +  $
$ C \langle \ell \rangle \leq \upsilon_0^{-3}  \langle \ell \rangle^{\tau_0} $.

As a consequence, by Lemma
\ref{lemma:inWave} we deduce that 
\begin{align}
\bigcup_{(\ell,n,n') \neq (0,n,n)} R_{\ell,n,n'}^{(-)} (\upsilon, \tau) 
\subset
\Big(\bigcup_{(\ell , j) \neq (0,0)} R_{\ell, j}^{(T)} (\upsilon_0, \tau_0)  \Big) 
\bigcup
 \Big(\bigcup_{(\ell,n,n') \in {\cal I}_\ell} 
 R_{\ell,n,n'}^{(-)} (\upsilon, \tau)  \Big)   \, .   \label{primainc}
 \end{align}

\begin{lemma}\label{espotau0}
 Let $ \tau > k_0 (2 \tau_0 + |\ST| ) $,  $ \tau_0 := 1+ k_0 (|\ST|+1)  $, 
$ \upsilon_0 = \upsilon^{\frak{a}} $ with $ \frak a = 1 / (6 k_0 + 1 )$.
Then 
\begin{equation} \label{tI2}
\Big| \bigcup_{(\ell,n,n') \neq (0,n,n)} R_{\ell,n,n'}^{(-)} (\upsilon, \tau)  \Big|\lesssim 
\upsilon^{\frac{\frak a}{k_0}} \, . 
\end{equation}
 \end{lemma}
 
 \begin{pf}
 By \eqref{secoriga} (applied with $ \upsilon_0, \tau_0 $ instead of $ \upsilon, \tau $), 
 and since $ \tau_0 := 1+ k_0 (|\ST|+1) $, the measure of 
\begin{equation}\label{stima.wtR}
\Big| \bigcup_{(\ell,j) \neq (0,0)} R^{(T)}_{\ell,j} (\upsilon_0, \tau_0)  \Big|
\lesssim
\sum_{\ell} \ell \big( \upsilon_0 \langle \ell \rangle^{- \tau_0 }\big)^{\frac{1}{k_0}}
 \lesssim \upsilon_0^{\frac{1}{k_0}}  \lesssim \upsilon^{\frac{\frak a}{k_0}} \, . 
\end{equation}
Moreover, by \eqref{secoriga}, 
Lemma \ref{restrizione indici risonanti} and  \eqref{defsetIl}
and the  choice of $ \tau, \upsilon_0 $, 
\be\label{stima.altra}
\Big| \bigcup_{(\ell,n,n') \in {\cal I}_\ell} 
 R_{\ell,n,n'}^{(-)} (\upsilon, \tau) \Big|  	
  \lesssim \sum_{\ell \in \Z^{|\ST|},\atop{ n,n'\leq \upsilon_0^{-3} \braket{\ell}^{\tau_0} }} \Big( \frac{\upsilon \la n- n' \ra}{\braket{\ell}^{\tau+1}} \Big)^{\frac{1}{k_0}}  \lesssim \sum_{\ell\in\Z^{|\ST|}} \frac{\upsilon^{\frac{1}{k_0}}\upsilon_0^{-6}}{\braket{\ell}^{\frac{\tau}{k_0}-2 \tau_0}} \leq C \upsilon^{\frac{1}{k_0} - 6 \frak a } 
\lesssim \upsilon^{\frac{\frak a}{k_0}} 
	\ee
	by the choice of $ \frak a $. 
	The bound 
	\eqref{tI2} follows by \eqref{stima.wtR} and \eqref{stima.altra}.
	\end{pf}

\noindent 
{\sc Proof of Theorem \ref{Teorema stima in misura} completed.}
By \eqref{relazione tau k0}, \eqref{meas1}, \eqref{meas2}, \eqref{meas3} and 
Lemma \ref{espotau0} we deduce that 
the measure of the set ${\cal G}_\e^c$ 
in \eqref{complementare insieme di cantor} is estimated by
$ |{\cal G}_\e^c|  \lesssim 
\upsilon^{  \frac{{\frak a}}{k_0}} \lesssim
\e^{  \frac{a}{k_0(1+ 6 k_0)}}   $ 
since $ \upsilon = \e^a $.  
The proof of Theorem \ref{Teorema stima in misura} is concluded. 
\\[1mm]
{\bf Proof of Theorem \ref{thm:VP}.} Fix $ {\cal J}_0 = 0 $. 
By Theorems \ref{main theorem} and \ref{Teorema stima in misura}, 
for any $ \gamma \in {\cal G}_\e $ defined in \eqref{defG-ep}, 
there exists a  quasi-periodic solution   $  \breve u_\e ( \vec \om_\e (\gamma) t, \theta) $ 
of the equation \eqref{HS:mu} with $ {\cal J}_0 = 0 $ and a diophantine frequency vector 
$ \vec \om_\e (\gamma) =(\om_{\e,n} (\gamma))_{n \in {\mathbb S}} $ 
as in \eqref{omega epsilon kappa}, 
of the form (recall  also \eqref{aacoordinates})  
\be\label{lapQP}
\breve u_\e ( \vphi, \theta) 
 =  
  {\mathtt v}^\intercal ( \vphi  + \Theta_\infty (
 \vphi  ) ,I_\infty (\vphi  ))+ 
  z_\infty (\vphi, \theta )   
  = {\mathtt v}^\intercal (\vphi,0 )   + \breve{\rm r}_\e( \vphi, \theta) \, .
  \ee
The remainder 
$ \breve{\rm r}_\e(\vphi, \theta) :=  z_\infty (\vphi, \theta ) + 
  {\mathtt v}^\intercal ( \vphi  + \Theta_\infty (
 \vphi  ) ,I_\infty (\vphi  )) - {\mathtt v}^\intercal ( \vphi ,0) $
 satisfies, using \eqref{stima toro finale},
the estimate 
$ \|   \breve{\rm r}_\e   \|_{s_0} = O(\e \upsilon^{-1})  $ where 
$ \upsilon = \e^{a} $ and $ a > 0 $ satisfies \eqref{relazione tau k0}.
By \eqref{aacoordinates} and \eqref{defcnsn},  
we have 
 \be\label{chiuev}
{\mathtt v}^\intercal (\vec \om_\e (\gamma) t ,0 ) = 
{\mathop \sum}_{n \in \ST} \symm_n  
\ta_n \cos (\om_{\e,n} (\gamma) t) \cos  (n \theta)  + \symm_n^{-1} 
   \ta_n \sin (\om_{\e,n} (\gamma) t) \sin (n \theta) 
\ee
with $ \ta_n :=  \sqrt{2n \zeta_n / \pi} $. 
Recalling the rescaling \eqref{rescavar} and \eqref{HSKt0res}, the function 
$  \wtilde u_\e (\vec \om_\e (\gamma) t, \theta) := \e 
 \breve u_\e (\vec \om_\e (\gamma) t, \theta)  $ is a quasi-periodic solution of \eqref{HSredu}
with $ \underline{\cal J}  = 0 $, and 
\be\label{defwtilde}
\begin{aligned}
\wtilde  {\xi}_\e (\vec \om_\e (\gamma) t , \theta) 
& := 
   \widetilde\beta_{2,\e}  (\vec \om_\e (\gamma) t) \ts_2 (  \theta)+ 
\e 
 \breve u_\e (\vec \om_\e (\gamma) t, \theta)  \, , 
 \end{aligned}
\ee
with $ \widetilde\beta_{2,\e} (\varphi) $ defined  in \eqref{defmubeta2},
is a quasi-periodic solution of the Hamiltonian system 
$ \pa_t \wtilde \xi = X_{K - \underline{\mu}_\e \wtilde \alpha_2} (\wtilde \xi) $,
see \eqref{linxi1}, with $  \underline{\mu}_\e $  defined  by the second line of 
\eqref{defmubeta2}. 
Recalling \eqref{lapQP}, and since 
 $ \| \widetilde\beta_{2,\e} \|_{s_0} = o(\e ) $ (by renaming $ s_0 $), 
the quasi-periodic function  \eqref{defwtilde} has the expansion
\be\label{formaxit}
\wtilde  {\xi}_\e (\vec \om_\e (\gamma) t , \theta)  =   
\e {\mathtt v}^\intercal (\vec \om_\e (\gamma) t,0 ) + 
\widetilde{\rm r}_\varepsilon (\vec \om_\e (\gamma) t , \theta ) 
\ee
where the remainder
$ \widetilde{\rm r}_\varepsilon (\vphi , \theta ) := \e \breve{\rm r}_\e(\vphi , \theta) + 
\widetilde\beta_{2,\e}  (\vphi) \ts_2 (  \theta)  $
satisfies
$  \|  \widetilde{\rm r}_\varepsilon   \|_{s_0} = o(\e)  $. 
In conclusion 
\be\label{defxiep}
\xi_\e (\vec \om_\e (\gamma) t , \theta) 
:=  \Phi^{-1} (\wtilde{\xi}_\e (\vec \om_\e (\gamma)  t )) 
\ee
is a  quasi-periodic solution  of the Hamiltonian system  
$ \pa_t \xi = X_{H -   {\underline \mu}_\e {\cal J}} ( \xi )  $ in  \eqref{linxi2} and, by \eqref{Homfin},   
of 
equation 
\eqref{Evera} with  
$ \Omega =  \Omega_\gamma -  {\underline \mu}_\e  \tfrac{2\sqrt{2}}{\sqrt{\pi}(\g -\g^{-1})} $. 
 We finally remark that, by \eqref{formaxit}, \eqref{chiuev},  
 and the properties of $ \Phi $ in Theorem 
\ref{thm:FM},  the quasi-periodic 
function  $ \xi_\e (\vec \om_\e (\gamma) t , \theta)  $ in \eqref{defxiep} has the form
\eqref{xiexpa}.  
This proves Theorem \ref{thm:VP} with ${\cal G} = {\cal G}_\e $,
$ \wtilde \omega(\gamma) = \vec \omega_\e (\gamma) $ and $\mu(\gamma)= -  
{\underline \mu}_\e  \tfrac{2\sqrt{2}}{\sqrt{\pi}(\g -\g^{-1})}$.

The rest of the paper will be devoted to  the proof of Theorem \ref{main theorem}.

\section{Tame estimates}\label{sec:regularity}

In this section we prove tame estimates 
for the composition of the nonlinear vector fields 
$X_{H} $, $ X_K $ and the   
rectification map $ \Phi $  with functions $ \xi (\lambda, \vphi, \theta) $ 
in the norm $ \| \ \|_s^{k_0,\upsilon} $ defined in \eqref{weinorm}.

\smallskip

In the analysis of the nonlinearity  we encounter integral operators acting on a 
$ 2 \pi $-periodic function $ \xi (\theta) $ as in  \eqref{def:int-op}, 
where the Kernel function  $ K(\varphi, \theta, \theta' ) $ may be  smooth or 
singular at the diagonal $ \theta = \theta' $. 
An example of integral operator with singular Kernel is
 \eqref{int14}. 
An integral operator with smooth Kernel is infinitely many times  
regularizing, see Lemma \ref{lem:Int}.
In order to quantify it conveniently
we first introduce the following  basic definition 
 of pseudo-differential operators. 
\begin{definition} \label{def:Ps2} 
A function $ a (\theta, \eta ) $ which is $ {\cal C}^\infty $-smooth on $ \R \times \R $,
 $ 2 \pi $-periodic in $ \theta $, and satisfies, for some $ m \in \R $,  the inequalities
$ \big| \pa_\theta^{n_1} \pa_\eta^{n_2} a (\theta,\eta ) \big| \leq C_{n_1,n_2} 
\langle \eta \rangle^{m - n_2} $, $  \forall n_1, n_2 \in \N_0 $, 
is called a symbol of order $ m $. 
Given a $ 2 \pi$-periodic function  $ u (\theta) = {\mathop \sum}_{j \in \Z} u_j e^{\ii j \theta } $, 
we associate to $ a(\theta, \eta) $, the pseudo-differential  operator of order $  m $, 
$ (Au) (\theta) 
:= {\mathop \sum}_{j \in \Z} a(\theta,j) u_j e^{\ii j \theta }  $. 
We denote 
$ A = {\rm Op} (a)  = a(\theta, D) $ with $  D := D_\theta := \tfrac{1}{\ii} \pa_\theta $.  
We denote by $ S^m $ 
the class of  
symbols of order $ m $, and by $ {\rm OPS}^m $ 
the set of pseudo-differential  operators of order $ m $.  The set 
$ {\rm OPS}^{-\infty} := \cap_{m  \in \R} {\rm OPS}^{m} $ are the
infinitely many times regularizing operators. 
\end{definition}

When the symbol $ a(\theta)$ is independent of $ \eta $, the operator  $ A = {\rm Op} (a) $
is the multiplication operator by the function $ a(\theta ) $.
In such a case we also  denote 
$ A = {\rm Op} (a) = a(\theta ) $. If the symbol $ a(\eta ) $ is independent of $ \theta $ 
then $ {\rm Op} (a)  $ is a Fourier multiplier. 

In the paper we encounter $ \vphi $-dependent 
symbols  $ a(\lambda, \vphi,  \theta, \eta ) $ which are $ {\cal C}^\infty $-smooth  in $ \vphi $ and $ k_0 $-times differentiable with respect to
  $ \lambda := (\omega, \gamma) \in {\mathtt \Lambda}_0 \subset \R^{{|\ST|}+1}$. 
The following norm (Definition 2.11 of \cite{BertiMontalto})
controls the regularity in $ (\vphi, \theta)$  and $ \lambda$, and  the decay in $ \eta $,  of a symbol $ a (\lambda, \vphi, \theta, \eta) $. 

\begin{definition}\label{def:pseudo-norm} 
Let $ A(\lambda) := a(\lambda, \vphi, \theta, D) \in {\rm OPS}^m $, $ m \in \R $,  
be 
pseudo-differential operators 
$k_0$-times differentiable with respect to $ \lambda \in \mathtt \Lambda_0 \subset \R^{{|\ST|} + 1} $. 
For $ \upsilon \in (0,1) $, $ \a \in \N_0 $, $ s \geq 0 $, we define  the 
norm 
\begin{equation}\label{norm1 parameter}
| A |_{m, s, \alpha}^{k_0, \upsilon} := 
{\mathop \sum}_{|k| \leq k_0} \upsilon^{|k|} 
{\rm sup}_{\lambda \in {\mathtt \Lambda}_0}|\partial_\lambda^k A(\lambda) |_{m, s, \alpha}  
\end{equation}
where  
$
| A( \lambda )  |_{m, s, \a} := {\rm max}_{0 \leq \beta  \leq \a} \sup_{\eta \in \R} \|  \partial_\eta^\beta 
a(\lambda, \cdot, \cdot, \eta )  \|_{s} \langle \eta \rangle^{-m + \beta} $. 
\end{definition}

We shall first prove that 
the difference between 
the integral operator $ \W (\xi) $  defined  in   \eqref{int-op} and the unperturbed operator 
$ \W_0 $  computed in Lemma \ref{intop0},   
is infinitely many times regularizing.

\begin{lemma}\label{Prop:csm}
Let  $ \| \xi \|_{s_0+1}^{k_0,\upsilon} \leq	 \d $ small enough. 
Then
the integral  operator $ \W (\xi) $ defined in  \eqref{int-op} 
decomposes as 
\be\label{decoWxiW0}
\W (\xi) = \W_0 + {\cal R}(\xi) 
\ee
where $ \W_0 $, defined in \eqref{defK0y}, is computed in \eqref{defK0}, 
and $ {\cal R} (\xi) $  
is  in $ {\rm OPS}^{-\infty} $, and
for any $ m,  \a \in \N_0 $,  for some constant $ \sigma (m,\a) > 0 $, 
for any $ s \geq s_0 $, 
\be\label{stimaRep}
|  {\cal R}(\xi)   |_{-m, s, \a}^{k_0, \upsilon} \lesssim_{m, s, \alpha, k_0}   
  \| \xi \|_{s + \sigma (m,\a)}^{k_0, \upsilon} \,.
\ee
\end{lemma}

\begin{pf}
In order to prove the decomposition \eqref{decoWxiW0}
we put in evidence the contribution at $ \xi = 0 $ of
the function $ \Ke(\xi) (\theta, \theta')  $ defined in \eqref{expR}, writing 
\be\label{decoker}
 \Ke(\xi) (\theta, \theta')   =  \Ke(0)(\theta, \theta')  + G(\xi)(\theta, \theta') =
 \Ke(0)(\theta, \theta') \Big(1 + \frac{G(\xi)(\theta, \theta')}{\Ke(0)(\theta, \theta')} \Big) 
\ee
where $ \Ke(0)(\theta, \theta')  $ is given in \eqref{rel-dis0},  
and 
\begin{align}\label{Gxi}
G(\xi)( \theta, \theta') & = \gamma \, G_1(\xi)( \theta, \theta')+ \gamma^{-1} 
G_2(\xi)( \theta, \theta')\, ,\\
  G_1(\xi)( \theta, \theta')&:=   \big[ 
 \sqrt{1 + 2 \xi (\theta)} \cos (\theta) -
  \sqrt{1 + 2 \xi (\theta')} \cos (\theta')
  \big]^2 - \big[ 
  \cos (\theta) - \cos (\theta')
  \big]^2\, , \label{l2}  \\ 
    G_2(\xi)( \theta, \theta')&:=
  \big[ 
 \sqrt{1 + 2 \xi (\theta)} \sin (\theta) -
  \sqrt{1 + 2 \xi (\theta')} \sin (\theta')
  \big]^2 -
  \big[ 
 \sin (\theta) -
 \sin (\theta')
  \big]^2 \, . \label{l3}
  \end{align}
Note that  the  function 
\be\label{defP}
p (\gamma, \theta, \theta' ) := \tfrac{\g^2+1}{\g^2-1} -  
 \cos (\theta + \theta' ) \geq  \tfrac{2}{\g^2-1} > 0 \, , \quad \forall \theta, \theta' \, , 
\ee
is strictly positive. 
Then by \eqref{decoker}-\eqref{rel-dis0} we write
the kernel of the integral operator $ \W (\xi) $ in \eqref{int-op} as  
\be\label{kerunp} 
 \tfrac{1}{4 \pi} \ln \big( \Ke(\xi)( \theta, \theta') \big)  = 
\tfrac{1}{4 \pi} \ln \big( \Ke(0)( \theta, \theta') \big) + R(\xi) (\theta, \theta')
\ee
where the first term is the kernel of the unperturbed operator $ \W_0 $ 
in Lemma \ref{intop0} and
\be\label{kerper}
R(\xi) (\theta, \theta') :=
 \frac{1}{4 \pi} 
 \ln   \Big( 1 + \frac{G(\xi)( \theta, \theta')}{\sin^2 \big( \frac{\theta'- \theta}{2} \big)} 
  \frac{\gamma}{2 (\gamma^2-1)p (\gamma, \theta, \theta') }\Big) \, .
\ee
This proves the decomposition \eqref{decoWxiW0} with integral operator 
\be\label{defcalR}
 {\cal R}(\xi) [q]:=\int_{\T} R(\xi)(\theta, \theta')q(\theta')d\theta'\, . 
\ee
We now show that  the Kernel $R(\xi) (\theta, \theta')  $ in 
 \eqref{kerper} extends to a $ \mathcal{C}^\infty (\T^2) $-function 
 and thus the integral operator $ {\cal R}(\xi)  $ 
 is in $ {\rm OPS}^{-\infty} $, see Lemma \ref{lem:Int}.  Since 
the $ \mathcal{C}^\infty $ function $ p (\gamma, \theta, \theta') > 0 $ is  strictly positive by \eqref{defP}, 
it is sufficient to prove that the function 
\be\label{contocrux}
(\theta, \theta')  \mapsto \frac{G(\xi)( \theta, \theta')}{\sin^2 \big( \frac{\theta'- \theta}{2} \big)}  
= \gamma \frac{G_1(\xi)( \theta, \theta')}{\sin^2 \big( \frac{\theta'- \theta}{2} \big)}   +
\gamma^{-1} \frac{G_2(\xi)( \theta, \theta')}{\sin^2 \big( \frac{\theta'- \theta}{2} \big)}
\ee
extends to a $  \mathcal{C}^\infty $ function on $ \T^2 $, which is small in $ \xi $. 
Note that the 
denominator in \eqref{contocrux}  vanishes for $ \theta' - \theta = 2 \pi k $ for any $ k \in \Z $,
and thus, due to the $ 2 \pi $-periodicity of $ G(\xi) (\theta, \theta')$ in $ \theta, \theta' $,  
it is sufficient to prove that 
the function in \eqref{contocrux}  extends to a $  \mathcal{C}^\infty $ function
outside the diagonal $ \{ \theta = \theta' \} $.  
We now prove this property  for the terms with  $ G_1 $ and $ G_2 $ separately.    
The function $ G_1 $   in \eqref{l2} can be written as 
\be\label{tt1}
\begin{aligned}
  G_1(\xi)( \theta, \theta') &= \big(a^+(\theta) -a^+(\theta') \big) \big( a^-(\theta)-a^-(\theta') \big)\, , \\
  a^+(\theta)&:=(\sqrt{1 + 2 \xi (\theta)} + 1 ) \cos (\theta)\, ,\quad
  a^-(\theta):=(\sqrt{1 + 2 \xi (\theta)} - 1 ) \cos (\theta)\, .
  \end{aligned}
  \ee
By the mean value theorem we write 
 \be\label{tt12}
 \frac{ G_1(\xi)( \theta, \theta')}{\sin^2 \big( \frac{\theta'- \theta}{2} \big)} = \Big( 
\int_0^1a^+_\theta(\theta'+ \tau (\theta-\theta'))d \tau  \Big) \Big( 
\int_0^1a^-_\theta(\theta' + \tau (\theta-\theta'))d \tau  \Big)\,  
 \frac{ ( \theta'- \theta)^2}{\sin^2 \big( \frac{\theta'- \theta}{2} \big)}  \, ,
  \ee
  where the last term admits a smooth extension at $ \theta' = \theta $, that for simplicity we denote 
  in the same way. 
  For any $ \xi (\theta) \in {\cal C}^\infty $ the functions $ a^\pm (\theta )$ are in $ \mathcal{C}^\infty $ and so is the function
  in \eqref{tt12} in the variables $ (\theta, \theta' )$.  
 We now estimate the norm 
  $   \| \  \|^{k_0,\upsilon}_{{\cal C}^s} :=    \| \ \|^{k_0,\upsilon}_{{\cal C}^s (\T^{|\ST|} \times \T \times \T)}  $
 of the function
   $
   f (\lambda; \varphi, \theta, \theta') := $ $ \tfrac{ G_1(\xi)( \theta, \theta')}{\sin^2 \big( \frac{\theta'- \theta}{2} \big)} $
with $ \xi =  \xi (\lambda, \varphi, \theta ) $. 
For any $ |k| \leq k_0 $  we estimate  
the norm $\| \ \|_{{\cal C}^s (\T^{|\ST|} \times \T \times \T)} $ of 
$ \pa_\lambda^k f (\lambda; \varphi, \theta, \theta') $ separately outside the diagonal 
$ \{ |\theta'- \theta|\geq \d_0 \} $ on $ \T^{2d}  $ and  close to the diagonal 
$ \{ |\theta'- \theta| <  \d_0 \} $ on $ \R^{2d}  $. In the first case the smooth  
function $ \sin^2 {((\theta'- \theta)/2)} \geq c_0 > 0 $ is strictly positive and then  
  \begin{align}
  \| \pa_\lambda^k f (\lambda; \varphi, \theta, \theta') \|_{{\cal C}^s(\T^{|\ST|} \times \{ |\theta'- \theta|\geq \d_0 \})}
 &  \lesssim  \upsilon^{-|k|} 
 \|  G_1 (\lambda; \varphi, \theta, \theta') \|_{ {\cal C}^s (\T^{|\ST|} \times \T \times \T) }^{k_0, \upsilon} \, . \label{outdia0}
 \end{align} 
  By interpolation estimates \eqref{prod} and 
  Lemma \ref{compo_moser} we deduce, for any 
  $ \| \xi \|_{s_0}^{k_0,\upsilon} $ small enough, 
  \be\label{outdia}
\|   G_1 (\lambda; \varphi, \theta, \theta') \|^{k_0,\upsilon}_{ {\cal C}^s (\T^{|\ST|} \times \T \times \T) }
\lesssim \| \xi \|_{{\cal C}^s}^{k_0,\upsilon} \, .
 \ee
Then we estimate the norm 
$ \| \pa_\lambda^k f (\lambda; \varphi, \theta, \theta') \|_{{\cal C}^s(\T^{|\ST|} \times \{ |\theta'- \theta|< \d_0 \})}$.
By \eqref{tt12} it is sufficient to bound, for any $ \tau_1, \tau_2 \in [0, 1]$,  
  \begin{align}
 &  \| \pa_\lambda^k \big( a^+_\theta (\lambda; \varphi, \theta + \tau_1 (\theta - \theta')) 
  a^-_\theta (\lambda; \varphi, \theta + \tau_1 (\theta - \theta'))   \big) \|_{{\cal C}^s(\T^{|\ST|} \times \{ |\theta'- \theta|< \d_0 \})} 	\notag \\
  & \lesssim 
  \| \pa_\lambda^k \big( a^+_\theta (\lambda; \varphi, x) 
  a^-_\theta (\lambda; \varphi, y)   \big) \|_{{\cal C}^s(\T^{|\ST|} \times \T \times \T \})} 
 \lesssim \upsilon^{-|k|}  \| \xi \|_{{\cal C}^{s+1}}^{k_0,\upsilon}  \label{closedia}
 \end{align} 
for any   $ \| \xi \|_{s_0+1}^{k_0,\upsilon} $ small enough, 
by using \eqref{prod} and 
  Lemma \ref{compo_moser}. 
In conclusion, by \eqref{outdia} and \eqref{closedia}, we get 
   $    \| f \|^{k_0,\upsilon}_{{\cal C}^s (\T^{|\ST|} \times \T \times \T)}   
    \lesssim   \| \xi \|_{{\cal C}^{s+1}}^{k_0,\upsilon} $. 
The  second term in  \eqref{Gxi} satisfies the same bound and we deduce that
\be\label{quasfin}
\Big\| \tfrac{G(\xi)( \theta, \theta')}{\sin^2 \big( \frac{\theta'- \theta}{2} \big)}  \Big\|_{{\cal C}^{s}}^{k_0,\upsilon} 
\lesssim  \| \xi \|_{{\cal C}^{s+1}}^{k_0,\upsilon} \lesssim 
 \| \xi \|_{s+s_0+1}^{k_0,\upsilon}  \, .
\ee
Finally, by \eqref{prod} and Lemma \ref{compo_moser}, 
 \eqref{defP} and \eqref{quasfin}, we conclude that 
 $R ( \xi (\lambda; \vphi, \cdot)) (\theta, \theta')  $ 
in \eqref{kerper} satisfies 
$ \| R (\xi(\lambda; \varphi, \cdot)) (\theta, \theta')\|_{{\cal C}^s} \lesssim 
 \| \xi \|_{s+s_0+1}^{k_0,\upsilon} $. 
Recalling  Lemma   \ref{lem:Int}, the operator $ {\cal R}$ in \eqref{defcalR} satisfies \eqref{stimaRep}. 
\end{pf}

\smallskip

We now provide tame estimates for the composition of the vector field 
$X_{H_{\geq 3}}(\xi ) $ defined in  \eqref{XHsv}  with 
functions $ \xi (\vphi, \lambda) $. 
Note that 
$ d X_{H_{\geq 3}} (0)= 0 $ and  $ d^2_\xi X_{H}(\xi) = d^2_\xi X_{H_{\geq 3}}(\xi) $.

\begin{lemma}\label{lemma:xp} 
 Assume that  $  \| \xi\|_{s_0+1}^{k_0, \upsilon}\leq \delta  $ is small enough. Then 
   $X_{H_{\geq 3}}(\xi (\vphi) ) $  satisfies the following tame estimates,
  for some $ \sigma >0 $, for any $ s \geq s_0 $, 
\begin{align}
\| X_{H_{\geq 3}} (\xi)\|_s^{k_0, \upsilon} & \lesssim_s 
 \| \xi\|_{s +\sigma}^{k_0, \upsilon} \| \xi\|_{s_0 +\sigma}^{k_0, \upsilon}\, ,\label{esti:xp} \\
 \| d_\xi X_{H_{\geq 3}}  (\xi)[\widehat \xi]\|_s^{k_0, \upsilon} & \lesssim_s 
 \| \xi \|_{s_0 + \sigma}^{k_0, \upsilon}  \| \widehat \xi \|_{s + \sigma}^{k_0, \upsilon} + \| \xi\|_{s + \sigma}^{k_0, \upsilon} \| \widehat \xi \|_{s_0 + \sigma}^{k_0, \upsilon} \, ,\label{esti:xp1} \\
\| d^2_\xi X_{H_{\geq 3}}(\xi)[\widehat \xi_1, \widehat \xi_2]\|_s^{k_0, \upsilon} & \lesssim_s 
\| \widehat \xi_1 \|_{s_0 + \sigma}^{k_0, \upsilon}  \| \widehat \xi_2\|_{s + \sigma}^{k_0, \upsilon} + \| \widehat\xi_2\|_{s_0 + \sigma}^{k_0, \upsilon} \big(\|  \widehat\xi_1\|_{s + \sigma}^{k_0, \upsilon} +\| \xi\|_{s + \sigma}^{k_0, \upsilon} \| \widehat \xi_1 \|_{s_0 + \sigma}^{k_0, \upsilon}\big) \,.\label{esti:xp2}
\end{align}
\end{lemma}

\begin{pf}  
For simplicity along the proof we denote $ P := H_{\geq 3}  $. We prove the estimate \eqref{esti:xp2}. The estimates \eqref{esti:xp} and \eqref{esti:xp1} then follow by Taylor theorem since $ X_{H_{\geq 3}} (0) = 0 $ and $ d X_{H_{\geq 3}} (0) = 0 $.
By  \eqref{linVF1} we have 
$ d X_H (\xi) [\widehat \xi ] = \pa_\theta
 \big(   \big(\Omega\,  g_\gamma(\theta)+v(\xi) (\theta ) \big) \widehat \xi ( \theta)-
 \W (\xi)\, [ \widehat \xi  ] (\theta)\big)$ and therefore 
 \begin{equation}\label{2nd-diff}
d^2_\xi X_{P}(\xi)[\widehat \xi_1, \widehat \xi_2](\theta) =\pa_\theta\Big(\big(d_\xi v(\xi)[\widehat \xi_2](\theta)\big)\, \widehat \xi_1(\theta)-d_\xi \big(\W (\xi) [\widehat \xi_1]\big)[\widehat \xi_2](\theta)\Big) 
 \end{equation}
where, by \eqref{defV} and \eqref{int-op}, 
  \begin{align}
&4\pi d_\xi \big(\W (\xi) [\widehat \xi_1]\big)[\widehat \xi_2](\theta)=\int_{\T}\tfrac{ d_\xi\Ke(\xi)[\widehat \xi_2](\theta,\theta')}{\Ke(\xi)(\theta, \theta')}
\widehat \xi_1(\theta')d\theta'\, , \label{dwxi}
\\
&4\pi  d_\xi v(\xi)[\widehat \xi_2](\theta)=\int_{\T}\tfrac{  d_\xi\Ke(\xi)[\widehat \xi_2 ](\theta, \theta')}{\Ke(\xi)(\theta, \theta')}
g_1(\xi) ( \theta,\theta')d\theta'  +
\int_{\T} \ln ( \Ke(\xi)(\theta, \theta'))g_2(\xi,\widehat \xi_2)(\theta,\theta') d \theta'\,  , \label{dvxi}
\\
& g_1(\xi)(\theta,\theta'):= \pa_{\theta'}\Big[  \Big(\tfrac{1 + 2 \xi (\theta')}{1 + 2 \xi (\theta)}\Big)^{1/2} \sin (\theta' - \theta) \Big],\label{g1}\\ 
&g_2(\xi,\widehat \xi_2)(\theta,\theta'):=\pa_{\theta'} 
\Big[\tfrac{\widehat \xi_2(\theta') (1 + 2 \xi (\theta))-\widehat \xi_2(\theta) (1 + 2 \xi (\theta'))}{ (1 + 2 \xi (\theta))^{\frac32} (1 + 2 \xi (\theta'))^{\frac12} }\sin (\theta' - \theta) \Big] \label{g2}\, .
  \end{align}
Using \eqref{decoker}, \eqref{rel-dis0} and  \eqref{defP} we may write
     \begin{align}\label{frac-dxi-M}
\tfrac{ d_\xi\Ke(\xi)[\widehat \xi_2](\theta,\theta')}{\Ke(\xi)(\theta, \theta')}&=
\tfrac{  d_\xi \Ke(\xi)[\widehat \xi_2 ](\theta, \theta')}{2\sin^2 \big( \tfrac{\theta'- \theta}{2} \big) }
\Big({p (\gamma, \theta, \theta' ) }+ \tfrac{G(\xi)(\theta, \theta')}{2\sin^2 \big( \tfrac{\theta'- \theta}{2} \big) }  \Big)^{-1}\,. 
\, 
  \end{align}
Differentiating \eqref{expR} with respect to $\xi$ gives 
$$
  d_\xi \Ke(\xi)[\widehat \xi_2 ](\theta, \theta') =
2\gamma\big(f_1(\theta)-f_1(\theta')\big)\big(h_1(\theta)-h_1(\theta')\big)+2\gamma^{-1}\big(f_2(\theta)-f_2(\theta')\big)\big(h_2(\theta)-h_2(\theta')\big) 
$$
with  $ f_1(\theta) :=\widehat\xi_2(\theta){(1 + 2 \xi(\theta))^{-\frac12}}\cos\theta $, 
$ h_1(\theta):= (1 + 2 \xi(\theta))^{\frac12}\cos\theta $,
$ f_2(\theta) :=\widehat\xi_2(\theta){(1 + 2 \xi(\theta))^{-\frac12}}\sin\theta $ and 
$ h_2(\theta):= (1 + 2 \xi(\theta))^{\frac12}\sin\theta $.
By the mean value theorem we  write 
\begin{align*}
F_1(\xi,\widehat \xi_2)(\theta, \theta')&:= \tfrac{  d_\xi \Ke(\xi)[\widehat \xi_2 ](\theta, \theta')}{2\sin^2 \big( \frac{\theta'- \theta}{2} \big)} \\ &= \gamma\Big( 
\int_0^1\partial_\theta f_1(\theta'+ \tau (\theta-\theta'))d \tau  \Big) \Big( 
\int_0^1\partial_\theta h_1(\theta' + \tau (\theta-\theta'))d \tau  \Big)\,  
 \tfrac{ ( \theta'- \theta)^2}{\sin^2 \big( \frac{\theta'- \theta}{2} \big)}\\
 &+\gamma^{-1} \Big( 
\int_0^1\partial_\theta f_2(\theta'+ \tau (\theta-\theta'))d \tau  \Big) \Big( 
\int_0^1\partial_\theta h_2(\theta' + \tau (\theta-\theta'))d \tau  \Big)\,  
 \tfrac{ ( \theta'- \theta)^2}{\sin^2 \big( \frac{\theta'- \theta}{2} \big)}  \, .
\end{align*}
Arguing as for \eqref{quasfin} we conclude that 
\be\label{quasfin2}
\big\| F_1 ( \xi (\lambda; \vphi, \cdot), \widehat  \xi_2 (\lambda; \vphi, \cdot)) \big\|_{{\cal C}^{s}}^{k_0,\upsilon} 
  \lesssim 
 \| \widehat \xi_2 \|_{s_0+\sigma}^{k_0,\upsilon}\big(1+ \| \xi \|_{s+\sigma}^{k_0,\upsilon} \big)\ + \| \widehat \xi_2 \|_{s+\sigma}^{k_0,\upsilon}\big(1+ \| \xi \|_{s_0+\sigma}^{k_0,\upsilon} \big)  \, .
\ee
Moreover, according to  \eqref{defP} and \eqref{quasfin} one has
\begin{equation}\label{est:rxi2}
\big\| \big({p (\gamma, \theta, \theta' ) }+ \tfrac{G(\xi (\lambda; \vphi, \cdot)) (\theta, \theta')}{2\sin^2 \big( \frac{\theta'- \theta}{2} \big) }  \big)^{-1}\big\|_{{\cal C}^{s}}^{k_0,\upsilon} \lesssim_s  1+ \| \xi \|_{s+\sigma}^{k_0,\upsilon} \, .
\end{equation}
By \eqref{prod},  \eqref{quasfin2}-\eqref{est:rxi2}, we conclude that 
the function 
in \eqref{frac-dxi-M} satisfies 
\begin{equation}\label{est:rxi3}
\Big\| \tfrac{ d_\xi\Ke\big(\xi(\lambda; \vphi, \cdot)\big)[\widehat \xi_2(\lambda; \vphi, \cdot)](\theta,\theta')}{\Ke\big(\xi(\lambda; \vphi, \cdot)\big)(\theta, \theta')} \Big\|_{{\cal C}^s}^{k_0,\upsilon}  \lesssim_s
 \| \widehat \xi_2 \|_{s_0+\sigma}^{k_0,\upsilon}\big(1+ \| \xi \|_{s+\sigma}^{k_0,\upsilon} \big)+ \| \widehat \xi_2 \|_{s+\sigma}^{k_0,\upsilon}\big(1+ \| \xi \|_{s_0+\sigma}^{k_0,\upsilon} \big) \, .
\end{equation}
Recalling  Lemma   \ref{lem:Int}, the operator 
\be\label{defcalW}
 {\cal W}\big(\xi(\lambda; \vphi, \cdot),\widehat \xi_2(\lambda; \vphi, \cdot)\big) [q]:=\int_{\T} \tfrac{ d_\xi\Ke\big(\xi(\lambda; \vphi, \cdot)\big)[\widehat \xi_2(\lambda; \vphi, \cdot)](\theta,\theta')}{\Ke\big(\xi(\lambda; \vphi, \cdot)\big)(\theta, \theta')}q(\theta')d\theta'\,  
\ee
is in $ {\rm OPS}^{-\infty} $  and satisfies,
for any $ m,  \a \in \N_0 $,  for some constant $ \sigma (m,\a) > 0 $, 
for any $ s \geq s_0 $, 
\be\label{stimaRepW}
\begin{aligned}
\big|  {\cal W}\big(\xi(\lambda; \vphi, \cdot),\widehat \xi_2(\lambda; \vphi, \cdot)\big)  \big|_{-m, s, \a}^{k_0, \upsilon} &\lesssim_{m, s, \alpha, k_0}   
 \| \widehat \xi_2 \|_{s_0+\sigma (m,\a)}^{k_0,\upsilon}\big(1+ \| \xi \|_{s+\sigma (m,\a)}^{k_0,\upsilon} \big)\\ 
 &\qquad\quad+ \| \widehat \xi_2 \|_{s+\sigma (m,\a)}^{k_0,\upsilon}\big(1+ \| \xi \|_{s_0+\sigma (m,\a)}^{k_0,\upsilon} \big)   \,.
 \end{aligned}
\ee
On the other hand, by \eqref{prod} and Lemma \ref{compo_moser},  the functions $g_1\big(\xi(\lambda; \vphi, \cdot)\big)$, $g_2\big(\xi(\lambda; \vphi, \cdot),\widehat \xi_2(\lambda; \vphi, \cdot)\big)$ in \eqref{g1}, \eqref{g2} satisfy
\begin{align}
\big\| g_1\big(\xi(\lambda; \vphi, \cdot)\big) \big\|_{s}^{k_0,\upsilon}&\lesssim 1+ \| \xi \|_{s+\sigma}^{k_0,\upsilon} \label{est-g1}
 ,\\
\big\| g_2\big(\xi(\lambda; \vphi, \cdot),\widehat \xi_2(\lambda; \vphi, \cdot)\big)\big\|_{s}^{k_0,\upsilon}& \lesssim  \| \widehat \xi_2 \|_{s_0+\sigma}^{k_0,\upsilon}\big(1+ \| \xi \|_{s+\sigma}^{k_0,\upsilon} \big)+ \| \widehat \xi_2 \|_{s+\sigma}^{k_0,\upsilon}\big(1+ \| \xi \|_{s_0+\sigma}^{k_0,\upsilon} \big)\, .\label{est-g2}
\end{align}
By \eqref{2nd-diff}, \eqref{dwxi}, \eqref{dvxi}, \eqref{stimaRepW}, \eqref{est-g1}, \eqref{est-g2}, Lemma \ref{Prop:csm} and using \eqref{prod} we deduce  \eqref{esti:xp2}. 
\end{pf}

Now we consider the symplectic rectification map $\Phi$.

\begin{lemma}\label{lem:dphi-phi-1}
Assume that $  \|  \xi\|_{s_0+2}^{k_0, \upsilon}\leq \delta$ is small enough.
Thus,  for any $s\geq s_0$,
\begin{align}\label{bartsmo}
& \| \bar t (\xi (\vphi))\|_{s}^{k_0,\upsilon} \lesssim_s \| \xi \|_s^{k_0,s} \, , \\
&  \label{est-phi}
  \| d\Phi(\xi)[\widehat\xi]\|_{s}^{k_0, \upsilon}  \lesssim_s \| \widehat \xi \|_{s}^{k_0, \upsilon}+\| \xi \|_{s+1}^{k_0, \upsilon} \| \widehat\xi \|_{s_0}^{k_0, \upsilon}\, , 
   \\
&  \| d^2\Phi(\xi)[\widehat\xi_1,\widehat\xi_2]\|_{s}^{k_0, \upsilon} 
 \lesssim_s \| \widehat \xi_1 \|_{s_0}^{k_0, \upsilon}\| \widehat \xi_2 \|_{s}^{k_0, \upsilon}+\| \widehat \xi_2 \|_{s_0}^{k_0, \upsilon}\big(\| \widehat \xi_1 \|_{s}^{k_0, \upsilon}+\| \xi \|_{s+2}^{k_0, \upsilon} \| \widehat\xi_1 \|_{s_0}^{k_0, \upsilon}\big)\, .\label{est-phi2}
\end{align}
Assume that $ \| \wtilde \xi\|_{s_0+1}^{k_0, \upsilon}\leq \delta$ is small enough. Then,
 for any $s\geq s_0$,
\be\label{est-phi-1}
\| \Phi^{-1}( \wtilde \xi)\|_{s}^{k_0, \upsilon}  \lesssim_s \| \wtilde \xi \|_{s}^{k_0, \upsilon}\, , 
\quad 
  \| d\Phi^{-1}(\wtilde\xi)[\widehat\eta]\|_{s}^{k_0, \upsilon}  \lesssim_s \| \widehat \eta \|_{s}^{k_0, \upsilon}+\| \wtilde\xi \|_{s+1}^{k_0, \upsilon} \| \widehat\eta \|_{s_0}^{k_0, \upsilon}\, . 
\ee
\end{lemma}

\begin{pf}
The function $  \bar t $ defined in Theorem \ref{thm:FM} is  $ C^\infty $
on a small ball of $ L^2 $ with  $ \bar t (0) = 0  $ and 
bounded derivatives. Then the composition estimate 
\eqref{bartsmo} follows as in \cite{BKM1}.
Then \eqref{est-phi}-\eqref{est-phi2} follow
 by \eqref{diffePhi}, \eqref{formaJ},   \eqref{gradbart},  \eqref{bartsmo}, Lemmata \ref{lem:FlowJ}- 
\ref{flowJ2old},  Lemma \ref{lem:prod0} and the identity
\be\label{Hess-Phi}
\begin{aligned}
 d^2 \Phi(\xi)[\widehat \xi_1,\widehat \xi_2] &=d^2  {\cal J} (\xi) [\widehat \xi_1, \widehat \xi_2] \tc_2+d^2 \bar t (\xi) [\widehat \xi_1,\widehat \xi_2] \ts_2+\Pi_{2}^\bot X_{{\cal J}_2}\big(\Phi^{  \overline t (\xi) }_{{\cal J}_2} [\widehat \xi_2]\big)d \bar t (\xi) [\widehat \xi_1]\\ &\quad +\Pi_{2}^\bot X_{{\cal J}_2} (\Phi_{\cal J}^{\overline t (\xi)} (\xi))  d^2 {\bar t} (\xi)[\widehat \xi_1,\widehat \xi_2]+\Pi_{2}^\bot X_{{\cal J}_2}\big(\Phi^{  \overline t (\xi) }_{{\cal J}_2} [\widehat \xi_1]\big)
d \bar t (\xi) [\widehat \xi_2] 
\\ 
&\quad + \Pi_{2}^\bot X_{{\cal J}_2}\big(X_{{\cal J}}  (\Phi_{\cal J}^{\overline t (\xi)} (\xi)) \big) \, d \bar t (\xi) [\widehat \xi_1] d \bar t (\xi) [\widehat \xi_2] 
 \, .
 \end{aligned}
 \ee
The estimates 
\eqref{est-phi-1} follow similarly by \eqref{formainv}. 
\end{pf}

We finally  provide  tame estimates for the composition operator induced
 by the Hamiltonian vector field $ X_{\sP} = ( - \pa_I \sP,  \pa_{\vartheta} \sP, \pa_\theta \nabla_{z} \sP) $ 
 in \eqref{operatorF}. 
 
 \begin{lemma}\label{lemma quantitativo forma normale}
Let $ i (\vphi) = \varphi + \fracchi (\vphi )$ a torus embedding satisfying 
\eqref{ansatz 0}. There exists $ \s > 0$ such that, for any $ s \geq s_0 $,  
\begin{equation}\label{stime XP}
\| X_{\sP} (i)\|_s^{k_0, \upsilon}  \lesssim_s 
 1 + \| {\mathfrak I}\|_{s + \s}^{k_0, \upsilon} \, , 
\end{equation}
and for all $\widehat \imath := (\widehat \vartheta, \widehat I, \widehat z)$, for any $ s \geq s_0 $,   
\begin{align}\label{stima derivata XP}
 \| d_i X_{\sP}  (i)[\widehat \imath]\|_s^{k_0, \upsilon} & \lesssim_s 
 \| \widehat \imath \|_{s +  \s}^{k_0, \upsilon} + \| \mathfrak I\|_{s +  \s}^{k_0, \upsilon} \| \widehat \imath \|_{s_0 +  \s}^{k_0, \upsilon} \, , \\
 \label{stima derivata seconda XP}
\| d^2_i X_{\sP}(i)[\widehat \imath, \widehat \imath]\|_s^{k_0, \upsilon} & \lesssim_s 
 \| \widehat \imath\|_{s +  \s}^{k_0, \upsilon} \| \widehat \imath \|_{s_0 +  \s}^{k_0, \upsilon} + \| \mathfrak I\|_{s +  \s}^{k_0, \upsilon} (\| \widehat \imath \|_{s_0 + \sigma}^{k_0, \upsilon})^2 \,.
\end{align}
\end{lemma}
\begin{pf}
By  \eqref{cNP2} and \eqref{formaHep} we have 
\be\label{def:sp}
\sP(   \vartheta,I, z) =\e^{-3} \cK_{\geq 3}( \e {\cal J}_0,\e \Xi  (   \vartheta,I, z) ) \, ,
\ee 
and, by  \eqref{XHemu}, \eqref{aacoordinates},  we get
$$
X_{\sP}(   \vartheta,I, z) = \e^{-2} \left(
\begin{array}{c}
   -  [\partial_I  {\mathtt v}^\intercal (\vartheta, I)]^\top\nabla_{\wtilde u}  \cK_{\geq 3}( \e {\cal J}_0,\e \Xi  (   \vartheta,I, z) )    \\
   \big[\partial_\vartheta {\mathtt v}^\intercal (\vartheta, I)\big]^\top \nabla_{\wtilde u} \cK_{\geq 3}( \e {\cal J}_0,\e \Xi  (   \vartheta,I, z) )    \\   \Pi_{\ST,2}^\bot \pa_\theta \nabla_{\wtilde u} \cK_{\geq 3}( \e {\cal J}_0,\e \Xi  (   \vartheta,I, z) )   
\end{array}
\right) 
$$
where  $ \Pi_{\ST,2}^\bot $ is the $ L^2 $-projector on $  \acca_{\ST,2}^\bot $ defined in \eqref{Hbot}. By \eqref{nablatuK}  we have 
\be\label{Kgeq3}
\nabla_{\wtilde u} \cK_{\geq 3}( \e {\cal J}_0,\e \Xi  (   \vartheta,I, z) ) = \Pi_2^\perp \nabla_{\wtilde \xi}   K_{\geq 3}\big( \e {\cal J}_0\tc_2+\e \Xi  (   \vartheta,I, z)\big) 
\ee
and 
$\nabla_{\wtilde \xi}   K_{\geq 3} (\wtilde \xi ) = \partial_\theta^{-1}X_{K_{\geq 3}} (\wtilde \xi ) $
where $X_{K_{\geq 3}} (\wtilde \xi )  $ is expanded as in \eqref{XKXP}.  Moreover,   using  \eqref{XKXP}   applying Taylor formula yields
\be\label{Id2-ham}
 \begin{aligned}
  X_{ K_{\geq 3}}( \wtilde \xi )  & = 
 d\Phi(\Phi^{-1}( \wtilde \xi))\Big[ X_{H_{\geq 3}}(\Phi^{-1}( \wtilde \xi) )
+ \partial_\theta {\bf \Omega}(\gamma) \int_0^1 d \Phi^{-1}\big(\tau\wtilde\xi\big)[\widetilde \xi]d\tau\Big]\\  &\quad+ \int_0^1 d^2 \Phi\big(\tau\Phi^{-1}( \wtilde \xi) \big)\big[\Phi^{-1}( \wtilde \xi), \partial_\theta {\bf \Omega}(\gamma) \wtilde \xi \, \big]d\tau \, . 
\end{aligned}
\ee
By the Moser composition Lemma \ref{compo_moser} we have,  for any $ s \geq s_0 $,  
\be\label{est:v}
\|\partial_\vartheta^\alpha\partial_I^\beta {\mathtt v}^\intercal (\vartheta( \cdot ), I( \cdot )) \|_s^{k_0,\upsilon} \lesssim_s 1+ \| {\mathfrak I}\|_{s }^{k_0, \upsilon}, \quad \forall \alpha,\beta\in \mathbb{N}_0^{|\ST|}\, , \quad |\alpha|+|\beta|\leq 3 \, ,
\ee
and 
the function $  \wtilde \xi_\e (\vphi, \theta ) := \e {\cal J}_0\tc_2+\e \Xi  ( i(\varphi)) $ satisfies 
$ \| \wtilde \xi_\e \|_{s}^{k_0, \upsilon} \lesssim_s \e |{\cal J}_0| 
+\e (1 + \| \fracchi \|_{s}^{k_0,\upsilon} )$. 
This bound and   \eqref{Kgeq3},  \eqref{Id2-ham},  Lemmata \ref{lemma:xp}-\ref{lem:dphi-phi-1}, \eqref{est:v}, 
imply  \eqref{stime XP}-\eqref{stima derivata seconda XP}.
\end{pf}

\section{Almost approximate inverse}\label{sezione approximate inverse}

In order to prove Theorem \ref{main theorem} we  implement a convergent Nash-Moser scheme to construct a solution of 
$ {\cal F}(i, \tg) = 0 $, where ${\cal F}(i, \tg) $ is the nonlinear operator defined in 
\eqref{operatorF}. For this aim we need to construct an \emph{almost-approximate right inverse} 
of the linearized operator 
$$
d_{i, \tg} {\cal F}(i_0, \tg_0 )[\widehat \imath \,, \widehat \tg ] =
\Dom \widehat \imath - d_i X_{\sK_{\tg}} ( i_0 (\vphi) ) [\widehat \imath ] - (\widehat \tg,0, 0 )
$$
where 
the torus  $ i_0 (\vphi) = (\vartheta_0 (\vphi), I_0 (\vphi), z_0 (\vphi)) $ is reversible, i.e. satisfies \eqref{parity solution}.
Here,  following \cite{BertiMontalto}, 
the adjective `approximate"  refers
to the presence of a remainder which is zero at an exact solution and  the adjective 
`almost" 
 refers to remainders which are small as $O(N_n^{-a})$  for some
$a>0$  (in suitable norms) at the $n$-th step of the Nash-Moser iteration.

\smallskip

 We  closely follow 
 the  strategy in \cite{BB13,BBM-auto,BertiMontalto,BBHM,BKM1,BFM,BFM1}, to 
reduce the problem to 
almost-invert a quasi-periodic operator acting on the normal subspace, see
\eqref{Lomega def} and the 
assumption (AI) below it. 
Thus we will be short, referring to the above papers for details.  
We consider the pull-back $ 1$-form $ i_0^* \Lambda $,  
where $ \Lambda $ is the 1-form in \eqref{Lambda 1 form},   
\begin{equation}\label{coefficienti pull back di Lambda}
 i_0^* \Lambda = {\mathop \sum}_{k = 1}^{|\ST|} a_k (\vphi) d \vphi_k \,, 
 \ \   a_k(\vphi) :=  \big( [\pa_\vphi \vartheta_0 (\vphi)]^\top I_0 (\vphi)  \big)_k 
+\tfrac12 ( \partial_\theta^{-1} z_0(\vphi), \partial_{\vphi_k} z_0(\ph) )_{L^2(\T)} 
\end{equation}
and 
\begin{equation} \label{def Akj} 
 i_0^* {\cal W} = d \, i_0^* \Lambda = {\mathop\sum}_{1 \leq k < j \leq {|\ST|}} A_{k j}(\vphi) d \vphi_k \wedge d \vphi_j\,, \quad  A_{k j} (\vphi) := 
\partial_{\vphi_k} a_j(\ph) - \partial_{\vphi_j} a_k(\ph) \, . 
\end{equation} 
Let define the ``error function''
\begin{equation} \label{def Zetone}
Z (\vphi) :=  (Z_1, Z_2, Z_3) (\vphi) := {\cal F}(i_0, \tg_0) (\vphi) =
\om \cdot \pa_\vphi i_0(\vphi) - X_{\sK_{\tg_0}}(i_0(\vphi)) \, .
\end{equation}
Along this section we  assume the following hypothesis, 
which will be verified by the approximate solutions  of the Nash-Moser iteration. 

\begin{itemize}
\item {\sc Ansatz.} 
The map $ (\omega, \gamma) \mapsto \fracchi_0 (\omega, \gamma) :=  i_0(\ph; \om, \gamma) - (\ph,0,0) $ 
is 
defined for all 
the parameters $(\omega, \gamma) \in \R^{|\ST|} \times [\g_1, \g_2] $, 
and for some $ \perd := \perd (\t, \ST) >  0 $,  $\upsilon \in (0, 1)$,   
\begin{equation}\label{ansatz 0}
\| {\mathfrak I}_0  \|_{s_0+ \perd}^{k_0, \upsilon} +  | \tg_0 - \omega|^{k_0, \upsilon} \leq C
\e \upsilon^{- 1} \, .
\end{equation}
\end{itemize}

In the next lemma,  following \cite{BB13,BBM-auto,BertiMontalto},  
 we first modify  the approximate torus $ i_0 $ to obtain a nearby isotropic torus $ i_\d $, namely  the pull-back $ 1$-form $ i_\delta^* \Lambda $ is closed. 
 We denote   
 $ \Delta_\vphi := {\mathop \sum}_{k=1}^{|\ST|} \partial_{\vphi_k}^2 $. 
 In the sequel
$ \s := \s(\ST, \tau,k_0) > 0 $ will denote possibly different larger ``loss of derivatives"  constants.

\begin{lemma}\label{toro isotropico modificato} {\bf (Isotropic torus)} 
The torus $ i_\delta(\vphi) := (\vartheta_0(\vphi), I_\delta(\vphi), z_0(\vphi) ) $ defined by 
\begin{equation}\label{y 0 - y delta}
I_\d (\varphi) := I_0 (\varphi) -  [\pa_\ph \vartheta_0(\vphi)]^{-\top}  \rho(\vphi) \, , \qquad 
\rho_j(\vphi) := \Delta_\vphi^{-1} {\mathop\sum}_{ k = 1}^{|\ST|} \partial_{\vphi_k}  A_{k j}(\vphi)  \, , 
\end{equation}
 is {\it isotropic}. 
There is $ \s  > 0 $ such that, for any $ s \geq s_0 $,  
\begin{align} \label{2015-2}
\| I_\delta - I_0 \|_s^{k_0, \upsilon} & \lesssim_s \| \fracchi_0 \|_{s+1}^{k_0,\upsilon} \\
\label{stima y - y delta}
\| I_\delta - I_0 \|_s^{k_0, \upsilon} 
& \lesssim_s  \upsilon^{-1} \big(\| Z \|_{s + \s}^{k_0, \upsilon} + 
\| Z \|_{s_0 + \s}^{k_0, \upsilon} \|  {\mathfrak I}_0 \|_{s + \s}^{k_0, \upsilon} \big) 
\\
\label{stima toro modificato}
\| {\cal F}(i_\delta, \tg_0) \|_s^{k_0, \upsilon}
& \lesssim_s  \| Z \|_{s + \s}^{k_0, \upsilon}  +  \| Z \|_{s_0 + \s}^{k_0, \upsilon} \|  {\mathfrak I}_0 \|_{s + \s}^{k_0, \upsilon} \\
\label{derivata i delta}
\| \pa_i [ i_\d][ \widehat \imath ] \|_s^{k_0, \upsilon} & \lesssim_s 
\| \widehat \imath \|_{s+1}^{k_0, \upsilon} +  
\| {\mathfrak I}_0\|_{s + \s}^{k_0, \upsilon} \| \widehat \imath  \|_{s_0}^{k_0, \upsilon} \, .
\end{align}
 \end{lemma}

In order to find an approximate inverse of the linearized operator $d_{i, \tg} {\cal F}(i_\delta )$, 
we introduce  the symplectic diffeomorpshim 
$ G_\delta : (\phi, y, w) \to (\vartheta, I, z)$ of the phase space $\T^{|\ST|} \times \R^{|\ST|} \times {\frak H}_{\ST,2}^\bot$ defined by
\begin{equation}\label{trasformazione modificata simplettica}
\begin{pmatrix}
\vartheta \\
I \\
z
\end{pmatrix} = G_\delta \begin{pmatrix}
\phi \\
y \\
w
\end{pmatrix} := 
\begin{pmatrix}
\!\!\!\!\!\!\!\!\!\!\!\!\!\!\!\!\!\!\!\!\!\!\!\!\!\!\!\!\!\!\!\!\!
\!\!\!\!\!\!\!\!\!\!\!\!\!\!\!\!\!\!\!\!\!\!\!\!\!\!\!\!\!\!\!\!\!
\!\!\!\!\!\!\!\!\!\!\!\!\!\!\!\!\!\!\!\!\!\!\!\!\!\!\!\!\!\!\!\! \! \vartheta_0(\phi) \\
\quad I_\delta (\phi) + [\pa_\phi \vartheta_0(\phi)]^{-\top} y - \big[ (\pa_\vartheta \wtilde{z}_0) (\vartheta_0(\phi)) \big]^\top \partial_\theta^{-1} w \\
\!\!\!\!\!\!\!\!\!\!\!\!\!\!\!\!\!\!\!\!\!\!\!\!\!\!\!\!\!
\!\!\!\!\!\!\!\!\!\!\!\!\!\!\!\!\!\!\!\!\!\!\!\!\!\!\!\!\!\!
\!\!\!\!\!\!\!\!\!\!\!\!\!\!\!\!\!\!\!\!\!\!\!\!\!\!\!\!\!  z_0(\phi) + w
\end{pmatrix} 
\end{equation}
where $ \wtilde{z}_0 (\vartheta) := z_0 (\vartheta_0^{-1} (\theta))$. 
It is proved in \cite{BB13} that $ G_\delta $ is symplectic, because  the torus $ i_\d $ is isotropic 
(Lemma \ref{toro isotropico modificato}).
In the new coordinates,  $ i_\delta $ is the trivial embedded torus
$ (\phi , y , w ) = (\phi , 0, 0 ) $.  Under the symplectic change of variables $ G_\d $ the Hamiltonian 
vector field $ X_{\sK_{\tg}} $ (the Hamiltonian $  \sK_{\tg}  $ is defined in \eqref{H alpha}) changes into 
\be\label{new-Hamilt-K}
X_{\tK_{\tg}} = (D G_\d)^{-1} X_{\sK_{\tg}} \circ G_\d \qquad {\rm where} \qquad \tK_\tg := 
\sK_{\tg} \circ G_\d  \, .
\ee
By \eqref{parity solution} the transformation $ G_\d $ is also reversibility preserving and so $ \tK_\tg $ is reversible, i.e. 
$ \tK_\tg \circ \vec {\cal S} = \tK_\tg $.   
The Taylor expansion of $ \tK_\tg $ at the trivial torus $ (\phi , 0, 0 ) $ is 
\be
\begin{aligned} 
\tK_\tg (\phi, y , w)
&  =  \tK_{00}(\phi, \tg) + \tK_{10}(\phi, \tg) \cdot y + (\tK_{0 1}(\phi, \tg), w)_{L^2(\T)} + 
\tfrac12 \tK_{2 0}(\phi) y \cdot y 
 \\ & 
\quad +  \big( \tK_{11}(\phi) y , w \big)_{L^2(\T)} 
+ \tfrac12 \big(\tK_{02}(\phi) w, w \big)_{L^2(\T)} + \tK_{\geq 3}(\phi, y, w)  
\label{KHG}
\end{aligned}
\ee
where $ \tK_{\geq 3} $ collects the terms at least cubic in the variables $ (y, w )$.
The Taylor coefficient $\tK_{00}(\phi, \tg) \in \R $,  
$\tK_{10}(\phi, \tg) \in \R^{|\ST|} $,  
$\tK_{01}(\phi, \tg) \in {\frak H}_{\ST,2}^\bot$, 
$\tK_{20}(\phi) $ is a ${|\ST|} \times {|\ST|}$ real matrix, 
$\tK_{02}(\phi)$ is a linear self-adjoint operator of $ {\frak H}_{\ST,2}^\bot $ and 
$\tK_{11}(\phi) \in {\cal L}(\R^{|\ST|}, {\frak H}_{\ST,2}^\bot )$. 
The coefficients $ \tK_{00} $, $ \tK_{10} $, $\tK_{01} $  in the Taylor expansion \eqref{KHG} vanish on an exact solution (i.e. $ Z  = 0 $), and 
$ \partial_\tg \tK_{10} \approx {\rm Id} $. 
The following lemma holds.

\begin{lemma} \label{lemma:Kapponi vari}
(Lemmata 5.6 and 5.7 in \cite{BertiMontalto})
There is $ \sigma > 0 $ such that for any $ s \geq s_0 $,  
\begin{align} \notag 
& \|  \partial_\phi \tK_{00}(\cdot, \tg_0) \|_s^{k_0, \upsilon} 
+ \| \tK_{10}(\cdot, \tg_0) - \om  \|_s^{k_0, \upsilon} +  \| \tK_{0 1}(\cdot, \tg_0) \|_s^{k_0, \upsilon} 
\lesssim_s  \| Z \|_{s + \s}^{k_0, \upsilon} +  \| Z \|_{s_0 + \s}^{k_0, \upsilon} \| {\mathfrak I}_0 \|_{s + \s}^{k_0, \upsilon} \\
& \| \partial_\tg \tK_{00}\|_s^{k_0, \upsilon} + \| \partial_\tg \tK_{10} - {\rm Id} \|_s^{k_0, \upsilon} + \| \partial_\tg \tK_{0 1}\|_s^{k_0, \upsilon} \lesssim_s   \| \fracchi_0\|_{s + \sigma}^{k_0, \upsilon} \, , 
\quad   \| \tK_{20}  \|_s^{k_0, \upsilon} \lesssim_s \e \big( 1 + \| \fracchi_0\|_{s + \s}^{k_0, \upsilon} \big) \, , \notag  \\ 
& \| \tK_{11} y \|_s^{k_0, \upsilon} 
\lesssim_s \e \big(\| y \|_s^{k_0, \upsilon}
+ \| \fracchi_0 \|_{s + \sigma}^{k_0, \upsilon}  
\| y \|_{s_0}^{k_0, \upsilon} \big) \, , \quad  \| \tK_{11}^\top w \|_s^{k_0, \upsilon}
\lesssim_s \e \big(\| w \|_{s + 2}^{k_0, \upsilon}
+  \| \fracchi_0 \|_{s + \sigma}^{k_0, \upsilon}
\| w \|_{s_0 + 2}^{k_0, \upsilon} \big)\, . \notag
\end{align}
\end{lemma}

Under the linear change of variables 
\begin{equation}\label{DGdelta}
D G_\delta(\vphi, 0, 0) 
\begin{pmatrix}
\widehat \phi \, \\
\widehat y \\
\widehat w
\end{pmatrix} 
:= 
\begin{pmatrix}
\pa_\phi \vartheta_0(\vphi) & 0 & 0 \\
\pa_\phi I_\delta(\vphi) & [\pa_\phi \vartheta_0(\vphi)]^{-\top} & 
- [(\pa_\theta \wtilde{z}_0)(\vartheta_0(\vphi))]^\top \partial_\theta^{-1} \\
\pa_\phi z_0(\vphi) & 0 & I
\end{pmatrix}
\begin{pmatrix}
\widehat \phi \, \\
\widehat y \\
\widehat w
\end{pmatrix} 
\end{equation}
the linearized operator  $d_{i, \tg}{\cal F}(i_\delta )$ is transformed
 (approximately) 
  into the one  obtained when one linearizes 
the Hamiltonian system  associated to $ \tK_\tg (\phi, y , w) $ in \eqref{KHG} 
at $(\phi, y , w ) = (\vphi, 0, 0 )$,
differentiating also in $ \tg $ at $ \tg_0 $, and changing $ \partial_t \rightsquigarrow \Dom $, 
namely 
\begin{equation}\label{lin idelta}
\begin{pmatrix}
\widehat \phi  \\
\widehat y    \\ 
\widehat w \\
\widehat \tg 
\end{pmatrix} \mapsto
\begin{pmatrix}
\Dom \widehat \phi +\partial_\phi \tK_{10}(\vphi)[\widehat \phi \, ] + \partial_\tg \tK_{10}(\vphi)[\widehat \tg] +
\tK_{2 0}(\vphi)\widehat y + \tK_{11}^\top (\vphi) \widehat w \\
 \Dom  \widehat y - \partial_{\phi\phi} \tK_{00}(\vphi)[\widehat \phi] - 
 \partial_\phi \partial_\tg  \tK_{00}(\vphi)[\widehat \tg] - 
[\partial_\phi \tK_{10}(\vphi)]^\top \widehat y - 
[\partial_\phi  \tK_{01}(\vphi)]^\top \widehat w   \\ 
\Dom  \widehat w -\partial_\theta
\{ \partial_\phi \tK_{01}(\vphi)[\widehat \phi] + \partial_\tg \tK_{01}(\vphi)[\widehat \tg] + \tK_{11}(\vphi) \widehat y + \tK_{02}(\vphi) \widehat w \}
\end{pmatrix} \! .  \hspace{-5pt}
\end{equation}
Here the transposed  operators
$ \partial_{\phi}\tK_{01}^\top $,  $ \tK_{11}^\top  : {{\frak H}_{\ST,2}^\bot \to \R^{|\ST|}} $ are defined by the 
duality relation $ ( \partial_{\phi} \tK_{01} [\hat \phi ],  w)_{L^2_x}  = \hat \phi \cdot [\partial_{\phi}\tK_{01}]^\top w  $,
$ \forall \hat \phi \in \R^{|\ST|}, w \in {\frak H}_{\ST,2}^\bot $, 
and similarly for $ \tK_{11} $. 

 By \eqref{DGdelta}, \eqref{ansatz 0}, \eqref{2015-2}, 
 the induced composition operator satisfies:   for all $ \widehat \imath := (\widehat \phi, \widehat y, \widehat w) $, $ s \geq s_0 $, 
\begin{gather} \label{DG delta}
\|DG_\delta(\vphi,0,0) [\widehat \imath] \|_s^{k_0, \upsilon} + \|DG_\delta(\vphi,0,0)^{-1} [\widehat \imath] \|_s^{k_0, \upsilon} 
\lesssim_s \| \widehat \imath \|_{s}^{k_0, \upsilon} +  \| {\mathfrak I}_0 \|_{s + \s}^{k_0, \upsilon}  \| \widehat \imath \|_{s_0}^{k_0, \upsilon}\,,
\\ 
\!\!\!\!\! \| D^2 G_\delta(\vphi,0,0)[\widehat \imath_1, \widehat \imath_2] \|_s^{k_0, \upsilon} 
\! \lesssim_s \! \| \widehat \imath_1\|_s^{k_0, \upsilon}  \| \widehat \imath_2 \|_{s_0}^{k_0, \upsilon} 
+ \| \widehat \imath_1\|_{s_0}^{k_0, \upsilon}  \| \widehat \imath_2 \|_{s}^{k_0, \upsilon} 
+  \| {\mathfrak I}_0  \|_{s + \s}^{k_0, \upsilon} \|\widehat \imath_1 \|_{s_0}^{k_0, \upsilon}  \| \widehat \imath_2\|_{s_0}^{k_0, \upsilon}.  \label{DG2 delta} 
\end{gather}
In order to construct an "approximate" inverse of \eqref{lin idelta} it is sufficient to prove 
that the operator  
\be\label{Lomega def}
{\cal L}_\omega := \Pi_{\ST,2}^\bot \big(\Dom   - 
\partial_\theta \tK_{02}(\vphi) \big)_{|{{\frak H}_{\ST,2}^\bot}} 
\ee
is "almost-invertible" up to remainders of size $O(N_{\tn - 1}^{- \mathtt a})$ (see precisely \eqref{stima R omega corsivo}) where 
 \be\label{NnKn}
 N_\tn := K_\tn^p \, , \ p > 0 \, , 
  \quad \forall \tn \geq 0 \quad \text{and} \quad 
K_\tn := K_0^{\chi^{\tn}} \, , \quad \chi := 3/ 2  \, , 
\ee
are the scales of the 
Nash-Moser iteration in Section \ref{sec:NM}.
Set $ {\frak H}_{\ST,2}^{s,\bot} (\T^{{|\ST|}+1}) :=  H^{s}(\T^{{|\ST|}+1}) 
\cap {\frak H}_{\ST,2}^\bot $.  

\begin{itemize}
\item {\sc Almost-invertibility assumption (AI).} 
{\it There exists a subset $ \tLm_o  \subset \mathtt{DC}(\upsilon, \tau) \times [\g_1, \g_2] $ such that, 
for all $ (\omega, \g) \in  \tLm_o  $ the operator $ {\cal L}_\omega $ in \eqref{Lomega def} may be decomposed as
\be\label{inversion assumption}
{\cal L}_\omega  
= {\cal L}_\omega^< + {\cal R}_\omega + {\cal R}_\omega^\bot+{\cal R}_\omega^Z
\ee
where $ {\cal L}_\omega^< $ is invertible. More precisely, there exist constants 
$  K_0, \sigma, \perd(\mathtt b), \mathtt a > 0$ such that for any $s_0 \leq s \leq S$,  the operators ${\cal R}_\omega$, ${\cal R}_\omega^\bot $ and ${\cal R}_\omega^Z$ satisfy the estimates  
\begin{align}  \label{stima R omega corsivo}
\|{\cal R}_\omega h \|_s^{k_0, \upsilon} & 
\lesssim_S  \e {\upsilon^{- 1}} N_{\tn - 1}^{- {\mathtt a}}\big( \|  h \|_{s + \sigma}^{k_0, \upsilon} +  \| \fracchi_0 \|_ {s + \perd (\mathtt b)  + \sigma }^{k_0, \upsilon} \| h \|_{s_0 + \sigma}^{k_0, \upsilon} \big)\, , \\
\label{stima R omega bot corsivo bassa}
\| {\cal R}_\omega^\bot h \|_{s_0}^{k_0, \upsilon} & \lesssim_S
K_\tn^{- b} \big( \| h \|_{s_0 + b 
+ \sigma}^{k_0, \upsilon} + 
\| \fracchi_0 \|_ { s_0  + \perd (\mathtt b)  + \sigma +b }^{k_0, \upsilon}  \|  h \|_{s_0 + \sigma}^{k_0, \upsilon}\big)\,, 
\qquad \forall b > 0\,, \\
\label{stima R omega bot corsivo alta}
\| {\cal R}_\omega^\bot h \|_s^{k_0, \upsilon} & \lesssim_S 
 \|  h \|_{s + \sigma}^{k_0, \upsilon} + \| \fracchi_0 \|_ {s  + \perd (\mathtt b) + \sigma}^{k_0, \upsilon} \| h \|_{s_0 + \sigma}^{k_0, \upsilon}
 \,, \\
\label{stima R omega Z}
\| {\cal R}_\omega^Z h \|_{s}^{k_0, \upsilon} & \lesssim_S  \big( \| {Z}  \|_{s+\sigma}^{k_0,\upsilon} 
+ 
\| {Z}  \|_{s_0+\sigma}^{k_0,\upsilon}
 \| \fracchi_0 \|_{s+\sigma}^{k_0,\upsilon}\big)\|  h \|_{s_0+\sigma}^{k_0,\upsilon}  
+ \| {Z}  \|_{s_0+\sigma}^{k_0,\upsilon}\|  h \|_{s+\sigma}^{k_0,\upsilon}
  \,.
\end{align}
Moreover, for any antireversible function
$ g \in {\frak H}_{\ST,2}^{s+\sigma,\bot} (\T^{{|\ST|}+1}) $, 
any $(\omega, \gamma) \in \Lambda_o$, there 
is a reversible 
solution $ h :=  ({\cal L}_\om^<)^{- 1} g  \in {\frak H}_{\ST,2}^{s,\bot} (\T^{{|\ST|}+1}) $  
of the linear equation $ {\cal L}_\om^< h = g $. There is an extension of the inverse operator (that we denote in the same way) to the whole 
$ \R^{|\ST|} \times [\g_1, \g_2]$ satisfying for all $s_0 \leq s \leq S$ the tame estimate
\begin{equation}\label{tame inverse}
\| ({\cal L}_\om^<)^{- 1} g \|_s^{k_0, \upsilon} \lesssim_S  \upsilon^{-1} 
\big(  \| g \|_{s + \sigma}^{k_0, \upsilon} + 
 \| {\mathfrak I}_0 \|_{s + \perd({\mathtt b}) + \sigma}^{k_0, \upsilon}  \|g \|_{s_0 + \sigma}^{k_0, \upsilon}  \big) \,.
\end{equation}}
\end{itemize}

The above assumption will be  verified in Sections \ref{reduction} and \ref{sec:redu}.

In order to find an almost-approximate inverse of the linear operator in \eqref{lin idelta} 
(and so of $ d_{i, \tg} {\cal F}(i_\d) $),
it is sufficient to invert the operator
\begin{equation}\label{operatore inverso approssimato} 
{\mathbb D} [\widehat \phi, \widehat y, \widehat w, \widehat \tg ] := 
  \begin{pmatrix}
\Dom \widehat \phi + \partial_\tg \tK_{10}(\vphi)[\widehat \tg] + 
\tK_{20}(\vphi) \widehat y  + \tK_{11}^\top(\vphi) \widehat w\\
\Dom  \widehat y - \partial_\phi \partial_\tg  \tK_{00}(\vphi)[\widehat \tg] \\
({\cal L}_\omega^{<}) \widehat w  -\partial_\theta\partial_\tg \tK_{01}(\vphi)[\widehat \tg] 
-\partial_\theta\tK_{11}(\vphi)\widehat y  
\end{pmatrix}
\end{equation}
which is obtained by neglecting in \eqref{lin idelta} 
the terms $ \partial_\phi \tK_{10} $, $ \partial_{\phi \phi} \tK_{00} $, $ \partial_\phi \tK_{00} $, 
$ \partial_\phi \tK_{01} $ 
(which vanish at an exact solution by Lemma \ref{lemma:Kapponi vari})
and replacing $ {\cal L}_\omega  $ with ${\cal L}_\omega^< $ , cfr. \eqref{inversion assumption}. 
 The following result is proved in a standard way in \cite{BBHM}.
\begin{proposition}\label{prop: ai}
Assume \eqref{ansatz 0} (with $ \perd = \perd(\mathtt b) + \sigma$) and \eqref{tame inverse}. 
Then, for all $(\om, \g) \in {\mathtt \Lambda}_o $, for all $ g := (g_1, g_2, g_3) $ 
satisfying  
the reversibility property 
\begin{equation}\label{parita g1 g2 g3}
g_1(\vphi) = g_1(- \vphi)\,,\quad g_2(\vphi) = - g_2(- \vphi)\,,\quad g_3(\vphi) = - ({\cal S} g_3)(- \vphi) 
\end{equation}
there exists a unique solution 
$ {\mathbb D}^{-1} g := (\widehat \phi, \widehat y, \widehat w, \widehat \tg ) $ of 
$ {\mathbb D} [\widehat \phi, \widehat y, \widehat w, \widehat \tg]  
= (g_1, g_2, g_3)^\top $
which satisfies \eqref{parity solution} and,  for any $s_0 \leq s \leq S $, 
$ \| {\mathbb D}^{-1} g \|_s^{k_0, \upsilon}
\lesssim_S \upsilon^{-1} \big( \| g \|_{s + \sigma }^{k_0, \upsilon} 
+  \| {\mathfrak I}_0  \|_{s + \perd(\mathtt b) + \sigma}^{k_0, \upsilon}
 \| g \|_{s_0 + \sigma}^{k_0, \upsilon}  \big) $.
\end{proposition}
Finally we deduce (see e.g.  Theorem 5.6 in \cite{BBHM}) that the operator 
\begin{equation}\label{definizione T} 
{\bf T}_0 := {\bf T}_0(i_0) := (D { \widetilde G}_\delta)(\vphi,0,0) \circ {\mathbb D}^{-1} \circ (D G_\delta) (\vphi,0,0)^{-1}
\end{equation}
is an almost-approximate right  inverse for $d_{i,\tg} {\cal F}(i_0)$ where
$ \widetilde{G}_\delta (\phi, y, w, \tg) := $  $ \big( G_\delta (\phi, y, w), \tg \big) $ 
 is the identity on the $ \tg $-component. 
Let $ \| (\phi, y, w, \tg) \|_s^{k_0, \upsilon} := $ $  \max \{  \| (\phi, y, w) \|_s^{k_0, \upsilon}, 
$ $ | \tg |^{k_0, \upsilon}  \} $.

\begin{theorem}  \label{thm:stima inverso approssimato}
{\bf (Almost-approximate inverse)}
Assume the inversion assumption (AI), that is \eqref{inversion assumption}-\eqref{tame inverse}.
Then there exists $ \bar \sigma := \bar \sigma(\tau,\ST, k_0) > 0 $ such that, 
if \eqref{ansatz 0} holds with $ \perd = \perd (\mathtt b) + \bar \sigma $, then for all $ (\om, \g) \in {\mathtt \Lambda}_o $, 
for all $ g := (g_1, g_2, g_3) $ satisfying \eqref{parita g1 g2 g3},  
the operator $ {\bf T}_0 $ defined in \eqref{definizione T} satisfies,  for all $s_0 \leq s \leq S $, 
\begin{equation}\label{stima inverso approssimato 1}
\| {\bf T}_0 g \|_{s}^{k_0, \upsilon} 
\lesssim_S  \upsilon^{-1}  \big(\| g \|_{s + \bar \sigma}^{k_0, \upsilon}  
+  \| {\mathfrak I}_0 \|_{s + \perd({\mathtt b}) +  \bar \sigma }^{k_0, \upsilon}
\| g \|_{s_0 + \bar \sigma}^{k_0, \upsilon}  \big)\, .
\end{equation}
 Moreover  ${\bf T}_0(i) $ is an almost-approximate  
 right 
 inverse of $d_{i, \tg} 
 {\cal F}(i_0)$. More precisely,
\begin{equation}\label{splitting per approximate inverse}
d_{i , \tg} {\cal F} (i_0) \circ {\bf T}_0(i_0) 
- {\rm Id} = {\cal P} + {\cal P}_\omega  + {\cal P}_\omega^\bot
\end{equation}
where the operators ${\cal P}$, ${\cal P}_\omega$,  ${\cal P}_\omega^\bot$  satisfy the following estimates 
\begin{align}
\| {\cal P} g \|_{s}^{k_0,\upsilon} & \lesssim_S \upsilon^{-1 } \big( \| {Z}  \|_{s+ \bar \sigma}^{k_0,\upsilon} 
+ 
\| {Z}  \|_{s_0+ \bar \sigma}^{k_0,\upsilon}
 \| \fracchi_0 \|_{s+\perd(\tb)+\bar \sigma}^{k_0,\upsilon}\big)\|  h \|_{s_0+
 \bar \sigma}^{k_0,\upsilon}  
+ \upsilon^{-1 } \| {Z}  \|_{s_0+ \bar \sigma}^{k_0,\upsilon}\|  h \|_{s+ \bar \sigma}^{k_0,\upsilon}  \, ,\label{stima inverso approssimato 2} \\
\| {\cal P}_\omega g \|_{s}^{k_0,\upsilon} & 
\lesssim_S  \e \upsilon^{- 2} N_{\tn - 1}^{- {\mathtt a}} \big( \|  g \|_{s + \bar \sigma}^{k_0,\upsilon} + \| \fracchi_0 \|_{s  + \perd({\mathtt b}) + \bar \sigma}^{k_0,\upsilon} \|  g \|_{s_0 + \bar \sigma}^{k_0,\upsilon}\big)    \,, \label{stima cal G omega}  \\
\| {\cal P}_\omega^\bot g\|_{s_0}^{k_0,\upsilon} & \lesssim_{S, b} 
\upsilon^{- 1} K_\tn^{- b } \big( \| g \|_{s_0 + \bar \sigma + b }^{k_0,\upsilon} +
\| \fracchi_0 \|_{s_0
+  \perd({\mathtt b}) + \bar \sigma  +b    }^{k_0,\upsilon} \big \| g \|_{s_0 +\bar \sigma}^{k_0,\upsilon} \big)\,,\,\,
\quad\forall b > 0 \, ,    \label{stima cal G omega bot bassa} 
 \\
\| {\cal P}_\omega^\bot g\|_{s}^{k_0,\upsilon} & \lesssim_{S} 
\upsilon^{- 1}  \big( \| g \|_{s + \bar \sigma  }^{k_0,\upsilon} +
\| \fracchi_0 \|_{s
+  \perd({\mathtt b}) + \bar \sigma  }^{k_0,\upsilon} \big \| g \|_{s_0 +\bar \sigma}^{k_0,\upsilon} \big) \, .    \label{stima cal G omega-bot} 
\end{align}
\end{theorem}

The next two sections are devoted to prove the  assumption (AI), namely \eqref{inversion assumption}-\eqref{tame inverse}.

\section{First reduction of the normal operator}\label{reduction}

We now write an explicit  expression of the linear operator 
$ {\cal L}_\om$ defined in \eqref{Lomega def}.  

\begin{lemma} \label{thm:Lin+FBR}
The Hamiltonian operator  $ {\cal L}_\om$ defined in \eqref{Lomega def}, acting in the normal 
subspace $ {\frak H}_{\ST,2}^\bot $,  has the form  
\begin{equation}\label{K 02}
{\cal L}_\om = \Pi_{\ST,2}^\bot ( {\cal L}_K + \e \pa_\theta {\cal R} )_{|{\frak H}_{\ST,2}^\bot}
\end{equation}
where:  
\\[1mm]
1. $  {\cal L}_K  $ is the Hamiltonian operator 
\begin{equation}\label{cL000K}
{\cal L}_K  := \omega\cdot \pa_\vphi  - d X_K (\wtilde \xi )_{| \wtilde \xi = \wtilde \xi_\delta (\vphi, \theta)}
 =  
 \omega\cdot \pa_\vphi  - \pa_\theta \pa_{\wtilde \xi}\nabla_{\wtilde \xi} K
( \wtilde \xi ))_{| \wtilde \xi = \wtilde \xi_\delta (\vphi, \theta)} 
\end{equation}
where   $ K $  is the  Hamiltonian  in \eqref{defKPhi} and 
\be\label{approqp}
\wtilde \xi_\delta (\vphi) 
:= \e {\cal J}_0 \tc_2  + 
 \widetilde {\beta}_{2,\d}   (\varphi)\ts_2 + \e \Xi  (  i_\delta (\vphi))
\ee
where $\Xi  $ is defined in \eqref{aacoordinates}, 
$i_\delta (\vphi) $ 
 is given in Lemma \ref{toro isotropico modificato},  
and 
\begin{align}
& \widetilde {\beta}_{2,\d} (\varphi) 
:=-
(\omega\cdot\pa_\vphi)_{\rm ext}^{-1}
\big[ (\pa_{\wtilde \alpha_2}  \cK) 
(\e {\cal J}_0, \e \Xi  (  i_\delta(\varphi) ))
- {\mu}_\delta   \big] \label{defbeta2delta} \, , \\
&   {\mu}_\delta 
:= 
\big\langle (\pa_{\wtilde \alpha_2}  \cK) (\e {\cal J}_0, \e \Xi  (  i_\delta(\varphi))) 
\big\rangle \, , \label{defmu}
\end{align} 
and $ \cK$ is the Hamiltonian in \eqref{tildeKeq}, 
the operator $ (\omega\cdot\pa_\vphi)_{\rm ext}^{-1} $ is defined in 
\eqref{paext0}.

Note that by \eqref{tildeKeq} the value of 
$ d X_K (\wtilde \xi_\delta )$ is actually independent of $ \wtilde {\beta}_{2,\d} (\vphi) $.
\\[1mm]
2. 
 $ {\cal R}  $ is a  self-adjoint operator with  the ``finite rank" form 
	\begin{equation}\label{finite_rank_R}
	{\cal R}(\vphi)[h] = {\mathop \sum}_{j=1}^{|\ST|} ( h,g_j)_{L^2} \chi_j \,, \quad 
	\forall\, h\in\acca_{\ST,2}^\bot \,,
	\end{equation}
	for functions $g_j,\chi_j \in \acca_{\ST,2}^\bot $ 
	which satisfy, for some $\sigma:= \sigma(\tau,\ST, k_0) > 0 $,  for all $ j = 1, \ldots, {|\ST|} $, for all $s\geq s_0$, 
\begin{equation}
\label{gjchij_est0}
	\begin{aligned}
\| g_j \|_s^{k_0,\upsilon} + \| \chi_j \|_s^{k_0,\upsilon} & \lesssim_s 1 + 
\| \fracchi_\delta \|_{s+\sigma}^{k_0,\upsilon} \,, \\
	\| d_i g_j [\widehat \imath] \|_s + \| d_i \chi_j [\widehat \imath ] \|_s & \lesssim_s 
	\| \widehat \imath \|_{s+\sigma} + \| \widehat \imath \|_{s_0+\sigma} 
	\| \fracchi_\delta \|_{s+\sigma} \,.
\end{aligned}
	\end{equation}
Moreover the operator $ {\cal L}_\omega $ is reversible.
\end{lemma}

\begin{pf}
By  \eqref{KHG} we have 
$\tK_{02} (\vphi) = \partial_w \nabla_w\tK_{\tg}(\vphi,0,0)  $
with $\tK_{\tg}$ defined in  \eqref{new-Hamilt-K},
and by \eqref{H alpha} 
\begin{align}\label{tK02}
\partial_w \nabla_w\tK_{\tg}(\vphi,0,0)  &=\partial_w \nabla_w (\sK_{\tg}\circ G_\delta) (\vphi,0,0)
={\bf \Omega}(\gamma)_{|{\frak H}_{\ST,2}^\bot}+\varepsilon \partial_w \nabla_w (\sP\circ G_\delta) (\vphi,0,0)
\end{align}
where ${\bf \Omega}(\gamma)$ is defined in \eqref{defQH}. From \eqref{trasformazione modificata simplettica} we write
$ (\sP\circ G_\delta) (\phi,y,w)
= $ $ \sP\big( \vartheta_0(\phi),
 I_\delta (\phi) +  L_1(\phi) y +   L_2(\phi)w,   z_0(\phi) + w\big) $
with
$ L_1(\phi):=[\pa_\phi \vartheta_0(\phi)]^{-\top} $ 
and 
$ L_2(\phi):= - \big[ (\pa_\vartheta \wtilde{z}_0) (\vartheta_0(\phi)) \big]^\top \partial_\theta^{-1} $. 
Hence by the chain rule, 
$ \nabla_w(\sP\circ G_\delta) (\phi,y,w)
= $ $ L_2(\phi)^\top(\partial_I  \sP)(G_\delta (\phi,y,w))+(\nabla_z  \sP)(G_\delta (\phi,y,w)) $
and, differentiating  with respect to $w$ at the point  $(\vphi,0,0)$, we obtain
\begin{align}\label{pa w nabla p}
&\partial_w\nabla_w(\sP\circ G_\delta) (\vphi,0,0)
= \partial_z\nabla_z  \sP(i_\delta (\vphi)) - {\cal R}(\varphi), 
\end{align}
with
$  - {\cal R}(\varphi):={\cal R}_1(\varphi)+{\cal R}_2(\varphi)+{\cal R}_3(\varphi) $ and
$$ 
{\cal R}_1(\varphi):=L_2(\vphi)^\top\partial_{II}  \sP(i_\delta (\vphi))L_2(\vphi) \, , \
 {\cal R}_2(\varphi):=L_2(\vphi)^\top\partial_z\partial_I  \sP(i_\delta (\vphi)) \, , \ 
 {\cal R}_3(\varphi):=\partial_I\nabla_z  \sP(i_\delta (\vphi))L_2(\vphi) \, . 
 $$ 
Writing  $L_2(\vphi):\acca_{\ST,2}^\bot\to \mathbb{R}^{|\ST|}$ as 
$ L_2(\vphi)[h] 
 = {\mathop \sum}_{j=1}^{|\ST|} ( h, \big[L_2(\vphi)\big]^\top \underline{e}_j )_{L^2} \underline{e}_j $,
for any $  h\in\acca_{\ST,2}^\bot $,	we deduce that  
  each of  ${\cal R}_i (\vphi) $, $ i = 1, 2, 3, $
has the finite rank form  \eqref{finite_rank_R}, thus   $ {\cal R} (\vphi) $. For example  
$
{\cal R}_1(\varphi) = {\mathop \sum}_{j=1}^{|\ST|} ( h,\big[L_2(\vphi)\big]^\top \underline{e}_j )_{L^2} A_1(\vphi)\underline{e}_j  $ with
$ A_1(\vphi) :=
	\big[L_2(\vphi)\big]^\top\partial_{II}  \sP(i_\delta (\vphi)) $. 
The estimates \eqref{gjchij_est0} then  follow from \eqref{ansatz 0}, \eqref{2015-2},  \eqref{stima y - y delta}, \eqref{stima toro modificato}, \eqref{DG delta}  and \eqref{DG2 delta}.		
By  \eqref{tK02} and \eqref{pa w nabla p}, we obtain
\begin{align}
\tK_{02} (\vphi)  
& = \partial_w \nabla_w\tK_{\tg}(\vphi,0,0) 
 ={\bf \Omega}(\gamma)_{|{\frak H}_{\ST,2}^\bot}+\varepsilon \partial_z\nabla_z  \sP(i_\delta (\vphi)) - \varepsilon {\cal R}(\varphi)
\notag \\ 
&\stackrel{\eqref{aacoordinates},\eqref{cNP2}}  = 
 \Pi_{\ST,2}^\bot
\Big( \partial_{\breve u}\nabla_{\breve u} 
  H_L(\breve u) +\varepsilon \partial_{\breve u}\nabla_{\breve u}  
{\cal P}( {\cal J}_0, \breve u)_{|\breve u = \Xi  (  i_\delta (\vphi))} \Big)_{|{\frak H}_{\ST,2}^\bot}
- \varepsilon {\cal R}(\varphi) \label{passpr}
\end{align}
Thus, by \eqref{passpr}, \eqref{aacoordinates}, \eqref{cNP2}, \eqref{formaHep}, \eqref{HS:mu-breve}, \eqref{tildeKeq}, 
and \eqref{tildeKeq},  
($ K (\wtilde \xi_\delta )$ is  independent of $ \wtilde {\beta}_2 $) 
we get 
\be\label{esprK02} 
\tK_{02} (\vphi)   = 
\Pi_{\ST,2}^\bot \big( \partial_{\wtilde \xi}\nabla_{\wtilde \xi}  
K( \wtilde \xi )_{| \wtilde \xi = \e {\cal J}_0 \tc_2 + \wtilde {\beta}_{2,\d} (\vphi) \ts_2 +
\e \Xi  (  i_\delta (\vphi))} \big)_{|{\frak H}_{\ST,2}^\bot} - \varepsilon {\cal R}(\varphi) 
\ee
By  \eqref{esprK02}, the operator in \eqref{Lomega def} has the form  \eqref{K 02}
with ${\cal L}_K$ defined in \eqref{cL000K}
\end{pf}

Notice  that the linear Hamiltonian  operator $ {\cal L}_K $ in \eqref{cL000K} 
is obtained linearizing  the Hamiltonian $ K $ in \eqref{defKPhi}
at the function $ \wtilde \xi_\delta (\vphi) $ in \eqref{approqp} and
that, by Lemma \ref{lem:REV}, the operator $ {\cal L}_K $ is reversible, according to 
Definition \ref{def:R-AR}. 

In the sequel we assume the following ansatz (satisfied by the approximate solutions obtained in the nonlinear Nash-Moser iteration): for some constant
 $ \perd_0 :=\perd_0(\tau,\ST)>0$, $\upsilon\in (0,1)$, 
\be\label{ansatz_I0_s0}
\| \fracchi_0 \|_{s_0+ \perd_0} ^{k_0,\upsilon} \, ,  \ 
\| \fracchi_\delta \|_{s_0+ \perd_0}^{k_0,\upsilon}  
\leq 1 \,.
\ee 
The constant $ \perd_0$ 
represents the loss of derivatives accumulated along the reduction procedure of the next sections. It is  independent of the Sobolev index $s$.

\begin{lemma}
The  function 
 $ \wtilde \xi_\delta (\vphi, \theta) $ in \eqref{approqp}, 
 resp. $\wtilde\beta_{2,\delta}(\varphi)$ in \eqref{defbeta2delta}, 
 is even in $  (\vphi, \theta) $, resp. odd  in $ \vphi $, 
and, for any  $ s \geq s_0 $, satisfy 
\be\label{esxidelta}
\| \wtilde \xi_\delta \|_{s}^{k_0, \upsilon}, \|\wtilde\beta_{2,\delta}\|_s^{k_0,\upsilon}  \lesssim_s 
 \e   (1 + \| \fracchi_0 \|_{s+\sigma}^{k_0,\upsilon} ) \, . 
\ee
The constant  $ \mu_\delta $ in \eqref{defmu}  satisfies 
\be\label{stima-defmudelta}
|\mu_\delta|^{k_0,\upsilon} \lesssim \varepsilon  \, . 
\ee
\end{lemma}

\begin{pf}
Since $ i_\delta(\vphi) = ( \vartheta_0(\vphi),I_\delta(\vphi),z_0(\vphi) )$
is reversible according to \eqref{parity solution0}-\eqref{parity solution}, the function 
 $ \wtilde \xi_\delta (\vphi, \theta) $ is even in $  (\vphi, \theta) $.
 By recalling \eqref{tildeKeq} and  \eqref{fre:uncha} we have that
$$
 (\pa_{\wtilde \alpha_2}  \cK) (\e {\cal J}_0, \e \Xi  (  i_\delta(\varphi))) =-\e \Omega_2 {\cal J}_0+(\pa_{\wtilde \alpha_2}  \cK_{\geq 3}) (\e {\cal J}_0, \e \Xi  (  i_\delta(\varphi)))  
$$
is even in $ \vphi $
and the function  $\wtilde\beta_{2,\delta}(\varphi)$ in \eqref{defbeta2delta} is odd in $ \vphi $.
The estimate \eqref{esxidelta} follows as in Lemma \ref{lemma quantitativo forma normale}.
Finally the constant $ \mu_\delta $ in \eqref{defmu}  satisfies
the estimate \eqref{stima-defmudelta}, by 
\eqref{def:sp}, \eqref{2015-2}, \eqref{stime XP}.
\end{pf}

Next we define
\be\label{defxid}
\xi_\delta (\vphi ) := \xi_\delta (\vphi, \theta) :=\Phi^{-1} (\wtilde \xi_\delta(\vphi))  \, ,
\ee 
that, by  \eqref{est-phi-1} and \eqref{esxidelta}, satisfies 
for any $ s \geq s_0 $, 
\be\label{esxidelta2}
\|  \xi_\delta \|_{s}^{k_0, \upsilon} \lesssim_s 
 \e (1 + \| \fracchi_0 \|_{s+\sigma}^{k_0,\upsilon} ) \, . 
\ee
Moreover, since $  \bar t (\xi_\delta (\vphi))  = \wtilde \beta_{2,\d} (\vphi) $ by \eqref{defxid}, 
we deduce by \eqref{esxidelta}, for any $ s \geq s_0 $, 
\be\label{bartesxidelta}
| \overline t (\xi_\d) |^{k_0, \upsilon}_s  \lesssim_s 
\e  (1 + \| \fracchi_0 \|_{s+\sigma}^{k_0,\upsilon} ) \, . 
\ee
In order to estimate the variation of the eigenvalues with respect to the approximate invariant
torus, one needs also to estimate the variation with respect to the torus $i(\varphi)$ in another low norm $\|\;\|_{s_1}$
for all Sobolev indexes $s_1$ such that
\be\label{def:s1}
s_1+\sigma_0\leq s_0+ \perd_0\, , \quad \textnormal{for some}\quad \sigma_0:=\sigma_0(\tau,\ST) \, . 
\ee
Hence, by \eqref{ansatz_I0_s0} we have
$ \| \fracchi_0 \|_{s_1+\sigma_0} ^{k_0,\upsilon} \, ,  \ 
\| \fracchi_\delta \|_{s_1+\sigma_0}^{k_0,\upsilon}  
\lesssim_{s_1} 1 $. 
For a quantity $g(i)$ (an operator, a map, a scalar function)
depending on the torus $i$, 
we denote the difference
$ \Delta_{12}g:=g(i_2)-g(i_1) $. 
Using \eqref{derivata i delta} we get 
$ \|\Delta_{12}  \wtilde \xi_\delta\|_{s_1}\lesssim \varepsilon \|i_1-i_2 \|_{s_1+\sigma} $. 
In the sequel we shall not insist much on this standard point. 

\subsection{Linearized operator after the 
symplectic reduction}
\label{Sec:KAMnew}

The goal of this section is to  obtain the expression  \eqref{bellaconiug} 
of 
$ \Pi_{\ST,2}^\bot ({\cal L}_K)_{|  {\frak H}_{\ST,2}^\bot} $ where  
${\cal L}_K $ is the operator defined  in \eqref{cL000K}. 
We need the following  lemma.
\begin{lemma}\label{phi-1}
We have
\begin{align}\label{flowpi}
\Pi_{\ST,2}^\bot\, d \Phi (\xi_\delta) &= \Pi_{\ST,2}^\bot\, \Phi^{\bar t (\xi_\delta)}_{{\cal J}_2}  +   {\mathscr R}_1\, , 
 \\
\big[ d \Phi (\xi_\delta) \big]^{-1}  \Pi_{\ST,2}^\bot 
&=   \Phi^{- \bar t (\xi_\delta)}_{{\cal J}_2} \Pi_{\ST,2}^\bot+   {\mathscr R}_2 \, ,  \label{flowpiinv} \\
L (\xi_\delta)^*  
& = \Pi_{\ST,2}^\bot (\Phi^{- \bar t (\xi_\delta)}_{{\cal J}_2})^*  +  {\mathscr R}_3
\, , \qquad 
L (\xi_\delta) := \big[ d \Phi (\xi_\delta)\big]^{-1}  \Pi_{\ST,2}^\bot \, , 
\end{align}
 where the operators 
 $ {\mathscr R}_1$, $ {\mathscr R}_2 $, $ {\mathscr R}_3 $ have the  finite rank form
 $ {\mathscr R} := (g , \cdot )  \chi    $  for functions $g,\chi \in \acca_{\ST,2}^\bot $ which satisfy, for some $\sigma:= \sigma(\tau,\ST, k_0) > 0 $ for all $s\geq s_0$, 
\begin{equation}
\label{gjchij_est}
	\begin{aligned}
\max \{ \| g \|_s^{k_0,\upsilon}, \| \chi \|_s^{k_0,\upsilon}  \} & \lesssim_s 
\e (1 + 
 \| \fracchi_0 \|_{s+\sigma}^{k_0,\upsilon}) \,, \\
\max \{ \| d_i g [\widehat \imath] \|_s, \| d_i \chi [\widehat \imath ] \|_s \} 
& \lesssim_s 
\e 
\big( \| \widehat \imath \|_{s+\sigma} +  \| \widehat \imath \|_{s_0+\sigma} 
\| \fracchi_0 \|_{s+\sigma} \big) \,.
\end{aligned}
\ee
\end{lemma}

\begin{pf}
Follows by Lemmata \ref{lem:dphi}-\ref{lem:adj} and   
\eqref{gradbart}, 
\eqref{flowaff}, \eqref{pr-comp1} and Section \ref{sec:regularity}.  
\end{pf} 

\noindent
In the sequel $ \sigma := \sigma(\tau,\ST, k_0) > 0 $  denotes a constant, possibly larger from lemma to lemma. 

\begin{lemma}\label{lem:AI} 
The operator   $\Pi_{\ST,2}^\bot ({\cal L}_K)_{|{\frak H}_{\ST,2}^\bot}$, 
where $ {\cal L}_K $ is defined  in \eqref{cL000K},  is, for 
any $\omega\in \mathtt {DC} (\upsilon, \tau)$, 
\be
\Pi_{\ST,2}^\bot ({\cal L}_K)_{| {\frak H}_{\ST,2}^\bot} 
 =
 \Pi_{\ST,2}^\bot \Big( 
  \Phi^{  \overline t (\xi_\d) }_{{\cal J}_2} \, {\cal L}_{H_\Omega}\, \Phi^{  - \overline t (\xi_\d) }_{{\cal J}_2} 
   - \breve{\mu}_\e \pa_\theta \circ g_\gamma (\theta) \Big)_{| {\frak H}_{\ST,2}^\bot} + {\cal R}_{Z}+  {\mathscr R}  \label{bellaconiug}
\ee
where each term in the right hand side of \eqref{bellaconiug} 
is defined for all $ (\omega,\gamma) \in 
\R^{|\ST|} \times [\g_1,\g_2] $ and
\\[1mm]
$1$. ${\cal L}_{H_\Omega}$ is the quasi-periodic Hamiltonian reversible operator 
  \be\label{defLH}
{\cal L}_{H_\Omega} := \om \cdot \pa_\vphi - d X_{H_\Omega} (\xi_\delta (\vphi)) 
\ee
with Hamiltonian
\be\label{defHOmega}
 H_\Omega:= -\frac12E+\frac{\Omega}{2}J \, , \quad  \Omega=\Omega_\gamma - \frac{2\sqrt{2}\gamma\mu_\delta}{\sqrt{\pi}(\g-\g^{-1})} \, , 
\ee 
where the constant $ \mu_\delta $ is defined in \eqref{defmu}.
\\[1mm]
2. The function $ g_\gamma (\theta) $ is defined in \eqref{fg0}  and 
the function $ \breve{\mu}_\e (\vphi ) $ satisfies, 
for some $ \sigma > 0 $,  
\be\label{stima-defmuep}
|\breve{\mu}_\e|^{k_0,\upsilon}_{s} \lesssim \e (1 + \| \fracchi_0 \|_{s+\sigma}^{k_0,\upsilon} ) \, , \quad \forall  s \geq s_0 \,   ;  
\ee
3.  the operator $ {\cal R}_{Z} $ acts in $ {\frak H}_{\ST,2}^\bot $ and 
satisfies, for some $ \sigma > 0 $, for any $ s \geq s_0 $,  the estimate
\be\label{estimate:RZ}
\| {\cal R}_{Z} h \|_s^{k_0,\upsilon} \lesssim_s  
\big( \| {Z}  \|_{s+\sigma}^{k_0,\upsilon} 
+ 
\| {Z}  \|_{s_0+\sigma}^{k_0,\upsilon}
 \| \fracchi_0 \|_{s+\sigma}^{k_0,\upsilon}\big)\|  h \|_{s_0+\sigma}^{k_0,\upsilon}  
+ \| {Z}  \|_{s_0+\sigma}^{k_0,\upsilon}\|  h \|_{s+\sigma}^{k_0,\upsilon}\, ; 
\ee
4. the operator  $ {\mathscr R} $ acts in $ {\frak H}_{\ST,2}^\bot $ and 
has the finite rank form 
\eqref{finite_rank_R} 
 	for functions $g_j,\chi_j \in \acca_{\ST,2}^\bot $ 
	which satisfy, for some $\sigma > 0 $,  for all $ j = 1, \ldots, {|\ST|} $, for all $s\geq s_0$,
\begin{equation}
\label{gjchij_est1}
	\begin{aligned}
\max\{ \| g_j \|_s^{k_0,\upsilon}, \| \chi_j \|_s^{k_0,\upsilon} \} & \lesssim_s 
\e  (1 + 
 \| \fracchi_0 \|_{s+\sigma}^{k_0,\upsilon}) \,, \\
	\max\{ \| d_i g_j [\widehat \imath] \|_s,  \| d_i \chi_j [\widehat \imath ] \|_s \} & \lesssim_s 
	\e  \big( \| \widehat \imath \|_{s+\sigma} +  \| \widehat \imath \|_{s_0+\sigma} 
	\| \fracchi_0 \|_{s+\sigma} \big) \,.
\end{aligned}
\ee
 \end{lemma}

\begin{pf}
 For any $  \wtilde \xi (\vphi) := \wtilde \xi (\vphi, \theta) $, 
 one has 
  \begin{align}\label{relKH3}
 \omega \cdot \partial_\vphi \wtilde \xi  -  X_{K -  \mu_\delta \wtilde \alpha_2} 
(\wtilde \xi ) 
  & = (d \Phi) ( \Phi^{-1} (\wtilde \xi))\Big(\omega \cdot \partial_\vphi  
 \big( \Phi^{-1} (\wtilde \xi )\big)- X_{H  - \mu_\delta {\cal J}}  ( \Phi^{-1}(\wtilde \xi) ) \Big)\, . 
  \end{align}
 Indeed, according to \eqref{phisymdiffeo} and  \eqref{defKPhi} the Hamiltonian vector field $ X_{K- \mu_\delta \wtilde \a_2}   $ is the push-forward of $X_{H  -  \mu_\delta {\cal J}} $ under the map $ \wtilde \xi =  \Phi (\xi)$, namely
 \be\label{prpf1}
 X_{K-\mu_\delta \wtilde \a_2} (\wtilde \xi )  = (d \Phi) ( \Phi^{-1} (\wtilde \xi)) 
 X_{H  -  \mu_\delta {\cal J}}  ( \Phi^{-1}(\wtilde \xi) ) \, , 
 \ee
 and, differentiating the identity $ \wtilde \xi (\vphi) =  \Phi (\Phi^{-1} (\wtilde \xi (\vphi)))  $,
 we get   
 \begin{align}
  \om \cdot \pa_\vphi \wtilde \xi (\vphi) 
  = d_\vphi  \wtilde \xi (\vphi) [\omega] 
  & =
(d \Phi) (\Phi^{-1} (\wtilde \xi (\vphi)))
\, d_\vphi \big( \Phi^{-1} (\wtilde \xi (\vphi) \big) [\omega] \notag \\
& =
(d \Phi) (\Phi^{-1} (\wtilde \xi (\vphi))) \, \om \cdot \pa_\vphi  \big( \Phi^{-1} (\wtilde \xi (\vphi))\big) \, . \label{prpf2}
\end{align}
The identities \eqref{prpf1}-\eqref{prpf2}  imply  \eqref{relKH3}. 
Next, differentiating  \eqref{relKH3} with respect to $\wtilde \xi$ at  $ \wtilde \xi_\delta (\vphi) $,
and since 
 \be \label{dphi-1-1}
 d \Phi^{-1} ( \wtilde \xi_\delta ) = \big[d \Phi (  \xi_\delta  )\big]^{-1} \quad {\rm with}\quad  \xi_\delta =  \Phi^{-1} ( \wtilde \xi_\d)\, , 
 \ee 
  we get
  \begin{align}
   \omega \cdot \partial_\vphi &-d X_{K -  \mu_\delta \wtilde \alpha_2}( \wtilde \xi_\delta) 
  =
d \Phi   ( \xi_\delta) \big(\omega \cdot \partial_\vphi  - d X_{H  -  \mu_\delta {\cal J}} (   \xi_\delta ) \big) \big[d \Phi   ( \xi_\delta) \big]^{-1}  \notag \\
&\quad  +
 d^2 \Phi ( \xi_\delta ) 
\Big[  
 \omega \cdot \partial_\vphi  ( \Phi^{-1} (\wtilde \xi_\delta))-X_{H  -  \mu_\delta {\cal J}} (  \Phi^{-1} (\wtilde \xi_\delta)  ),  d \Phi^{-1} ( \wtilde \xi_\delta )  [ \,]  \Big] \notag \\
 &  \stackrel{\eqref{defLH}, \eqref{defHOmega}, \eqref{relKH3},\eqref{dphi-1-1}} =
d \Phi   ( \xi_\delta) \, {\cal L}_{H_\Omega}  \, 
\big[d \Phi   ( \xi_\delta) \big]^{-1}  \notag  \\
& \qquad \qquad  +  d^2 \Phi ( \xi_\delta ) 
\Big[  
d \Phi^{-1} ( \wtilde \xi_\delta )
 \big( \omega \cdot \partial_\vphi \wtilde \xi_\delta -  X_{K - \mu_\delta \wtilde \alpha_2} (\wtilde \xi_\delta) \big),  d \Phi^{-1} ( \wtilde \xi_\delta )   [ \,]  \Big]  \, .  \label{terzoquasi}
\end{align}
The next goal is to prove \eqref{dxi-omg}.
In view of  \eqref{tildeKeq}, 
the vector field $ X_{K- \mu_\delta \wtilde \a_2}   $ is 
\be\label{primoPF}
X_{K- \mu_\delta \wtilde \a_2} (\wtilde \xi )  =
X_{K} (\wtilde \xi ) + \mu_\delta \ts_2  =\mu_\delta \ts_2  - \pa_{\wtilde \a_2} {\cal K} ( \wtilde \a_2,  \wtilde u ) \ts_2 + 
\pa_\theta \nabla_{\wtilde u} {\cal K} ( \wtilde \a_2,  \wtilde u  \, ) \, .
\ee
Differentiating  \eqref{approqp} and subtracting \eqref{primoPF} we get 
  \begin{align*}
\omega \cdot \partial_\vphi \wtilde \xi_\delta -  X_{K -  \mu_\delta \wtilde \alpha_2} (\wtilde \xi_\delta) 
 &=\big(\omega \cdot \partial_\vphi  \widetilde \beta_{2,\delta}+(\pa_{\wtilde \alpha_2}  \cK) 
(\e {\cal J}_0, \e \Xi  (  i_\delta(\varphi) )) -\mu_\delta \big)\ts_2 \\
& \quad  +\omega \cdot \partial_\vphi  \big(\e \Xi  (  i_\delta (\vphi))\big) - \pa_\theta (\nabla_{\wtilde u} \cK) \big(\varepsilon {\cal J}_0, \wtilde u \big)_{| \wtilde u = \e \Xi  (  i_\delta (\vphi))}\,  .
\end{align*}  
Then, in view of \eqref{defbeta2delta} and \eqref{paext0}, for all $\omega\in \mathtt {DC} (\upsilon, \tau)$  one has
  \begin{align}
\omega \cdot \partial_\vphi \wtilde \xi_\delta -  X_{K -  \mu_\delta \wtilde \alpha_2} (\wtilde \xi_\delta) 
&
 =
\e d \big( \Xi  (  i_\delta (\vphi)) \big) [\omega] - \pa_\theta (\nabla_{\wtilde u} \cK) \big(\varepsilon {\cal J}_0, \wtilde u \big)_{| \wtilde u = \e \Xi  (  i_\delta (\vphi))}\nonumber
\\
& \stackrel{\eqref{HSKt0res}} = \e 
d \Xi  (  i_\delta (\vphi)) \, \omega \cdot  \pa_\vphi i_\delta (\vphi)   -\e \pa_\theta (\nabla_{\breve u}{\breve \cK}) ( {\cal J}_0, \breve u )_{| \breve u =  \Xi  (  i_\delta (\vphi))}\, . 
\label{quasgo}
\end{align}  
According to   \eqref{Hepsilon}, the Hamiltonian 
vector field $ \pa_\theta (\nabla_{\breve u}{\breve \cK}) ( {\cal J}_0, \breve u ) $ is the push-forward of $ X_{\sK_\e}$ under the map $\breve u = \Xi (\vartheta, I, z) $ 
in \eqref{aacoordinates}, namely 
\be\label{secopuf}
 \pa_\theta (\nabla_{\breve u}{\breve \cK}) ( {\cal J}_0, \breve u )_{| \breve u =  \Xi  (  i_\delta (\vphi))}  =
d \Xi (i_\delta (\vphi))  X_{\sK_\e} (i_\delta (\vphi))\, .
\ee
Combining \eqref{quasgo} and \eqref{secopuf} we get
\begin{align}\label{ef3}
\omega \cdot \partial_\vphi \wtilde \xi_\delta -  X_{K -  \mu_\delta \wtilde \alpha_2} (\wtilde \xi_\delta)  = \e d\Xi  (  i_\delta (\vphi))\big(\omega \cdot \partial_\vphi   i_\delta (\vphi) -X_{\sK_\e} ( i_\delta (\vphi))\big)\,. 
\end{align}
By  \eqref{cNP}-\eqref{cNP2} and \eqref{H alpha} we have 
$ \sK_\e(\vartheta,I, z) = \sK_{\tg_0} (\vartheta,I, z)+(\tg_0-\vec{\omega}(\g))\cdot I $, 
so that, by \eqref{XHemu},  the corresponding vector field is 
\be\label{pvfsvf}
X_{\sK_\e} = X_{\sK_{\tg_0}} - \big(
  \tg_0- \vec{\omega}(\g),
   0,
   0
\big)^\top \, , 
\ee
Plugging \eqref{pvfsvf} into \eqref{ef3}  we finally find
\begin{align}
\omega \cdot \partial_\vphi \wtilde \xi_\delta -  X_{K -  \mu_\delta \wtilde \alpha_2} (\wtilde \xi_\delta) & =  \e d\Xi  (  i_\delta (\vphi))\big(\omega \cdot \partial_\vphi   i_\delta (\vphi) -X_{\sK_{\tg_0}}( i_\delta (\vphi))\big) \nonumber
 + \e d\Xi  (  i_\delta (\vphi))  \big(
  \tg_0- \vec{\omega}(\g),
   0,
   0
\big)^\top \notag \\
& 
\stackrel{\eqref{operatorF}, \eqref{aacoordinates}} = \e d\Xi  (  i_\delta (\vphi)) {\cal F}(i_\delta, \tg_0) 
 + \e \partial_\vartheta {\mathtt v}^\intercal (\vartheta_0,I_\delta) ( \tg_0- \vec{\omega}(\g))   \, , \label{dxi-omg}
\end{align}
and inserting \eqref{dxi-omg} into  \eqref{terzoquasi} we deduce
 \begin{align}
   \omega \cdot \partial_\vphi -d X_{K -  \mu_\delta \wtilde \alpha_2}( \wtilde \xi_\delta) 
&  =
d \Phi   ( \xi_\delta) \, {\cal L}_{H_\Omega} \,  \big[d \Phi   ( \xi_\delta) \big]^{-1} \label{primacong1} 
\\
&\quad  +\e  d^2 \Phi ( \xi_\delta ) 
\Big[  
d \Phi^{-1} ( \wtilde \xi_\delta )  d\Xi  (  i_\delta (\vphi)) {\cal F}(i_\delta, \tg_0),  d \Phi^{-1} ( \wtilde \xi_\delta )   [ \,]  \Big]    \notag  
 \\
&\quad   +\e d^2 \Phi ( \xi_\delta ) \Big[ d \Phi^{-1} ( \wtilde \xi_\delta ) 
  \partial_\vartheta {\mathtt v}^\intercal (\vartheta_0,I_\delta)( \tg_0- \vec{\omega}(\g)) ,  d \Phi^{-1} ( \wtilde \xi_\delta )   [ \,]  \Big]  \, .  \notag
\end{align}
On the other hand, by \eqref{primoPF}, 
the Hamiltonian vector field $ X_{ \wtilde \alpha_2 } $ has a component only  on the mode $ \ts_2  $. Hence, applying the projector $ \Pi_{\ST,2}^\bot $ to \eqref{primacong1}  and using  \eqref{flowpi}-\eqref{flowpiinv}  we obtain
\be\label{quasifia}
\begin{aligned}
\Pi_{\ST,2}^\bot \big( \omega\cdot \pa_\vphi  - d X_K ( \wtilde \xi_\delta   ) \big)_{| {\frak H}_{\ST,2}^\bot} 
 & =
 \Pi_{\ST,2}^\bot \Big( 
  \Phi^{  \overline t (\xi_\d) }_{{\cal J}_2} \, {\cal L}_{H_\Omega}\, \Phi^{  - \overline t (\xi_\d) }_{{\cal J}_2} \Big)_{| {\frak H}_{\ST,2}^\bot} + {\cal R}_{Z} +   {\cal E}+   \widetilde{\cal E} 
\end{aligned}
\ee
with 
\begin{align}
{\cal R}_{Z} &:= \e \, \Pi_{\ST,2}^\bot \, d^2 \Phi (\xi_\delta) 
\Big[  d \Phi^{-1} ( \wtilde \xi_\delta ) 
d\Xi  (  i_\delta (\vphi)) {\cal F}(i_\delta, \tg_0)  ,  d \Phi^{-1} ( \wtilde \xi_\delta ) _{| {\frak H}_{\ST,2}^\bot}   [ \,]  \Big], \label{RZ}
\\
 {\cal E} &:= \e\,  
 \Pi_{\ST,2}^\bot \, d^2 \Phi (\xi_\delta) 
\Big[  d \Phi^{-1} ( \wtilde \xi_\delta ) 
\partial_\vartheta {\mathtt v}^\intercal (\vartheta_0,I_\delta)( \tg_0- \vec{\omega}(\g)), d \Phi^{-1} ( \wtilde \xi_\delta ) _{| {\frak H}_{\ST,2}^\bot}   [ \,]  \Big]\, ,  \label{RE}
\\  
\widetilde{\cal E} &:= \Pi_{\ST,2}^\bot   \Phi^{  \overline t (\xi_\d) }_{{\cal J}_2}  {\cal L}_{H_\Omega} {\mathscr R} _2 + {\mathscr R} _1  {\cal L}_{H_\Omega} \big[d \Phi   ( \xi_\delta) \big]^{-1}_{| {\frak H}_{\ST,2}^\bot} \, .\label{REE}
 \end{align}
By \eqref{RZ},   Lemma \ref{lem:dphi-phi-1}, \eqref{ansatz 0}, \eqref{est:v},   \eqref{esxidelta2},   \eqref{stima toro modificato} we obtain the estimate \eqref{estimate:RZ}.
Now we shall expand the  term ${\cal E}$ in \eqref{RE}. 
Then, by setting  
\be\label{deff}
f:=d \Phi^{-1} ( \wtilde \xi_\delta ) \partial_\vartheta {\mathtt v}^\intercal (\vartheta_0,I_\delta)( \tg_0- \vec{\omega}(\g))\, ,
\ee
and using the identity \eqref{Hess-Phi} we deduce that the operator in \eqref{RE} writes
\begin{align} 
&{\cal E}=\e\,  d \bar t (\xi_\delta) [f] \, \Pi_{\ST,2}^\bot X_{{\cal J}_2}\big(\Phi^{  \overline t (\xi_\delta) }_{{\cal J}_2} 
d \Phi^{-1} ( \wtilde \xi_\delta )_{| {\frak H}_{\ST,2}^\bot}  [\,] \big)+{\cal E}_1+{\cal E}_2+{\cal E}_3 \qquad\textnormal{with} \label{svicalE0} \\
&{\cal E}_1:=\e\,\Pi_{\ST,2}^\bot X_{{\cal J}_2} (\Phi_{\cal J}^{\overline t (\xi_\delta)} (\xi_\delta))  d^2 {\bar t} (\xi_\delta)\big[  f, d \Phi^{-1} ( \wtilde \xi_\delta )_{| {\frak H}_{\ST,2}^\bot} [\, ]  \big], \notag \\ 
&{\cal E}_2:=\e\, d \bar t (\xi_\delta)[f]\,
\Pi_{\ST,2}^\bot X_{{\cal J}_2}\big(X_{{\cal J}}  (\Phi_{\cal J}^{\overline t (\xi_\delta)} (\xi_\delta)) \big) \,  d \bar t (\xi_\delta) d \Phi^{-1} ( \wtilde \xi_\delta )_{| {\frak H}_{\ST,2}^\bot} \, ,
\notag \\ 
&{\cal E}_3 := \e\,
\Pi_{\ST,2}^\bot X_{{\cal J}_2}\big(\Phi^{  \overline t (\xi_\delta) }_{{\cal J}_2} [f]  \big)d \bar t (\xi_\delta) d \Phi^{-1} ( \wtilde \xi_\delta ) _{| {\frak H}_{\ST,2}^\bot}  \notag
 \, .
\end{align}
 Using 
 \eqref{XJcal}, \eqref{dphi-1-1} and Lemma \ref{phi-1}, we have   
\begin{align*}
\Pi_{\ST,2}^\bot  X_{{\cal J}_2}\big(\Phi^{  \overline t (\xi_\d) }_{{\cal J}_2} d \Phi^{-1} ( \wtilde \xi_\delta )_{| {\frak H}_{\ST,2}^\bot}   [\, ] \big)
&  = 2 \aleph\,  \Pi_{\ST,2}^\bot  \big( \pa_\teta \circ g_\gamma (\theta)\big)  +  \Pi_{\ST,2}^\bot  X_{{\cal J}_2}\big(\Phi^{  \overline t (\xi_\d) }_{{\cal J}_2} {\mathscr R}_2 [\, ] \big) 
\end{align*}
and thus the first operator in \eqref{svicalE0} is 
\be\label{svicalE}
{\cal E}=\breve{\mu}_\e \, \Pi_{\ST,2}^\bot  \big( \pa_\teta \circ g_\gamma (\theta)\big)+{\cal E}_1+{\cal E}_2+{\cal E}_3+{\cal E}_4  
\ee
with
\be\label{defmuep}
\breve{\mu}_\e (\vphi ) :=  2  \aleph\,  \e\, 	\, d \bar t (\xi_\delta) [f]  \, , \quad
{\cal E}_4:= \e\,  d \bar t (\xi_\delta) [f] \, 
\Pi_{\ST,2}^\bot  X_{{\cal J}_2}\big(\Phi^{  \overline t (\xi_\d) }_{{\cal J}_2} {\mathscr R}_2 [\, ] \big) \, . 
\ee
The function $ \breve{\mu}_\e $ 
is defined for all 
$ (\omega,\gamma) \in \R^{|\ST|} \times [\g_1,\g_2] $
and satisfies \eqref{stima-defmuep}.

Inserting   \eqref{svicalE} into   \eqref{quasifia} 
we deduce \eqref{bellaconiug} with 
\be\label{defR1ap}
 {\mathscr R}  := {\cal E}_1 + {\cal E}_2 + {\cal E}_3+ {\cal E}_4
+ \widetilde{\cal E}  . 
\ee
Note that  each $ \mathcal{E}_1 $, $ \mathcal{E}_2   $, $ \mathcal{E}_3   $
has the finite rank form  \eqref{finite_rank_R}, since 
\begin{align*} 
{\cal E}_1&=(  \cdot , g_1)_{L^2} \chi_1,\quad  g_1 := L^*  d\nabla \bar t (\xi_\delta)[f],\quad\chi_1 :=\e\,\Pi_{\ST,2}^\bot X_{{\cal J}_2} (\Phi_{\cal J}^{\overline t (\xi_\delta)} (\xi_\delta))  \, , \notag \\
{\cal E}_2&=( \cdot ,  g_2 )_{L^2} \chi_2,\quad\quad\quad g_2 := L^*  \nabla \bar t (\xi_\delta),\quad  \chi_2 := \frac{\breve{\mu}_\e}{2  \aleph}  \, \Pi_{\ST,2}^\bot X_{{\cal J}_2}\big(X_{{\cal J}}  (\Phi_{\cal J}^{\overline t (\xi_\delta)} (\xi_\delta)) \big)  \, , \notag \\
{\cal E}_3&= ( \cdot ,  g_3 )_{L^2} \chi_3,\quad\quad\quad g_3:= 
L^*  \nabla \bar t (\xi_\delta),\quad  \chi_3 :=\e\, \Pi_{\ST,2}^\bot X_{{\cal J}_2}\big(\Phi^{  \overline t (\xi_\delta) }_{{\cal J}_2} [f]\big)\, , 
\end{align*}
with  $ L =[d \Phi(\xi_\delta)]^{-1}_{| {\frak H}_{\ST,2}^\bot} $. Moreover the operators ${\cal E}_4 $ and $\widetilde{\cal E} $, defined  in \eqref{defmuep} and  \eqref{REE},
are finite rank as well, being composition of the finite rank operators 
$ {\mathscr R}_1$ and $ {\mathscr R} _2 $, and the finite rank operator  
$ {\mathscr R} $ in \eqref{defR1ap} satisfies \eqref{gjchij_est1}.
The operators $ {\cal R}_Z $ in \eqref{RZ} and  ${\cal R}$ in \eqref{defR1ap}
are defined  for all 
$ (\omega,\gamma) \in \R^{|\ST|} \times [\g_1,\g_2] $
and satisfy the estimates \eqref{estimate:RZ}, resp. \eqref{gjchij_est1}.
 \end{pf}

\subsection{Conjugation by the flow  $ \Phi^{  \overline t (\xi_\d) }_{{\cal J}_2} $
}\label{sec:LH}

 The 
 following proposition 
 shows that the linearized operator 
 $ \Pi_{\ST,2}^\bot ({\cal L}_K)_{| \acca_{\ST,2}^\bot} $ 
 is very similar to the restriction of 
 $ {\cal L}_{H_\Omega} $  defined in \eqref{defLH}, i.e. 
  before the symplectic rectification $ \Phi $. 
  We first recall  that the Hilbert transform $ {\cal H} $, 
acting on the $ 2 \pi $-periodic functions, is the Fourier multiplier 
\be\label{Hilbert-transf}
{\cal H} ( e^{\ii j \theta} ) := - \ii \, \sign(j) e^{\ii j \theta} \, , \    \forall j \neq 0 \, ,  \quad {\cal H} (1) := 0\,.
\ee
  
\begin{proposition}\label{lem:Lorto} 
{\bf (Structure of $\Pi_{\ST,2}^\bot ({\cal L}_K)_{| \acca_{\ST,2}^\bot}$)} 
The linear operator $\Pi_{\ST,2}^\bot ({\cal L}_K)_{| \acca_{\ST,2}^\bot}$ in \eqref{bellaconiug} has the form, for 
any $\omega\in \mathtt {DC} (\upsilon, \tau)$, 
\be\label{linKOK}
\Pi_{\ST,2}^\bot ({\cal L}_K)_{| \acca_{\ST,2}^\bot}
 =  \Pi_{\ST,2}^\bot ( 
\omega \cdot \partial_\vphi  -  
\pa_\theta \,  {\cal V}+ \pa_\theta W_0  +{\rm R}_\e
 )_{| \acca_{\ST,2}^\bot}
+{\cal R}_{Z}+ {\mathscr R} 
\ee
where each term in the right hand side of \eqref{linKOK} 
is defined for all $ (\omega,\gamma) \in 
\R^{|\ST|} \times [\g_1,\g_2] $ and
\\[1mm] 
1. ${\cal V} (\vphi, \theta) $ is a real even function of the form 
\be\label{new V}
{\cal V} (\vphi, \theta) = -  \Omega_\gamma + b (\vphi, \theta) 
\ee
where $ b (\vphi, \theta) $ satisfies, for some $ \sigma > 0 $, for any $ s \geq s_0 $, 
\begin{equation}\label{new V1} \| b  \|_s^{k_0,\upsilon} \lesssim_s 
\e (1 + \| \fracchi_0 \|_{s + \sigma}^{k_0, \upsilon}) \, ; 
\end{equation}
2.
$ \W_0 $ is the self-adjoint, reversibility preserving real operator computed in Lemma \ref{intop0}. The Hamiltonian operator $\pa_\theta \W_0 $ has the form 
\be\label{laH}
\pa_\theta \W_0 = \tfrac12 {\cal H} + {\cal Q}_\infty
\ee
where $ {\cal H} $ is the Hilbert transform defined in \eqref{Hilbert-transf}
and ${\cal Q}_{\infty} $ is the operator in $ {\rm OPS}^{-\infty}$ 
\be\label{defQinfty}
{\cal Q}_{\infty} (q) := 
 - \ii \sum_{j \in \Z \setminus \{0\}} \frac{j}{2|j|}   \ka_j q_{-j}  e^{ \ii j \theta}  \, , 
 \quad  q (\theta) = \sum_{j \in \Z \setminus \{0\}} q_j e^{\ii j \theta} \, , \quad
 \ka_j =   \Big( \frac{\gamma-1}{\gamma + 1} \Big)^{|j|}  \, ;
\ee   
3.   ${\rm R_\e}$ is  a reversibility preserving real operator in $  {\rm OPS}^{-\infty} $ and satisfies, for any $ m, \alpha \geq 0 $, $ s \geq s_0 $,
\be\label{stimaRep4}
| {\rm R}_\e  |_{-m, s, \a}^{k_0, \upsilon} \lesssim_{m, s, \alpha, k_0} \e   
 (1 + \| \fracchi_0 \|_{s + \sigma (m,\a)}^{k_0, \upsilon}) \, ;
\ee
4.  $ {\cal R}_{Z} $ satisfies the estimate \eqref{estimate:RZ};
\\[1mm] 
5. $ {\mathscr R} $ has the finite rank form 
\eqref{finite_rank_R} 
 	for functions $g_j,\chi_j \in \acca_{\ST,2}^\bot $ 
	 satisfying  \eqref{gjchij_est1}.
\end{proposition} 

The rest of this section is devoted to prove  Proposition \ref{lem:Lorto}.
In view of \eqref{bellaconiug}, in the next lemma
we  first  provide the structure of the
linear operator $ {\cal L}_{H_\Omega} $  defined by \eqref{defLH}.
\begin{lemma}\label{lem11.4}
{\bf (Structure of  $ {\cal L}_{H_\Omega} $)}
 The  Hamiltonian, reversible, real linear operator 
$ {\cal L}_{H_\Omega} $ in \eqref{defLH} 
 has the form 
\be\label{QPLIN}
{\cal L}_{H_\Omega} = 
\om \cdot \pa_\vphi - \pa_\theta\circ  V_{\mu_\delta}
 (\vphi, \theta) + \pa_\theta  \W_0 + \pa_\theta { R}_\e 
\ee
where 
\begin{enumerate}
\item 
 $ V_{\mu_\delta} (\vphi, \theta) $ is a  real even function, 
 $ V_{\mu_\delta} ( \vphi, \theta ) = V_{\mu_\delta} ( - \vphi, - \theta) $,   with the form  
 \be \label{defg0AA}
 V_{\mu_\delta} (\vphi, \theta)   =  - \Omega_\gamma  - \frac{2\sqrt{2}\gamma\mu_\delta}{\sqrt{\pi}(\g-\g^{-1})} g_\gamma (\theta)+  b_\delta (\vphi, \theta)
 \ee
where  the function 
$ g_\gamma (\theta) $  is introduced in \eqref{fg0}, the constant 
$ \mu_\delta $ in \eqref{defmu}, and the function $ b_\delta (\vphi, \theta) $ satisfies,
for some $ \s > 0 $, for any $ s \geq s_0 $, 
\be\label{bpicco}
\| b_\delta \|_{s}^{k_0, \upsilon} \lesssim_s 
\e( 1 + \| \fracchi_0 \|_{s+\s}^{k_0, \upsilon}) \, ; 
\ee
\item  The Hamiltonian operator $\pa_\theta \W_0 $ has the form   \eqref{laH}; 
\item \label{lite3}
$ { R}_\e $ 
is a self-adjoint reversibility preserving real operator in $ {\rm OPS}^{-\infty }$ satisfying  for all $  m,  \a \in \N_0 $,   for some constant $ \sigma (m,\a) > 0 $,  for any  $ s \geq s_0 $, 
an estimate as \eqref{stimaRep4}.
\end{enumerate}
\end{lemma}
\begin{pf}
In view of  Proposition \ref{lin-eq}, the   operator ${\cal L}_{H_\Omega} $ in \eqref{bellaconiug}  is 
\be\label{linop}
{\cal L}_{H_\Omega}     = \om \cdot \pa_\vphi     -   \pa_\theta \circ
V_{\mu_\delta} (\vphi, \theta)    +
\pa_\theta \circ \W ({\xi_\delta}(\vphi))  
\ee
where the operator $  \W({\xi})$ is defined in \eqref{int-op}  and
$V_{\mu_\delta} (\vphi, \theta) $ is the real function  
\be\label{Vdeltapt} 
V_{\mu_\delta} (\vphi, \theta) := \Omega g_\gamma (\theta) + v(\vphi, \theta) \, , 
\quad 
\Omega=\Omega_\gamma - \tfrac{2\sqrt{2}\gamma\mu_\delta}{\sqrt{\pi}(\g-\g^{-1})} \, ,  
\ee
with $ \Omega $  in \eqref{defLH},   and
$ v (\vphi, \theta)  := v(\xi_\delta(\vphi))(\theta )   $ 
is obtained evaluating $ v(\xi) (\theta) $ in \eqref{defV}  at $ \xi_\delta (\vphi, \theta)  $. 
By \eqref{Vdeltapt},
writing $v(\vphi, \theta)  = v_0( \theta) + b (\vphi, \theta)  $,  and 
using the identity \eqref{id-g-v0}, 
 we get  \eqref{defg0AA}. The function 
$ b (\vphi, \theta)  := v (\vphi, \theta ) - v_0 (\theta) $
satisfies the estimate \eqref{bpicco} by   \eqref{esxidelta2} and Lemma~\ref{Prop:csm}.
From Lemma \ref{intop0} we deduce the expansion \eqref{laH}.
Finally, by  \eqref{linop}, we get \eqref{QPLIN} with $ R_\e := W(\xi_\delta) - W_0 $,
which satisfies an estimate as \eqref{stimaRep4}  by Lemma \ref{Prop:csm} and 
\eqref{esxidelta2}. 
\end{pf}

Next, we study 
the conjugated operator $ \Phi^{\bar t (\xi_\delta)}_{{\cal J}_2} {\cal L}_{H_\Omega}  \Phi^{-\bar t (\xi_\delta)}_{{\cal J}_2} $  in \eqref{bellaconiug}. 
We  use the representation of the linear 
 symplectic flow $ \Phi^{t}_{{\cal J}_2} $ 
  provided by Lemma \ref{flowJ2old}. We anticipate the following lemma. 

\begin{lemma}\label{lem:diffest}
Let ${\cal B} $ be the composition operator in \eqref{phidin} induced by the 
${\cal C}^\infty $ diffeomorphism \eqref{invediff} of the torus $ \T $. 
Then 
\be\label{lemdifp}
{\cal B} \circ  T_{2\aleph \bar t(\xi_\delta)} \circ  {\cal B}^{-1} = 
 {\cal B}_{( \bar t (\xi_\delta))} 
\ee
where  ${\cal B}_{( \bar t (\xi_\delta))} $ is the composition operator 
\be \label{compdif}
{\cal B}_{(\bar t (\xi_\delta))} q (\theta):= q ( a_\e (\varphi,\theta))
\ee
induced by the $ \vphi$-dependent  family of diffeomorphisms of $ \T $ close to the identity, 
\begin{align}
a_\e (\varphi,\theta)  := 
\theta +  \alpha_\e (\varphi, \theta)  \, ,\quad  \alpha_\e (\varphi, \theta)
 &  := 
  \beta (\theta) +  2\aleph \bar t (\xi_\delta(\varphi) ) +
\breve \beta \big(\theta + \beta(\theta) +2\aleph \bar t (\xi_\delta (\varphi)) \big) \label{ath1} \\
&  = 2\aleph \bar t (\xi_\delta (\varphi) ) +
 \breve \beta \big(\theta + \beta(\theta) +2\aleph \bar t (\xi_\delta(\varphi)) \big) -
\breve \beta (\theta + \beta(\theta))   \notag \, .
\end{align}
The function $ \alpha_\e (\varphi, \theta) $ is odd$ (\vphi, \theta)  $
and, for some $ \s > 0 $, 
for any $ s \geq s_0 $, 
 \be\label{estim:a}
\|  \alpha_\e   \|_{s}^{k_0,\upsilon }\lesssim_s \, \e (1 + \| \fracchi_0 \|_{s+\sigma}^{k_0,\upsilon}) \, . 
 \ee
\end{lemma}

 \begin{pf}
By a direct computation we deduce 
\eqref{lemdifp}-\eqref{compdif} with $ \alpha_\e (\vphi, \theta) $ 
defined in \eqref{ath1}. 
Then note that 
 $ \beta (\theta) + \breve \beta (\theta + \beta (\theta)) = 0 $, 
 $  \forall \theta \in \R $. In addition, by \eqref{defxid}, the function 
 $ \bar t (\xi_\delta (\vphi)) = \wtilde \beta_{2,\d} (\vphi)$  is odd in $ \vphi $ 
as well as  $ \alpha_\e (\vphi, \theta) $. Finally
\eqref{estim:a} follows by \eqref{bartesxidelta}.
 \end{pf}

We now conjugate the transport operator  $ \omega \cdot \pa_\vphi - 
\pa_\theta \circ V_{\mu_\delta}(\vphi, \theta)  $.

\begin{lemma}\label{lemC5}
The conjugation of the Hamiltonian operator 
$  \omega \cdot \pa_\vphi -  
\pa_\theta \circ V_{\mu_\delta}(\vphi, \theta)  $ under the symplectic map 
$ \Phi^{\overline t (\xi_\d)}_{{\cal J}_2} $ is 
 \be\label{primaconc1}
\Phi^{\overline t (\xi_\d)}_{{\cal J}_2}  
\big(  \omega \cdot \pa_\vphi -  
\pa_\theta \circ V_{\mu_\delta}(\vphi, \theta) \big)  \Phi^{-\overline t (\xi_\d)}_{{\cal J}_2}  =
  \omega \cdot \pa_\vphi -  \pa_\theta \circ V_{(\overline t (\xi_\d))}  (\vphi, \theta)
\ee
where $V_{(\overline t (\xi_\d))} (\vphi, \theta)  $ is the even, real function 
\be\label{defVt}
V_{(\overline t (\xi_\d))} (\vphi, \theta) := 
\frac{2 \aleph  \big( \omega \cdot \pa_\vphi \overline t (\xi_\d)) \big)}{1+\beta_\theta}
  + 
  \frac{{\cal B}_{(\overline t (\xi_\d))}  \big(  
V_{\mu_\delta} (\vphi, \theta )
(1+ \b_\theta ( \theta) )  \big)}{1+ \b_\theta ( \theta )}
\ee
and ${\cal B}_{(\overline t (\xi_\d))} $ is the composition operator defined in \eqref{compdif}. 
\end{lemma}

\begin{pf}
By  \eqref{reprePhi2}  and since $ \Psi  = (1+ \beta_\theta) {\cal B} $ 
 is independent of $ \vphi $ we have
$$
\Phi^{\overline t (\xi_\d) }_{{\cal J}_2}  \circ 
\omega \cdot \pa_\vphi \circ \Phi^{- \overline t (\xi_\d) }_{{\cal J}_2} 
 = \Psi \circ T_{2\aleph \overline t (\xi_\d) }   \circ 
\omega \cdot \pa_\vphi  \circ T_{-2\aleph \overline t (\xi_\d) } \circ \Psi^{-1}   \, . 
$$
Moreover
$ T_{2 \aleph \overline t (\xi_\d)}   \circ 
\omega \cdot \pa_\vphi  \circ T_{-2 \aleph \overline t (\xi_\d)} = \omega \cdot \pa_\vphi - 
2 \aleph 
\big( \omega \cdot \pa_\vphi \overline t (\xi_\d)) \big) \pa_y $
and then
\begin{align}
\Phi^{\overline t (\xi_\d) }_{{\cal J}_2}  \circ 
\omega \cdot \pa_\vphi \circ \Phi^{- \overline t (\xi_\d) }_{{\cal J}_2} &  =
\Psi \circ \Big( \omega \cdot \pa_\vphi - 
2 \aleph 
\big( \omega \cdot \pa_\vphi \overline t (\xi_\d)) \big) \pa_y \Big) \circ \Psi^{-1} \notag \\
 &  = 
\omega \cdot \pa_\vphi - 
\pa_\theta \, 
 \frac{2 \aleph  \big( \omega \cdot \pa_\vphi \overline t (\xi_\d)) \big)}{1+\beta_\theta}  \, ,  
 \label{pezzmiss}
\end{align}
using that
\be\label{diffeo1st} 
\tfrac{1}{1+ \beta_\theta (\theta)} = (1 + \breve \beta_y ( y))_{|y 
= \theta + \beta ( \theta)} = {\cal B} (1 + \breve \beta_y ( y))) \, , 
\ee
and $ \Psi \circ \pa_y \circ  \Psi^{-1}  =  $
$ (1+ \beta_\theta ) \circ {\cal B} \circ \pa_y \circ {\cal B}^{-1}  \circ \frac{1}{1+\beta_\theta}
= $ $  \pa_\theta \circ \frac{1}{1+\beta_\theta} $. 
In the sequel to simplify notation we write $ t  := \overline t (\xi_\d) $.
By  \eqref{reprePhi2} we write 
\be\label{comp-transport}
\Phi^{t}_{{\cal J}_2} 
\circ \pa_\theta \circ V_{\mu_\delta}(\vphi, \theta) \circ \Phi^{-t}_{{\cal J}_2} 
=
\Psi \circ T_{2\aleph t}  \circ \Psi ^{-1} \circ \pa_\theta \circ V_{\mu_\delta}(\vphi, \theta) \circ
 \Psi  \circ T_{-2\aleph t} \circ \Psi^{-1} \, . 
\ee
{\sc Step 1.}  {\it It results} 
\be\label{defV1}
\Psi^{-1}  \circ  \pa_\theta \circ V_{\mu_\delta} (\vphi, \theta)  \circ \Psi  = 
\pa_y\circ V_1(\vphi, y) \quad \text{with} \quad 
V_1(\vphi, y) :=  
{\cal B}^{-1}  \big(
V_{\mu_\delta} (\vphi, \theta )(1+ \b_\theta ( \theta) )\big) \, . 
\ee
Indeed, using the conjugation rules for the multiplication operator for a function $ f $, 
and for the derivative operator $ \pa_\theta $, 
$ {\cal B}^{-1} \circ
f \circ 
{\cal B}   = \big( {\cal B}^{-1} f \big) $, $
{\cal B}^{-1} \circ \pa_\theta \circ 
{\cal B}   = \big( {\cal B}^{-1} (1+ \beta_\theta) \big) \circ \pa_y $, 
we get 
$$
\begin{aligned}
\Psi^{-1} \circ \pa_\theta \circ V_{\mu_\delta} \circ \Psi 
& =
{\cal B}^{-1} \circ \frac{1}{1+\beta_\theta} \circ {\cal B} \circ {\cal B}^{-1} \circ 
\pa_\theta \circ {\cal B} \circ {\cal B}^{-1} \circ V_{\mu_\delta}
\circ {\cal B} \circ {\cal B}^{-1}  \circ (1+ \beta_\theta ) \circ {\cal B} \\
& =
\pa_y  \circ 
\big( {\cal B}^{-1} \big(V_{\mu_\delta} (1+ \beta_\theta)\big) \big) \, .
\end{aligned}
$$
{\sc Step 2.}  
Since $ T_{2\aleph t}  $ and $ \pa_y $ commute, 
we have that
\be\label{defV2}
T_{2\aleph t}  \circ \pa_y \circ  
V_1 (\vphi, y)    \circ T_{-2\aleph t} = \pa_y \circ  V_2(\vphi, y)  
\quad \text{where} \quad V_2 (\vphi, y)  :=  T_{2\aleph t} V_1 \, .
\ee
{\sc Step 3.} 
By Step 1, applied with $ \Psi $ instead of $ \Psi^{-1} $, we get 
\be\label{defV3}
 \Psi \circ \pa_y \circ V_2 (\vphi, y) \circ \Psi^{-1}
= \pa_\theta \circ V_3 (\vphi, \theta) \quad
\text{with} \quad  V_3 (\vphi, \theta) :=  
{\cal B}  (    V_2 (\vphi,y ) (1+ \breve \b_y ( y))) \, . 
\ee
By \eqref{defV2}, \eqref{defV1}, \eqref{diffeo1st}, \eqref{lemdifp}  we obtain 
$ V_3 (\vphi, \theta)  = 
  \tfrac{ {\cal B}_{(t)}  ( V_{\mu_\delta} (1+ \beta_\theta))}{1+ \b_\theta ( \theta )}  $
and by  \eqref{pezzmiss}, \eqref{comp-transport}, \eqref{defV1}, \eqref{defV2}, \eqref{defV3}
we deduce \eqref{primaconc1}-\eqref{defVt}. 
 The function 
$ V_{(\overline t (\xi_\d))}  $ in \eqref{defVt} is even since $ \bar t (\xi_\delta) 
= \wtilde \beta_{2,\d} (\vphi) $ 
is odd, since  $ V_{\mu_\delta}$ is 
even, $ \beta  $ is odd, and 
the function $ \alpha_\e (\vphi,\theta)$ in Lemma \ref{lem:diffest}  is odd in $ (\vphi, \theta)$.
\end{pf}

We now conjugate $   \pa_\theta  \W_0 $. 

\begin{lemma}\label{lem:HRe}
 The conjugation of the Hamiltonian operator 
$   \pa_\theta  \W_0 $ 
under 
 $ \Phi^{ \overline t (\xi_\d) }_{ {\cal J}_2} $ is 
 \be\label{conH}
 \begin{aligned}
 \Phi^{ \overline t (\xi_\d) }_{{\cal J}_2}  \circ \pa_\theta  \W_0  \circ \Phi^{ - \overline t (\xi_\d) }_{{\cal J}_2} 
 & = 
 \pa_\theta  \W_0 + {\cal E}_{ \overline t (\xi_\d) } \, 
 \end{aligned}
\ee
where $ {\cal E}_{ \overline t (\xi_\d) } \in {\rm OPS}^{-\infty} $ satisfies, 
for all $  m,  \a \in \N_0 $,  for any $ s \geq s_0 $,
an estimate as \eqref{stimaRep4}. 
\end{lemma}

\begin{pf}
To simplify notation we write $ t $ for $\overline t (\xi_\d) $.
We recall that $ \partial_\theta W_0 = \tfrac12 {\cal H} + {\cal Q}_\infty $  
see \eqref{laH}, and we conjugate separately these two operators. 
\\[1mm]
{\sc Step 1.} {\it Conjugation of $ {\cal H }$.} By  \eqref{reprePhi2}  we write 
$ \Phi^{t}_{{\cal J}_2}  \circ {\cal H}  \circ \Phi^{-t}_{{\cal J}_2}= $ 
 \begin{align} 
 (1+ \beta_\theta) \circ  {\cal B} \circ T_{2\aleph t}  \circ {\cal B}^{-1} \circ
\frac{1}{1+ \beta_\theta} \circ {\cal B} \circ {\cal B}^{-1} \circ {\cal H} \circ {\cal B} 
\circ {\cal B}^{-1} \circ 
(1+ \beta_\theta) \circ  {\cal B} \circ T_{-2\aleph t}    \circ {\cal B}^{-1} \circ
\frac{1}{1+ \beta_\theta}  \, . \notag
\end{align}
Using the conjugation rules  
$ (1 + \beta_\theta) \circ {\cal B} = {\cal B} \circ \big({\cal B}^{-1}(1 + \beta_\theta)\big) $
and
$   {\cal B}^{-1} \circ
\tfrac{1}{1+ \beta_\theta}   = \big( {\cal B}^{-1} \tfrac{1}{1+ \beta_\theta} \big) 
 \circ {\cal B}^{-1} $, 
we deduce that
\begin{align}
\Phi^{t}_{{\cal J}_2}  \circ {\cal H}  \circ \Phi^{-t}_{{\cal J}_2} 
 = {\cal B} \circ f    
\circ T_{2\aleph t}  \circ  \frac{1}{f}  \circ {\cal B}^{-1}\circ {\cal H}  \circ
{\cal B} 
\circ f  \circ  T_{-2\aleph t}    \circ 
\frac{1}{f}
 \circ {\cal B}^{-1} \label{coniu2}
 \end{align}
where  $ f := {\cal B}^{-1} ( 1 + \beta_\theta) $.
Moreover, since
$ f \circ T_{2\aleph t} \circ \frac{1}{f} = T_{2\aleph t} \circ  \frac{f(\theta-2\aleph t)}{f(\theta )}$, and $ f \circ T_{-2\aleph t} \circ \frac{1}{f} = \frac{f(\theta)}{f(\theta - 2\aleph t)} \circ T_{-2\aleph t} $,
then
 \begin{align}
  \Phi^{t}_{{\cal J}_2}  \circ {\cal H}  \circ \Phi^{-t}_{{\cal J}_2} 
  & = \big({\cal B} \circ T_{2\aleph t} \circ {\cal B}^{-1}\big) \circ \frac{1}{h(t)}  \circ {\cal H}   \circ 
h(t)  \circ \big({\cal B} \circ T_{2\aleph t} \circ {\cal B}^{-1}\big)^{-1} \label{coniu3}
 \end{align}
where $ h  $ is the  $ {\cal C}^\infty $ function
$ h := {\cal B} \big( \frac{f}{f (\cdot  -2\aleph \bar t (\xi_\delta)) } \big) $, $ 
f := {\cal B}^{-1} ( 1 + \beta_\theta) $. 
In view of \eqref{invediff} and \eqref{bartesxidelta}, we have 
 \be\label{estim:h}
\|  h  - 1 \|_{s}^{k_0,\upsilon }\lesssim_s  \e 
(1 + \| \fracchi_0\|_{s+\sigma}^{k_, \upsilon})\, .
 \ee
Then, by \eqref{coniu3} and  \eqref{lemdifp} we have 
\be
 \Phi^{t}_{{\cal J}_2}  \circ  {\cal H} \circ
   \Phi^{-t}_{{\cal J}_2}
   = {\cal B}_{(t)} \circ \Big({\cal H} +\frac{1}{h} \circ [{\cal H}, h-1 ] \Big)\circ{\cal B}_{(t)}^{-1} = 
   {\cal H} + {\cal E}_1 + {\cal E}_2 \label{conH00}
\ee
with
 $ {\cal E}_1 := 
{\cal B}_{(\bar t (\xi_\delta))} 
\circ  {\cal H}  \circ {\cal B}_{(\bar t (\xi_\delta))}^{-1}-{\cal H}  $ and $
{\cal E}_2 := 
{\cal B}_{(\bar t (\xi_\delta))} \circ 
\frac{1}{h} \circ [{\cal H}, h-1 ] \circ 
{\cal B}_{(\bar t (\xi_\delta))}^{-1} $. 
By  Lemma \ref{lemma:conjug-Hilbert} and \eqref{estim:a} we deduce that 
$ {\cal E}_1 \in  {\rm OPS}^{-\infty}   $ 
 and it satisfies an estimate as \eqref{stimaRep4}.
Lemma \ref{lemma:commutator-Hilbert}, \eqref{estim:h}, 
and  Lemma \ref{lemma:conj-integ-op} imply that $ {\cal E}_2 
$  is in $  \text{OPS}^{-\infty } $ and 
 it satisfies  \eqref{stimaRep4}.
 \\[1mm]
{\sc Step 2.} {\it Conjugation of $ {\cal Q}_\infty  $.}
In a similar way to \eqref{conH00} we get
 \be \label{conH2}
 \Phi^{t}_{{\cal J}_2} \circ  {\cal Q}_\infty  \circ
  \Phi^{-t}_{{\cal J}_2}
 = {\cal Q}_\infty+ {\cal E}_3  \quad
\text{ with} \quad  
  {\cal E}_3:={\cal B}_{(\bar t (\xi_\delta))} \circ \tfrac{1}{h }  \circ {\cal Q}_\infty  \circ 
h  \circ {\cal B}_{(\bar t (\xi_\delta))}^{-1} - {\cal Q}_\infty \, . 
 \ee
Since $ {\cal Q}_\infty\in  {\rm OPS}^{-\infty} $ then, by  Lemma \ref{lemma:conj-integ-op}, the term ${\cal E}_3 \in  \text{OPS}^{-\infty } $ 
and satisfies an estimate as \eqref{stimaRep4}.

In conclusion, by \eqref{conH00}, \eqref{conH2} we get \eqref{conH}
with  $ {\cal E}_{\bar t (\xi_\delta)} := 
\tfrac12 ({\cal E}_1 + {\cal E}_2) + {\cal E}_3 $
that satisfies an estimate as \eqref{stimaRep4}  as
each ${\cal E}_i $, $ i = 1,2,3 $. 
This concludes the proof of Lemma \ref{lem:HRe}. 
\end{pf}

\begin{pfn}{\sc of Proposition \ref{lem:Lorto}.} 
By Lemma \ref{lem:AI} and Lemmata \ref{lem11.4}, \ref{lemC5}, \ref{lem:HRe} 
we obtain Proposition \ref{lem:Lorto} with the function 
\be\label{deflav}
{\cal V} (\vphi, \theta):=  
V_{(\overline t (\xi_\d))} (\vphi, \theta) +\breve{\mu}_\e (\vphi) \, g_\gamma (\theta)
\, , 
\ee
where $V_{(\overline t (\xi_\d))} (\vphi, \theta)  $ is defined in 
\eqref{defVt}, the function $\breve{\mu}_\e (\vphi) $ in Lemma \ref{lem:AI},
and the remainder 
$ {\rm R}_\e:=  $ $ \Phi^{ \overline{t}(\xi_\delta)}_{{\cal J}_2}  \circ  \partial_\theta R_\e  \circ \Phi^{ \overline{t}(\xi_\delta)}_{-{\cal J}_2}  +{\cal E}_{\overline{t}(\xi_\delta)} $. 
The function $ {\cal V}$ in \eqref{deflav}  has the form \eqref{new V}-\eqref{new V1} 
by \eqref{defVt}, \eqref{bartesxidelta}, \eqref{defg0AA}, \eqref{bpicco}, \eqref{stima-defmudelta} and \eqref{stima-defmuep}. 
Moreover 
\eqref{stimaRep4}  follows by  Lemmata \ref{lem11.4}-\ref{lite3} 
and \ref{lem:HRe}. 
\end{pfn}

\subsection{Almost 
approximate reduction  up to smoothing remainders}\label{sec:diffeo}

We now conjugate the quasi-periodic linear operator 
$ \omega \cdot \partial_\vphi  -  
\pa_\theta \circ {\cal V}+ \pa_\theta W_0  +{\rm R}_\e $
in \eqref{linKOK}
by a $ \vphi $-dependent family of symplectic transformations 
(Definition \ref{def:sympl})
\be\label{def:Phi}
{\mathfrak S} := {\mathfrak S} (\vphi)  := ( 1 + {\mathfrak b}_\theta ) 
{\fB}  (\vphi)
\, , 
\ee
induced by a  $ \vphi $-dependent family of diffeomorphisms of $ \T $, 
$ y = \theta +  {\mb}  (\vphi, \theta ) $
where $ \mb (\vphi, \theta )  $ is a small periodic function  
chosen in Lemma \ref{lem:traspo} and $ {\fB} := {\fB} (\vphi)  $ 
is the  induced composition operator 
\be\label{defB}
({\fB}  u )(\vphi, \theta) := u( \vphi, \theta + \mb (\vphi, \theta) )\, .  
\ee
Let
$ N_{\bar{\mathtt n}}:= N_0^{\chi^{\bar{\mathtt n}}} $, $  \chi=3/2 \, $, $ N_{-1}:=1 $, 
and define the constants 
\be\label{choiceab}
\ta := 3( \tau_1+1)  \, , \quad \tau_1 :=(  k_0 + 1)\tau_0+ k_0 \, , \quad \tb := [\ta] + 2 \, .
\ee

\begin{lemma}\label{lem:traspo}
{\bf (Almost straightening of the transport)}
There exists $\tau_2:=\tau_2(|\mathbb{S}|,\tau)> \tau_1+{\mathtt a}+1$ such that, for all $S>s_0+ k_0$, there are $N_0:=N_0(S, {\mathtt b})\in \mathbb{N}$ and $\delta:=\delta(S, {\mathtt b})\in(0,1)$ such that, if $N_0^{\tau_2}\varepsilon \upsilon^{-1}<\delta$,  the following holds true. For any $ \bar {\mathtt n}\in \mathbb{N}_0 $ there exist
\begin{enumerate}
\item 
a constant $ {\mathtt m}_{\bar{\mathtt n}}: ={\mathtt m}_{\bar{\mathtt n}} (\omega, \gamma)\in \R $, with $ {\mathtt m}_{0} := \Omega_\gamma $, 
 defined for any $(\omega, \gamma) \in \R^{|\ST|} \times [\g_1, \g_2]$, of the form 
\be\label{c1-picco}
{\mathtt m}_{\bar{\mathtt n}} = \Omega_\gamma + 
{\mathtt r}_{\e,\bar{\mathtt n}}  
 \, , \quad 
 |{\mathtt r}_{\e,\bar{\mathtt n}}|^{k_0, \upsilon} \lesssim_{k_0} \e  \, ,
\ee
and 
$ | {\mathtt m}_{\bar \tn}- {\mathtt m}_{\bar \tn-1}|^{k_0,\upsilon}  
		\lesssim_{k_0} \e N_{\bar \tn-2}^{-\ta} $, 
		for any $ \bar \tn \geq 2 $;
\item
a  function $ \mb (\vphi, \theta ):= \mb_{\bar{\mathtt n}} (\vphi, \theta ) $, odd in $ (\vphi, \theta )$, defined for all $ (\omega, \gamma )\in \R^{|\ST|} \times [\g_1, \g_2] $,  
satisfying for some $ \s:=\s(\tau,\ST,k_0) > 0 $, the estimate (independent of $ \bar \tn $)  
\be\label{betainfty} 
\| \mb \|^{k_0,\upsilon}_s \lesssim_S \e \upsilon^{-1}  
(1 + \| \fracchi_0 \|_{s+\sigma}^{k_0, \upsilon}) \, ,  \quad \forall s_0\leq s\leq S \, ; 
\ee
\end{enumerate}
such that, for any $ (\omega, \gamma)$ in 
\be\label{1-tras}
\begin{aligned}
\mathtt{TC}_{\bar{\mathtt n}+1}(2\upsilon,\tau) := \Big\{  (\omega, \gamma) \in \R^{|\ST|} \times &[\g_1, \g_2] 
\, : \, 
| \om \cdot \ell  + {\mathtt m}_{\bar{\mathtt n}}  j  | \geq {2 \upsilon }{\langle \ell \rangle^{-\tau}} \, \\ & \forall (\ell,j)\in (\mathbb{Z}^{|\ST|}\times \mathbb{Z})\setminus \{0\} \, ,
\, |(\ell,j) |\leq N_{\bar{\mathtt n}} \Big\} \,  ,
\end{aligned}
\ee
we have
\be\label{NMKH}
\bphi^{-1}  \big( \om \cdot \pa_\vphi - \pa_\theta\,  {\cal  V} +
\pa_\theta \, W_0 +{\rm R}_\e \big) 
\bphi   = \om \cdot \pa_\vphi + 
 {\mathtt m}_{\bar{\mathtt n}} \pa_y + \pa_y \, W_0  
 +  {\mathfrak R}+  {\mathfrak P}_{\bar{\mathtt n}} \, , 
\ee
where the transformation
$\bphi =\bphi_{ \bar{\mathtt n}}  $ is defined in \eqref{def:Phi}. Each term in the right hand side of \eqref{NMKH}
is defined for all $ (\omega,\gamma) \in 
\R^{|\ST|} \times [\g_1,\g_2] $ and
\\[1mm]
3. $ {\mathfrak R} $ is  a Hamiltonian reversible real
operator in $ {\rm OPS}^{-\infty}$ 
satisfying,  
 for all $  m,  \a \in \N_0 $,  for some $ \sigma (m, \alpha) > 0 $,
\be\label{fina-Re1}
| {\mathfrak R}   |_{-m, s, \a}^{k_0, \upsilon} \lesssim_{m, S, \alpha, k_0} \e   \upsilon^{-1}
 (1 + \| \fracchi_0 \|_{s + \sigma (m,\a)}^{k_0, \upsilon}) \, ,  \quad \forall s_0\leq s\leq S \, ; 
\ee
4. 
$ {\mathfrak P}_{\bar{\mathtt n}} := \pa_y  \circ {\mathtt p}_{\bar{\mathtt n}}$
with a real even function ${\mathtt p}_{\bar{\mathtt n}} (\vphi, y ) $ satisfying 
\be\label{petainfty} 
\| {\mathtt p}_{\bar{\mathtt n}} \|^{k_0,\upsilon}_s \lesssim_s 
\e  N_{\bar{\mathtt n}-1}^{-{\mathtt a}}  
(1 + \| \fracchi_0 \|_{s+\sigma+{\mathtt b}}^{k_0, \upsilon}) \, ,  \quad \forall s_0\leq s\leq S \, ; 
\ee  
5. for any $s_1$ as in \eqref{def:s1},
$ |\Delta_{12} {\mathtt m}_{\bar\tn} |\lesssim \varepsilon \|i_1-i_2 \|_{s_1+\sigma}$, 
$ |\Delta_{12} \mb |\lesssim \varepsilon \upsilon^{-1} \|i_1-i_2 \|_{s_1+\sigma} $. 
Moreover, for  $0 < 2\upsilon-\rho<2\upsilon<1$,
\be\label{inclusion-transport}
\varepsilon C(s_1)N_{\bar{\mathtt n}}^{\tau+1}\|i_1-i_2 \|_{s_1+\sigma}\leq \rho\quad\Rightarrow\quad \mathtt{TC}_{\bar{\mathtt n}+1}(2\upsilon,\tau) (i_1)\subseteq \mathtt{TC}_{\bar{\mathtt n}+1}(2\upsilon-\rho,\tau)(i_2) \, .
\ee
\end{lemma}

\begin{pf}
By \eqref{new V}-\eqref{new V1} and (a small variant of) the almost straightening 
Theorem A.2 and Corollary A.4 of \cite{BFM}, cfr. \cite{FGMP},
there exist a constant $ {\mathtt m}_{ \bar {\mathtt n}} (\omega, \gamma) $ satisfying \eqref{c1-picco}, 
 an odd  function $\mb_{\bar{\mathtt n}} (\vphi, \theta )  $ satisfying \eqref{betainfty}, such that, for any $ (\omega, \gamma)$ in  $ \mathtt{TC}_{\bar{\mathtt n}+1}(2\upsilon,\tau) $
it results
\be\label{B-1Btra}
{\fB}^{-1} \big( 
\om \cdot \pa_\vphi -  {\cal  V} (\vphi, \theta ) \pa_\theta\,  \big) {\fB} =
\om \cdot \pa_\vphi + ({\mathtt m}_{ \bar {\mathtt n}} + {\mathtt p}_{\bar{\mathtt n}}
(\vphi, y ))  \pa_y  
\ee
for a function $ {\mathtt p}_{\bar{\mathtt n}} (\vphi, y ) $ satisfying \eqref{petainfty}. 
The function $\mb_{\bar{\mathtt n}} (\vphi, \theta )  $ 
is odd in $ (\vphi, \theta )$ since the function 
${\cal V} (\vphi, \theta )  $ in Proposition \ref{lem:Lorto} is even. 
Since (see \cite{BFM1}) 
$$
{\fB}^{-1} (\om \cdot \pa_\vphi -  {\cal  V}\pa_\theta\,  ){\fB} =
\om \cdot \pa_\vphi + {\cal  V}_1 \pa_y \, , \quad 
{\cal  V}_1 := \fB^{-1} ( \omega\cdot \pa_\vphi \mb - {\cal V}(1+ \mb_\theta))  \, , 
$$
if and only if 
$ {\bphi}^{-1} (\om \cdot \pa_\vphi -  \pa_\theta \circ {\cal  V} ){\bphi} =
\om \cdot \pa_\vphi + \pa_y \circ {\cal  V}_1 $,   
we deduce, by \eqref{B-1Btra}, for any $ (\omega, \gamma)$ in  $ \mathtt{TC}_{\bar{\mathtt n}+1}(2\upsilon,\tau) $, the conjugation 
$ \bphi^{-1} (\om \cdot \pa_\vphi - \pa_\theta\,  {\cal  V} ) 
\bphi   = \om \cdot \pa_\vphi + 
 {\mathtt m}_{\bar{\mathtt n}} \pa_y + 
 \pa_y   \circ  {\mathtt p}_{\bar{\mathtt n}} $. 
 
We now conjugate  the whole Hamiltonian 
operator 
 in \eqref{NMKH}  
with the symplectic and reversibility preserving 
map $\bphi  $ 
(since $ \mb $ is odd $ (\vphi, \theta)$).
Recalling \eqref{laH}, \eqref{def:Phi} and \eqref{defB}, 
we have 
\begin{align}
\bphi^{-1} \pa_\theta \W_0 \bphi 
&   = 
\tfrac12 \fB^{-1} \tfrac{1}{1+\mb_\theta } {\cal H} (1 +\mb_\theta  )
\fB  + 
\bphi^{-1} {\cal Q}_\infty\bphi  \, .  \label{terza-con}
\end{align}
We expand
$  \fB^{-1} \tfrac{1}{1+\mb_\theta } {\cal H} (1 +\mb_\theta ) \fB  
 =   $ $ {\cal H} + \big(\fB^{-1} {\cal H} \fB  - {\cal H} \big) +  $
$ \fB^{-1} 
\tfrac{1}{1+\mb_\theta } [{\cal H},\mb_\theta  ]  \fB $ where,   
by Lemmata 
\ref{lemma:conjug-Hilbert}, \ref{lem:Int}  and \ref{lemma:commutator-Hilbert}, 
 the operators $ \fB^{-1} {\cal H} \fB  - {\cal H} $ and $ [{\cal H},\mb_\theta  ] $   
are in $ {\rm OPS}^{-\infty } $. 
Moreover
$ \bphi^{-1} {\cal Q}_\infty \bphi  =  
{\cal Q}_{\infty} + \big(\fB^{-1} \tfrac{1}{1+\mb_\theta} {\cal Q}_{\infty} 
(1+\mb_\theta) \fB  - {\cal Q}_{\infty} \big) $
and, since ${\cal Q}_{\infty}  $ is in $ {\rm OPS}^{-\infty} $, 
  the last term  is in $ {\rm OPS}^{-\infty} $ 
(see Lemma \ref{lemma:conj-integ-op}). In conclusion, by \eqref{QPLIN},
\eqref{NMKH}, \eqref{laH}, \eqref{terza-con} and the above decomposition 
we deduce \eqref{NMKH} with the remainder 
\be\label{defR1expli}
{\mathfrak R}  :=  \bphi^{-1}  {\rm R}_\e  \bphi +
\tfrac12 
\fB^{-1} \tfrac{1}{1+\mb_\theta }   [{\cal H},\mb_\theta  ] \fB +
\tfrac12 \big(\fB^{-1} {\cal H} \fB  - {\cal H} \big) 
+\fB^{-1} \tfrac{1}{1+\mb_\theta} {\cal Q}_{\infty} 
(1+\mb_\theta) \fB  - {\cal Q}_{\infty}  \, . 
\ee
Note that $ {\fB} $ and $ \bphi $ are defined for all the values of
$ (\omega, \gamma) \in \R^{|\ST|} \times [\g_1, \g_2]$
because the function $ \mb = \mb_{\bar{\mathtt n}} $ is defined for all 
$ (\omega, \gamma) \in \R^{|\ST|} \times [\g_1, \g_2] $. Thus ${\mathfrak R}   $ in 
\eqref{defR1expli} is defined on the whole $ \R^{|\ST|} \times [\g_1, \g_2] $. 
By 
\eqref{stimaRep} and 
applying  Lemmata  \ref{lemma:conj-integ-op}, \ref{lemma:commutator-Hilbert}, 
\ref{lem:Int}, \ref{lemma:conjug-Hilbert}  
and  \eqref{betainfty} we deduce that 
$ {\mathfrak R}  $ 
satisfies \eqref{fina-Re1}.
\end{pf}

Now we deduce a similar conjugation result for the operator 
$ {\cal L}_\omega $ defined
in \eqref{Lomega def}, which  acts on the normal subspace $ {\frak H}_{\ST,2}^\bot $.
The operator $ {\cal L}_\om$ is expressed in terms of $ {\cal L}_K $ as in \eqref{K 02}. 
The operator $ {\cal L}_\om$ is reversible as stated in Lemma \ref{thm:Lin+FBR}.
We conjugate $ {\cal L}_\om $ via the map
\begin{equation}\label{Phi 1 Phi 2 proiettate}
\bphi_\bot  := \Pi_{\ST,2}^\bot   \, \bphi \,  \Pi_{\ST,2}^\bot 
\end{equation}
where  $\bphi$ is defined in \eqref{def:Phi} and 
$ \Pi_{\ST,2}^\bot $ is the $ L^2 $ projector on 
$ {\frak H}_{\ST,2}^\bot $. We first have the following lemma.

\begin{lemma} (\cite{BBM-auto}) \label{Lemma:trasformazioni finali W}
There exists a constant $\sigma > 0$  such that, 
assuming \eqref{ansatz 0} with $\perd \geq \sigma $, 
 for any $ S > s_0$ there exists a constant 
$\delta(S) > 0$ such that, if $\e \upsilon^{-1} \leq \delta(S)$, then the operator
$ \bphi_\bot $ defined in \eqref{Phi 1 Phi 2 proiettate} is invertible and 
for all $ s \in [s_0,S] $, 
for all  $h : = h(\lambda) \in H^{s + \sigma}_\bot $,  
\begin{align} \label{stima Phi 1 Phi 2 proiettate}
\| \bphi_\bot^{\pm 1} h \|_s^{k_0, \upsilon} 
& \lesssim_{S} \| h \|_{s + \sigma }^{k_0, \upsilon} 
+ \| \fracchi_0 \|_{s + \sigma }^{k_0, \upsilon} \| h \|_{s_0 + \sigma }^{k_0, \upsilon}\, . 
\end{align}
Moreover the operator $  \bphi_\bot $ is reversibility preserving. 
\end{lemma}

In the sequel we do not keep further the Hamiltonian structure 
of the conjugated operator
and we preserve just  the reversible one. The main of conclusion of this section is the following proposition:

\begin{proposition}\label{prop: sintesi linearized}
{\bf (Almost approximate reduction of $ {\cal L}_\om$ up to smoothing remainders)}
For any $  \bar \tn \in \N_0 $
and for all $  (\om, \gamma) \in \mathtt{TC}_{ \bar {\mathtt n}+1}(2\upsilon,\tau) $ defined in \eqref{1-tras}, 
the operator ${\cal L}_\om$ in \eqref{Lomega def}, i.e. \eqref{K 02}, is conjugated  
to the real reversible  operator ${\cal L}_\bot $, namely 
\begin{equation}\label{forma finale operatore pre riducibilita}
 \bphi_\bot^{-1} {\cal L}_\om  \bphi_\bot = 
\om \cdot \pa_\vphi \mathbb{I}_\bot +
\Pi_{\ST,2}^\bot  \big(  {\mathtt m}_{ \bar{\mathtt n}}   \pa_\theta  
+  \pa_\theta \W_0  \big)_{|\acca_{\ST,2}^\bot} + {\mathfrak R}_\bot+{\mathfrak P}_{\bot,\bar{\mathtt n}} + {\mathfrak R}_{Z} :={\cal L}_\bot  
\end{equation}
where  $\mathbb{I}_\bot$ denotes the identity map of $ {\frak H}_{\ST,2}^\bot $,  each term in the right hand side of \eqref{forma finale operatore pre riducibilita} 
is defined for all $ (\omega,\gamma) \in 
\R^{|\ST|} \times [\g_1,\g_2] $  and
\\[1mm]
1. the constant ${\mathtt m}_{\bar \tn}   : \R^{|\ST|} \times [\g_1, \g_2] \to \R $ 
satisfies  \eqref{c1-picco};
\\[1mm]
2.
the reversible real Hamiltonian operator  $ \pa_\theta \W_0 $ 
as the form \eqref{laH};
\\[1mm]
3.
$ {\mathfrak R}_\bot $ is  a   reversible real 
operator in $ {\rm OPS}^{-\infty}$
satisfying, 
 for all $  m \in \N $,  for some $ \sigma (m)$,   
\be\label{fina-Re}
|  {\mathfrak R}_\bot  |_{-m, s,0}^{k_0, \upsilon} \lesssim_{m, S, k_0} \e  \upsilon^{-1} 
 (1 + \| \fracchi_0 \|_{s + \sigma (m)}^{k_0, \upsilon}) \, , \quad \forall s \in [s_0,S]  \, ; 
\ee
4. the operator ${\mathfrak P}_{\bot,\bar{\mathtt n}}  $ satisfies, for some
$ \sigma > 0  $,   for all $  s_0\leq s\leq S $, 
\be\label{petainfty2} 
\| {\mathfrak P}_{\bot,\bar{\mathtt n}}  h \|^{k_0,\upsilon}_s \lesssim_s 
\e  N_{\bar{\mathtt n}-1}^{-{\mathtt a}} 
\big( \| h \|_{s+\sigma}^{k_0, \upsilon} + \| \fracchi_0 \|_{s+\sigma+{\mathtt b}}^{k_0, \upsilon} \| h \|_{s_0+\sigma}^{k_0, \upsilon}\big) \, ;
\ee
5. the operator $ {\mathfrak R}_{Z}$ satisfies, for some $ \s > 0 $,
for any $ s \in [s_0,S] $,
\be\label{estimate:RZ2}
\|  {\mathfrak R}_{Z} h \|_s^{k_0,\upsilon} \lesssim_S 
\big( \| {Z}  \|_{s+\sigma}^{k_0,\upsilon} 
+ 
\| {Z}  \|_{s_0+\sigma}^{k_0,\upsilon}
 \| \fracchi_0 \|_{s+\sigma}^{k_0,\upsilon}\big)\|  h \|_{s_0+\sigma}^{k_0,\upsilon}  
+ \| {Z}  \|_{s_0+\sigma}^{k_0,\upsilon}\|  h \|_{s+\sigma}^{k_0,\upsilon}
 \, .
\ee
\end{proposition}

\begin{pf}
Set for brevity 
$ {\cal L}_1 := \omega \cdot \partial_\vphi  -  
\pa_\theta \,  {\cal V}+ \pa_\theta W_0  +{\rm R}_\e $. By \eqref{K 02} we have 
$$
 \bphi_\bot^{-1} {\cal L}_\om  \bphi_\bot 
   =  \bphi_\bot^{-1} \Pi_{\ST,2}^\bot  {\cal L}_K   \bphi_\bot +  {\cal R}_I  \qquad \text{where}
 \qquad {\cal R}_I :=  
 \e \bphi_\bot^{-1}  \Pi_{\ST,2}^\bot  \pa_\theta {\cal R} \bphi_\bot  \, .  
$$
 Thus by \eqref{linKOK}, \eqref{Phi 1 Phi 2 proiettate} 
and using the identity 
 $ \Pi_{\ST,2}^\bot    = {\rm Id}  -  \Pi_{\ST,2,0} $, 
   where  
$  \Pi_{\ST,2,0} := \Pi_{\ST,2} + \pi_0 $ with the projector 
$\pi_0 $,  
defined on a $ 2 \pi $-periodic function $ u(\theta)$,  as
$ \pi_0 u := \frac{1}{2\pi} \int_\T u(\theta)\, d \theta $, 
we get 
 \begin{align}\label{R1}
  \bphi_\bot^{-1} {\cal L}_\om  \bphi_\bot  & = 
\bphi_\bot^{-1} \Pi_{\ST,2}^\bot  {\cal L}_1  \bphi \,  \Pi_{\ST,2}^\bot+{\mathfrak R}_{Z}  + {\cal R}_{II} + {\cal R}_{I}   
\end{align}
with
$ {\mathfrak R}_{Z} :=\bphi_\bot^{-1}  {\cal R}_{Z}   \, \bphi_\bot $ and 
$ {\cal R}_{II} := - \bphi_\bot^{-1}  \Pi_{\ST,2}^\bot {\cal L}_1  \Pi_{\ST,2,0}   \, \bphi \,  \Pi_{\ST,2}^\bot+\bphi_\bot^{-1}  {\mathscr R}    \, \bphi_\bot $. 
Moreover, setting 
$ {\cal L}_2 :=  \bphi^{-1}  {\cal L}_1 \bphi  $, 
we have 
\be
\bphi_\bot^{-1} \Pi_{\ST,2}^\bot  {\cal L}_1  \bphi \,  \Pi_{\ST,2}^\bot  = \bphi_\bot^{-1}   \Pi_{\ST,2}^\bot  \bphi  {\cal L}_2\,  \Pi_{\ST,2}^\bot     = 
\bphi_\bot^{-1}  
 \Pi_{\ST,2}^\bot   \bphi  \Pi_{\ST,2}^\bot  
{\cal L}_2\,  \Pi_{\ST,2}^\bot + {\cal R}_{III} \label{R3}  
\ee
with
$ {\cal R}_{III} := \bphi_\bot^{-1}   
\Pi_{\ST,2}^\bot   \bphi  \Pi_{\ST,2,0}  
{\cal L}_2\,   \Pi_{\ST,2}^\bot $, and  the identity 
\eqref{R1}-\eqref{R3}, \eqref{Phi 1 Phi 2 proiettate} imply  that 
$ \bphi_\bot^{-1} {\cal L}_\om  \bphi_\bot  =  $ $  \Pi_{\ST,2}^\bot
{\cal L}_2\,  \Pi_{\ST,2}^\bot +{\mathfrak R}_{Z} + {\cal R}_f $ with 
$  {\cal R}_f :=   {\cal R}_I + {\cal R}_{II} + {\cal R}_{III} $. 
Finally, by \eqref{NMKH}, 
we get \eqref{forma finale operatore pre riducibilita} with
remainders
$$
 {\mathfrak P}_{\bot,\bar{\mathtt n}}  :=\Pi_{\ST,2}^\bot \,  {\mathfrak P}_{\bar{\mathtt n}}  \,\Pi_{\ST,2}^\bot \quad \text{where} \quad 
  {\mathfrak P}_{\bar{\mathtt n}} = \pa_y  \circ {\mathtt p}_{\bar{\mathtt n}} \, ,  
\qquad {\mathfrak R}_\bot :=\Pi_{\ST,2}^\bot   {\mathfrak R}\, \Pi_{\ST,2}^\bot+ {\cal R}_f \, . 
$$
By \eqref{fina-Re1}  the operator $\Pi_{\ST,2}^\bot   {\mathfrak R}\, \Pi_{\ST,2}^\bot$ 
satisfies an estimate like \eqref{fina-Re}. 
 In view of  \eqref{R1}-\eqref{R3}, the operator $ {\cal R}_f $ has the finite rank form 
\eqref{finite_rank_R} with functions $ g_j, \chi_j $ satisfying 
$ \max \{ \| g_j \|_s^{k_0,\upsilon}, \| \chi_j \|_s^{k_0,\upsilon}  \}  \lesssim_s \e
(1 + \| \fracchi_0 \|_{s+\sigma}^{k_0,\upsilon})  $ and 
thus, by Lemma \ref{lem:Int}, 
the integral operator $ {\cal R}_f $ satisfies  an estimate like \eqref{fina-Re}. Finally \eqref{estimate:RZ2} follows by \eqref{estimate:RZ}, \eqref{stima Phi 1 Phi 2 proiettate},
\eqref{ansatz_I0_s0}. 
\end{pf}

\section{Reducibility and inversion} 
\label{sec:redu}

In this section we almost-diagonalize 
the quasi-periodic  real reversible operator 
\be\label{inizioKAMRe}
{\mathfrak L}_\bot    := \om \cdot \pa_\vphi  \mathbb{I}_\bot
+
\Pi_{\ST,2}^\bot  \big(  {\mathtt m}_{\bar {\mathtt n}}   \pa_\theta  
+  \pa_\theta \W_0  \big)_{|\acca_{\ST,2}^\bot} + {\mathfrak R}_\bot \, , 
\ee
obtained by  neglecting from $ {\cal L}_\bot $ 
in \eqref{forma finale operatore pre riducibilita} the remainders  
${\mathfrak P}_{\bot, \bar {\mathtt n}}   $ and $ {\mathfrak R}_{Z}  $,  
by a KAM iterative scheme.
The operator $ {\mathfrak L}_\bot  $   acts on 
$ {\frak H}_{\ST,2}^\bot =  \oplus_{n \in \ST^c} V_n  $, $ 
\ST^c := \N \setminus (\ST \cup \{2\}) $, 
where 
\be\label{def:Vnprima}
V_n := \Big\{ q(\theta ) =  \a_n {\mathtt c}_n (\theta) + 
\b_n {\mathtt s}_n (\theta) \, , \ ( \a_n,  \b_n ) \in \R^2   \Big\} 
\ee 
and the functions $ {\mathtt c}_n $, 
$ {\mathtt s}_n $ are defined in \eqref{defcnsn}.
We  represent $ {\mathfrak L}_\bot  $ 
as a 
 matrix of $ 2 \times 2 $ matrices,  
in the basis
$ \{ {\mathtt c}_n , {\mathtt s}_n  \}_{n \in \ST^c} $. 
We now present such matrix representation and its main  properties. 
\\[1mm]
{\bf Matrices with decay.} We consider  $ \vphi $-dependent real linear operators 
 $ A(\vphi) $ acting on
$  \oplus_{n \in \ST^c} V_n  $, 
where $ V_n $ are the $ 2 $-dimensional $ L^2$ pair-wise orthogonal  subspaces 
$ V_n $  in \eqref{def:Vnprima}. 
The action of  
$ A(\vphi) $ is represented, with respect to the basis $ \{ {\mathtt c}_n  \, , {\mathtt s}_n  \} $ 
 in each  $ V_n $, 
by the infinite dimensional matrix 
$ \big( [A]_n^{n'} (\vphi) \big)_{n,n' \in \ST^c}   $ 
where 
\be\label{Ann'}
 [A]_n^{n'} (\vphi) = 
 \begin{pmatrix}
 (A(\vphi)  \tc_{n'}  , \tc_{n} )_{L^2}  
&  (A(\vphi) \ts_{n'} , \tc_{n})_{L^2}   \\
 (A(\vphi) \tc_{n'}  , \ts_{n} )_{L^2}  
& (A(\vphi) \ts_{n'} , \ts_{n} )_{L^2} 
\end{pmatrix} \, . 
\ee
We also Fourier expand with respect to  $ \vphi $ each 
\be\label{Ann'-Fourier}
[A]_n^{n'} (\vphi) = {\mathop \sum}_{\ell \in \Z^{|\ST|}} \widehat {[A]_n^{n'}} (\ell) e^{\ii \ell \cdot \vphi} \, ,
\quad \overline{\widehat {[A]_n^{n'}} (\ell)} = \widehat {[A]_n^{n'}} (-\ell)  \, , 
\ee
where  $ \widehat {[A]_n^{n'}} (\ell) \in \text{Mat}_2 (\C) $.
We identify a function
\be\label{defqthetaFou}
q(\vphi, \theta) 
 = {\mathop \sum}_{n \in \ST^c, \ell \in \Z^{|\ST|}} 
\big( \widehat \alpha_{n} (\ell) \, \tc_n (\theta ) +  \widehat{\b}_{n}(\ell) \, \ts_n (\theta ) \big)
e^{\ii \ell \cdot \vphi}
\ee
with the sequence of Fourier coefficients 
$ \big\{ (\widehat \alpha_{n} (\ell), \widehat \beta_{n} (\ell) )\big\}_{n\in \ST^c, \ell \in \Z^{|\ST|}} $, 
and 
the function $ \Pi_K q $ with the truncated sequence 
\be\label{PiK2co}
\Pi_K  \big\{ (\widehat \alpha_{n} (\ell), \widehat \beta_{n} (\ell) )\big\}_{n\in \ST^c, \ell \in \Z^{|\ST|}} = \big\{ (\widehat \alpha_{n} (\ell), \widehat \beta_{n} (\ell) )\big\}_{n\in \ST^c, \ell \in \Z^{|\ST|}, |(n,\ell)| \leq K}  \, .
\ee
 Moreover we identify
the operator $ A(\vphi) $ with the matrix 
\be\label{matricesAnn}
 \big( {[A]_{n,\ell}^{n',\ell'}}  \big)_{(n,\ell),(n',\ell')\in \ST^c \times \Z^{|\ST|} } \, , \quad
  {[A]_{n,\ell}^{n',\ell'}} := \widehat {[A]_n^{n'}} (\ell- \ell') \, ,  
  \ee 
of $ 2 \times 2 $ complex matrices $ \widehat {[A]_n^{n'}} (\ell- \ell')  $, T\"oplitz in time.

\begin{definition} {\bf ($s$-decay norm)} \label{def:sdecay}
We define the $ s $-decay norm 
\be\label{s-decay}
| A |_s^2  := 
  {\mathop \sum}_{m \in \Z, L \in \Z^{|\ST|}}  
\Big(\sup_{n-n' = m} \| \widehat {[A]_n^{n'}} (L) \| \Big)^2 \langle m, L \rangle^{2s} 
\ee
where  $ \| \widehat {[A]_n^{n'}} (L) \| $ denotes the operator  
norm on $ {\rm Mat}_2(\C) $. 
For a family of operators $ \lambda \mapsto 
A(\lambda) $ which are $ k_0 $-times differentiable in $ \lambda = (\omega, \gamma ) \in \R^{|\ST|} \times [\g_1, \g_2 ]$, we set
$ |A|_s^{k_0,\upsilon} := {\mathop \sum}_{|k| \leq k_0} 
\upsilon^k | \pa_\lambda^k A |_s $. 
\end{definition}

In view of \eqref{Ann'} we clearly have 
$ | A |_s^{k_0,\upsilon} \leq | A \langle D \rangle^M |_s^{k_0,\upsilon} $.

The decay norm \eqref{s-decay} satisfies interpolation tame estimates, cfr. e.g. 
 \cite{BBM-Airy,BKM,BB20}:  
 for any $ s \geq  s_0 $,  
\be\label{tame-prod}
|A B |_s^{k_0,\upsilon} \leq C(s,k_0) |A |_s^{k_0,\upsilon} | B |_{s_0}^{k_0,\upsilon}  + 
C(s_0,k_0) |A |_{s_0}^{k_0,\upsilon} | B |_{s}^{k_0,\upsilon} \, . 
\ee
{\bf (Action)} For any $ s \geq s_0 $, $ q \in H^s $ 
\be\label{normapseudo-action}
\| A q \|_s^{k_0,\upsilon} \leq C(s,k_0) |A |_s^{k_0,\upsilon} \| q \|_{s_0}^{k_0,\upsilon} 
+  C(s,k_0) |A |_{s_0}^{k_0,\upsilon} \| q \|_{s}^{k_0,\upsilon} \, . 
\ee
{\bf (Exponential map)}
If $ |A \langle D \rangle^M |_{s_0}^{k_0,\upsilon} \leq \delta (s_0,k_0) $ is small enough,  
then, for any $s \geq s_0 $, 
\be \label{Neumann pseudo diff}
| e^A - {\rm Id}  |_s^{k_0,\upsilon} \lesssim_s |A |_s^{k_0,\upsilon}  \, , \quad 
| (e^A - {\rm Id}) \langle D \rangle^M |_s^{k_0,\upsilon} \lesssim_s 
|A \langle D \rangle^M |_s^{k_0,\upsilon} \, .
\ee
Given a linear operator $A(\vphi) $  we define the 
\emph{smoothed operator}
 ${\it \Pi_N} A$, $N\in\N $, with matrix entries 
		\begin{equation}\label{PiNA}
		\widehat{[{\it \Pi_N }A]_n^{n'}}(\ell-\ell') := \begin{cases}
		\widehat  {[A]_n^{n'}} (\ell-\ell') & \text{if } |\ell-\ell'| \leq N \\
		0 & \text{otherwise} \, . 
		\end{cases}
		\end{equation}
		We also denote ${\it \Pi^\perp_N} := \mathbb{I}_\bot -{\it \Pi_N }$.
For any $ \tb \geq 0 $, 
$ s \in \R $, 
\be\label{smo-pro}
| {\it \Pi_N^\bot} A |_s^{k_0,\upsilon} \leq N^{-\tb} | A |_{s+\tb}^{k_0,\upsilon} \, . 
\ee
The next lemma  embeds pseudo-differential operators into 
matrices with off-diagonal decay.

\begin{lemma}\label{lem:embed}
Let $ A = A(\vphi) = {\rm Op} (a (\vphi, \theta, \xi)) $ be a 
$ \vphi  $-dependent family of 
pseudo-differential operators in $  {\rm OPS}^{-M}  $, $ M \in \R $.
Then the decay norm of  the operator
$ {\mathtt A} := \Pi_{\ST,2}^\bot A_{|\acca_{\ST,2}^\bot} $ satisfies  
$ | {\mathtt A} \langle D \rangle^M |_s^{k_0,\upsilon} \lesssim 
| A  |_{-M,s+ s_0,0}^{k_0,\upsilon}   $, for any $ s \geq s_0 $,
where the norm $ | \  |_{-M,s+ s_0,0}^{k_0,\upsilon} $ is defined in  \eqref{norm1 parameter}. 
\end{lemma}

\begin{pf}
We represent  the operator 
$ {\mathtt A} \langle D \rangle^M$ as the matrix of $ 2 \times 2 $-matrices
$ [{\mathtt A} \langle D \rangle^M ]_n^{n'}  $ as in \eqref{Ann'}, whose 
elements are linear combinations
of the matrix elements 
$ (A \langle D \rangle^M)_{\pm n}^{\pm n'} (\ell) = 
(A)_{\pm n}^{\pm n'} (\ell) \langle n' \rangle^M $ 
with respect to the exponential basis, for any $ j, j' \in \Z $, 
\be\label{Aop-Aps}
\begin{aligned}
\widehat{A^{j'}_j} (\ell ) 
& := \int_{\T^{|\ST|}} e^{- \ii \ell \cdot \vphi} 
( A(\vphi) e^{\ii j' \theta}, e^{\ii j \theta} )_{L^2(\T)} d \vphi \\
&  =
\int_{\T^{{|\ST|}+1}} e^{- \ii \ell \cdot \vphi} 
a(\vphi, \theta, j') e^{\ii j' \theta} e^{- \ii j \theta}  d \vphi  d \theta 
 = \widehat a (\ell, j-j', j') \, . 
 \end{aligned}
\ee
We claim that
\be\label{stima-pre}
| \widehat{A^{j'}_j} (\ell )| \langle j' \rangle^M \lesssim  
\frac{|A|_{-M,s,0}}{\langle j - j', \ell \rangle^{s}} \, . 
\ee
Indeed, recalling  the Definition  \ref{def:pseudo-norm}, 
$ \sum_{\ell,J}  |\widehat a (\ell, J, \xi)|^2 \langle \ell,J \rangle^{2s}
= $ $ \| a (\cdot, \cdot, \xi) \|_s^2 \leq  |A|_{-M,s,0}^2 \langle \xi \rangle^{-2M} $, 
$ \forall \xi \in \R $, 
and therefore
$ |\widehat a (\ell, J, \xi)| \langle \xi \rangle^{M}
\leq  \frac{|A|_{-M,s,0}}{\langle \ell,J \rangle^{s} } $, $  \forall \ell, J, \xi $. 
As a consequence, recalling \eqref{Aop-Aps}, 
we deduce 
 \eqref{stima-pre}. Thus, for any $ n, n' \in \N $,  
\be\label{Ann'de}
\| \widehat{ [ {\mathtt A} \langle D \rangle^M ]_n^{n'}} (\ell) \| \lesssim
\frac{|A|_{-M,s,0}}{\langle n - n', \ell \rangle^{s}} +
\frac{|A|_{-M,s,0}}{\langle n + n', \ell \rangle^{s}}
\lesssim 
\frac{|A|_{-M,s,0}}{\langle n - n', \ell \rangle^{s}} \, .
\ee
Therefore, by  \eqref{s-decay} and \eqref{Ann'de},  
$ | {\mathtt A} \langle D \rangle^M |_s^2 \leq $ 
$ |A|_{-M,s+s_0,0}^2  $ for $ 2 s_0 > {|\ST|} + 1 $ and, for any $ |k| \leq k_0 $, we deduce 
$ | \pa_\l^k {\mathtt A} \langle D \rangle^M |_s \lesssim $ 
$ |\pa_\l^k A|_{-M,s+s_0,0} \lesssim  $ $ \upsilon^{-k} |A|_{-M,s+s_0,0}^{k_0,\upsilon} $ 
proving the lemma. 
\end{pf}

Finally 
we characterize 
the reality and reversibility properties of an operator by  its matrix entries.  
 
\begin{lemma}\label{lem:rrar}
An operator $ A(\vphi) \equiv (\widehat{[A]_n^{n'}}(\ell-\ell'))_{n,n' \in \ST^c, \ell, \ell' \in \Z^{|\ST|}} $ is {\sc real}, if and only if $ \overline{\widehat{[A]_n^{n'}}(L)} = \widehat{[A]_n^{n'}}(-L) $; 
{\sc reversible},  if and only if 
\be\label{Defs2}
\widehat{[A]_n^{n'}}(-L) {\cal S}_2
= - {\cal S}_2  \widehat{[A]_n^{n'}}(L) \, , \quad {\cal S}_2 := \begin{psmallmatrix}
1 & 0 \\
0 & - 1 
\end{psmallmatrix}; 
\ee
{\sc reversibility preserving},  if and only if 
$
\widehat{[A]_n^{n'}}(-L) {\cal S}_2
=  {\cal S}_2 \widehat{[A]_n^{n'}}(L) $. 
\end{lemma}

\begin{pf}
The reality condition is given in \eqref{Ann'-Fourier}. 
By Definition \ref{def:R-AR}, 
recalling \eqref{matricesAnn}  and  the form of  $ {\cal S}  $ in \eqref{Sanbn}
in the coordinates $ (\a_n, \b_n ) $,  the lemma follows. 
\end{pf}

\noindent
{\bf Preparation of the reducibility scheme.}
Using \eqref{defK0} we represent 
the quasi-periodic real reversible linear operator   $ \frak L_\bot $
in \eqref{inizioKAMRe}
as the infinite dimensional matrix (see \eqref{Ann'},  \eqref{matricesAnn} and  \eqref{anbn})
\begin{align}\label{defL0}
{\mathtt L} &= \om \cdot \pa_\vphi \mathbb{I}_\bot  + {\mathtt D} + {\mathtt R}\, ,  \\
\nonumber {\mathtt D} &:= {\rm diag}_{n \in \ST^c} {\mathtt D}_n \, , \quad 
{\mathtt D}_n  :=  \begin{pmatrix}
0 &  {\mathtt m}_{\bar \tn} n - \tfrac12 + \tfrac{\ka_n}{2}  \\
 - {\mathtt m}_{\bar \tn} n + \tfrac12 + \tfrac{\ka_n}{2} & 0 
\end{pmatrix} \, , \quad 
 {\mathtt R} := ( \widehat{[ {\mathfrak R}_\bot ]_n^{n'}}(\ell-\ell'))_{n, n' \in \ST^c} \, . 
\end{align}
By  Lemma \ref{lem:embed} and \eqref{fina-Re} we deduce that, 
for any $ M \in \N  $, there is $\s_M $ such that, for any $ s \in [s_0,S] $, 
\be\label{small-inizi}
| {\mathtt R} \langle D \rangle^M |_s^{k_0,\upsilon} \lesssim 
| {\mathfrak R}_\bot  |_{-M,s+ s_0,0}^{k_0,\upsilon} \lesssim_S 
\e  \upsilon^{-1} (1 + \| \fracchi_0 \|_{s + \sigma_M}^{k_0, \upsilon}) \, .
\ee
Recalling  that $  {\mathtt m}_{\bar \tn}  = \Omega_\gamma + 
{\mathtt r}_{\e, \bar \tn}  $ (see  \eqref{c1-picco}) and  
 the definition on $\mu_n^\pm $  in \eqref{anbn}, 
we also write (cfr. \eqref{anbn})  
$
{\mathtt D}_n = 
\begin{psmallmatrix}
0 & \mu_n^+ +  n {\mathtt r}_{\e, \bar \tn}   \\
- (\mu_n^{-} +  n {\mathtt r}_{\e, \bar \tn}  )  & 0 
\end{psmallmatrix} $.
We then conjugate the operator $ {\mathtt L} $ in \eqref{defL0} with the 
 bounded, real and reversibility preserving map
 (Lemma \ref{lem:rrar}),  see  \eqref{def:Lan},
 $ {\bf \symm}_\e := {\rm Diag}_{n \in \ST^c} [{\bf \symm}_\e]_{n}^n  $  where 
\be\label{def:Lanp}
[{\bf \symm}_\e]_{n}^n := 
\begin{cases}
\begin{psmallmatrix}
\symm_1(\e) & 0 \\
0 & - \symm_1^{-1}(\e)
\end{psmallmatrix} \, , \ \text{if} \ n = 1 \cr
\begin{psmallmatrix}
\symm_n (\e) & 0 \\
0 & \symm_n^{-1}(\e)  
\end{psmallmatrix} \, , \  \forall n \geq 3 \, , 
\end{cases} \qquad
\symm_n (\e) := \Big( \frac{|\mu_n^+ + n {\mathtt r}_{\e,\bar \tn}|}{|\mu_n^- +  
n {\mathtt r}_{\e, \bar \tn}|} \Big)^{\frac14} \, ,  
\ee
obtaining  the quasi-periodic real reversible linear operator 
\be\label{OpL1}
{\mathtt L}_0 := {\bf \symm}_\e^{-1} {\mathtt L} {\bf \symm}_\e = 
\om \cdot \pa_\vphi {\mathbb I}_\bot + {\mathtt D}_0 +  {\mathtt R}_0 \, , \quad
{\mathtt D}_0 := {\rm diag}_{n \in \ST^c} {\cal D}_n^{(0)} (\e) \, ,  
\ee
where, recalling \eqref{anbn-ri0}, 
\be\label{anbn-ri}
\begin{aligned}
& 
{\cal D}_n^{(0)} (\e) := 
\begin{cases}
& \begin{psmallmatrix}
0 & \Omega_n^{(0)} (\e)    \\
\Omega_n^{(0)} (\e) & 0 
\end{psmallmatrix} 
\, , \ \ 3 \leq n \leq \bar n \, ,  \\
& 
\begin{psmallmatrix}
0 & \Om_n^{(0)} (\e)  \\
- \Om_n^{(0)} (\e) & 0 
\end{psmallmatrix} 
\, , \ \forall  n \in \{ 1, \bar n + 1, \bar n + 2 \ldots \} \, ,  
\end{cases}
\end{aligned} 
\ee
with
\be\label{defomega0}
\Om_n^{(0)} (\e) := 
\big| (\mu_n^+ + n {\mathtt r}_{\e, \bar \tn}) 
( \mu_n^- +  n {\mathtt r}_{\e,  \bar \tn} ) \big|^{\frac12}  
\ee
and ${\mathtt R}_0 := {\bf \symm}_\e^{-1} {\mathtt R} {\bf \symm}_\e $. 
Since the decay norms of the map ${\bf \symm}_\e $ in \eqref{def:Lanp} satisfies
(recall also that the constant $ {\mathtt r}_{\e,\bar \tn} $ satisfies \eqref{c1-picco}),   
\be\label{MM-1decay} 
|{\bf \symm}_\e^{\pm 1}|_{s}^{k_0,\upsilon} \, , \quad
| \la  {\bf \symm}_\e^{\pm 1} \la D \ra^{M}|_{s}^{k_0,\upsilon} \leq C  \, ,
\ee
we deduce by \eqref{small-inizi} and \eqref{tame-prod} that, for any $ s \in [s_0,S] $, 
\be\label{small-inizi1}
| {\mathtt R}_0 \langle D \rangle^M |_s^{k_0,\upsilon} \lesssim_S 
\e \upsilon^{-1} (1 + \| \fracchi_0 \|_{s + \sigma_M}^{k_0, \upsilon}) \, . 
\ee

\noindent
{\bf KAM almost-diagonalization theorem.}
We now almost-diagonalize the operator $ {\mathtt L}_0 $ 
in \eqref{OpL1} along the scales
$ N_{-1} := 1 $, $ N_\tn = N_{\tn-1}^\chi $, $ \forall \tn \geq 1 $, $\chi = 3/2 $.
We fix  the constants
$ \ta := 3( \tau_1+1)  $, $  \tau_1 :=(  k_0 + 1)\tau_0+ k_0 $, $ \tb := [\ta] + 2 $,
as in \eqref{choiceab}.  
\begin{theorem}\label{iterative_KAM}
{\bf (Almost-diagonalization of $ {\mathtt L}_0 $)}
Let $ M \in \N $. 
There is  $ \tau_2 \geq \tau_1 +\ta$ such that, for any $ S > s_0 $, there is 
$N_0:=N_0(S, \tb,M)\in\N$ such that, if
	\begin{equation}\label{small_KAM_con}
		N_0^{\tau_2} 
		| {\mathtt R}_0 \langle D \rangle^M |_{s_0 + \tb}^{k_0,\upsilon}  \upsilon^{-1} 
		\leq 1 \,,
	\end{equation}
	then, for all $ \tn\in\N_0$, $ \nu =0,1,\ldots,\tn$:
	\\[1mm]
	$({\bf S1})_\nu $ There exists a real and reversible operator
\be\label{def:Ln}	
{\mathtt L}_\nu = \omega \cdot \pa_\vphi {\mathbb I}_\bot 
+ {\mathtt D}_\nu + {\mathtt R}_\nu \, , \quad 
{\mathtt D}_\nu := {\rm diag}_{n \in \ST^c} {\cal D}_n^{(\nu)} \, ,  
\ee
defined for all $ (\omega, \gamma) \in \R^{|\ST} \times [\g_1, \g_2] $, 
where ${\cal D}_n^{(\nu)} := {\cal D}_n^{(\nu)} (\e) $ are $ k_0 $-times differentiable $2 \times 2 $ real matrices of the form 
\begin{align}\label{anbn-rin}
& 
{\cal D}_n^{(\nu)}  := 
\begin{cases} 
\begin{psmallmatrix}
0 & \Omega_{n,h}^{(\nu)}     \\
\Omega_{n,h}^{(\nu)} & 0 
\end{psmallmatrix}, \ 3 \leq n \leq \bar n \, , 
\, \ \Om_{n,h}^{(\nu)}  = \Om_n^{(0)} (\e) 
+  {\frak r}_{n,h}^{(\nu)} \, ,  
  \\ 
\begin{psmallmatrix}
0 & \Om_{n,e}^{(\nu)}  \\
- \Om_{n,e}^{(\nu)}  & 0 
\end{psmallmatrix}, \, n \in \ST^c \setminus \{3, \ldots, \bar n \} \, ,  \ 
\Om_{n,e}^{(\nu)}  = \Om_n^{(0)} (\e) 
+  {\frak r}_{n,e}^{(\nu)}  \, , 
\end{cases}
\end{align} 
where  $ \Om_n^{(0)} (\e)$ are defined in  \eqref{defomega0}, 
$ {\frak r}_{n,h}^{(0)} = {\frak r}_{n,e}^{(0)} =0 $, 
 and, for any $ \nu \geq 1 $, 
\begin{align}
&
 \sup_{3 \leq n \leq \bar n}
 |  {\frak r}_{n,h}^{(\nu)}|^{k_0,\upsilon}  \leq C(S,\tb) \e \upsilon^{-1} \, , \quad
 \sup_{n \in \ST^c \setminus \{3, \ldots, \bar n \}}  
 |  {\frak r}_{n,e}^{(\nu)}|^{k_0,\upsilon} \leq C(S, \tb) 
n^{-M}  \e  \upsilon^{-1} \, , \label{eintot} \\
&\sup_{3 \leq n \leq \bar n}  |
{\frak r}_{n,h}^{(\nu)}   - {\frak r}_{n,h}^{(\nu-1)} |^{k_0,\upsilon} +
\sup_{\ST^c \setminus \{3, \ldots, \bar n \}} |{\frak r}_{n,e}^{(\nu)}   - {\frak r}_{n,e}^{(\nu-1)} |^{k_0,\upsilon}n^{M}  \leq C(s_0, \tb) 
 | {\mathtt R}_{\nu-1} \langle D \rangle^M |_{s_0}^{k_0,\upsilon}   \, . \label{einen}
 \end{align} 
The remainder $ {\mathtt R}_\nu $  satisfies, for any $ s \in [s_0,S] $,  
\begin{align}
 | {\mathtt R}_\nu \langle D \rangle^M |_{s}^{k_0,\upsilon} 
\leq 
 | {\mathtt R}_0 \langle D \rangle^M |_{s+\tb}^{k_0,\upsilon} 
N_{\nu-1}^{-\mathtt a} \, , \label{estRn:iterativa} \\
 | {\mathtt R}_\nu \langle D \rangle^M |_{s+\tb}^{k_0,\upsilon}  \leq 
| {\mathtt R}_0 \langle D \rangle^M |_{s+\tb}^{k_0,\upsilon}  N_{\nu-1}\,.
\label{estRn:iterativahigh}
\end{align}
Define the sets $\Lambda_0^\upsilon = \mathtt{DC}(2\upsilon,\tau)\times [\g_1, \g_2] $ and, for $ \nu \geq 1 $, 
\begin{align}
\Lambda_\nu^\upsilon :=  
\Big\{ & \lambda = (\omega, \gamma) \in   \Lambda_{\nu-1}^\upsilon  \ : 
\big| \omega \cdot \ell + \Omega_{n,e}^{(\nu-1)}  - \Omega_{n',e}^{(\nu-1)}   
\big| \geq { \upsilon \langle n - n' \rangle}{\langle \ell \rangle^{-\tau}} \, , 
\label{sec-Meln-p} \\
& \  
\forall (\ell, n, n') \neq (0, n, n) \, , \ n, n' \in \ST^c \setminus \{3, \ldots, \bar n \}  \, , \ 
\ell \in \Z^{|\ST|}  \, , \  | \ell | \leq N_{\nu-1}  \, , \nonumber   \\ 
& 
\big| \omega \cdot \ell + \Omega_{n,e}^{(\nu-1)}  + \Omega_{n',e}^{(\nu-1)}  
\big| \geq { \upsilon (n+n')}{\langle \ell \rangle^{-\tau}} \, ,  \label{sec-Meln-p-22} \\
& \ 
 \forall (\ell, n, n') \, , \ n, n' \in \ST^c \setminus \{3, \ldots, \bar n \} \, , \ \ell \in \Z^{|\ST|}  \, , \ | \ell | \leq N_{\nu-1}   \nonumber 
\Big\}  \, . 
 \end{align} 
For  $ \nu \geq 1 $ there exists a real and reversibility preserving operator 
$\Phi_\nu $ acting in 
$ {\frak H}_{\ST}^\bot $, defined for all
$(\omega, \gamma) \in \R^{|\ST|} \times [\g_1, \g_2 ] $, of the form
$ \Phi_\nu = e^{\Psi_{\nu-1}} $ satisfying, for all $ s \in [s_0,S] $, 
\be\label{smallnessPsin-1}
|\Psi_{\nu-1} \langle D \rangle^M |_s^{k_0,\upsilon} \leq C(S,\tb) 
\upsilon^{-1}  N_{\nu-1}^{\tau_1}N_{\nu-2}^{-\ta}
  | {\mathtt R}_0\langle D \rangle^M |_{s+\tb}^{k_0,\upsilon}   \, , 
\ee
such that, for all $ \l \in \Lambda_\nu^\upsilon $
the  conjugation formula 
$ {\mathtt L}_\nu = \Phi_{\nu-1}^{-1} {\mathtt L}_{\nu-1} \Phi_{\nu-1} $ holds.
 \item[${\bf(S2)_{\nu}}$] Let $ i_1 $, $ i_2 $ be such that 
${\mathtt R}_0(i_1)$,  ${\mathtt R}_0(i_2 )$ satisfy \eqref{small_KAM_con}. 
Then  for all $(\omega, \g) \in \R^{|\ST|} \times [\g_1, \g_2] $, 
$\nu \geq 1$, 
\begin{align}
&  \sup_{n \in \ST^c \setminus \{3, \ldots, \bar n \}} \ |\Delta_{12} {\frak r}_{n,e}^{\nu}|  n^{ M}  \lesssim_{S, \mathtt b} \e \upsilon^{- 1}\| i_1 - i_2  \|_{ s_0  + \perd(\mathtt b)}\, ,  \label{r nu i1 - r nu i2} \\
& \label{r nu - 1 r nu i1 i2}
\sup_{n \in \ST^c \setminus \{3, \ldots, \bar n \}}\big| \Delta_{12}( {\frak r}_{n,e}^\nu - {\frak r}_{n,e}^{\nu - 1}) \big| n^M \lesssim_{S, \mathtt b} \e 
\upsilon^{- 1}    N_{\nu - 2}^{- \mathtt a}  \| i_1 - i_2  \|_{ s_0  + \perd(\mathtt b)}\, .
\end{align}
\item[${\bf(S3)_{\nu}}$]  Let $i_1$, $i_2$ be like in ${\bf(S2)_{\nu}}$ 
and $0 < \rho \leq \upsilon/2$. 
Then 
\begin{equation}\label{inclusione cantor riducibilita S4 nu} 
\varepsilon \upsilon^{- 1 }  C(S) N_{\nu-1}^{ \tau + 1 } 
\|i_2 - i_1 \|_{s_0 + \perd(\mathtt b)} 
\leq \rho \quad \Longrightarrow \quad 
\tLm_\nu^\upsilon(i_1) \subseteq \tLm_\nu^{\upsilon - \rho}(i_2) \, . 
\end{equation}
\end{theorem}

Theorem \ref{iterative_KAM} implies that the invertible operator
$ U_\tn := \Phi_0 \circ \ldots \circ  \Phi_{\tn-1} $, $  \tn \geq 1 $, 
has almost diagonalized $ {\mathtt L}_0 $ for any $(\om,\gamma) 
  \in \cap_{\nu = 0}^{\tn } \Lambda_\nu^\upsilon = \Lambda_\tn^\upsilon $, meaning that \eqref{Linf} below holds.
Arguing as in Corollary 4.1 in \cite{BBM-Airy} we deduce the following result. 

\begin{theorem} {\bf (Almost-reducibility of $ {\mathtt L}_0  $)}\label{thm:riduKAM}
For all $  S > s_0 $ there is $ N_0 (S, \tb), \d_0 (S,\tb ) > 0 $ such that, if the 
 the smallness condition
\be\label{KAM-small}
N_0^{\tau_2}\e \upsilon^{-2} \leq \delta_0
\ee
hold, where $\tau_2$ is defined in Theorem \ref{iterative_KAM}, then, for any $ \tn \in \N $, 
\be \label{Unvicina}
|U_\tn - {\rm Id} |_s^{k_0,\upsilon} + |U_\tn^{-1} - {\rm Id} |_s^{k_0,\upsilon}
\leq \varepsilon\upsilon^{-2} C(S) | {\mathtt R}_0 \langle D \rangle^M |_{s+\tb}^{k_0, \upsilon} \, ,  
\ee
and, for any $ (\omega, \g) \in   \Lambda_{\tn}^\upsilon $ 
defined in \eqref{sec-Meln-p}-\eqref{sec-Meln-p-22} (with $ \nu = \tn $), 
we have
\be\label{Linf}
U_{ \tn}^{-1}  {\mathtt L}_0 U_{\tn} = {\mathtt L}_{ \tn} =  
\omega \cdot \pa_\vphi  {\mathbb I}_\bot  + {\mathtt D}_{ \tn} + {\mathtt R}_{ \tn}
\ee
where ${\mathtt D}_{ \tn} $ has the form \eqref{def:Ln}-\eqref{anbn-rin}.
\end{theorem}

\noindent
{\bf Proof of almost-diagonalization Theorem \ref{iterative_KAM}.} 
The proof is inductive.  
The key of the convergence  are 
the inductive relations
\eqref{lownorm1}-\eqref{highnorm1}.
In order to prove \eqref{lownorm1}-\eqref{highnorm1}  we use  that the class of 
matrices with finite $ s $-decay norm satisfies
\eqref{tame-prod}, \eqref{smo-pro} and it is closed
for the solution of the  homological equation obtained in Lemma
\ref{lemma:homological}, see \eqref{psi-v-r}. We will prove in detail only $ ({\bf S1})_\nu $.
\\[1mm]
{\bf Initialization.}
{\sc Proof of} $({\bf S1})_0 $. The operator $ {\mathtt L}_0  $
in \eqref{OpL1}-\eqref{anbn-ri} has the form \eqref{def:Ln}-\eqref{anbn-rin}
with $ {\frak r}_{n,h}^{(0)} = {\frak r}_{n,e}^{(0)} =0 $.
\\[1mm]
{\bf The reducibility step.}
We now describe the generic inductive step, showing how to define 
 $ {\mathtt L}_{\nu+1}  $, $ \Phi_\nu $, $ \Psi_\nu $.  
To simplify notation in this section we drop the index $ \nu $ and we write
$ + $ for $ \nu + 1 $. 

We conjugate $ {\mathtt L} $ in \eqref{def:Ln}-\eqref{anbn-rin} by the flow map $ \Phi = e^{\Psi}  $ where
\be\label{decoAnn'}
\Psi = (\Psi_n^{n'}(\vphi))_{n,n' \in \ST^c} \, , \quad \Psi_n^{n'} (\vphi) = 
{\mathop \sum}_{\ell \in \Z^{|\ST|}, |\ell| \leq N} \widehat {\Psi_n^{n'}} (\ell) e^{\ii \ell \cdot \vphi} 
\in {\rm Mat}_2(\R) \, . 
\ee
By  \eqref{def:Ln} and a Lie expansion (e.g. formulas (3.17)-(3.18) in \cite{BFM})
\be\label{coniugn+1}
\begin{aligned}
e^{-\Psi} {\mathtt L} e^{\Psi}  & = 
\omega \cdot \pa_\vphi + {\mathtt D} 
+ (\om \cdot \pa_\vphi \Psi) + [{\mathtt D}, \Psi] + {\it \Pi_N } {\mathtt R} +  
{\it \Pi^\bot_N} {\mathtt R} \\
& 
- \int_0^1 e^{- \tau \Psi} [\Psi,{\mathtt R}] e^{ \tau \Psi} d \tau 
- \int_0^1 (1- \tau) e^{- \tau \Psi} 
[\Psi,  (\om \cdot \pa_\vphi \Psi) + [{\mathtt D},\Psi ] ] e^{ \tau \Psi} d \tau  
\end{aligned}
\ee 
where ${\it \Pi_N }$ is defined in \eqref{PiNA} and $ {\it \Pi_N^\bot } = \mathbb{I}_\bot  - 
{\it \Pi_N } $. 
\\[1mm]{\bf Homological equation.}
We look for a solution $ \Psi $ of  the homological equation
\be\label{homo1}
\om \cdot \pa_\vphi \Psi + [{\mathtt D}, \Psi] +  {\it \Pi_N } {\mathtt R} =  \bl{\mathtt R} \br
\ee
where $ \bl {\mathtt R} \br  =  (\bl {\mathtt R} \br^{n'}_n(\ell))_{n,n' \in \ST^c, \ell \in \Z^{|\ST|}} $ is the $ \vphi $-independent 
normal form  
\be\label{forma-nor}
\bl {\mathtt R} \br^{n'}_n (\ell)
:= 
\begin{cases}
r_{n,h} \begin{psmallmatrix}
0 & 1  \\
1 & 0 
\end{psmallmatrix} \, , \ n = n'  \, , \ n \in \{3, \ldots, \bar n \} \, , \ \ell = 0 \\
r_{n,e} \begin{psmallmatrix}
0 & 1  \\
-1 & 0 
\end{psmallmatrix} \, , \ n = n'  \, , \ n \notin \{3, \ldots, \bar n \} \, ,  \ \ell = 0 \\
0 \qquad \qquad \qquad \text{otherwise} \, , \\
\end{cases}
\ee 
where 
\be\label{Rev-NF-gene}
r_{n,h} := \tfrac12 (b_n+c_n) \, , \ 
r_{n,e} := \tfrac12 (b_n-c_n) \, , \  
\widehat{[{\mathtt R}]_n^n} (0) =  \begin{psmallmatrix}
0 & b_n  \\
c_n & 0 
\end{psmallmatrix} \, , \quad b_n, c_n \in \R \, . 
\ee  
Note that 
$ \widehat{[{\mathtt R}]_n^n}(0) $  has the form \eqref{Rev-NF-gene}, by Lemma \ref{lem:rrar}.

\begin{lemma}\label{lemma:homological}
{\bf (Homological equation)}
There exists a real and reversibility preserving 
linear operator $ \Psi \equiv (\Psi_n^{n'} (\ell-\ell'))_{n,n' \in \ST^c, \ell, \ell' \in \Z^{|\ST|}}$
as in \eqref{decoAnn'}, 
defined for any $ (\omega, \gamma) \in \R^{|\ST|} \times [\gamma_1, \gamma_2 ]$, 
which is a solution of the homological equation \eqref{homo1}  
for any $ \Lambda_{\nu + 1}^{\upsilon} $ 
(see \eqref{sec-Meln-p}-\eqref{sec-Meln-p-22}), 
satisfying, for any $ s \in [s_0,S] $,  $ \tau_1 := \tau (k_0+1) + k_0  $,  
\be\label{psi-v-r}
 | \Psi |_s^{k_0,\upsilon} \leq C N^{\tau_1} \upsilon^{-1} | {\mathtt R} |_s^{k_0,\upsilon} \, ,
 \quad
| \Psi \langle D \rangle^M |_s^{k_0,\upsilon} \leq C N^{\tau_1} \upsilon^{-1} 
| {\mathtt R} \langle D \rangle^M |_s^{k_0,\upsilon} \,.
\ee
\end{lemma}

The rest of this paragraph is devoted to the proof  of Lemma \ref{lemma:homological}.

Recalling 
\eqref{def:Ln}-\eqref{anbn-rin} and 
\eqref{decoAnn'}, 
the homological equation \eqref{homo1} 
reduces to the set of equations, for any 
$  \ell \in \Z^{|\ST|} $, $ |\ell| \leq N $,  $ n, n' \in \ST^c \setminus \{3, \ldots, \bar n \} $,  
\be\label{homo2-deve}
\ii \omega \cdot \ell \, \widehat{\Psi_n^{n'}} (\ell) + {\cal D}_n 
\widehat{\Psi_n^{n'}} (\ell) -  
\widehat{\Psi_n^{n'}} (\ell)  {\cal D}_{n'} +
 \widehat{[{\mathtt R}]_n^{n'}}(\ell) =  \bl{\mathtt R} \br_n^{n'} (\ell) \, .   
\ee
In order to solve the linear equations  \eqref{homo2-deve}, 
we have to study the eigenvalues of the map 
$$
T_{\ell,n,n'} : {\rm Mat}_2(\C) \to {\rm Mat}_2(\C) \, , \quad 
T_{\ell,n,n'} (B)  :=  \ii \omega \cdot \ell B + {\cal D}_n  B -
B  {\cal D}_{n'} \, . 
$$
The spectrum of 
$ T_{\ell,n,n'} $ is the sum of  $ \ii \omega \cdot \ell $ plus  the spectrum of the linear map 
\be\label{Dnn'B}
{\cal C}_{n,n'} : {\rm Mat}_2(\C) \to {\rm Mat}_2(\C)   \, , \quad 
{\cal C}_{n,n'} (B) :=   {\cal D}_n B - B  {\cal D}_{n'} \, . 
\ee
All the eigenvalues of $ {\cal C}_{n,n'} (B) $
are the differences between the eigenvalues 
of ${\cal D}_n $ and 
those of  $ {\cal D}_{n'} $.
Due to the different form of  ${\cal D}_n $ in 
\eqref{anbn-rin} we distinguish the following cases. 
\\[1mm]
{\bf Case Hyperbolic/Elliptic:}
$ n, n' \in \ST^c  $, $ n \in \{ 3, \ldots, \bar n \} $, $ n' \notin \{ 3, \ldots, \bar n \} $, or viceversa.
In  the case $ n \in \{ 3, \ldots, \bar n \} $, $ n' \notin \{ 3, \ldots, \bar n \} $,  
by  \eqref{anbn-rin},
all the eigenvalues of 
$ T_{\ell,n,n'}$ are 
$ \ii (\omega \cdot \ell \pm \Omega_{n',e} ) \pm \Omega_{n,h}  $. 
Thus they are all different from zero with modulus larger than 
\be\label{lobou}
|\omega \cdot  \ell \pm \Omega_{n',e} | + | \Omega_{n,h}  | \geq 
\max \{ |\omega \cdot  \ell \pm \Omega_{n',e} |, | \Omega_{n,h}  |  \}
\geq c > 0  \, .
\ee
Thus 
$ T_{\ell,n,n'} $ is invertible for all the parameters $ (\omega, \gamma) \in \R^{|\ST|} \times [\gamma_1, \gamma_2] $ and
\be\label{defPsinn'}
\widehat{\Psi_n^{n'}} (\ell) := -  T_{\ell,n,n'}^{-1} \widehat{[{\mathtt R}]_n^{n'}} (\ell)  \, , \
\forall ( \ell, n, n' ) \in \{|\ell | \leq N \}\times \{ 3, \ldots, \bar n \} \times 
(\ST^c \setminus \{ 3, \ldots, \bar n \})  \, , 
\ee
is  the unique solution  of the homological equation \eqref{homo2-deve} (note that 
by \eqref{forma-nor} we have $ \bl {\mathtt R} \br_n^{n'} ( \ell )= 0 $).  
Since $ \overline{\widehat{[{\mathtt R}]_n^{n'}} (\ell)} = \widehat{[{\mathtt R}]_n^{n'}} (- \ell) $ (reality condition), taking the complex conjugate in \eqref{homo2-deve} and by uniqueness we deduce that
$ \widehat{\Psi_n^{n'}} (\ell) $ satisfies the reality condition 
as well. Moreover, multiplying \eqref{homo2-deve} from the left and the right by
$ {\cal S}_2$ we get
$ \ii \omega \cdot \ell {\cal S}_2 \widehat{\Psi_n^{n'}} (\ell) {\cal S}_2 - 
 {\cal S}_2 {\cal D}_n  
\widehat{\Psi_n^{n'}} (\ell)  {\cal S}_2 +  
 {\cal S}_2 \widehat{\Psi_n^{n'}} (\ell)  {\cal D}_{n'}  {\cal S}_2 +
 {\cal S}_2  \widehat{[{\mathtt R}]_n^{n'}}(\ell)  {\cal S}_2 =  0 $. 
Since ${\cal D}_{n}, {\cal D}_{n'} $ anti-commute with ${\cal S}_2 $ 
defined in \eqref{Defs2}
and  $ [{\mathtt R}]_n^{n'} (\ell)$ satisfies the reversibility condition in Lemma \ref{lem:rrar} we deduce that 
$ {\cal S}_2 \widehat{\Psi_n^{n'}} (\ell) {\cal S}_2 =
\widehat{\Psi_n^{n'}}(-\ell)   $, namely 
$ \widehat{\Psi_n^{n'}} (\ell)   $ satisfies the 
anti-reversibility condition. 
\\[1mm]
{\bf Case Hyperbolic/Hyperbolic:  
$ n, n' \in \ST^c  $, $ n, n' \in \{ 3, \ldots, \bar n \} $.} 
Let us consider the basis of $ {\rm Mat}_2(\C) $ defined by  
(it is constructed by the eigenvectors of $ {\cal D}_n $ and $ {\cal D}_{n'} $, see e.g  \cite{BKM}, Lemma 7.3)
\be \label{F1F4}
F_{1} := 
 \begin{psmallmatrix}
1 & 1  \\
1 & 1 
\end{psmallmatrix} \, , \
F_2 := \begin{psmallmatrix}
1 & -1 \\
- 1 &  1 
\end{psmallmatrix} \, , \ 
F_{3} := 
 \begin{psmallmatrix}
1 & - 1  \\
1 & - 1 
\end{psmallmatrix} \, , \
F_4 := \begin{psmallmatrix}
1 & 1 \\
- 1 & - 1 
\end{psmallmatrix} 
 \, .
\ee 
In the basis $ \{ F_1, F_2, F_3, F_4 \} $ 
 any matrix in $ {\rm Mat}_2 (\C)$ decomposes as 
\be\label{proiez1}
\begin{psmallmatrix}
\alpha & \beta \\
\gamma & \delta 
\end{psmallmatrix} =
x_1  F_1 + x_2 F_2 + x_3 F_3 + x_4 F_4  
\ee
with
\be\label{proiez2}
\begin{aligned}
&  x_1 := \tfrac14 (\a+\b+\g+\d) \, , \ 
x_2 := \tfrac14 (\a-\b-\g+\d) \, , \\
& x_3 := \tfrac14 (\a-\b+\g-\d) \, , \
x_4 := \tfrac14 (\a+\b-\g-\d) \, .
\end{aligned}
\ee
Recalling \eqref{anbn-rin},  
the operator $ {\cal C}_{n,n'} $ in \eqref{Dnn'B} is represented, in the basis
\eqref{F1F4}, by the diagonal matrix
\be\label{nF1nF3} \small
\begin{psmallmatrix} 
\Om_{n} - \Om_{n'} & 0 & 0 & 0  \\
0 & -  (\Om_{n} - \Om_{n'})  & 0 & 0  \\
0 & 0 & \Om_{n} + \Om_{n'}  & 0  \\
0 & 0 & 0 & - (\Om_{n} + \Om_{n'})   
\end{psmallmatrix} \, , \  \Om_{n} \equiv \Om_{n,h}  \, . 
\ee
Therefore all the eigenvalues of  $ T_{\ell,n,n'}$ are 
$ \ii \omega \cdot \ell \pm \Omega_{n',h}  \pm \Omega_{n,h } $ 
and, since the $ \Omega_{n,h}  $
are all simple, their modulus is  bounded by 
\be\label{lbhh}
|\omega \cdot  \ell | + | \Omega_{n,h}  - \Omega_{n',h} |  \geq c > 0 \, , 
\quad \forall (\ell,n,n') \neq (0,n,n) \, . 
\ee
Thus, for any $ (\ell,n, n') \neq (0,n,n)$  the operator 
$ T_{\ell,n,n'} $ is invertible for all the parameters $ (\omega, \gamma) \in 
\R^{|\ST|} \times [\gamma_1, \gamma_2] $ and 
\be\label{defPsinn'h}
\widehat{\Psi_n^{n'}} (\ell) := -  T_{\ell,n,n'}^{-1} [{\mathtt R}]_n^{n'} (\ell)  \, , 
\quad  (\ell,n, n') \neq (0,n,n) \, ,  |\ell| \leq N \, , \ n, n' \in \{ 3, \ldots, \bar n \} \, , 
\ee
is the unique solution 
of the homological equation 
\eqref{homo2-deve} (note that 
by \eqref{forma-nor} we have $ \bl {\mathtt R} \br_n^{n'} ( \ell )= 0 $).  
The reality and the anti-reversibility condition follow as well. 
Next we consider the case $ (\ell,n,n') = (0,n,n)$. 
By \eqref{nF1nF3} the operator $ T_{0,n,n} $ has range the span 
$ \la F_3, F_4 \ra  $. 
By \eqref{Rev-NF-gene}, \eqref{proiez1},
\eqref{proiez2},  the projection of 
$  \widehat{[ {\mathtt R}]_n^n}(0) = \begin{psmallmatrix}
0 & b_n \\
c_n & 0 
\end{psmallmatrix}
$ on $ \la F_1, F_2 \ra $ 
 is 
$ x_1  F_1 + x_2 F_2 = \tfrac12 (b_n+c_n) 
\begin{psmallmatrix}
0 & 1 \\
1 & 0  
\end{psmallmatrix} $ with 
$ x_1 = \tfrac14 (b_n+c_n) = - x_2  $. 
By \eqref{forma-nor}-\eqref{Rev-NF-gene}, 
the homological equation \eqref{homo2-deve} reduces to  
$  {\cal C}_{n,n} \big[ \widehat { \Psi_n^{n} } (0)\big] =  - (x_3 F_3 + x_4 F_4) $
where $  x_3 = - \tfrac14 (b_n-c_n) = - x_4 $, 
and, by \eqref{nF1nF3}, we define, for any $ n = 3, \ldots, \bar n $, its  
 solution 
\be\label{psinnh0}
\widehat{\Psi_n^{n}} (0) 
:= -\frac{x_3}{2 \Omega_n} F_3 + \frac{x_4}{2 \Omega_n} F_4 
=  \frac{b_n - c_n}{4 \Omega_n} 
\begin{pmatrix}
1 & 0 \\
0 & - 1 
\end{pmatrix} \, , \quad \forall n \in \{3, \ldots,  \bar n \} 
\, . 
\ee
The real matrix
$ \widehat{\Psi_n^{n}}(0)  $ commutes with $ {\cal S}_2 $
and thus 
is reversibility preserving 
by  Lemma \ref{lem:rrar}. 
\\[1mm]
{\bf Case Elliptic/Elliptic:  
$ n, n' \in \ST^c  $, $ n, n' \notin \{ 3, \ldots, \bar n \} $.} 
Consider the basis of $ {\rm Mat}_2(\C) $
(constructed by the eigenvectors of $ {\cal D}_n $ and $ {\cal D}_{n'} $, see e.g  \cite{BKM}, Lemma 7.3)   
\be \label{M1M4}
M_{1} := 
 \begin{psmallmatrix}
\ii & 1  \\
-1 & \ii 
\end{psmallmatrix} \, , \
M_{2} := 
 \begin{psmallmatrix}
-\ii & 1  \\
-1 & - \ii 
\end{psmallmatrix} \, , \
M_3 := \begin{psmallmatrix}
1 & \ii  \\
\ii & - 1 
\end{psmallmatrix} \, , \
M_4 := \begin{psmallmatrix}
1 & - \ii  \\
- \ii & - 1 
\end{psmallmatrix} \, .
\ee 
In the basis $ \{ M_1, M_2, M_3, M_4 \} $ 
 any matrix in $ {\rm Mat}_2 (\C)$ decomposes as 
\be\label{proiezM1}
\begin{psmallmatrix}
\alpha & \beta \\
\gamma & \delta 
\end{psmallmatrix} =
x_1  M_1 + x_2 M_2 + x_3 M_3 + x_4 M_4  
\ee
with
\be\label{proiezM2}
\begin{aligned}
&  x_1 := \tfrac14 ( - \ii \a+\b-\g- \ii \d) \, , \ 
x_2 := \tfrac14 (\ii \a+\b-\g+ \ii \d) \, , \\
& x_3 := \tfrac14 (\a- \ii \b- \ii \g-\d) \, , \
x_4 := \tfrac14 (\a+\ii \b+ \ii \g-\d) \, .
\end{aligned}
\ee
Recalling \eqref{anbn-rin} the operator $ {\cal C}_{n,n'} $ in \eqref{Dnn'B} is represented, in the basis
\eqref{M1M4}, by the diagonal matrix
\be\label{n1n3} \small
\begin{psmallmatrix}
\ii (\Om_{n} - \Om_{n'}) & 0 & 0 & 0  \\
& - \ii (\Om_{n} - \Om_{n'})  & 0 & 0  \\
0 & 0 & \ii(\Om_{n} + \Om_{n'})  & 0  \\
0 & 0 & 0 & - \ii(\Om_{n} + \Om_{n'})   
\end{psmallmatrix} \, , \quad \Om_{n} \equiv \Om_{n,e}  \, . 
\ee
Therefore  the eigenvalues of $ T_{\ell,n, n'} $ 
are
$ \ii ( \omega \cdot \ell \pm \Om_{n}  \pm \Om_{n'}  )  $.
We first consider the case $ (\ell, n, n') \neq (0,n,n) $.  We decompose,
for any $ (\ell, n, n') \neq (0,n,n)$, $ |\ell| \leq N $, 
 $ n, n' \in \ST^c \setminus \{3, \ldots, \bar n \} $, 
$ [{\mathtt R}]_n^{n'} (\ell) = x_1 M_1 + x_2 M_2 + x_3 M_3 + x_4 M_4  $,   
 and 
we define 
\be\label{extellell}
\widehat{\Psi_n^{n'}} (\ell) =  y_1 M_1 + y_2 M_2 + y_3 M_3 +   y_4 M_4 \, ,
\ee
where, denoting $ \chi $  the cut-off function defined in \eqref{cutoff},
\be\label{defest1}
\begin{aligned}
& y_1  :=
- \tfrac{\chi ((\omega \cdot \ell + \Om_{n}  -  \Om_{n'} ) \upsilon^{-1} \langle \ell \rangle^\tau)}{\ii (\omega \cdot \ell + \Om_{n}  - \Om_{n'} )  } x_1 \, ,  \quad 
y_2  := - 
\tfrac{\chi ((\omega \cdot \ell - \Om_{n}  + \Om_{n'} ) \upsilon^{-1} \langle \ell \rangle^\tau)}{\ii (\omega \cdot \ell - \Om_{n}  + \Om_{n'} )  } x_2 \, ,  \\
& y_3  := - 
\tfrac{\chi ((\omega \cdot \ell + \Om_{n}  + \Om_{n'}) \upsilon^{-1} \langle \ell \rangle^\tau)}{\ii (\omega \cdot \ell + \Om_{n}  + \Om_{n'})  } x_3 \, ,  \quad 
y_4  := - 
\tfrac{\chi ((\omega \cdot \ell - \Om_{n}  - \Om_{n'} ) \upsilon^{-1} \langle \ell \rangle^\tau)}{\ii (\omega \cdot \ell - \Om_{n} - \Om_{n'})  } x_4 \, , 
\end{aligned}
\ee
in such a way that
$ \widehat{\Psi_n^{n'}} (\ell) $  
solves the homological equation 
\eqref{homo2-deve} for any  $ \l  \in \Lambda_{\nu+1}^\upsilon $, see 
\eqref{sec-Meln-p}-\eqref{sec-Meln-p-22}.

Next we consider the case $ (\ell,n,n') = (0,n,n) $. 
By \eqref{n1n3} the operator $ T_{0,n,n} $ has range $ \la M_3, M_4 \ra  $. 
Then by \eqref{Rev-NF-gene}, \eqref{proiezM1},
\eqref{proiezM2}, the projection of
$ \widehat{[{\mathtt R}]_n^n} (0) = \begin{psmallmatrix}
0 & b_n \\
c_n & 0 
\end{psmallmatrix} $ on $ \la M_1, M_2 \ra $ 
is 
$
x_1  M_1 + x_2 M_2  = \tfrac{b_n-c_n}{2} 
\begin{psmallmatrix}
0 & 1 \\
-1 & 0 
\end{psmallmatrix}$ where
$ x_1 = \tfrac14 (b_n-c_n) = x_2 $. 
Thus, by \eqref{forma-nor}-\eqref{Rev-NF-gene}, 
the homological equation \eqref{homo2-deve} reduces to
$  {\cal C}_{n,n} \big[ \widehat { \Psi_n^{n} } (0)\big] =  - \big( x_3 M_3 + x_4 M_4\big) $, 
$ x_3 = - \ii \tfrac14 (b_n + c_n) = - x_4 $, 
and, by \eqref{n1n3} we define the 
solution
\be\label{solpsi}
\widehat{\Psi_n^{n}} (0) := 
- \frac{x_3}{2 \ii \Omega_n} M_3 + \frac{x_4}{2 \ii \Omega_n} M_4 
=  \frac{b_n + c_n}{4 \Omega_n} \begin{pmatrix}
1 & 0 \\
0 & -1 
\end{pmatrix} \, , \quad 
\forall n \in \ST^c \setminus \{ 3, \ldots, \bar n \} 
 \, . 
\ee 
The real matrix
$ \widehat{\Psi_n^{n}} (0) $ 
 commutes with $ {\cal S}_2 $, 
and thus satisfies the condition of being reversibility preserving as required by Lemma \ref{lem:rrar}. 
The fact that map $ \Psi (\vphi) \equiv \big( \Psi_n^{n'} (\ell))_{\ell \in \Z^{|\ST|}, |\ell| \leq N, n,n' \in \ST^c} $ 
defined in \eqref{defPsinn'}, \eqref{defPsinn'h}, \eqref{psinnh0},
\eqref{extellell}, \eqref{defest1}, \eqref{solpsi} satisfies 
the estimates   \eqref{psi-v-r} 
follows by the lower bounds of the 
eigenvalues in  \eqref{lobou}, \eqref{lbhh}, \eqref{sec-Meln-p}-\eqref{sec-Meln-p-22}
and the properties of the cut-off function $ \chi $ in \eqref{cutoff}.
Lemma \ref{lemma:homological} is proved. \qed
\\[1mm]
{\bf Conclusion of the iterative step.}
By \eqref{coniugn+1}, \eqref{homo1}
  and Lemma \ref{lemma:homological}, for any 
$ (\omega, \gamma) $ in $ \Lambda_{\nu + 1}^{\upsilon} $, we have 
\begin{align}
\label{Ln+1KAM}
&{\mathtt L}_+ = \Phi^{-1} {\mathtt L} \Phi =
\omega \cdot \pa_\vphi  {\mathbb I}_\bot  + {\mathtt D}_{+} + {\mathtt R}_+,\, \\
\label{Rn+1KAMd}
& {\mathtt D}_{+} := {\mathtt D} +  \bl {\mathtt R} \br  \\
& {\mathtt R}_+ := 
{\it \Pi_N^\bot } {\mathtt R}  
- \int_0^1 e^{- \tau \Psi} [\Psi,{\mathtt R}] e^{ \tau \Psi} d \tau 
- \int_0^1 (1- \tau) e^{- \tau \Psi} 
\big[ \Psi,  \bl{\mathtt R} \br - {\it \Pi_N} {\mathtt R}  \big] e^{ \tau \Psi} d \tau \, . \label{Rn+1KAM}
\end{align}
The right hand side in \eqref{Rn+1KAM} defines an extension of
$ {\mathtt L}_+ $ to the whole parameter space $ \R^{|\ST|} \times [\gamma_1, \gamma_2 ] $
since $ {\mathtt R} $ and $ \Psi $ are defined on $ \R^{|\ST|} \times [\gamma_1, \gamma_2 ] $.  
The operator $ {\mathtt L}_+ $ in \eqref{Ln+1KAM} has the same form of $ \mathtt L $
in \eqref{def:Ln}-\eqref{anbn-rin}.
The following lemma follows by 
\eqref{Rev-NF-gene} and Definition \ref{def:sdecay}. 

\begin{lemma}\label{lem:NNF}
{\bf (New normal form)}
For all $ (\omega, \gamma) \in \R^{|\ST|} \times [\gamma_1, \gamma_2] $
the operator  ${\mathtt D}_{+} $ in \eqref{Rn+1KAMd}
has the same form of   ${\mathtt D} $ in \eqref{def:Ln}-\eqref{anbn-rin} 
 with $ \Om_{n,h}^{(\nu)} $ and $  \Om_{n,e}^{(\nu)}   $ replaced by 
$ \Om_{n,h}^{(\nu+1)}  := \Om_{n,h}^{(\nu)}  + r_{n,h}^{(\nu)} 
:= \Om_n^{(0)} (\e)  +  {\frak r}_{n,h}^{(\nu+1)}  \in \R $
and  
$ \Om_{n,e}^{(\nu+1)}  = \Om_{n,e}^{(\nu)}  + r_{n,e}^{(\nu)}  
= \Om_n^{(0)} (\e)  +  {\frak r}_{n,e}^{(\nu+1)} \in \R $, and 
$\sup_{3 \leq n \leq \bar n} |r_{n,h}^{(\nu)}|^{k_0, \upsilon} + 
\sup_{n \in \ST^c \setminus \{3, \ldots, \bar n \}} 
|r_{n,e}^{(\nu)}|^{k_0, \upsilon} n^{M}   \lesssim 
 |{\mathtt R}_\nu \la D\ra^M |^{k_0,\upsilon}_{s_0} \, .$ 
\end{lemma}

\noindent{\bf The iteration.}
Suppose $ ({\bf S1})_\nu $ is true. Now we prove $({\bf S1})_{\nu+1}$.

\begin{lemma}
The remainder  $ {\mathtt R}_{\nu+1}  $ defined in \eqref{Rn+1KAM}  satisfies, for any 
$ s \geq s_0 $, 
\begin{align}\label{lownorm1}
& | {\mathtt R}_{\nu+1} \langle D \rangle^M |_s^{k_0,\upsilon} \leq 
N_\nu^{-\tb} |  {\mathtt R}_\nu \langle D \rangle^M  |_{s+\tb}^{k_0,\upsilon} + C(s) 
N_\nu^{\tau_1} \upsilon^{-1} 
| {\mathtt R}_\nu \langle D \rangle^M |_{s_0}^{k_0,\upsilon} 
|  {\mathtt R}_\nu \langle D \rangle^M  |_{s}^{k_0,\upsilon} \, , \\
& \label{highnorm1}
| {\mathtt R}_{\nu+1} \langle D \rangle^M |_{s+\tb}^{k_0,\upsilon} \leq 
|  {\mathtt R}_\nu \langle D \rangle^M  |_{s+\tb}^{k_0,\upsilon} + C(s+\b)
N_\nu^{\tau_1} \upsilon^{-1} 
 |  {\mathtt R}_\nu \langle D \rangle^M  |_{s_0}^{k_0,\upsilon} 
 | {\mathtt R}_\nu \langle D \rangle^M |_{s+\tb}^{k_0,\upsilon} \, . 
\end{align}
\end{lemma}

\begin{pf}
The remainder  $ {\mathtt R}_{\nu+1}  $ defined in \eqref{Rn+1KAM}  satisfies, 
by \eqref{tame-prod}, \eqref{Neumann pseudo diff}, \eqref{psi-v-r}, 
\eqref{smallnessPsin-1} and  \eqref{small_KAM_con}, the estimate
$ |  {\mathtt R}_{\nu+1} \langle D \rangle^M |_s^{k_0,\upsilon} 
 \leq  $ $  | {\it \Pi_{N_\nu}^\bot}  {\mathtt R}_\nu \langle D \rangle^M  |_s^{k_0,\upsilon} + C(s)
|  {\mathtt R}_\nu \langle D \rangle^M |_{s_0}^{k_0,\upsilon} 
|   {\mathtt R}_\nu \langle D \rangle^M  |_{s}^{k_0,\upsilon} $. 
Using \eqref{smo-pro}  we obtain \eqref{lownorm1}-\eqref{highnorm1}. 
\end{pf}

As a corollary we deduce that 
\begin{lemma}
\eqref{estRn:iterativa} and \eqref{estRn:iterativahigh}
hold at the step $ \nu + 1 $.
\end{lemma}

\begin{pf}
By \eqref{lownorm1} and the inductive assumption \eqref{estRn:iterativa} we have 
$$
\begin{aligned}
| {\mathtt R}_{\nu+1} \langle D \rangle^M |_s^{k_0,\upsilon} 
& \leq 
N_\nu^{-\tb}  | {\mathtt R}_0 \langle D \rangle^M |_{s+\tb}^{k_0,\upsilon} 
N_{\nu-1} + C(s) 
N_\nu^{\tau_1} N_{\nu-1}^{-2\ta}  \upsilon^{-1} 
 | {\mathtt R}_0 \langle D \rangle^M |_{s+\tb}^{k_0,\upsilon} 
 | {\mathtt R}_0 \langle D \rangle^M |_{s_0+\tb}^{k_0,\upsilon}  \\
 & \leq | {\mathtt R}_0 \langle D \rangle^M |_{s+\tb}^{k_0,\upsilon} N_\nu^{-\ta} 
 \end{aligned}
$$
if 
$ N_\nu^{-\tb} N_{\nu-1} N_\nu^{\ta} \leq \tfrac12 $ and 
$ C(s) 
N_\nu^{\tau_1} N_{\nu-1}^{-2\ta} N_\nu^{\ta}  \upsilon^{-1} 
 | {\mathtt R}_0 \langle D \rangle^M |_{s_0+\tb}^{k_0,\upsilon} \leq \tfrac12 $. 
 These conditions are fulfilled by \eqref{choiceab}, \eqref{small_KAM_con} 
 and taking $ N_0 (s) $ large enough. The estimate
 \eqref{estRn:iterativa} at the step $ \nu + 1 $ is proved. 
 Moreover, by \eqref{highnorm1}, 
 \eqref{estRn:iterativa} at $ s = s_0 $, \eqref{choiceab} and \eqref{small_KAM_con}, we get
$ | {\mathtt R}_{\nu+1} \langle D \rangle^M |_{s+\tb}^{k_0,\upsilon}
\leq  $ $ 
 |  {\mathtt R}_\nu \langle D \rangle^M  |_{s+\tb}^{k_0,\upsilon} C \leq $
$  | {\mathtt R}_0 \langle D \rangle^M |_{s+\tb}^{k_0,\upsilon} N_\nu $
proving \eqref{estRn:iterativahigh} at the step $ \nu + 1 $. 
\end{pf}

The bound \eqref{einen} at $ \nu + 1 $ follows by Lemma \ref{lem:NNF}, and 
 \eqref{eintot} at $ \nu + 1 $ 
using also   \eqref{estRn:iterativa} and \eqref{small_KAM_con}. 
The statement $ ({\bf S1})_\nu $  in  Theorem \ref{iterative_KAM} is complete.
The proof of  $ ({\bf S2})_\nu $-$ ({\bf S3})_\nu $ follows similarly. 
\\[1mm]
{\bf Almost-approximate-invertibility  of  $ {\cal L}_\om $.}
We are now able to verify the almost-invertibility assumption (AI) in 
 \eqref{inversion assumption}-\eqref{tame inverse}. 
In Sections \ref{reduction} and \ref{sec:redu}  we have 
``almost-approximately" transformed the linear operator 
$ {\cal L}_\omega $ defined in \eqref{Lomega def} 
into the diagonal constant coefficient reversible operator 
$ {\mathtt L}_{ \tn} $ in \eqref{Linf}, through a  sequence of  transformations 
which map high norm Sobolev spaces into themselves. 
More precisely, for any $  \bar \tn \in \N $, by \eqref{forma finale operatore pre riducibilita}, \eqref{inizioKAMRe}, \eqref{defL0}, \eqref{OpL1} and \eqref{Linf} with $ \tn = \bar \tn $, we have that, 
for any $ (\omega, \gamma) \in  \mathtt{TC}_{\bar {\mathtt n}+1}(2\upsilon,\tau)\cap \Lambda^\upsilon_{\bar \tn} $,
\begin{equation}\label{final semi conjugation}
{\mathcal L}_\omega =  {\cal V}_{\bar \tn}  {\mathtt L}_{\bar \tn}  {\cal V}_{\bar \tn}^{-1} +\bphi_\bot {\mathfrak P}_{\bot,\bar \tn} \bphi_\bot^{-1} + \bphi_\bot {\mathfrak R}_{Z} \bphi_\bot^{-1}   \, ,  
\qquad {\cal V}_{\bar \tn} :=  \bphi_\bot {\bf M}_\e U_{\bar \tn}  \, ,
\end{equation}
where $ {\mathfrak P}_{\bot,\bar \tn}$ and $ {\mathfrak R}_{Z}$  are defined in Proposition \ref{prop: sintesi linearized}. Moreover,  by \eqref{stima Phi 1 Phi 2 proiettate}, 
\eqref{MM-1decay}, 
\eqref{Unvicina}, \eqref{small-inizi1} 
and \eqref{normapseudo-action} we deduce that the maps 
$ {\cal V}_{\bar \tn}^{\pm 1}  $ satisfy the tame estimates 
\begin{equation}\label{stime W1 W2}
\| {\cal V}_{\bar \tn}^{\pm 1} h \|_s^{k_0,\upsilon} \lesssim_S \| h \|_{s + \sigma }^{k_0,\upsilon}  + 
\| \fracchi_0 \|_{s + \perd(\mathtt b) + \sigma}^{k_0,\upsilon} 
\| h \|_{s_0 + \sigma}^{k_0,\upsilon} \,.
\end{equation}
We now  decompose 
$ {\mathtt L}_{\bar \tn} = \omega \cdot \pa_\vphi  {\mathbb I}_\bot  + {\mathtt D}_{\bar \tn} + {\mathtt R}_{ \bar \tn}
$ in \eqref{Linf} as 
\begin{align}\label{decomposizione bf Ln}
& {\mathtt L}_{\bar \tn} = {\mathtt L}_{ \bar \tn}^{<}  + {\mathtt R}_{\bar \tn} + {\mathtt R}_{ 
\bar \tn}^\bot    \qquad \text{with} \\
& \label{Rn-bot}
{\mathtt L}_{\bar \tn}^{<} := \Pi_{K_{\bar \tn}} \big( \Dom   {\mathbb I}_\bot  + {\mathtt D}_{ \bar \tn} \big) \Pi_{K_{\bar \tn}} + \Pi_{K_{\bar \tn}}^{\bot,\text{rev}}  \, , 
\quad {\mathtt R}_{\bar \tn}^\bot := 
\Pi_{K_{\bar \tn}}^\bot \big( \Dom  {\mathbb I}_\bot  + {\mathtt D}_{ \bar 
\tn} \big) \Pi_{K_{\bar \tn}}^\bot 
- \Pi_{K_{\bar \tn}}^{\bot,\text{rev}}  \, , 
\end{align}
where 
the operator $ {\mathtt D}_{\bar \tn} $ has the form 
\eqref{def:Ln}-\eqref{anbn-rin}, $ \Pi_K $ is in  \eqref{PiK2co},
 $ K_{\bar \tn} := K_0^{\chi^{ \bar\tn}} $
is the scale of the nonlinear Nash-Moser iterative scheme and, in accordance with the notation introduced in \eqref{defqthetaFou}-\eqref{matricesAnn}, 
$ \Pi_K^{\bot,\text{rev}}  $ is the reversible operator
$ 
  \Pi_K^{\bot,\text{rev}} := 
\text{diag}_{|(n,\ell)| > K \atop
{n\in \ST^c, \ell \in \Z^{|\ST|}} }
\begin{psmallmatrix}
0 & 1\\
1 & 0 
\end{psmallmatrix}.$

\begin{lemma} {\bf (Inverse operator)} 
\label{lem:first-Mel}
For all $ \lambda = (\om, \g)  $ in 
\begin{equation}\label{prime di melnikov}
 {\mathtt \Lambda}_{\bar \tn + 1}^{\upsilon, I}  
 :=  {\mathtt \Lambda}_{\bar \tn + 1}^{\upsilon, I} ( i ) 
 := 
\big\{ \lambda \in  \R^{|\ST|} \times [\g_1, \g_2] : |\omega \cdot \ell  
+  \Omega_{n,e}^{(\bar \tn)}| \geq  \tfrac{2\upsilon n}{  \langle \ell  \rangle^{ \tau}}, \,
 \forall | \ell  | \leq K_{\bar \tn},\, n \in  \ST^c \setminus \{3, \ldots, \bar n  \} \big\} 
\end{equation}
for any antireversible function $ g $, there is a unique reversible solution $ h $
of  $ {\mathtt L}_{\bar \tn}^< h = g $.
There exists an extension of the inverse operator (that we denote in the same way) to the whole $\R^{|\ST|} \times [\g_1, \g_2]$ satisfying the estimate, 
for $\mu =  k_0 + \t (k_0 +1)$, for any $ s \geq s_0 $,  
\begin{equation}\label{stima tilde cal Dn}
\| ({\mathtt L}_{\bar \tn}^<)^{- 1} g \|_s^{k_0, \upsilon} 
\lesssim_{k_0} \upsilon^{- 1} \| g \|_{s + \mu}^{k_0, \upsilon}\, .
\end{equation}
\end{lemma}

\begin{pf}
We identify  $ g(\vphi )$, resp. $ h(\vphi) $,  with the sequence 
$ \{ ( \widehat{g}_n^+ (\ell), \widehat{g}_n^- (\ell)  )\}_{n \in \ST^c, \ell \in \Z^{|\ST|}} $
of Fourier coefficients, resp.
$ \{ ( \widehat{\alpha}_n (\ell), \widehat{\beta}_n (\ell)  )\}_{n \in \ST^c, \ell \in \Z^{|\ST|}} $,  as in \eqref{defqthetaFou}. 
We denote $ \tn := {\bar \tn}$ for simplicity. 
Recalling \eqref{anbn-rin} and \eqref{Rn-bot}, 
the equation $ {\mathtt L}_\tn^< h = g $ amounts to 
\be\label{sys-inv}
\begin{split}
{\cal T}_{n,\ell,h}^{(\tn)}\begin{psmallmatrix}
  \widehat   \a_n  (\ell) \\
 \widehat  \b_n  (\ell)
\end{psmallmatrix}  &=  \widehat   g_n  (\ell)    \, ,  \quad  3\leq n \leq \bar n   \, ,
\quad  \widehat   g_n  (\ell) := 
\begin{psmallmatrix}
 \widehat   g_n^+  (\ell) \\
 \widehat   g_n^-  (\ell)
\end{psmallmatrix}  \\
{\cal T}_{n,\ell,e}^{(\tn)}\begin{psmallmatrix}
  \widehat   \a_n  (\ell) \\
 \widehat  \b_n  (\ell)
\end{psmallmatrix}  &= \widehat   g_n  (\ell)    \, ,  \quad  n \in \ST^c \setminus \{3,\ldots, \bar n \} \, , \quad \forall | (\ell,n)| \leq K_\tn  \, , 
\end{split}
\ee
where
$ {\cal T}_{n,\ell,h}^{(\tn)}:= \begin{psmallmatrix}
 \ii \om \cdot \ell   & \Omega_{n,h}^{( \tn)}    \\
  \Omega_{n,h}^{( \tn)}    & \ii \om \cdot \ell 
\end{psmallmatrix} $, 
$ {\cal T}_{n,\ell,e}^{(\tn)}:= \begin{psmallmatrix}
 \ii \om \cdot \ell   &   \Omega_{n,e}^{( \tn)}    \\
 - \Omega_{n,e}^{( \tn)}     & \ii \om \cdot \ell 
\end{psmallmatrix}  $
and $ (  \widehat   \a_n  (\ell),
 \widehat  \b_n  (\ell) ) = (g_n^-  (\ell),g_n^+  (\ell)) $ for any $ |(\ell,n)| > K_\tn $. 
The eigenvalues of ${\cal T}_{n,\ell,h}^{(\tn)}$ are $ \ii \om \cdot \ell\pm  \Omega_{n,h}^{( \tn)} (\e)$ and they are all different from zero with 
\be\label{deriv1hyp}
|\ii \om \cdot \ell\pm  \Omega_{n,h}^{( \tn)} (\e) | = 
|\omega \cdot  \ell  | + | \Omega_{n,h}^{( \tn)}  | \geq 
\max \{ |\omega \cdot  \ell  |, | \Omega_{n,h}^{( \tn)}  |  \}
\geq c > 0  \, .
\ee
Thus ${\cal T}_{n,\ell,h}^{(\tn)}$ is invertible for all the parameters $ (\omega, \gamma) \in \R^{|\ST|} \times [\gamma_1, \gamma_2] $.
 and 
$ \big({\cal T}_{n,\ell,h}^{(\tn)}\big)^{-1} \widehat   g_n  (\ell)$ is the unique solution 
of the first equation in  \eqref{sys-inv}. 
On the other hand, each matrix ${\cal T}_{n,\ell,e}^{(\tn)}  $ 
is diagonalized by 
$$ 
M := \frac{1}{\sqrt{2}} 
\begin{pmatrix} 
\ii & 1 \\
1 & \ii
\end{pmatrix} \, , 
\ \ 
M^{-1} {\cal T}_{n,\ell,e}^{(\tn)} M = D_{n,\ell,e} \, , \ \  
D_{n,\ell,e}  := 
\begin{psmallmatrix} 
\ii (\omega \cdot \ell - \Omega_{n,e}^{( \tn)} ) & 0  \\
0 & \ii (\omega \cdot \ell + \Omega_{n,e}^{( \tn)} ) 
\end{psmallmatrix} \, . 
$$
Let 
\be\label{deriv1ell}
(D_{n,\ell,e}^{-1})_{\text{ext}}  := 
\begin{psmallmatrix} 
I_{n,\ell}^-  & 0  \\
0 &  I_{n,\ell}^+
\end{psmallmatrix}, \quad 
I_{n,\ell}^\pm  := I_{n,\ell}^\pm  (\omega, \gamma ):=
\tfrac{\chi ((\omega \cdot \ell \pm \Omega_{n,e}^{( \tn)} ) \upsilon^{-1} \langle \ell \rangle^\tau)}{\ii (\omega \cdot \ell \pm  \Omega_{n,e}^{( \tn)}) } 
\ee
where $ \chi $  is the cut-off function defined in \eqref{cutoff}. 
Thus, defining
\be\label{deriv2ell}
\big({\cal T}_{n,\ell,e}^{(\tn)}\big)^{-1}_{\rm ext} 
:= M (D_{n,\ell,e}^{-1})_{\text{ext}}  M^{-1} 
= -  \tfrac{\ii}{2}  I_{n,\ell}^+ M_1 +  \tfrac{\ii}{2}  I_{n,\ell}^- M_2  \, ,
\ee
where  $M_1$, $M_2$ are the matrices  in \eqref{M1M4}, we deduce that 
$\big({\cal T}_{n,\ell,e}^{(\tn)}\big)^{-1}_{\rm ext} \widehat   g_n  (\ell)$ solves the second equation in  \eqref{sys-inv} 
 for any $ (\omega, \g)\in{\mathtt \Lambda}_{\tn + 1}^{\upsilon, I} $, see \eqref{prime di melnikov}. 
The estimate \eqref{stima tilde cal Dn} follows by the lower bounds of the 
small divisors in
\eqref{deriv1hyp}-\eqref{deriv2ell} 
and the properties of the cut-off function $ \chi $ in \eqref{cutoff}.
\end{pf}

By \eqref{final semi conjugation}, \eqref{decomposizione bf Ln}, 
Theorem \ref{thm:riduKAM}, 
 estimates \eqref{stime W1 W2}, \eqref{stima tilde cal Dn},   and using that, 
for all $ b  > 0$,   
\begin{equation}\label{stima tilde cal Rn}
\| {\mathtt R}_{\bar \tn}^\bot h \|_{s_0}^{k_0,\upsilon} \lesssim K_{\bar \tn}^{- b} \| h \|_{s_0 + b + 1}^{k_0,\upsilon}\,,
\qquad \| {\mathtt R}_{\bar \tn}^\bot h\|_s^{k_0,\upsilon} \lesssim \| h \|_{s + 1}^{k_0,\upsilon} \, , 
\end{equation}
we deduce the following result, stating the assumption (AI) on the almost-invertibility of  ${\cal L}_\omega$. 
\begin{theorem}\label{inversione parziale cal L omega}
{\bf (Almost-approximate-invertibility of $ {\cal L}_\omega $)} 
Let $ {\mathtt a}, {\mathtt b} $ as in \eqref{choiceab}, $ M \geq 1 $. 
There exists $\perd (\mathtt b) := \perd (\mathtt b,k_0,\tau)  > 0 $ such that, assuming 
\eqref{ansatz 0}  with $ \perd \geq \perd (\mathtt b)$ and,  for all 
$S > s_0$, the smallness condition
\eqref{KAM-small}, then 
 there exists $\sigma = \sigma(\tau, \mathbb S) > 0$ so that,  
 for any $ \bar \tn \ge 0$ and any
\begin{equation}\label{Melnikov-invert}
\lambda \in  {\bf \Lambda}_{\bar \tn + 1}^{\upsilon}:=\lambda \in  {\bf \Lambda}_{\bar \tn + 1}^{\upsilon} (i)  
:=\mathtt{TC}_{\bar {\mathtt n}+1}(2\upsilon,\tau) \cap  \Lambda_{\bar \tn + 1}^\upsilon   
\cap {\mathtt \Lambda}_{\bar \tn + 1}^{\upsilon, I} 
\end{equation}
(see \eqref{1-tras}, \eqref{sec-Meln-p}, \eqref{sec-Meln-p-22}, \eqref{prime di melnikov}),  
the operator $ {\mathcal L}_\omega$, defined in \eqref{Lomega def}, 
can be decomposed as 
\be\label{splitting cal L omega}
 {\mathcal L}_\omega  = {\mathcal L}_\omega^{<} + {\cal R}_\omega +  {\cal R}_\omega^\bot+ {\cal R}_\omega^Z  \quad \text{ with} \quad \begin{aligned}
&  {\mathcal L}_\omega^{<} := {\cal V}_{\bar \tn}^{-1} {\mathtt L}_{\bar \tn}^{<} {\cal V}_{\bar \tn} \,,\quad 
{\cal R}_\omega := {\cal V}_{\bar \tn}^{- 1} {\mathtt R}_{\bar \tn} {\cal V}_{\bar \tn} 
+\bphi_\bot {\mathfrak P}_{\bot,\bar \tn} \bphi_\bot^{-1} 
\,,\\
& {\cal R}_\omega^\bot := {\cal V}_{\bar \tn}^{- 1} {\mathtt R}_{\bar \tn}^\bot {\cal V}_{\bar \tn} \, , \quad {\cal R}_\omega^Z := \bphi_\bot {\mathfrak R}_{Z} \bphi_\bot^{-1}\, , 
 \end{aligned}
\ee
where  ${\mathcal L}_\omega^{<} $ is invertible and satisfies \eqref{tame inverse} and the operators ${\cal R}_\omega$, ${\cal R}_\omega^\bot$ and  $ {\cal R}_\omega^Z$ satisfy \eqref{stima R omega corsivo}-\eqref{stima R omega Z}. 
\end{theorem}

\section{Proof of Theorem \ref{main theorem}}\label{sec:NM}

Theorem \ref{main theorem} 
is a by now standard consequence of Theorem \ref{iterazione-non-lineare} below
which  provides a sequence of better and better approximate solutions 
of 
$ {\cal F} ( i, \tg) = 0$, 
where $\mF(i,\tg)$ is the nonlinear operator  
in \eqref{operatorF}. 
Note that the estimates 
\eqref{stima inverso approssimato 2}-\eqref{stima cal G omega-bot}
coincide with (5.63)-(5.66)  in \cite{BertiMontalto}. Thus we shall be short. 

We consider the finite-dimensional subspaces 
$$
E_\tn := \Big\{ \fracchi (\vphi ) = (\Theta , I , z) (\vphi) , \ \  
\Theta = \Pi_\tn \Theta, \ I = \Pi_\tn I, \ z = \Pi_\tn z \Big\}
$$
where $ \Pi_\tn := \Pi_{K_\tn}  $ is the projector
$ \Pi_\tn z(\ph,\theta) := \sum_{|(\ell ,j)| \leq K_\tn} 
z_{\ell,  j} e^{\ii (\ell  \cdot \ph + j\theta)}  $
with $ K_\tn = K_0^{\chi^\tn} $, $ \chi = 3/ 2 $,  and 
we denote with the same symbol 
$ \Pi_\tn p(\ph) :=  {\mathop \sum}_{|\ell | \leq K_\tn}  p_\ell e^{\ii \ell  \cdot \ph} $. 
We define $ \Pi_\tn^\bot := {\rm Id} - \Pi_\tn $.  
The projectors $ \Pi_\tn $, $ \Pi_\tn^\bot$ satisfy the smoothing properties \eqref{p2-proi}. 

In view of the Nash-Moser Theorem \ref{iterazione-non-lineare} we introduce the following constants: 
\begin{align}
& {\mathtt a}_1 :=  {\rm max}\{6  \sigma_1 + 13, \chi p (\tau + 1) (4\perd + 1) + \chi(\perd(\mathtt b) + 2 \sigma_1) +1 \} \, ,  \\
& \mathtt a_2 := \chi^{- 1} \mathtt a_1  - \perd(\mathtt b) - 2 \sigma_1 \, , \label{costanti nash moser} \\
& \mu_1 :=  3( \perd({\mathtt b}) + 2\sigma_1 ) + 1, 
\qquad {\mathtt b}_1 := {\mathtt a}_1 + \perd({\mathtt b}) +  3 \sigma_1 + 3 + \chi^{-1} \mu_1 ,
\qquad \chi = 3/ 2 \, ,
\label{costanti nash moser 1} 
\\
& \sigma_1 := \max \{ \bar \sigma\,, \s \} \, , \quad  
S = s_0 + {\mathtt b}_1 + \perd ({\mathtt b}) + \bar \sigma  \, ,  
\label{costanti nash moser 2}
\end{align}
where $\bar \sigma := \bar \sigma(\tau,\ST, k_0) > 0$ is defined in Theorem \ref{thm:stima inverso approssimato}, 
$ \s $ is the largest loss of regularity in the estimates of the Hamiltonian vector field $X_{\sP} $ in Lemma \ref{lemma quantitativo forma normale}, 
$ \perd(\mathtt b)$  is defined in 
Theorem \ref{inversione parziale cal L omega}, 
$\mathtt b$ is the constant $ {\mathtt b} := [{\mathtt a}] + 2 \in \N $ 
where $ {\mathtt a} $ is defined in \eqref{choiceab}, 
and the exponent $ p $ in \eqref{NnKn} satisfies  
$  p {\mathtt a} > $ $ (\chi - 1 ) {\mathtt a}_1 + \chi \sigma_1 = $ $ \tfrac12 {\mathtt a}_1 + \tfrac32 \sigma_1 $. There exists 
$ p := p(\tau,|\ST, k_0) $ such that this condition is verified 
by the definition of $ {\mathtt a}_1 $ in \eqref{costanti nash moser}. For example we fix
$ p := 3 (\perd(\mathtt b) + 3 \sigma_1 + 1)/ \mathtt a $.

We remark that the constant $ {\mathtt a}_1 $ is the exponent in \eqref{P2n}. 
The constant $ {\mathtt a}_2 $ is the exponent in \eqref{Hn}. 
The constant $ \mu_1 $ is the exponent in $({\cal P}3)_{\tn}$.

Given
$  W = ( \fracchi, \beta ) $ where
$   \fracchi = \fracchi (\lambda) $ 
is the periodic component of a torus as in \eqref{componente periodica}, and $ \b = \b (\l) \in \R^{|\ST|} $
we denote $ \|  W \|_{s}^{k_0, \upsilon} := \max\{ \|  \fracchi \|_{s}^{k_0, \upsilon} ,  |  \beta |^{k_0, \upsilon} \} $, 
where  $ \|  \fracchi \|_{s}^{k_0, \upsilon}  $ is defined in \eqref{def:norma-cp}.

\begin{theorem}\label{iterazione-non-lineare} 
{\bf (Nash-Moser)} 
There exist $ \d_0$, $ C_* > 0 $, such that, if
\be\label{nash moser smallness condition}  
K_0^{\tau_3} \e \upsilon^{- 2} < \d_0, 
\quad \tau_3 := \max \{ p \tau_2, 2 \sigma_1 + {\mathtt a}_1 + 4 \}, \  K_0 := \upsilon^{- 1}, 
\quad \upsilon := \e^a \, , 
\quad 0 < a < \tfrac{1}{\tau_3  + 2}\,,
\ee
where  $M \geq 1  $  and $ \tau_2 := \tau_2(\tau, {|\ST|})$ is  defined in Theorem \ref{iterative_KAM}, 
then, for all $ \tn \geq 0 $: 

\begin{itemize}
\item[$({\cal P}1)_{\tn}$] 
there exists a $k_0$ times differentiable function $\tilde W_\tn : \R^{|\ST|}  \times [\g_1, \g_2] 
\to E_{\tn -1} \times \R^{|\ST|} $, $ \lambda = (\omega, \g) \mapsto \tilde W_\tn (\lambda) 
:=  (\tilde \fracchi_\tn, \tilde \alpha_\tn - \om ) $, for  $ \tn \geq 1$,  
and $\tilde W_0 := 0 $, satisfying 
$ \| \tilde W_\tn \|_{s_0 + \perd({\mathtt b}) + \sigma_1}^{k_0, \upsilon} \leq C_*   \e \upsilon^{-1} $.  
Let $\tilde U_\tn := U_0 + \tilde W_\tn$ where $ U_0 := (\vphi,0,0, \om )$.
The difference $\tilde H_\tn := \tilde U_{\tn} - \tilde U_{\tn-1}$, $ \tn \geq 1 $,  satisfies
\begin{equation}  \label{Hn}
\|\tilde H_1 \|_{s_0 + \perd({\mathtt b}) + \sigma_1}^{k_0, \upsilon} \leq	 
C_* \e \upsilon^{- 1} \,,\quad \| \tilde H_{\tn} \|_{ s_0 + \perd({\mathtt b}) + \sigma_1}^{k_0, \upsilon} \leq C_* \e \upsilon^{-1} K_{\tn - 1}^{- \mathtt a_2} \,,\quad \forall \tn \geq 2. 
\end{equation}

\item[$({\cal P}2)_{\tn}$]   
Setting $ {\tilde \imath}_\tn := (\vphi, 0, 0) + \tilde \fracchi_\tn $, we define 
\be\label{def:cal-Gn}
{\cal G}_{0} := \tOm \times [\g_1, \g_2]\,, \quad 
{\cal G}_{\tn+1} :=  {\cal G}_\tn \cap {\bf \Lambda}_{\tn + 1 }^{\upsilon}({\tilde \imath}_\tn)
\,, \quad \tn \geq 0 \, , 
\ee
where $  {\bf \Lambda}_{n + 1}^{ \upsilon}({\tilde \imath}_\tn) $ is defined in \eqref{Melnikov-invert}. 
Then, for all $\lambda \in {\cal G}_\tn$,  setting $ K_{-1} := 1 $, we have 
\be\label{P2n}
\| {\cal F}(\tilde U_\tn) \|_{ s_{0}}^{k_0, \upsilon}  \leq C_* \e K_{\tn - 1}^{- {\mathtt a}_1} \, .
\ee
\item[$({\cal P}3)_{\tn}$] \emph{(High norms).} 
$ \| \tilde W_\tn \|_{ s_{0}+ {\mathtt b}_1}^{k_0, \upsilon} 
\leq C_* \e \upsilon^{-1}  K_{\tn - 1}^{\mu_1}$ for all $\lambda \in {\cal G}_\tn$.
\end{itemize}
\end{theorem}

\begin{corollary}
Let $ \gamma = \e^a  $ with $ a \in (0, a_0) $ and $ a_0 := 1 / (2+ \tau_3 ) $ where $\tau_3$ is defined in \eqref{nash moser smallness condition} and $ K_0 = \upsilon^{-1} $.
Then there is $ \e_0 > 0 $ so that for any $ 0 < \e \leq \e_0 $ the following holds:
\begin{enumerate}
\item there exists a function 
$ U_\infty (\lambda) := (i_\infty(\lambda), \a_\infty(\lambda)) 
\in H^{\bar s}_\vphi  \times H^{\bar s}_\vphi \times H^{\bar s}_{\vphi, x}
\times \R^{|\ST|}  $, defined for all $ \l \in \R^{|\ST|} \times [\g_1, \g_2] $, 
 where $ \bar s := s_0 + \perd(\mathtt b) + \sigma_1 $,  
satisfying 
\be\label{U infty - U n}
\|  U_\infty -  U_0 \|_{\bar s}^{k_0, \upsilon} \leq C_* \e \upsilon^{- 1} \,, \quad \| U_\infty - {\tilde U}_\tn \|_{\bar s}^{k_0, \upsilon} \leq C \e \upsilon^{-1} K_{\tn }^{- \mathtt a_2}\,, \ \  \tn \geq 1 \, .
\end{equation}
\item for any $ \lambda $ in the set

$
\bigcap_{\tn \geq 0} {\cal G}_\tn = 
\mG_0 \cap \bigcap_{\tn \geq 1} 
 {\bf \Lambda}_{\tn}^{\upsilon}(\tilde \imath_{\tn-1}) 
\stackrel{\eqref{Melnikov-invert}}{=} 
\mG_0 \cap \Big( \bigcap_{\tn \geq 1}  \mathtt{TC}_{\tn}(2\upsilon,\tau)
(\tilde \imath_{\tn-1})\Big)  \cap \Big( \bigcap_{\tn \geq 1}  \Lambda_{\tn}^{\upsilon}(\tilde \imath_{\tn-1}) \Big) \cap 
 \Big( \bigcap_{\tn \geq 1}   \mathtt \Lambda_{\tn}^{\upsilon, I}(\tilde \imath_{\tn-1}) \Big),
$

with $\mG_0 =  \mathtt \Omega \times [\g_1, \g_2] $,
the torus embedding $ i_\infty (\lambda) $ solves $ {\cal F}(\lambda, U_\infty(\lambda)) = 0 $. 
\end{enumerate}
\end{corollary}

\begin{pf}
For any  $ 0 < \e < \e_0 $ small enough, the smallness condition  
in \eqref{nash moser smallness condition} holds  and Theorem \ref{iterazione-non-lineare} applies.   
By $({\cal P}1)_{\tn}$  
the  sequence of functions 
$  {\tilde U}_\tn $ converges as $ \tn \to + \infty $ to a function $ U_\infty (\omega) $
satisfying  \eqref{U infty - U n}.
By Theorem \ref{iterazione-non-lineare}-$ ({\cal P}2)_\tn$, we deduce that 
$ {\cal F}(\lambda, U_\infty(\lambda)) = 0 $ for any $ \lambda  $ in
$ \cap_{\tn \geq 0} {\cal G}_\tn $.  
\end{pf}

To conclude the proof of Theorem \ref{main theorem},
we define the ``final eigenvalues". By \eqref{einen}, the sequence $({\frak r}_{n,e}^{({\mathtt n})} (i_\infty))_{\mathtt n\in \mathbb{N}}$, with ${\frak r}_{n,e}^{({\mathtt n})} (i_\infty)$ given by Theorem \ref{iterative_KAM} (evaluated at $i_\infty$), is a Cauchy sequence in $|\;\;|^{k_0,\upsilon}$.
Then we define
$ {\frak r}_{n}^{\infty} := \lim_{\bar{\mathtt n} \to+ \infty } {\frak r}_{n,e}^{({\mathtt n})} (i_\infty)  \,, $ for any $ n \in \ST^c \setminus \{3, \ldots, \bar n\} $, and 
 \eqref{stime autovalori infiniti} holds, 
  and,  by $({\bf S1})_\nu$ and \eqref{small-inizi1}   (note that 
$ S $ is fixed in \eqref{costanti nash moser 2}). 
Similarly, recalling \eqref{c1-picco}, 
we define  
${\mathtt r}^\infty_\e$  in \eqref{form1infty} as
$ {\mathtt r}^\infty_\e:=\lim_{{\mathtt n} \to + \infty }{\mathtt r}_{\e,{\mathtt n}}(i_\infty)
$ 
and \eqref{rnrs} holds. 
Finally, arguing as in \cite{BertiMontalto,BBHM} we deduce that 
the Cantor set $ {\cal C}_\infty^{\upsilon} $ in \eqref{Cantor set infinito riccardo} 
is contained in $ \cap_{\tn \geq 0} \mG_\tn $.
Indeed, by the inclusion properties \eqref{inclusion-transport}, \eqref{inclusione cantor riducibilita S4 nu} and \eqref{U infty - U n},
we   have that 
$$ 
{\cal  G}_\infty 
:= \mG_0   \cap \Big[ \bigcap_{\tn \geq 1}  \mathtt{TC}_{ {\mathtt n}}(4\upsilon,\tau)
(i_\infty)\Big]  \cap \Big[ \bigcap_{\tn \geq 1} \Lambda_\tn^{2 \upsilon}( i_\infty) \Big] 
\cap \Big[ \bigcap_{\tn \geq 1} \mathtt \Lambda_\tn^{2 \upsilon, I}(i_\infty)  \Big]
\subset 
\bigcap_{\tn \geq 0 } {\cal G}_\tn 
$$
where $\mG_\tn$ is defined in \eqref{def:cal-Gn}
and 
$ {\cal C}_\infty^\upsilon \subseteq {\cal  G}_\infty $, cfr. Lemma 8.6 in \cite{BertiMontalto}. 

\appendix

\section{Contour dynamics radial equation \eqref{Evera}}\label{app:CDE}

In this Appendix we provide the derivation of the evolutionary  
equation 
for the radial variable $ \xi (t, \theta)$. 

{{\begin{lemma}
Let $ D(t) $ 
be a bounded simply connected region with smooth boundary parametrized by 
$ z\,:\,  \R_+\times \T\, \to \,  \mathbb{C}  $, $  (t,\theta)\, \mapsto \, z(t, \theta) $. 
 Assume that $ \bw(t)= {\bf 1}_{D(t)} $. 
On the boundary,  the stream function $\psi$ in   \eqref{stream-patch}   
  is given   by
\be\label{formapsi}
\psi(t,z(t, \theta)) 
=  \frac{1}{8  \pi} \int_{\T}
\big[\ln \big(\vert  z( t, \theta')-z(t,  \theta)\vert^2\big) -1\big]  
{\rm Im} \big[ (\overline{ z(t,  \theta')}-\overline{z( t, \theta)}) \partial_{\theta'}z(t,  \theta') \big] d \theta' \, .
\ee
Moreover
\be\label{formpapsi}
\partial_\theta \psi (t,z(t, \theta)) 
=  -\frac{ 1}{4\pi }\int_{\T}\ln \big( |z(t,\theta)- z(t,\theta')|^2 \big) \partial^2_{\theta\theta'}\textnormal{ Im}\big[ \, \overline{z(t,\theta)} z(t,\theta')\big]d\theta' \, .
\ee
\end{lemma}

\begin{pf} By the complex form of Green's formula 
\be\label{green}
\int_D \pa_{\bar \zeta} f(\zeta, \bar \zeta) d A(\zeta)= \frac{1}{2 \ii} \int_{\pa D} f(\zeta, \bar \zeta) d \zeta \, , \quad \pa_{\bar \zeta} := \frac12 ( \pa_x + \ii \pa_y ) \, ,
\ee
 we can  write  the stream function $ \psi $ in \eqref{stream-patch} as an integral over the boundary  $\partial D(t)$,
  \be\label{stream-c}
\psi(t, z) 
= \frac{ 1}{8 \ii \pi}\int_{\partial D(t)}(\overline{\zeta}-\overline{z})\big[\ln
(\vert \zeta-z \vert^2) -1\big] d\zeta\, ,\quad 
\forall z\in \mathbb{C} \, .
\ee
Using the change of variables $\zeta=z( t, \theta')$  we obtain
\begin{align*}
\psi(t,z)   
 & = \frac{ 1}{8 \ii \pi}\int_{\T}(\overline{ z(t, \theta')}-\overline{z})
 \Big[\ln \big( \vert  z(t, \theta')-z\vert^2 \big) -1\Big]  \partial_{\theta'}z(t, \theta') d \theta'\, , \quad 
\forall z\in \mathbb{C}\, .
\end{align*}
On the boundary, namely $z=z(t,\theta)$ with $ \theta\in \mathbb{T}$, one gets the formula 
$$
\psi(t,z(t,\theta))   
 = \frac{ 1}{8 \ii \pi}\int_{\T}(\overline{ z(t, \theta')}-\overline{z(t, \theta)})
 \Big[\ln \big( \vert  z(t, \theta')-z(t, \theta)\vert^2 \big) -1\Big]  \partial_{\theta'}z(t, \theta') d \theta'\, 
$$
and, since $ \psi $
 is real valued,  formula  \eqref{formapsi} follows. Next, differentiating  \eqref{formapsi} with respect to $ \theta $ gives
\be\label{vbar}
\begin{aligned}
\partial_\theta \psi(t,z(t,\theta))
 &=-\frac{ 1}{8\pi }\int_{\T}\textnormal{Im}\Big[\frac{\overline{z(t,\theta')}- \overline{z(t,\theta)}}{z(t,\theta')- z(t,\theta)}\partial_{\theta'} z(t,\theta')\partial_\theta z(t,\theta)\Big]d\theta' \\
 &\quad - \frac{ 1}{8\pi }\int_{\T}\ln \big( \vert  z( t, \theta')-z(t,  \theta)\vert^2 \big) \textnormal{Im}\big[\partial_{\theta'} z(t,\theta')\partial_\theta\overline{z(t,\theta)}\big]d\theta' \, . 
\end{aligned}
\ee
Inserting the identity 
$
\tfrac{\overline{z(\theta')}- \overline{z(\theta)}}{z(\theta')- z(\theta)}\partial_{\theta'} z(\theta')=\big(\overline{z(\theta')}- \overline{z(\theta)}\big)\partial_{\theta'}\ln\big(|z(\theta')- z(\theta)|^2\big) -\partial_{\theta'} \overline{z(\theta')}  
$
 into \eqref{vbar} gives 
\be\label{vbar2}
\begin{aligned}
\partial_\theta \psi(t,z(t,\theta))&=
- \frac{ 1}{8\pi}\int_{\T}\partial_{\theta'}\ln\big(|z(\theta')- z(\theta)|^2\big) \textnormal{Im}\big[\big(\overline{z(t,\theta')}- \overline{z(t,\theta)}\big)\partial_\theta z(t,\theta)\big]d\theta' 
  \\&\quad-\frac{ 1}{8\pi }\int_{\T}\ln \big( \vert  z( t, \theta')-z(t,  \theta)\vert^2 \big) \textnormal{Im}\big[\partial_{\theta'} z(t,\theta')\partial_\theta\overline{z(t,\theta)}\big]d\theta' \, .  
\end{aligned}
\ee
By integrating by parts the first term in \eqref{vbar2}  
we deduce   \eqref{formpapsi}.
\end{pf}

\smallskip

\noindent
{\bf Proof of Lemma \ref{lem:eq-rad-def}.} 
Differentiating $w (t, \theta)$ in  \eqref{time-para} with respect to $ t $ and $ \theta $, we deduce  
$ \pa_t w (t, \theta) =  (1 + 2 \xi (t, \theta))^{- \frac12}  \pa_t \xi (t, \theta) \we(\theta) $
and  
\be
 \begin{aligned}\label{der-eq11}
\pa_\theta w(t, \theta) 
& = 
 (1 + 2 \xi (t, \theta))^{- \frac12} \pa_\theta \xi (t, \theta) \we(\theta)
+ 
(1 + 2 \xi (t, \theta))^{\frac12} \pa_\theta \we(\theta) \, . 
\end{aligned}
\ee
Then, the left hand side of \eqref{sub-HS} writes
\be\label{pat}
{\rm Im}\big[ \pa_t w (t, \theta) \ov{\pa_\theta w(t, \theta)} \big] =  
\pa_t \xi (t, \theta) \, {\rm Im}\big[  \we(\theta) \ov{\pa_\theta \we( \theta)} \big]
 = - \pa_t \xi (t, \theta) \, ,
\ee
having used the identity \eqref{id:we}.
Moreover, by a direct calculus we obtain
\be \label{zq}
 \begin{aligned}
& | w(t, \theta)|^2 =
| \we(\theta)|^2 (1 + 2 \xi (t, \theta) )=
 g_\gamma(\theta)(1 + 2 \xi (t, \theta) ) \, ,\\
 &  | w(t, \theta )  - w (t, \theta')|^2 = \Ke(\xi)(\theta, \theta')\, ,
\end{aligned}
\ee
where   $ \Ke(\xi)(\theta, \theta') $ and $g_\gamma(\theta)$ are defined in \eqref{expR} and \eqref{fg0} respectively. 
By \eqref{time-para} we get
 \begin{align}
{\rm Im} \big[ \partial_{\theta'} w (t,\theta')   
\overline{\pa_\theta w(t, \theta) } \big]
& = \pa^2_{\theta \theta'}
\big[ (1+ 2 \xi (t, \theta))^\frac12 (1+ 2 \xi (t, \theta') ^\frac12   \sin (\theta' - \theta) \big] \, .\label{ide}
\end{align}
Combining  \eqref{formpapsi} with the identities \eqref{zq} and \eqref{ide} we  get that 
 \be
 \partial_\theta \psi(t, w(t, \theta))  
 = \frac{ -1}{4  \pi} \int_{\T}
\ln (\Ke(\xi) (\theta, \theta'))   \pa^2_{\theta\theta'}\big[(1 + 2 \xi(t,\theta))^{ \frac12}(1 + 2 \xi (t,\theta'))^{ \frac12}\sin (\theta' - \theta)\big]
d \theta' \, . \label{parpsi-new}
\ee
Inserting  the identities 
\eqref{pat}, \eqref{zq} and \eqref{parpsi-new}  into \eqref{sub-HS} we obtain \eqref{Evera}. This completes the proof of Lemma \ref{lem:eq-rad-def}. \qed

We now prove that the  vector field in the right hand side of \eqref{Evera} 
vanishes 
at $ \xi (\theta) = 0 $ if and only if the angular velocity $ \Omega = \Omega_\g $ as in \eqref{choiceOmega}, proving 
 Lemma \ref{lem:equilibrium}. For  this aim, let us compute the integral
$ \tfrac{1}{4 \pi} 
\int_{\T} \ln ( \Ke(0)(\theta, \theta') ) 
 \sin (\theta' - \theta ) d \theta' $, 
where,  by \eqref{expR}, 
\be\label{expR0}
\Ke(0)(\theta, \theta')  = 
 \g \big(  \cos \theta -  \cos \theta' \big)^2 + \g^{-1}
 \big(  \sin \theta - \sin \theta' \big)^2 \, . 
\ee
We shall use the following  singular integrals computed  in \cite{CCS}: for any $ k \in \Z \setminus \{0\}, $
\be\label{Iccs0}
\int_{\T} e^{\ii k \theta } \ln \Big( \sin^2 \Big( \frac{\theta}{2} \Big) \Big) d \theta = - \frac{2\pi}{|k|} \, ,  \quad 
\int_{\T} e^{\ii k \theta } 
\ln \Big( \frac{1+r^2}{1-r^2} - \cos (\theta) \Big) d \theta = 
- \frac{2\pi}{|k|} \Big( \frac{1-r}{1+r}  \Big)^{|k|}\, .
\ee

\begin{lemma}\label{lem:comp-integ} For all $\gamma \geq 1$ we have 
\begin{align}\label{int14}
 \frac{1}{4 \pi} 
\int_{\T} \ln ( \Ke(0)(\theta, \theta') ) 
 \sin (\theta' - \theta ) d \theta'
=-\frac{\Omega_\gamma}{2}\partial_\theta g_\gamma(\theta) \, , 
\end{align}
where $g_\gamma(\theta)$ is introduced in  \eqref{fg0}. 
\end{lemma}
\begin{pf}
Using 
$ \cos \theta -  \cos \theta'  =-2\sin \big( \tfrac{\theta- \theta'}{2} \big)\sin \big( \tfrac{\theta'+ \theta}{2} \big) $ and
$   \sin \theta -  \sin \theta' =2\sin \big( \tfrac{\theta- \theta'}{2} \big)\cos \big( \tfrac{\theta'+ \theta}{2} \big) $, 
 we get, by \eqref{expR0} and the identity
 $ \sin^2 (\alpha) = (1- \cos (2\a))/2 $, 
   \be\label{rel-dis0}
\Ke(0)(\theta, \theta')    = 2 
\frac{\gamma^2-1}{\gamma} \sin^2 \Big( \frac{\theta'- \theta}{2} \Big) \Big[  
\frac{\g^2+1}{\g^2-1} -
  \cos (\theta + \theta' ) \Big] \, . 
  \ee
By \eqref{rel-dis0}, we obtain
\be\label{pe1}
\begin{aligned}
&\int_{\T} \ln ( \Ke(0)(\theta, \theta') ) 
 \sin (\theta' - \theta ) d \theta'  =  
\int_{\T}  \ln \Big( 2 \frac{\gamma^2-1}{\gamma} \Big)   \sin (\theta' - \theta ) d \theta'  \\
&\quad + 
\int_{\T}  \ln \Big( \sin^2 \Big( \frac{\theta'- \theta}{2} \Big) \Big) \sin (\theta' - \theta ) d \theta' 
+
\int_{\T}   \ln \Big(  
\frac{\g^2+1}{\g^2-1} -
  \cos (\theta + \theta' ) \Big)
 \sin (\theta' - \theta ) d \theta'  \, .
 \end{aligned}
 \ee
The first term in \eqref{pe1} is  zero. The second  integral in \eqref{pe1} is zero
 by oddness.
In order to compute the third one, 
we write
$ \sin (\theta' - \theta) = \sin ( \theta' + \theta - 2 \theta )  = 
\sin ( \theta' + \theta) \cos  ( 2 \theta ) - \cos ( \theta' + \theta) \sin  ( 2 \theta ) $, 
and therefore 
\be \label{lic1}
\begin{aligned}
\int_{\T}   \ln \Big(  
\frac{\g^2+1}{\g^2-1} -
  \cos (\theta + \theta' ) \Big)
 \sin (\theta' - \theta ) d \theta'
& = 
\cos  ( 2 \theta ) \int_{\T}   \ln \Big(  
\frac{\g^2+1}{\g^2-1} -   \cos  \theta'  \Big)
\sin  \theta' d \theta'  \\
&  - \sin  ( 2 \theta ) \int_{\T}   \ln \Big(  
\frac{\g^2+1}{\g^2-1} -   \cos \theta'  \Big) 
\cos  \theta'  d \theta'  \, . 
\end{aligned}
\ee
The first integral in \eqref{lic1} is zero because the integrand is odd and the second one
is
\be\label{lastsi}
\int_{\T}   \ln \Big(  
\frac{\g^2+1}{\g^2-1} -   \cos \theta'  \Big) 
\cos  \theta'  d \theta' = 
\int_{\T}   \ln \Big(  
\frac{ 1 + (\g^{-1})^{2}}{ 1 - (\g^{-1})^2 } -   \cos ( \theta' ) \Big) 
\cos ( \theta' ) d \theta' \stackrel{\eqref{Iccs0}} = -
2 \pi  \frac{\g-1}{\g+1} \, . 
\ee
The identity \eqref{int14} follows by \eqref{lastsi}, \eqref{pe1} and \eqref{fg0}. The proof of Lemma \ref{lem:comp-integ}
is complete. 
\end{pf}

\smallskip

{\bf Proof of  Lemma \ref{lem:equilibrium}.} 
Follows by   \eqref{Evera}, \eqref{fg0}, \eqref{choiceOmega} and \eqref{int14}.

\section{Technical lemmata}\label{sec:B}

For completeness we report in this Appendix some  properties
of the norm $ \normk{\ }{s} $ in \eqref{weinorm},   
proved in \cite{BertiMontalto}.
For any $N>0$, we define the smoothing operator
\begin{equation}\label{def:smoothings}
\Pi_N :  
u(\vphi,\theta) =  
{\mathop \sum}_{\ell\in\Z^{|\ST|},j\in\Z}u_{\ell,j}e^{\ii (\ell\cdot \vphi +j \theta )} \mapsto
(\Pi_N u)(\ph, \theta) := {\mathop \sum}_{\la \ell,j \ra \leq N} u_{\ell, j} e^{\ii (\ell\cdot\ph + j \theta )} \, 
\end{equation}
and set $ \Pi^\perp_N := {\rm Id} - \Pi_N $. 

\begin{lemma} \label{lem:prod0}
The following  tame estimates for the product and composition 
 hold. 
\\[1mm]
1. For all $s\geq s_0 > (|\ST|+1)/2$, 
	\begin{equation}
	\label{prod}
	\normk{u v}{s} \leq C(s,k_0) \normk{u}{s}\normk{v}{s_0} + C(s_0,k_0)\normk{u}{s_0}\normk{v}{s}\,.
	\end{equation}
2. For any $N\geq 1$  the 
operators $\Pi_N, \Pi_N^\perp$ in \eqref{def:smoothings} 
satisfy 
\be\label{p2-proi} 
\| \Pi_N u \|_{s}^\zug 
 \leq N^\alpha \| u \|_{s-\alpha}^\zug \, , \quad 0 \leq \a \leq s \, , 
 \quad 
\| \Pi_N^\bot u \|_{s}^\zug  \leq N^{-\alpha} \| u \|_{s + \alpha}^\zug\, , \quad  \a \geq 0 \, .
\ee	
3. Let $ \| \b \|_{2s_0+k_0+1}^\zug \leq \d (s_0, k_0) $ small enough. 
Then the composition operator 
\be\label{change-variab-diff}
(\mB u)(\ph,\theta) := u(\ph, \theta + \b (\ph,\theta))  
\ee
satisfies the following tame estimates: for all $ s \geq s_0$, 
\be\label{pr-comp1}
\| {\cal B} u \|_{s}^\zug \leq C(s,k_0)\big( \| u \|_{s+k_0}^\zug 
+ \| \b \|_{s}^\zug \| u \|_{s_0+k_0+1}^\zug\big) \, .
\ee
4. Let $ \| \b \|_{2s_0+k_0+1}^\zug \leq \d (s_0, k_0) $ small enough.  The function $ \breve \beta $ defined by 
the inverse diffeomorphism 
$ y = \theta + \beta (\vphi, \theta) $ if and only if $ \theta = y + \breve \b ( \vphi, y ) $,  
satisfies 
$ \| \breve \beta \|_{s}^\zug \leq C(s,k_0) \| \b \|_{s+k_0}^\zug $. 
\end{lemma}	

We also state a standard Moser tame estimate for the nonlinear composition operator, see for instance Lemma 2.31 in \cite{BertiMontalto}, 
$ u(\vphi,\theta)\mapsto {\mathtt f}(u)(\vphi,\theta) :=f(\vphi,\theta,u(\vphi,\theta)) $. 
\begin{lemma}{\bf (Composition operator)}\label{compo_moser}
	Let $f\in  \mathcal{C}^\infty(\T^{|\ST|+1},\R)$. If $u(\lambda)\in H^s $ is a family of Sobolev functions satisfying $\normk{u}{s_0}\leq 1$, then, for all $s > s_0 >  (d+1)/2$,
	$	 \normk{{\mathtt f}(u)}{s}\leq C(s,k_0,f)\big( 1+\normk{u}{s} \big) $. 
	If $ f(\varphi, \theta, 0) = 0 $ then 
	$	  \normk{{\mathtt f}(u)}{s}\leq C(s,k_0,f) \normk{u}{s} $. 
\end{lemma}
An integral operator with smooth kernel is infinitely many times 
regularizing.  

\begin{lemma}[Lemma 2.32 of~\cite{BertiMontalto}]  
\label{lem:Int}
Let $ K := K( \lambda, \cdot ) \in {\cal C}^\infty(\T^{|\ST|} \times \T \times \T)$ for all
$ \l \in {\mathtt \Lambda}_0 $. Then the integral operator 
\be\label{integral operator}
({\cal R} u ) ( \vphi, \theta) := \int_\T K(\lambda, \vphi, \theta, \theta') u( \theta')\,d \theta' 
\ee
is in $ {\rm OPS}^{- \infty}$ and, for all $ m, s,  \a \in \N_0 $,  
$ | {\cal R}  |_{-m, s, \a}^{k_0, \upsilon} \leq C(m, s, \alpha, k_0)  \| K \|_{{\cal C}^{s + m + \alpha}}^{k_0, \upsilon} $. 
\end{lemma}
An integral operator transforms into another integral operator under the changes of variables \eqref{change-variab-diff}.
The following result is  Lemma 2.34 of \cite{BertiMontalto}.
\begin{lemma}\label{lemma:conj-integ-op}
Let $K(\lambda, \cdot ) \in {\mathcal C}^\infty(\T^{|\ST|} \times \T \times \T)$ 
and $p(\lambda, \cdot ) \in {\mathcal C}^\infty(\T^{|\ST|} \times \T, \R) $. 
There exists $\delta := \delta(s_0, k_0) > 0$ such that if $\| \beta\|_{2 s_0 + k_0 + 1}^{k_0, \upsilon} \leq \delta$, then 
the integral operator $ {\mathcal R}$ as in \eqref{integral operator} transforms into the integral operator
$ \big({\cal B}^{- 1} {\mathcal R} {\cal B} \big) u(\vphi,  \theta) = \int_\T \widetilde K(\lambda, \vphi, \theta, y) u(\vphi, y)\,d y $
with a  $ {\mathcal C}^\infty $ 
Kernel $\widetilde K(\lambda, \cdot, \cdot, \cdot)$ which satisfies
$
\| \widetilde K\|_{s}^{k_0, \upsilon} \leq $ $ C(s, k_0)  \big( \| K \|_{s + k_0}^{k_0, \upsilon} 
+ \| p \|_{{s + k_0 + 1}}^{k_0, \upsilon} \| K \|_{s_0 + k_0 + 1}^{k_0, \upsilon}\big) $,
for any  $ s \geq s_0$. 
\end{lemma}
The commutator between the Hilbert transform ${\cal H}$ in \eqref{Hilbert-transf} and the multiplication operator
for a smooth function is a regularizing operator. 

\begin{lemma}[Lemma 2.35 in \cite{BertiMontalto}]\label{lemma:commutator-Hilbert}
Let $ a( \lambda, \cdot, \cdot ) \in {\mathcal C}^{\infty} (\T^{|\ST|} \times \T, \R)$. Then the commutator $[a, {\cal H} ] \in $ $ \text{OPS}^{-\infty }$ and, for all $ m,  s,  \a \in \N_0 $, 
$ \big|  [a, {\cal H} ] \big|_{-m, s, \a}^{k_0, \upsilon} \leq C(m, s, \alpha, k_0) \| a \|_{{s + k_0 + 1+ m + \alpha}}^{k_0, \upsilon} $. 
\end{lemma}
\begin{lemma}[Lemma 2.36 in \cite{BertiMontalto}] \label{lemma:conjug-Hilbert}
Let $ \beta(\lambda, \cdot) \in {\mathcal C}^\infty(\T^{|\ST| + 1})$. There exists $\delta(s_0, k_0) > 0$ such that, if $\|\beta \|_{2 s_0 + k_0 + 1}^{k_0, \upsilon} \leq \delta(s_0, k_0)$, then ${\cal B}^{-1} {\cal H} {\cal B}^{-1} - {\cal H}$ is an integral operator of the form
$ ( {\cal B}^{-1} {\cal H} {\cal B} - {\cal H}) u (\vphi, \theta) 
= $ $ \int_\T  \,K(\lambda, \vphi, \theta,z) u(\ph, z)\,dz $
where $ K = K(\lambda, \cdot) \in {\mathcal C}^\infty (\T^{|\ST|} \times \T \times \T ) $ satisfies 
$ \| K \|_{s}^{k_0, \upsilon} \leq C(s, k_0) \| \beta \|_{s +  k_0 + 2}^{k_0, \upsilon} $, 
for all $  s \geq s_0 $. 
\end{lemma}

\noindent
{\bf Diophantine equation.} 
We recall basic facts about diophantine equations.
If  $\omega $ is a Diophantine vector in $\mathtt {DC} (\upsilon, \tau) $
 the equation $\omega\cdot \pa_\vphi v = u$, where $u(\vphi,\theta)$ has zero average with respect to $\vphi$, has the periodic solution
$
	(\omega\cdot\pa_\vphi)^{-1} u := \sum_{\ell\in\Z^{|\ST|}\setminus\{0\},j\in\Z} \frac{u_{\ell,j}}{\ii\,\omega\cdot \ell} e^{\ii(\ell\cdot\vphi+j\theta)}$. 
For all $\omega\in\R^{|\ST|}$, we define its extension
\begin{equation}\label{paext0}
	(\omega\cdot\pa_\vphi)_{\rm ext}^{-1} u(\vphi,\theta) := 
	{\mathop \sum}_{(\ell,j)\in\Z^{|\ST|+1}}\frac{\chi(\omega \cdot \ell \upsilon^{-1}\braket{\ell}^\tau)}{\ii \omega\cdot\ell} u_{\ell,j}e^{\ii(\ell\cdot\vphi+j\theta)}\,,
\end{equation}
where $\chi\in  \mathcal{C}^\infty(\R,\R)$ is an even positive $  \mathcal{C}^\infty$ cut-off function such that
\begin{equation}\label{cutoff}
	\chi(\eta) := \begin{cases}
	0 & \text{ if } \ |\eta|\leq 1/3 \\
	1 & \text{ if } \ |\eta| \geq 2/3 
	\end{cases}\,, \qquad \pa_\eta \chi(\eta) >0   \quad \forall\,\eta \in (\tfrac13,\tfrac23) \,.
\end{equation}
Note that $(\omega\cdot\pa_\vphi)_{\rm ext}^{-1} u = (\omega\cdot\pa_\vphi)^{-1}u$ for all $\omega\in \mathtt{DC}(\upsilon,\tau)$. 
The following estimate holds
$		\normk{(\omega\cdot\pa_\vphi)_{\rm ext}^{-1}u}{s} \leq C(k_0)\upsilon^{-1}\normk{u}{s+\tau_1} $, $  \tau_1 :=\tau(k_0+1)+k_0$


\footnotesize

\begin{thebibliography}{99}


\bibitem{Alaz-Bal} Alazard T., Baldi P.,
 {\it Gravity capillary standing water waves},
 {Arch. Rat. Mech.  Anal.},  217, 3, 741--830,  2015. 

\bibitem{Ar} Arnold V.I., 
{\it Proof of a theorem of A. N. Kolmogorov on the persistence of quasi-periodic motions under small perturbations of the Hamiltonian},
 Russ. Math. Surv., 18, 9--36, 1963.

\bibitem{BBHM} Baldi P., Berti M., Haus E., Montalto R., 
{\it Time quasi-periodic gravity water waves in finite depth}, 
Inventiones Math. 214 (2), 739--911, 2018.

\bibitem{BBM-Airy} {Baldi P.,  Berti M., Montalto R.},
 {\it KAM for quasi-linear and fully nonlinear forced perturbations of Airy equation}, 
Math. Annalen 359, 471--536, 2014. 

\bibitem{BBM-auto} {Baldi P., Berti M., Montalto R.}, 
{\it KAM for autonomous quasi-linear perturbations of KdV}, 
 Ann. Inst. H. Poincar\'e Analyse Non. Lin. 33, no. 6, 1589--1638, 2016. 


\bibitem{BM20} Baldi P., Montalto R.,
{\it Quasi-periodic incompressible Euler flows in 3D},
Advances in Mathematics, DOI: 10.1016/j.aim.2021.107730.

\bibitem{BaBM}
Bambusi D., Berti M., Magistrelli E.,
 {\it Degenerate KAM theory for partial differential equations}, 
Journal Diff. Equations, 250, 8, 3379--3397, 2011.


\bibitem{BBP2}   Berti M., Biasco L., Procesi M., 
{\it KAM for Reversible Derivative Wave Equations},
 Arch. Ration. Mech. Anal. 212(3), 905--955, 2014.


\bibitem{BBo1} Berti M., Bolle P., {\it Sobolev quasi periodic solutions of multidimensional 
wave equations with a multiplicative potential},  Nonlinearity, 25, 2579-2613, 2012.


\bibitem{BB10}  Berti M., Bolle P.,  
{\it Quasi-periodic solutions  with Sobolev regularity of NLS on $ \T^d $ with a multiplicative potential}, 
 J. Eur. Math. Soc., Vol. 15, 229-286, 2013.




\bibitem{BB13} Berti M., Bolle P.,  
{\it A Nash-Moser approach to KAM theory}, 
Fields Institute Communications, special volume ``Hamiltonian PDEs and Applications'', 
255--284,  2015. 


\bibitem{BB20}  Berti M.,  Bolle P., 
{\it ``Quasi-periodic solutions of nonlinear wave equations on $\T^d$"}, 
vii + 355,  { Monographs of the EMS}. ISBN print 978-3-03719-211-5, 2020. 



\bibitem{BFM1} Berti M., Franzoi L., Maspero A.,
{\it Traveling quasi-periodic water waves with constant vorticity},  
Arch. Ration. Mech. Anal. 240, 99--202, 2021. 

\bibitem{BFM}  Berti M.,  Franzoi L.,  Maspero A., 
{\it Pure gravity traveling quasi-periodic water waves with constant vorticity}, 
arXiv:2101.12006, to appear on Comm. Pure Applied Math. 


\bibitem{BKM} Berti M., Kappeler T., Montalto R., 
{\it Large KAM tori for perturbations of the dNLS equation}, 
Ast\'erisque, 403, viii + 148, 2018.


\bibitem{BKM1} Berti M.,  Kappeler T., Montalto R., 
{\it Large KAM tori for quasi-linear perturbations of KdV}, 
Archive for Rational Mechanics, 239, 1395--1500, 2021. 


\bibitem{BertiMontalto} 
Berti M., Montalto R., 
{\it Quasi-periodic standing wave solutions of gravity-capillary water waves}, 
Memoires AMS, Volume 263, 1273, ISSN 0065--9266, 2020.

\bibitem{BC} Bertozzi A., Constantin P., 
{\it Global regularity for vortex patches}, 
Comm. Math. Phys. 152, 19--28, 1993.


\bibitem{BM} Bertozzi A.,  Majda A., 
 {\it Vorticity and Incompressible Flow}, 
Cambridge Univ. Press, 2001.


\bibitem{B3}  Bourgain J., {\it Quasi-periodic solutions of Hamiltonian
perturbations of $2D$ linear Schr\"odinger equations},
Annals of Math. 148, 363-439, 1998.

\bibitem{B5}  Bourgain J., 
 {\it Green's function estimates for lattice Schr\"odinger operators and applications}, 
Annals of Mathematics Studies 158, Princeton University Press, Princeton, 2005.


\bibitem{Bur} Burbea J., 
 {\it Motions of vortex patches.} 
Lett. Math. Phys. 6, no. 1,  1--16, 1982.

\bibitem{CCS} Castro A., C\'ordoba D., G\'omez-Serrano J., 
{\it Uniformly rotating analytic global patch solutions for active scalars}, 
Ann. PDE  2, no. 1, 1--34, 2016.


\bibitem{Ch} Chemin J.-Y., 
{\it Persistance de structures geometriques dans les fluides incompressibles bidimensionnels},
 Ann. Ec. Norm. Sup. 26, 4, 1--16, 1993.

\bibitem{CP}
Chierchia L., Pinzari G.,
{\it The planetary N-body problem: symplectic foliation, reductions and invariant tori},
Inventiones Math. 186, no.1, 1--77, 2011.



\bibitem{DZ}  Deem G.S., Zabusky N. J.,
 {\it Vortex waves: Stationary "V-states", Interactions, Recurrence, and Breaking},
Phys. Rev. Lett.  40, no. 13, 859--862, 1978.

\bibitem{DH} Duistermaat J.J., H\"ormander L. { \it Fourier integral operators. II.} Acta Math. 128, 183-269, 1972. 


\bibitem{HFMV} de la Hoz F.,   Hmidi T.,    Mateu J.,  Verdera J., 
 {\it  Doubly connected V-states for the planar Euler equations},
 SIAM J. Math. Anal. 48, no. 3, 1892--1928, 2016. 

\bibitem{Eliasson-Greber-Kuksin} 
Eliasson L.H., Gr\'ebert B., Kuksin S., 
{\it KAM for the nonlinear beam equation},  Geom. Funct. Anal. 26, 1588-1715, 2016.

 
\bibitem{Eliasson-Kuksin} Eliasson L.H., Kuksin S.,
{\it KAM for nonlinear Schr\"odinger equation},
 Ann. Math.  172, 371--435, 2010.

 \bibitem{Marmi-Fasano} Fasano A., Marmi S.,
  {\it Analytical Mechanics, an introduction}, 
 Oxford graduate Texts, 2006.  
 

\bibitem{FG} Feola R., Giuliani F., 
 {\it Quasi-periodic traveling waves on an infinitely deep fluid under gravity}. 	ArXiv:2005.08280, to appear on Memoires American Mathematical Society. 


\bibitem{FGMP}
Feola R.,  Giuliani F.,  Montalto R., Procesi M.,
{\it Reducibility of first order linear operators on tori via Moser's theorem},  
J. Funct. Anal. 276, no. 3, 932--970, 2019.

\bibitem{FGP} 
Feola R., Giuliani F.,  Procesi M.,  
{\it Reducible KAM tori for the Degasperis-Procesi equation}, 
Comm. Math. Phys., 377(3):1681-1759, 2020.

\bibitem{Feola-Procesi} Feola R., Procesi M., 
{\it Quasi-periodic solutions for fully nonlinear forced reversible Schr\"odinger equations},
 J. Diff. Eq. 259, no.7, 3389--3447, 2015.

\bibitem{Giu}
Giuliani F.,
{\it Quasi-periodic solutions for quasi-linear generalized KdV equations}, 
J. Differential Equations, 262(10):5052-5132, 2017.


\bibitem{GPS}  G\'omez-Serrano J.,   Park J.,   Shi J.,  
{\it Existence of non-trivial non-concentrated compactly supported stationary solutions of the 2D Euler equation with finite energy},
arXiv:2112.03821. 


\bibitem{GHS}
Guo Y., Hallstrom C.,  Spirn D.,  
{\it Dynamics near an unstable Kirchhoff ellipse},  
Comm. Math. Phys. 245, no.2, 297--354, 2004. 

\bibitem{HH2} Hassainia  Z.,  Hmidi  T.,
 {\it Steady asymmetric vortex pairs for Euler equations},
 { Discrete Contin. Dyn. Syst.} 41, no. 4, 1939--1969, 2021.
 
 	\bibitem{HHM21} Hassainia Z., Hmidi T., Masmoudi N., 
	\textit{KAM theory for active scalar equations,} arXiv:2110.08615.

 \bibitem{Hass-Mass-Wheel}  Hassainia Z.,  Masmoudi N., Wheeler M. H., 
 {\it Global bifurcation of rotating vortex patches},
 Comm. Pure Appl Math., 73, no. 9, 1933--1980, 2020.

\bibitem{HaR} Hassainia Z., Roulley E., 
{\it Boundary effects on the emergence of quasi-periodic solutions for Euler equations},
arXiv 2202.10053.

 
 
    \bibitem{HM2}  Hmidi T.,   Mateu J.,
     {\it Bifurcation of rotating patches from Kirchhoff vortices}, 
     Discrete Contin. Dyn. Syst. 36, no. 10, 5401--5422, 2016. 
     
      \bibitem{HM3} Hmidi T.,   Mateu J.,
  {\it Degenerate bifurcation of the rotating patches}. 
  Adv. Math. 302 , 799--850, 2016. 
  
   \bibitem{HM} Hmidi T.,   Mateu J.,
   {\it  Existence of corotating and counter-rotating vortex pairs for active scalar equations}, 
    Comm. Math. Phys. 350, no. 2, 699--747, 2017.   
     
      \bibitem{HMV} Hmidi T.,   Mateu J.,  Verdera J.,
  {\it Boundary Regularity of Rotating Vortex Patches}, 
  Arch. Ration. Mech. Anal. 209, no. 1, 171--208, 2013.
  
   \bibitem{HR}  Hmidi T., Renault C.,
  {\it Existence of small loops in a bifurcation diagram near degenerate eigenvalues}. 
  Nonlinearity, 30, no. 10, 3821--3852, 2017.
  
  	\bibitem{HR21} Hmidi T., Roulley E., \textit{Time quasi-periodic vortex patches for quasi-geostrophic shallow-water equations}, arXiv2110.13751.

\bibitem{IP-Mem-2009}
Iooss G., Plotnikov P.,
\newblock {\it Small divisor problem in the theory of three-dimensional water
  gravity waves}, 
\newblock { Mem. Amer. Math. Soc.}, 200(940):viii+128, 2009.

\bibitem{Ioo-Plo-Tol}
Iooss G., Plotnikov P., Toland J.,
{\it Standing waves on an infinitely deep perfect fluid under gravity},
{Arch. Ration. Mech. Anal.}, 177, no.3, 367--478, 2005.


\bibitem{Kappeler-Poschel} Kappeler T., P\"oschel J.,
{\it KdV \& KAM}. 
 {Springer, Berlin}, 2003. 




\bibitem{Kir} Kirchhoff G., 
{\it Vorlesungen uber mathematische Physik}, 
Leipzig, 1874.

\bibitem{KP} Kuksin S., P\"oschel J., {\it
Invariant Cantor manifolds of quasi-periodic oscillations
for a nonlinear Schr\"{o}dinger equation}, Annals of Math. 2  143, , 149-179, 1996.

\bibitem{Kuksin} Kuksin S., 
 {\it Hamiltonian perturbations of in finite-dimensional linear systems with imaginary spectrum},
Funktsional. Anal. i Prilozhen. 21, no. 3, 22--37, 95, 1987.

\bibitem{K2} Kuksin S.,
 {\it A KAM theorem for equations of the Korteweg-de Vries type},
Rev. Math. Phys. 10, no. 3, 1--64, 1998.

\bibitem{Kuksin2} Kuksin S.,
 {\it  Analysis of Hamiltonian PDEs}. 
 {Oxford Lecture Series in Mathematics and its Applications}, vol. 19. Oxford University Press, Oxford, 2000.


\bibitem{Liu-Yuan} Liu J., Yuan X.,
{\it A KAM theorem for Hamiltonian partial differential equations with
unbounded perturbations}, 
 {Commun. Math. Phys}. 307, no. 3, 629--673, 2011.

\bibitem{Love}  Love A. E. H.,
 {\it On the Stability of certain Vortex Motions},
  Proc. London Math. Soc., 25(1):18--42, 1893.


\bibitem{PlTo}
Plotnikov P.,  Toland J.,
 {\it Nash-{M}oser theory for standing water waves},
 { Arch. Ration. Mech. Anal.} 159, no. 1, 1--83, 2001.


\bibitem{Po3} P\"oschel J.,
 {\it Quasi-periodic solutions for a nonlinear wave equation},
  Comment. Math. Helv.,  71, no. 2, 269--296, 1996.

\bibitem{PP1}
Procesi M., Procesi C., 
{\it A normal form for the Schr\"odinger equation with analytic non-linearities}, Comm. Math. Phys. 312, 501-557, 2012. 

\bibitem{PP}
Procesi M., Procesi C.,
{\it A  KAM algorithm for the resonant non-linear Schr\"odinger equation}, 
Advances in Math., 399--470, 2015.


 \bibitem{PP3} Procesi C., Procesi M., {\it Reducible quasi-periodic solutions for the Non Linear
Schr\"odinger equation}, Bollettino unione Matematica Italiana. 

\bibitem{Ru1} R\"ussmann H.,
 {\it Invariant tori in non-degenerate nearly integrable Hamiltonian systems},
 {Regul. Chaotic Dyn.} 6, no. 2, 199--204, 2001.



\bibitem{Tang} Tang Y., 
{\it Nonlinear stability of vortex patches}, 
Transactions of AMS 304, no. 2, 1987. 

\bibitem{TdL}  Torres de Lizaur F., 
{\it Chaos in the incompressible Euler equation on manifolds of high dimension}, 
Invent. Math., 228:687-715, 2022.



 
\bibitem{Wan}  Wan Y.H.,
 {\it The stability of rotating vortex patches}, 
 Comm. Math. Phys. 107, no.1, 1--20, 1986. 

\bibitem{WP} Wan Y.H., Pulvirenti M.,  
{\it Nonlinear stability of circular vortex patches}, 
Comm. Math. Phys. 99, 435--450, 1985.

\bibitem{Wayne} Wayne E.,
 {\it Periodic and quasi-periodic solutions of nonlinear wave equations via KAM theory},
  Comm. Math. Phys. 127, 479--528, 1990.

\bibitem{Y} Yudovitch VI.,
{\it Non-stationary flow of an ideal incompressible liquid}, 
USSR Comput.  Math. Math. Phys. 3, 1407--1456, 1963 [transl. from: 1963 Zh. Vych. Mat. 3, 1032--1066]. 




\end{thebibliography}
\end{document}